\documentclass[a4paper]{amsart}
\usepackage{amsmath,amsthm,amsfonts,amssymb,amscd,mathrsfs}
\usepackage{psfrag,graphicx}
\usepackage{bbding}
\usepackage{epic,eepic}
\usepackage{color}
\usepackage{ebezier}

\theoremstyle{plain}

\newtheorem{Thm}{Theorem}[section]
\newtheorem{Lem}[Thm]{Lemma}
\newtheorem{Prop}[Thm]{Proposition}
\newtheorem{Cor}[Thm]{Corollary}
\theoremstyle{remark}
\newtheorem{Def}[Thm] {Definition}
\newtheorem{Rem}[Thm] {Remark}

\newtheorem{Ex}[Thm] {Example}
\newtheorem{Que}[Thm] {Question}

\long\def\begcom#1\endcom{}

\newcommand{\length}{\operatorname{\length}}

\newcommand{\Diff}{\operatorname{Diff}}

\def\Diff{\operatorname{Diff}}

\def\length{\operatorname{length}}

\def\ln{\operatorname{ln}}

\def\e{{\mathcal E}}

\def\m{{\mathcal M}}

\def\In{\ln}

\begin{document}



\title[Diffeomorphisms with Liao-Pesin set  ]
      {Diffeomorphisms with Liao-Pesin set$^{1}$  }


\footnotetext[1]{This version puts    arXiv:1004.0486 \&  arXiv:1011.6011
      (year 2010) together and updates some other  new observation.}

\footnotetext[2]{Tian is the corresponding author.}

\author[W. Sun] {Wenxiang Sun}
\address[W. Sun] {LMAM, School of Mathematical Sciences, Peking University, Beijing 100871, People's Republic of China}
\email{sunwx@math.pku.edu.cn }

\author[X. Tian] {Xueting Tian$^2$}
\address[X. Tian]{School of Mathematical Science,  Fudan University\\Shanghai 200433, People's Republic of China}
\email{xuetingtian@fudan.edu.cn}

\maketitle


\begin{abstract}

\medskip

 {
  In this paper we mainly deal with an invariant  (ergodic) hyperbolic measure $\mu$ for a  diffeomorphism $f,$ assuming that $f$ is {\it just  $C^1$} and for $\mu$ a.e. $x$, the sum of   Oseledec spaces corresponding to   negative
Lyapunov exponents     (quasi-limit-)dominates  the sum of Oseledec spaces corresponding to   positive
Lyapunov exponents at $x$. We generalize  a certain of  results of  Pesin theory from  $C^{1+\alpha}$ to the $C^1$ system
$ (f,\mu)$, including  a sufficient condition for existence of horseshoe, Livshitz theorem, exponential growth of periodic points and entropy, distribution of periodic points, periodic measures,  horseshoes, nonuniform specification and
lower semi-continuity of entropy function
 etc. These results give us more information on ergodic theory of   $C^1$ non-uniformly hyperbolic systems. In particular, they are applied  for {\it $C^1$ partially hyperbolic}  systems  whose central bundle displays some non-uniform hyperbolicity, including some robust systems.  Moreover, for some {\it $C^1$ partially hyperbolic} (not necessarily volume-preserving)  systems,  we get some information of Lebesgue measure on ``Average-nonuniform hyperbolicity" and ``volume-non-expanding".


In this process a  constructed machinery is developed for   $C^1$ (not necessarily $C^{1+\alpha}$) diffeomorphisms:   a new filtration of   Pesin blocks  is established topologically ({\it independent   on measures})  such that every block  has {\it stable manifold}  theorem and simultaneously has {\it exponential  shadowing} lemma.
The new filtration construction,  different with classical construction of  $C^{1+\alpha}$ Pesin blocks in  \cite{P1,Katok}, is  mainly  inspired from Liao's    {\it quasi-hyperbolicity}  and  so  here we  call new blocks by {\it Liao-Pesin blocks} and call the new established $C^1$ Pesin theory by {\it $C^1$ Liao-Pesin Theory}.
 Liao-Pesin set   not only exists for  hyperbolic invariant measures, 
  but also exists for general probability measures. For example, Liao-Pesin set has full measure for Lebesgue measure (not assuming invariant) in some partially hyperbolic systems.


}


\newpage
\end{abstract}

\newpage

\tableofcontents

\section{Introduction}\label{Introduction}

The study of hyperbolic dynamics began with the study of uniformly hyperbolic dynamical systems. This study was tremendously useful in the development of technical tools and insights, and in the shaping of a body of concepts suitable for the description and study of complicated dynamics. Its development of the
so-called nonuniform hyperbolicity theory by Pesin, Katok and
others was an important breakthrough, see  \cite{Katok,BP,P1} etc. The theory of nonuniformly hyperbolic dynamical systems, which was also known as Pesin Theory (or dynamical systems with nonuniformly hyperbolic behaviour), builded on the notions and paradigms from the theory of uniformly hyperbolic dynamical systems. Recently there are some generalization for Hilbert spaces, see  \cite{LY-map,LY-flow} etc.

Nonuniform hyperbolicity allows the asymptotic expansion and contraction rates to depend on the point in a way that does not admit uniform bounds, which provides a generalization that is broad enough to include a wide range of applications.
However, unlike the uniform
hyperbolicity theory, one important requirement is that the nonuniform hyperbolicity theory assumes
 not only the differentiability of the given dynamics being of
class $C^1$ but also  the first derivative satisfying  an
$\alpha$-H$\ddot{o}$lder condition for some $\alpha>0$. Thus there
appears to be  a gap between these two theories. A natural question arises:\\

{\it \bf  Whether $C^{1+\alpha}$ nonuniform hyperbolicity theory can be
established only under the $C^1$ differentiability hypothesis?}

\bigskip

 The
general answer to this problem is negative. Pugh pointed out in
 \cite{Pugh} that the $\alpha$-H$\ddot{o}$lder condition for the
first derivative is necessary in the Pesin's stable manifold theorem
and thus the $C^1$ setting   is very different from the setting  of
$C^{1+\alpha}$. However, it is  still interesting to
investigate  $C^1$
nonuniformly hyperbolic systems plus some assumptions, for example, domination condition.
There are some   advance in recent days.
One such kind result, generalized from $C^{1+\alpha}$ to $C^1$ case with dominated splitting,  is Pesin's entropy formula which is an formula in Pesin theory that the entropy of a measure that is invariant under a dynamical system is given by the total asymptotic expansion rate present in the dynamical system \cite{SunTianEntropy,CCE}.
Another is  an analog of the
Pesin's stable manifold theorem in  $C^1$ nonuniformly hyperbolic
systems with dominated Oseledec splitting  \cite{ABC}.  Even for uniformly hyperbolic case, $C^1$ and $C^{1+\alpha}$ are different.  For example, it is known that $C^{1+\alpha}$ volume-preserving Anosov diffemorphism is ergodic but it is still  unknown whether it is true just assuming  $C^1.$

\bigskip

In  present paper, we reobtain
 Katok's   closing lemma  \cite{Katok} under the hypothesis of  $C^1$ nonuniformly hyperbolic systems with limit-dominated splitting  (Definition
\ref{def:limit-dom}), which is weaker than the usual domination from topological viewpoint. Moreover, we also rebuild  a {\it exponential   shadowing}  lemma which can be as a weaker version of  Katok's  shadowing \cite{P1}.
However, exponential shadowing   is still   enough  to generalize lots of classical results in   Pesin theory
from $C^{1+\alpha}$ Pesin theory to  the $C^1$ setting with  (limit) domination, including classical results (e.g.  \cite{Katok,BP}) and recent ones (e.g.  \cite{Hir,LLS,Uga,ACW,OliTian,LST}).
These results give us more information on ergodic theory of   $C^1$ non-uniformly hyperbolic systems. These results include:

(1) existence and density of 
  periodic points, existence of horseshoe, density of periodic measures, nonuniform specification;

(2) Livshitz Theorem;

(3) average-nonuniform hyperbolicity of Lebesgue measure etc.





  In particular, a certain of results are applied  for {\it $C^1$ partially hyperbolic }  systems  whose central bundle displays some non-uniform hyperbolicity, including some robust systems.  Moreover, for some {\it $C^1$ partially hyperbolic} (not necessarily volume-preserving)  systems,  we get some information of Lebesgue measure on hyperbolicity. In particular, we obtain hyperbolic SRB-like measures in some partially hyperbolic systems.

\medskip

In this process a  constructed machinery is developed for   $C^1$ (not necessarily $C^{1+\alpha}$) diffeomorphisms:   a new filtration of   Pesin blocks  is established topologically (independent   on measures)  such that every block  has {\it stable manifold}  theorem and simultaneously has {\it exponential  shadowing} lemma.
The new filtration construction,  different with classical construction of  $C^{1+\alpha}$ Pesin blocks in  \cite{P1,Katok}, is  mainly  inspired from Liao's    {\it quasi-hyperbolicity}  and  so  here we  call new blocks by {\it Liao-Pesin blocks} and call the new established $C^1$ Pesin theory by {\it $C^1$ Liao-Pesin Theory}.


Recall that for the partially hyperbolic (not uniformly hyperbolic) systems introduced in  \cite{D2009}, Leplaideur et al proved in  \cite{LOR} for the central direction, all ergodic invariant measures only have negative exponents, with the exception of a Dirac measure supported on a saddle with positive exponent.     Moreover, the examples in  \cite{D2009} and Example \ref{Ex:nonuni-hyp-sys} below  tell us that even though every  ergodic measure is hyperbolic with a dominated Oseledec hyperbolic splitting, the system is not necessary to be uniformly hyperbolic. This implies that  our assumption in present paper is a very weak version of non-uniform hyperbolicity, since we just assume one measure to be hyperbolic with  (quasi-limit-)dominated splitting.

\section{Preliminaries}\label{section-preliminary}

Let $M$ be a compact   smooth Riemannian manifold and
let $d$ denote the distance induced by the Riemannian metric.
Denote the tangent bundle of $M$  by $TM$.
For $r\geq 1,$ denote by $\Diff^r (M)$ the space of $C^{r}$
diffeomorphisms of M.  Let $m$ be a Lebesgue measure  and
denote the set of all $C^r$  volume-preserving diffeomorphisms by $\Diff^r_m (M),\,r\geq 1$.

Let $ \m(M)$ be the set of all probability measures supported on $M.$ Given a subset $\Delta\subseteq M,$ let $\overline{\Delta}$ denotes the closure of $\Delta$. Let $\m_{f} (\Delta),\m_e (\Delta),\m_p (\Delta)$ denote the
space of invariant measures, ergodic measures, periodic measures with full measure  on $\Delta$. Here a measure $\mu$ is called {\it periodic}, if there is a
periodic point $z$ with period $p$ such that
 $\mu=\frac1p\sum_{i=0}^{p-1}\delta_{f^i (z)},$  where $\delta_x$
denotes the Dirac measure at $x$. Denote by $Per (f)$ and $Per_h (f)$ the set of all  periodic points and the subset of  hyperbolic ones,
respectively.
Let $P_n (f)$ denote the set of all periodic points with period $n$.

For any $\mu\in\m_f (M)$, let $supp (\mu)$ denote the
support of $\mu,$ the minimal compact set $\Delta\subseteq M$ such that $\mu (\Delta)=1.$
We denote by $V_f (\nu)$ the set of accumulation measures of time
averages  ${\e}_N (\nu)=\frac 1{N}{\sum_{j=0}^{N-1}f_\ast^j\nu}.$  Then $V_f (\nu)$ is a nonempty, closed and connected subset of ${\m}_{f} (M)$. And we denote by
$V_f (x)$ the set of accumulation measures of time averages
$${\e}_N (x)=\frac 1{N}{\sum_{j=0}^{N-1}\delta_{f^jx}}.$$



\subsection{Uniform Hyperbolicity, Dominated splitting and Partial Hyperbolicity}

 Denote the minimal norm of an invertible linear map $A$ by
$m (A)=\|A^{-1}\|^{-1}$. Let $\Omega (f)$ denote the non-wandering set of $f$.
Let $\Delta\subseteq M$ be an $f$-invariant   set. Let $E\subseteq T_\Delta M$ be a $Df$-invariant bundle. We say
that $E$ is uniformly contracting, if there exist $C>0$ and $0<\lambda<1$ such that $$ {\|Df^n|_{E (x)}\|} \leq C \lambda^n, \forall x\in M,\,\, n\geq 1.$$ In parallel,  we say that $E$ is uniformly expanding, if  there exist $C>0$ and $0<\lambda<1$ such that $$ {\|Df^{-n}|_{E (x)}\|} \leq C \lambda^n, \forall x\in M,\,\, n\geq 1.$$

 We say that $\Delta$ to be a {\it hyperbolic set}, if there is  a $Df-$invariant splitting $T_{\Delta}M=E^s\oplus E^u$  on
$\Delta$ such that  $E^s$ is is uniformly contracting and $E^u$ is uniformly expanding.
Here the splitting $T_{\Delta}M=E^s\oplus E^u$ is called {\it hyperbolic splitting.}
If $\Delta=M$ is hyperbolic, then the system $f$ is called {\it Anosov}.  If the non-wandering set $\Omega (f)$ is hyperbolic and   the periodic points is dense in $\Omega (f),$ then $f$ is called {\it Axiom A.}
 Recall that a {\it horseshoe}  for a diffeomorphism $f$ is a transitive, locally maximal  (or called isolated)
hyperbolic set $\Lambda$, that is totally disconnected and not finite  (such a set must be
perfect, hence a Cantor set).

It is known that a hyperbolic splitting $T_{\Delta}M=E^s\oplus E^u$ is always unique, continuous and can be extended on the closure of $\Delta$ and even its neighborhood.

We
recall the notion of 
 dominated splitting.


\begin{Def}\label{def:dominated}
A $Df-$invariant splitting $T_{\Delta}M=E\oplus F$ is called to be  dominated on $\Delta$,  if there exists $C>0,0<\lambda<1$ such that
 \begin{eqnarray}\label{eq-domination-def}
\frac
{\|Df^n|_{E (x)}\|}{m (Df^n|_{F (x)})}\leq C \lambda^n,\,\,\,\,\forall \,\,n\geq 1,x\in \Delta.
\end{eqnarray}
We write   $E\prec F$.
\end{Def}


  This definition may be formulated, equivalently, as follows: there exists $S\in \mathbb{Z}^+$ and $\bar{\lambda}\in (0,1)$ such that
\begin{eqnarray}\label{eq-domination-equvialent1}
 \frac
{\|Df^n|_{E (x)}\|}{m (Df^n|_{F (x)})}\leq   \bar{\lambda}^n,\,\,\forall\,\, x \in
\Delta,\,\forall \,\,n\geq S.
\end{eqnarray}
Another equivalent statement of dominated splitting is that:
there exists $L\in \mathbb{Z}^+$ and $\bar{\bar{\lambda}}\in (0,1)$ such that
\begin{eqnarray}\label{eq-domination-equvialent2}
 \frac
{\|Df^L|_{E (x)}\|}{m (Df^L|_{F (x)})}\leq   \bar{\bar{\lambda}}^L,\,\,\forall\,\, x \in
\Delta.
\end{eqnarray}
Remark that Gourmelon  ( \cite{Gour}) proved that there always exists an adapted metric for which
$C = 1$ in  (\ref{eq-domination-def}). For this adapted metric,  we have $S=1$ in  (\ref{eq-domination-equvialent1}) and $L=1$ in   (\ref{eq-domination-equvialent2}).

Remark that for any hyperbolic set $\Delta\subseteq M$, the hyperbolic splitting $T_{\Delta}M=E^s\oplus E^u$ is always dominated.

It is known  that dominated splitting  is always continuous and  can be extended to its closure. Moreover, such a splitting can be extended in a dominated way to the maximal invariant set in a neighborhood of $\overline{\Delta}$ (for example, see  \cite [Appendix B.1 , P.287-P.290] {BLV}).    For convenience, we introduce a concept called dominated-$\epsilon$-neighborhood. We say that $\Delta$ has a {\it dominated-$\epsilon$-neighborhood}, if   there is   $\epsilon>0$  such that the dominated splitting is extended and dominated on $\Gamma:=\cap_{n\in \mathbb{Z}}f^n (U)$ where   $U=B ( \overline{\Delta},\epsilon)  $ is the $\epsilon-$neighborhood of  $\overline{\Delta}.$

Recall that  partial hyperbolicity usually means that there exists a splitting in three subbundles such that one  (which is called the unstable bundle) is uniformly expanding, one (which is called the stable bundle) is uniformly contracting and the other one  (which is called the center bundle) may have no hyperbolicity  but is dominated by the unstable bundle and dominates the stable one.
  Here we adopt a  more general notion of partial hyperbolicity, by using a splitting into two subbundles:

 \begin{Def}
 \label{DefinitionPartialHyperbolicity} {\bf  (Partial hyperbolicity)}  \em
 We call a dominated splitting $T_\Delta M=E\oplus F$  to be   {\it partially hyperbolic,}  if  $E$ is uniformly contracting
 or $F$ is uniformly expanding. Here $\Delta$ is called a partially hyperbolic set.  In particular, if $\Delta=M,$ we say $f$ to be partially hyperbolic.
   \end{Def}

 Let $A^r(M)$, $PH^r (M)$ and $D^r(M)$ denote the spaces of all $C^r$ Anosov diffeomorphisms,  partially hyperbolic diffeomorphisms  and diffeomorphisms with (global) dominated splitting respectively
    where $r\geq 1$.   Note that  $A^r(M)\subseteq PH^r (M)\subseteq D^r(M).$

Recall that a $Df$ invariant subbundle $G\subset TM$ is called to be {\it quasi-conformal}, if for any $\epsilon>0$, there exists $C_\epsilon>0$ such that
  for any $x\in M$ and $n\geq 1,$  $$C_\epsilon^{-1}e^{-n\epsilon} \leq \frac{\|Df^n|_{G (x)}\|}{m (Df^n|_{G (x)})}\leq C_\epsilon e^{n\epsilon}.$$
  The left inequality is trivial if  take $C_\epsilon \geq 1,$ since  $\frac{\|Df^n|_{G (x)}\|}{m (Df^n|_{G (x)})}\geq 1. $
In particular, quasi-conformal condition  implies that  for any $x\in M,$
$$\limsup_{n\rightarrow +\infty}\frac1n\log {\|Df^{\pm n}|_{G (x)}\|}=\limsup_{n\rightarrow +\infty}\frac1n\log m (Df^{\pm n}|_{G (x)}),$$
$$\liminf_{n\rightarrow +\infty}\frac1n\log {\|Df^{\pm n}|_{G (x)}\|}=\liminf_{n\rightarrow +\infty}\frac1n\log m (Df^{\pm n}|_{G (x)}).$$ It is obvious that every one dimensional
$Df$ invariant subbundle $G\subset TM$ is quasi-conformal.

It is known that for any hyperbolic set $\Lambda,$   \begin{eqnarray}\label{eq-growth-period-hyperbolic}
    h_{top} (f|_\Lambda)\geq\limsup_{n\rightarrow \infty }\frac 1n \log \# P_n (f)\cap \Lambda .
      \end{eqnarray}
      It is not difficult to prove because $f|_\Lambda$ is expansive and  $ \# P_n (f)\cap \Lambda$ is a particular separated set. In general for any hyperbolic set $\Lambda,$
       \begin{eqnarray*}
    h_{top} (f|_\Lambda)\leq\limsup_{n\rightarrow \infty }\frac 1n \log \# P_n (f).
      \end{eqnarray*}
      Moreover, for any horseshoe $\Lambda$ (or general isolated hyperbolic set, for example see \cite{Bowen-book}),   $$  h_{top} (f|_\Lambda)= \limsup_{n\rightarrow \infty }\frac 1n \log \# P_n (f)\cap \Lambda .$$
           From uniform hyperbolicity, there is some $C_0>0, \zeta_0\in (0,1)$ such that for any $C>C_0,\zeta\in  (0,\zeta_0)$  such that for all $n\geq 1,$ $P_n(f)= P_n (f,\zeta, C),$
            where   \begin{eqnarray*}
  & &P_n (f,\zeta, C):=\{x\in M\,|\,f^nx=x
   \text{ and for any }l\geq 1, j=0,1,2,\cdots, n-1,\,\,\,\,\,\,\,\,\,
   \,\,\,\,\,\,\,\,\,\,\\\nonumber
   & & \,\,\,\,\,\,\,\,\,\,\,\,\,\,\,\,\,\,\,\,\,\,\,\,\,\,\,\,\,\,\,\,
   \,\,\,\,\,\,\,\,\,\,\,\,\,\,
  \,\max\{\|Df^l|_{E (f^jx)}\|,\|Df^{-l}|_{F (f^jx)}\|\}\leq C \zeta^l \},
  \end{eqnarray*}
$E,F$ denote the stable bundles and the unstable bundles, respectively.
So for any horseshoe $\Lambda$ (or general isolated hyperbolic set),   there is some $C_0>0, \zeta_0\in (0,1)$ such that for any $C>C_0,\zeta\in  (0,\zeta_0)$
such that
        \begin{eqnarray} \label{eq-horseshoe-periodicgrowth}
           h_{top} (f|_\Lambda)= \limsup_{n\rightarrow \infty }\frac 1n \log \# P_n (f, \zeta, C)\cap \Lambda.
           \end{eqnarray}

\subsection{Oseledec Theorem, Hyperbolic Measures}


Let $f\in \Diff^1 (M)$ and $\mu\in\mathcal{M}_{f} (M).$
By the Oseledec Multiplicative Ergodic Theorem  \cite{Os}, there is  a Borel set $L (\mu)$ satisfying $fL (\mu)=L (\mu)$
and $\mu (L (\mu))=1$, called Oseledec basin of $\mu$, such that for every $x\in L (\mu),$ there exist

 (a)  real numbers, called Lyapunov exponents,
\begin {equation} \label{eq:Lya-exps}
 \lambda_{1} (x)
<\lambda_{2} (x)< \cdot\cdot\cdot< \lambda_{W (x)} (x)\,\, (W (x) \leq dim (M))
\end {equation}

 (b)  positive integers $m_{1} (x),\,\, \cdot\cdot\cdot,\,\,\,\,m_{W (x)} (x)$, satisfying $m_{1} (x)+
\cdot\cdot\cdot +m_{W (x)} (x)=dim (M)$;

 (c)  a measurable splitting
$T_{x}M=E_{x}^{1}\oplus\cdot\cdot\cdot\oplus E_{x}^{W (x)}$
   with dim$E_{x}^{i}=m_{i} (x)$ and $Df (E_{x}^{i})=E_{f (x)}^{i}$,
such    that
$$\lim_{n\rightarrow \pm\infty}\frac{\log\|Df^{n}v\|}{n}=\lambda_{i} (x),
   $$
   with  uniform convergence  on
    $\{ v \in E_{x}^{i}\,|\,\,\|v\|=1\},\,\, i=1,\, 2,\, \cdots,\,
   W (x)  $. \\




For convenience,  we say every point $x\in L (\mu)$ to be {\it Lyapunov-regular}.   For $x\in L (\mu),$ let  $E^s (x)$ (called {\it stable bundle}) denote the bundle composed by the Oseledec bundles whose Lyapunov exponents are negative,   let $E^u (x)$ (called {\it unstable bundle}) denote the bundle composed by the Oseledec bundles whose Lyapunov exponents are positive and respectively, $E^0 (x)$ (called {\it quasi-identity  bundle}) denote the bundle composed by the Oseledec bundles whose Lyapunov exponents are zero.  In general $dim (E^s (x))$ is called the index of $x$. In particular, if $\mu$ is ergodic, $\lambda_{i} (x),m_i (x),\,W (x),dim (E^s (x))$ are constants for $\mu$ a.e. $x.$
 Define the {\it integrable Lyapunov exponents} of $\mu$ as $$\lambda_i (\mu)=\int \lambda_i (x) d\mu.$$ If $\mu$ is ergodic, $\mu$ a.e. $x,$ $\lambda_{i} (x)=\lambda_i (\mu)$. In other words, for ergodic $\mu$, its Lyapunov exponents are same as its integrable Lyapunov exponents.


Now we recall the notion of hyperbolic invariant measure. Let $\mu\in\mathcal{M}_{f} (M).$

\begin{Def}\label{def:hyp-measFixedindex}  We call $\mu$ to be  hyperbolic (with a fixed index), if there exists an invariant set $\Delta\subseteq L (\mu)$ with $\mu (\Delta)=1$ and two positive integers $n_1,n_2$ with $n_1+n_2=dim (M)$ such that for any $x\in \Delta,$

 (1) none of the Lyapunov exponents of  $x$  are  zero (i.e., $dim (E^0 (x))=0$);

 (2) there exist Lyapunov exponents of $x$  with different signs (i.e., $dim (E^s (x))\cdot dim  (E^u (x))  \neq 0$);

 (3)    $dim (E^s (x))=n_1$ (called index of $\mu$, denoted by $ind (\mu)$) and $dim (E^u (x))=n_2$.
\end{Def}

Denote the space of hyperbolic invariant measures supported a set $\Delta$ by $\mathcal{M}^{h}_{f} (\Delta).$ For any $\mu\in \mathcal{M}^{h}_{f} (\Delta),$ we call $E^s\oplus E^u$ to be {\it Oseledec's hyperbolic splitting} (simply, hyperbolic splitting).
Recall that the original hyperbolic measure just requires  the conditions  (1) and  (2) but in present paper for convenience of statements, we further require fixed index. If considering the original hyperbolic measure, one can write it by a finite convex sum  of at most $dim (M)$ different hyperbolic measures with fixed index.  More precisely,  let $$\Delta_i=\{x\in \Delta|\, dim (E^s (x))=i\},\,\,\,\,\text{ and }\, \Gamma=\{i|\,1\leq i\leq dim (M),\, \mu (\Delta_i)>0\}.$$ For $i\in \Gamma,$ define $\mu_i=\mu|_{\Delta_i}$, then $$\mu=\sum_{i\in \Gamma} \mu (\Delta_i)\mu_i.$$


We say  a $\mu\in\m_f (M)$ to be {\it uniformly hyperbolic}, if $supp (\mu)$ is a hyperbolic set.  Denote the space of uniformly hyperbolic invariant measures by $\mathcal{M}^{uh}_{f} (M).$
 Denote the set of all hyperbolic periodic measures by $\m_p^{h} (M)$.   It is obvious that $$ \m_p^{h} (M)\subseteq \m_f^{uh} (M)\subseteq \m_f^{h} (M).$$
We say a measure $\mu$ is {\it nonatomic} (in  \cite{Katok1}, it is also called continuous  measure),  if for any point $x\in M$, $\mu (\{x\})=0.$ Denote by $\m_f^n (M)$ the set of all nonatomic invariant measures.
 Let $\m_f^{+} (M)$ denote the set of all invariant measures with positive metric entropy.  Remark that $\m_f^+ (M)\cap\m_e (M)\subseteq \m_e (M)\setminus \m_p (M)=\m_f^n (M)\cap \m_e (M).$

Remark that for any hyperbolic set $\Delta\subseteq M$ and any $\mu\in\m_f (\Delta)$,  $\mu$ is uniformly hyperbolic, its all   Lyapunov exponents are far from zero and its Oseledec hyperbolic splitting coincides with the corresponding hyperbolic splitting $T_{\Delta}M=E^s\oplus E^u$. In other words, for a hyperbolic set $\Delta\subseteq M$, we have $$\m_f (\Delta)=\m_f^h (\Delta)=\mathcal{M}^{uh}_{f} (\Delta). $$
\subsection{Dominated Splitting \& Limit-dominated Splitting}

Let $x\in M$ and $T_{Orb (x)}M=E\oplus F$ be a $Df-$invariant
splitting on the orbit  of $x$, denoted by $Orb (x)$.  We introduce dominated splitting at one point (or one orbit).

\begin{Def}\label{def-dominated-point}
We call $T_{Orb (x)}M=E\oplus F$ to be dominated at $x$,
if    $T_{Orb (x)}M=E\oplus F$ is  dominated on the set
 $\Delta=Orb (x).$
\end{Def}

Note that dominated splitting  at different orbits admit different  $S$ and $\lambda$. So we can introduce a following weaker notion to generalize  original  dominated splitting. Let $\Delta$ be an $f-$invariant set
and $T_{\Delta}M=E\oplus F$ be a $Df-$invariant splitting on
$\Delta$.

\begin{Def}\label{def-quasi-dominated}
$T_{\Delta}M=E\oplus F$ is called to be  quasi-dominated on $\Delta$, if   for any $x\in \Delta$, $T_{Orb (x)}M=E\oplus F$ is    dominated     at $x$.
\end{Def}


In the study to find dominated splitting, quasi-dominated splitting may be one medium step.  For example,  recall a result of   \cite{BV}  that for a $C^1$ generic volume
preserving diffeomorphism and Lebegue a.e. $x,$ its
Oseledec splittings is either trivial (i.e., all
Lyapunov exponents are zero) or dominated at $x$ (i.e., the sum of two subbundles $E_{x}^{i}\oplus E_x^{i+1}$  in Oseledec splitting  is dominated for all $i$).  On the other hand, this splitting is dependent on a measure. So we also  introduce dominated splitting and quasi-dominated splitting with respect to a  measure. Let  $\mu\in \mathcal{M} (M)$ and let $\Delta$ be an $f-$invariant set with $\mu (\Delta)=1$. Let $T_{\Delta}M=E\oplus F$ be a $Df-$invariant splitting on
$\Delta$.

\begin{Def}\label{def:dominated-meas}
We say that $T_{\Delta}M=E\oplus F$ is a $\mu-$dominated splitting (or $\mu-$quasi-dominated splitting), if there is   an $f-$invariant set $\Delta'\subseteq \Delta$ with  $\mu (\Delta')=1$   such that  $T_{\Delta'}M=E\oplus F$ is dominated (or quasi-dominated) on $\Delta'.$ \\
 In particular, if $\mu\in \mathcal{M}^{h}_{f} (M),$    we say that the Oseledec's hyperbolic splitting of $\mu$ is $\mu-$dominated (or $\mu-$quasi-dominated), if there is   an $f-$invariant set $\Delta'\subseteq L (\mu)$ with  $\mu (\Delta')=1$   such that  $T_{\Delta'}M=E^s\oplus E^u$ is dominated (or quasi-dominated) on $\Delta'.$
\end{Def}
Since dominated splitting     can be  extended to the closure,  then   there is   a $\mu-$dominated splitting on a set
 with full measure  $\Leftrightarrow$ there is a $\mu-$dominated splitting on $supp (\mu).$ However, it is unknown for the case of quasi-dominated splitting, except that $\mu$ is ergodic,  see section \ref{section-limitdomination}.

Now we start to introduce another similar
notion as domination, called limit domination. Let $\Delta$ be an $f-$invariant set  and $T_{\Delta}M=E\oplus F$ be
a $Df-$invariant splitting on $\Delta$.

\begin{Def}\label{def:limit-dom}
$T_{\Delta}M=E\oplus F$ is limit-dominated, if there
exists $S\in \mathbb{Z}^+,\,\lambda\in (0,1)$  such that
$$\limsup_{l\rightarrow+\infty}  \frac
{\|Df^S|_{E (f^{lS} (x))}\|}{m (Df^S|_{F (f^{lS} (x))})}\leq  \lambda^S,\,\,\forall
x \in \Delta.$$
We write   $E\prec^l F$.
\end{Def}

Since $\Delta$ is
$f-$invariant, one has
\begin{eqnarray}\label{def-equal-limit-domination}
\begin{split}&\limsup_{l\rightarrow+\infty} \frac
{\|Df^S|_{E (f^{lS} (x))}\|}{m (Df^S|_{F (f^{lS} (x))})}\leq  \lambda^S,\,
\forall\, x\in \Delta\, \\
&  \Leftrightarrow\,
 \limsup_{l\rightarrow+\infty} \frac
{\|Df^S|_{E (f^{l} (x))}\|}{m (Df^S|_{F (f^{l} (x))})}\leq  \lambda^S,\,
\forall\, x\in \Delta.
\end{split}
\end{eqnarray}
In another equivalent way, there is some $S\geq 1$ and $\zeta>0$ in the sense that
$$\limsup_{l\rightarrow+\infty}  \frac1S \log\frac
{\|Df^S|_{E (f^{l} (x))}\|}{m (Df^S|_{F (f^{l} (x))})}\leq  -2\zeta ,\,\,\forall
x \in \Delta.$$




In parallel, we introduce limit-domination for one point (or orbit).

  \begin{Def}\label{def-limitdominated-point}  We say that a $Df-$invariant splitting $T_{Orb (x)}M=E\oplus F$ to be limit-dominated at $x$, if    $T_{Orb (x)}M=E\oplus F$ is  limit-dominated on the set $\Delta=Orb (x).$

  \end{Def}

  By sub-multiplication of norms,  that $T_{Orb (x)}M=E\oplus F$ is  limit-dominated at $x$ is equivalent that
$$  \liminf_{S\rightarrow+\infty} \limsup_{l\rightarrow+\infty}  \frac1S \log\frac
{\|Df^S|_{E (f^{l} (x))}\|}{m (Df^S|_{F (f^{l} (x))})}<0.
$$
 Moreover, let us  introduce quasi-limit-domination.
 Let $\Delta$ be an $f-$invariant set
and $T_{\Delta}M=E\oplus F$ be a $Df-$invariant splitting on
$\Delta$.

\begin{Def}\label{def-quasi-limit-dominated}
$T_{\Delta}M=E\oplus F$ is called to be  quasi-limit-dominated on $\Delta$, if for any $x\in \Delta$, the   splitting  $T_{Orb (x)}M=E\oplus F$ is limit-dominated at $x$.

\end{Def}

Now let us introduce limit-domination and quasi-limit-domination for a measure  $\mu\in \mathcal{M} (M)$.
Let $\Delta$ be an $f-$invariant set with $\mu (\Delta)=1$
and $T_{\Delta}M=E\oplus F$ be a $Df-$invariant splitting on
$\Delta$.

\begin{Def}\label{def-limit-dominated-measure}
 We say that $T_{\Delta}M=E\oplus F$ is a $\mu-$limit-dominated splitting (resp., $\mu-$quasi-limit-dominated splitting), if there is   an $f-$invariant set $\Delta'\subseteq \Delta$ with  $\mu (\Delta')=1$   such that  $T_{\Delta'}M=E\oplus F$ is limit-dominated (resp., quasi-limit-dominated) on $\Delta'.$ \\
 In particular, if $\mu\in \mathcal{M}^{h}_{f} (M),$    we say that the Oseledec's hyperbolic splitting of $\mu$ is $\mu-$limit-dominated (resp., $\mu-$quasi-limit-dominated), if there is   an $f-$invariant set $\Delta'\subseteq L (\mu)$ with  $\mu (\Delta')=1$   such that  $T_{\Delta'}M=E^s\oplus E^u$ is limit-dominated (resp., quasi-limit-dominated) on $\Delta'.$
  \end{Def}

For convenience, for an invariant set $\Delta\subseteq M,$ let $$\mathcal{M}^{dh}_{f} (\Delta):=\{ \mu\in \mathcal{M}^{h}_{f} (\Delta)|  \text{   hyperbolic splitting of }\mu \text{ is } \mu-\text{dominated}\}$$ and $$\mathcal{M}^{qdh}_{f} (\Delta):=\{ \mu\in \mathcal{M}^{h}_{f} (\Delta)|  \text{   hyperbolic splitting of }\mu \text{ is } \mu-\text{quasi-dominated}\}.$$
     Define $$\mathcal{M}^{ldh}_{f} (\Delta):=\{ \mu\in \mathcal{M}^{h}_{f} (\Delta)|  \text{
    hyperbolic splitting of }\mu \text{ is } \mu-\text{limit-dominated}\}$$ and
   $\mathcal{M}^{qldh}_{f} (\Delta):= $
   $$\{ \mu\in \mathcal{M}^{h}_{f} (\Delta)|  \text{   hyperbolic splitting of } \mu \text{ is }  \mu-\text{quasi-limit-dominated}\}.$$ It   is easy to see that $  \mathcal{M}^{h}_{f} (\Delta)\supseteq \mathcal{M}^{qldh}_{f} (\Delta)  \supseteq \mathcal{M}^{ldh}_{f} (\Delta)\cup \mathcal{M}^{qdh}_{f} (\Delta)\supseteq \mathcal{M}^{ldh}_{f} (\Delta)\cap \mathcal{M}^{qdh}_{f} (\Delta)\supseteq \mathcal{M}^{dh}_{f} (\Delta)\supseteq\mathcal{M}^{uh}_{f} (\Delta).$
 Recall that for a  hyperbolic set
 $\Delta\subseteq M,$  the hyperbolic splitting $T_{\Delta}M=E^s\oplus E^u$ is always dominated,
  every invariant measure $\mu\in\m_f (\Delta)$ is hyperbolic,  its Oseledec hyperbolic splitting
  corresponds to the hyperbolic    splitting and so is dominated.
So for any hyperbolic set  $\Delta\subseteq M $,
 $$ \mathcal{M}^{qldh}_{f} (\Delta) = \mathcal{M}^{ldh}_{f} (\Delta)= \mathcal{M}^{qdh}_{f} (\Delta)= \mathcal{M}^{dh}_{f} (\Delta)=\mathcal{M}^{uh}_{f} (\Delta)=\mathcal{M}^{h}_{f} (\Delta)=\mathcal{M}_{f} (\Delta).$$



For any invariant measure $\mu$, we say its  Oseledec's   splitting  is dominated (or quasi-dominated, limit-dominated, quasi-limit-dominated), if there is an $f-$invariant set $\Delta\subseteq L (\mu)$ with  $\mu (\Delta)=1$   such that   for all $i$  the sum of two subbundles in Oseledec splitting $$E_{x}^{i}\oplus E_x^{i+1}$$    is dominated (or quasi-dominated, limit-dominated, quasi-limit-dominated) on  $\Delta$.
 It is obvious for any $\mu\in \mathcal{M}^{h}_{f} (M),$ if its Oseledec's splitting is dominated (or quasi-dominated, limit-dominated, quasi-limit-dominated), then  so does its Oseledec's hyperbolic splitting.


  Recall that for Anosov systems and  Axiom A systems,   they coincide on a set with totally full measure (i.e., being equivalent in the sense of probabilistic perspective) but they are different from geometric sense that every Anosov system carries a {\it global}  hyperbolic splitting.
  Similarly,  we will show  (quasi-)domination and  (quasi-)limit-domination are equivalent  in the sense of probabilistic perspective in last section. However, they are different from topological or dimensional  viewpoint.

  More precisely, by  (\ref{eq-domination-equvialent2})  (quasi-)dominated spitting always implies  (quasi-)limit-dominated splitting but the inverse is unknown.   On the other hand, recall  that dominated splitting is an open condition that it  can be extended to the closure, even neighborhoods and dominated splitting is always continuous.  However,
in Definition~\ref{def:limit-dom} it is unknown whether the limit-dominated  splitting $T_{x}M=E (x)\oplus F (x)$ is
  continuous on $\Delta$ and can be extended to the closure of $\Delta$ and neighborhoods.

    Limit domination only requires that $E$ can dominate $F$ for  large enough positive iterate  of the orbit (e.g., see a simple but extreme  example:
Example~\ref{Ex:nonuni-hyp-sys}) and it is enough to construct a topological definition of $C^1$ new Pesin set  (independent of measures) that carries shadowing lemma. This implies we can realize shadowing on a new Pesin set as large as possible. Recall that irregular set (the set of points that Birkhoff average does not converge) always has zero measure for any invariant measure but for many systems (including Axiom A systems and some non-hyperbolic ones), it carries full topological entropy (Bowen's dimensional definition) and so it is important from dimensional perspective. So  limit domination may admit   the new Pesin set to  contain more   points  (for example, irregular points)  which may be useless in probabilistic  perspective but may have other important information.

Furthermore, the left limit in Definition \ref{def:limit-dom}
 is more convenient to connect Lyapunov exponents of
Birhorff average (Proposition \ref{Prop:relation-limdom-Lyaexp}). Another observation is that a {\it global}  dominated splitting is important to obtain entropy formula for SRB-like measures but $\mu$-dominated splitting (w.r.t a measure $\mu$) is not enough, see  \cite{CCE} for more details.
So for a {\it global limit-domination},  more information from topological viewpoint than $\mu$-dominated splitting, we guess that maybe entropy formula still holds for SRB-like measures (in future work for consideration). So we prefer to introduce the concept of limit domination in present paper. More
discussion will appear in Section \ref{section-limitdomination}.

\subsection{Average-nonuniform Hyperbolicity}\label{section-average-hyperbolic}


 We introduce a notion of
\emph{degree of average-nonuniform hyperbolicity}. Let $K\in\mathbb{N},\,\,\zeta>0$. For a given  $f$-invariant subset
  $\Delta$, let
$T_{x}M=E (x)\oplus  F (x),\,\,x\in \Delta$ be a $Df$-invariant
splitting.
\begin{Def}\label{def:degree-Pesinset}
 We call $\Delta$ an average-nonuniformly hyperbolic set with
$ (K,\,\zeta)$-degree corresponding to $T_{x}M=E (x)\oplus  F (x)$, if for $ \forall x \in\Delta ,$ one has
$$\limsup_{l\rightarrow
+\infty}\sum_{j=0}^{l-1}\frac{\log\|Df^{K}|_{E ({f^{
jK} (x)})}\|}{lK}\leq- \zeta,\,\,\,\,\,\,and$$
$$ \liminf_{l\rightarrow +\infty}\sum_{j=-l}^{-1}\frac{\log
m (Df^{K}|_{F ({f^{ jK} (x)})})}{lK}\geq\zeta.$$
We say a probability measure $\mu$ (not necessarily invariant) to be average-nonuniformly hyperbolic, if there is an invariant set $\Delta$ with
$\mu (\Delta)=1$, $K\in\mathbb{N},\,\,\zeta>0$ such that $\Delta$ is  average-nonuniformly hyperbolic  with
$ (K,\,\zeta)$-degree corresponding to some $Df$-invariant
splitting $T_\Delta M=E\oplus F.$

\end{Def}

By the sub-multiplications of the norms, it is easy to see if an {\it invariant} measure is average-nonuniformly hyperbolic, then it is also hyperbolic. Moreover, the inverse is also true, see Lemma \ref{Lem:Generalized-Multi-erg-thm}.  For any Anosov diffeomorphism, the whole space is average-nonuniformly hyperbolic and so every probability measure is average-nonuniformly hyperbolic. In particular, Lebesgue measure is average-nonuniformly hyperbolic but Lebesgue measure
 may be not average-nonuniformly hyperbolic for general Axiom A systems.  Here we interest on average-nonuniform  hyperbolicity of  Lebesgue measure, which is possibly useful to find SRB (physical) measures.

\subsection{Entropy}


Let $(X,d)$ be a compact
metric space with Borel $\sigma-$algebra $\mathfrak{B}(X)$ and let
$f:X\rightarrow X$ be a continuous map. Let $\mathcal{M}(X)$ denote the space of all probability measures supported on $X$.  The set of all invariant
measures and the set of all ergodic invariant measures are denoted
by $\mathcal{M}_f(X)$ and $\mathcal{M}_f^e(X)$, respectively.
 For $x,y\in X$ and  $n\in \mathbb{N}$, let
$$d_n(x,y)=\max_{0\leq i\leq n-1}d(f^i(x),\,f^i(y)).$$   Let $x\in X$. The dynamical open ball $B_n (x,\varepsilon)$  and dynamical closed ball $\overline{B}_n (x,\varepsilon)$ are defined respectively  as
 $$B_n (x,\varepsilon):=\{y\in X|\, d_n(x,y) <\varepsilon\},\overline{B}_n (x,\varepsilon):=\{y\in X|\, d_n(x,y) \leq \varepsilon\}.$$  A set $S$ is $(n,\varepsilon)$-separated for $Z$ if $S\subset Z$ and $d_n(x,y)>\varepsilon$ for any $x,y\in S$ and $x\neq y$. A set $S\subset Z$ if $(n,\varepsilon)$-spanning for $Z$ if for any $x\in Z$, there exists $y\in S$ such that $d_n(x,y)\leq \varepsilon$.

 We have the following definition of entropy for compact set and thus the definition of entropy for a general subset.
\begin{Def}\label{compact-set}
For $E\subset X$ compact, we have the following Bowen's definition of topological entropy (c.f.  \cite{Walters}).
\begin{equation}\label{Bowen-entropy}
  h_{top}(f,E)=\lim_{\varepsilon\to 0}\lim_{n\to \infty}\frac{\log s_n(E,\varepsilon)}{n},
\end{equation}
where $s_n(E,\varepsilon)$ denotes the maximal cardinality of set which is  $(n,\varepsilon)$-separated for $E.$
For a general subset $Y\subset X$, define
\begin{equation}\label{entropy}
  h_{top}(f,Y)=\sup\{h_{top}(f,E):E\subset Y ~\textrm{is compact}\}.
\end{equation}
Finally, we put $h_{top}(f)=h_{top}(f,X)$. Since $X$ is a compact metric space, the definition depends only on the topology
on $X$, i.e. it is independent of the choice of metric defining the same topology on $X$.
\end{Def}

Let $\mu\in\mathcal{M}_f(X)$. Given $\xi=\{A_1,\cdots,A_k\}$ a finite measurable partition of $X$,
 i.e., a disjoint collection of elements of $\mathfrak{B}(X)$ whose
 union is $X$,  we define the entropy of $\xi$ by
 $$H_{\mu}(\xi)=-\sum_{i=1}^{k}\mu(A_i)\log \mu(A_i).$$
The metric entropy of $f$ with respect to $\xi$ is given by
 $$h_{\mu}(f,\,\xi)=\lim_{n\rightarrow\infty}\frac{1}{n}\log H_{\mu}(\bigvee_{i=0}^{n-1}f^{-i}\xi).$$
 The metric entropy of $f$ with respect to $\mu$ is given by
 $$h_{\mu}(f)=\sup_{\xi}h_{\mu}(f,\xi),$$
 where $\xi$ ranges over all finite measurable partitions of $X$.

Remark that entropy is   a classical concept to describe the dynamical complexity:  larger entropy denotes  more stronger complexity.


 We say $f$ is {\it entropy-hyperbolic,}  if for any $\epsilon>0,$ there is a horseshoe $H_\epsilon$ such that
  $$h_{top} (H_\epsilon) >h_{top} (f)-\epsilon.$$

\section{Results}\label{Section-Results}

\subsection{Periodic Orbits, Periodic  Measures and Horseshoe} \label{sub-density of periodic pts and meas}




In $C^{1+\alpha}$ Pesin theory, there are several  basic  results.
 One is about approximating properties of hyperbolic periodic points and measures.

\begin{Thm}\label{Thm:Horshoe-PosiStaMnfd}   Let $f\in \Diff^1 (M)$ and   $\mu\in \mathcal{M}^{qldh}_{f} (M)\cap \m_f^n (M).$     Then there exists a hyperbolic periodic point with homoclinic point  (implying the existence of  horseshoe and positive topological entropy) such that the closure of its global unstable manifold has positive measure. Moreover, the support of the measure $\mu$ is contained in the closure of such hyperbolic periodic points, that is, $$\textrm{supp} (\mu)\subseteq \overline{Per_h (f)}.$$
\end{Thm}

\begin{Rem}\label{Rem:Horshoe-PosiStaMnfd} Here $\mu$ is not necessary to be ergodic.  In the $C^{1+\alpha}$ case, the existence of horseshoes can be found in
 \cite{Katok,Katok1,P1} (for example, see  \cite{Katok1}, Theorem S.5.1).
 Some related results of existence of periodic points with homoclinic points  for an \textit{ergodic} hyperbolic measure with dominated splitting in  $C^1$ setting  also appeared in  \cite{Gan2,Yang}. In the $C^{1+\alpha}$ case, the result that the closure of unstable manifold of the shadowing periodic
point has positive measure can be found in  \cite{Uga}. Here Theorem \ref{Thm:Horshoe-PosiStaMnfd} generalizes the result of  \cite{Uga} to some $C^1$ case.
\end{Rem}

\begin{Rem}
For $C^{1+\alpha}$ case, there is a result  ( \cite{Katok1}, Corollary S.5.2 on Page 694) that if $\mu$ is an ergodic measure with all Lyapunov exponents negative, then $\mu$ is concentrated on the orbit of a periodic sink $p$, that is, $\exists m>0$ such that $supp (\mu)=\{p,f (p),\cdots,f^{m-1} (p)\}.$ Its proof relies on a technique called regular neighborhood which is established for $C^{1+\alpha}$ case. However, for our present $C^1$ case, it is still unknown.

\end{Rem}
\medskip

We will prove Theorem \ref{Thm:Horshoe-PosiStaMnfd} in section \ref{proof-horseshoe-periodic-dense}. By Theorem \ref{Thm:Horshoe-PosiStaMnfd}, clearly one has
\begin{Cor}\label{Cor:Horshoe-PosiStaMnfd-imply}
$$ \mathcal{M}^{qldh}_{f} (M)\cap \m_f^n (M)\neq \emptyset \Leftrightarrow \mathcal{M}^{uh}_{f} (M)\cap \m_f^n (M)\neq \emptyset$$$$ \Leftrightarrow \mathcal{M}^{uh}_{f} (M)\cap \m_f^n (M)\cap\m_f^+ (M)\neq \emptyset.$$
\end{Cor}







\bigskip

By (\ref{eq-horseshoe-periodicgrowth}) we have a direct corollary of Theorem \ref{Thm:Horshoe-PosiStaMnfd}.
 \begin{Cor}\label{Cor:1111111111Horshoe-numbergrowth-periodic} Let $f\in \Diff^1 (M)$ and   $\mu\in \mathcal{M}^{qldh}_{f} (M)$   be  a  nonatomic     measure. Then \\
  (1) $f$ has a compact invariant set $\Lambda$ such that $f|_\Lambda$ is a horseshoe.\\
  (2)  there is some $ C_0>0, \zeta_0\in (0,1)$ such that for any $\zeta\in  (0,\zeta_0),C>C_0,$ we have  $$    \limsup_{n\rightarrow \infty }\frac 1n \log \# P_n (f, \zeta, C) >0.$$ In particular, $$\limsup_{n\rightarrow \infty }\frac1n \log \# {P_n (f)}>0.$$
\end{Cor}

 Moreover, we have a following
  relation between metric entropy and the exponential
  growth of periodic points, which is firstly proved for $C^{1+\alpha}$ case \cite{Katok}.

\begin{Thm}\label{Thm:exponential-rate-per-pts}Let $f\in \Diff^1 (M)$ and   $\mu\in \mathcal{M}^{qldh}_{f} (M)$   be  an ergodic measure. Then
$$ \limsup_{n\rightarrow \infty }\frac1n \log \# {P_n (f)}\geq h_\mu (f).$$
\end{Thm}

\begin{Rem}\label{rem-equal-expoential-growth} For $C^{1+\alpha}$ case, there is a generalized version which states that the metric entropy equals to the exponential growth of periodic measures which approximate to the given measure in weak$^*$ topology \cite{LST1}. However, Lyapunov neighborhood
is an important technique 
in $C^{1+\alpha}$ case but in our present case, it is still unknown for the inverse inequality.
\end{Rem}

 Theorem \ref{Thm:exponential-rate-per-pts} will be a direct application of Corollary \ref{Cor:1Horshoe-numbergrowth-periodic} below.
\bigskip


Now let us generalize one classical result (see also Chapter 21 in
 \cite{DGS}) by Sigmund to $C^1$ non-uniformly hyperbolic case.  It was proved that for any topological dynamical system with specification, including Axiom A systems, the periodic measures are dense in the space of invariant measures. There are some generalizations for $C^{1+\alpha}$ non-uniformly hyperbolic case (see  \cite{Hir,LLS}). Here we study $C^1$ non-uniformly hyperbolic case.

\begin{Thm}\label{Thm:density-per-meas}
Let $f\in \Diff^1 (M)$ and   $\mu\in \mathcal{M}^{qldh}_{f} (M)$   be  an ergodic measure. Then there is a $\mu$ full-measured set
$\tilde{\Lambda}\subseteq supp (\mu)$ (corresponding to new established Liao-Pesin set as below) such that
for any $t>0,$ the set of periodic measures supported on $t-$neighborhood of $supp (\mu)$
is dense
in the set of all $f-$invariant measures supported on
$\tilde\Lambda$.
\end{Thm}

In particular, we have following.
\begin{Thm}\label{Thm:density-per-meas22222}
Let $f\in \Diff^1 (M)$ and   $\mu\in \mathcal{M}^{qldh}_{f} (M)$   be  an ergodic measure. Then     for any $t>0,$  $\mu$ can be approximated by periodic measures supported on $t-$neighborhood of $supp (\mu).$
\end{Thm}

In other words,   one has following corollary which suggests a ``weak stability" of hyperbolic measures.

\begin{Cor}\label{Cor:density-per-meas}
Let $f\in \Diff^1 (M)$ and   $\mu\in \mathcal{M}^{qldh}_{f} (M)$   be  an ergodic measure. If $f_n$ is a sequence converging to $f$ in the $C^1$ topology, then $f_n$ has an invariant hyperbolic probability measure $\mu_n$
 such that $\mu_n$ converges to $\mu$ in weak$^*$ topology. Furthermore $\mu_n$ may be assumed to be supported on hyperbolic periodic points.
 \end{Cor}




\subsection{Quasi-invariant and Quasi-ergodic}

For possible applications of more general dynamical systems, we want to introduce two new concepts called quasi-invariant and quasi-ergodic, which are inspired from Poincar$\acute{\text{e}}$ Recurrence Theorem.

\begin{Def}\label{Def-quasi-invariant-ergodic}
Let $\mu\in \m  (M).$ We say $\mu$ to be quasi-invariant, if for any Borel set $A$ with $\mu (A)>0$, there exists a sequence  $\{n_i\}\uparrow \infty$ such that $f^{n_i} (A)\cap A \neq \emptyset.$\\
 Moreover, we say  $\mu$ to be quasi-ergodic, if for any Borel set $A,B$ with $\mu (A)\cdot\mu (B)>0$, there exists a sequence  $\{n_i\}\uparrow \infty$ such that $f^{n_i} (A)\cap B \neq \emptyset.$

\end{Def}

It is obvious that every quasi-ergodic measure is quasi-invariant and by  Poincar$\acute{\text{e}}$ Recurrence Theorem every invariant measure is quasi-invariant.  If a measure $\mu$ is quasi-ergodic and also invariant, then it is not difficult to see $\mu$ should be ergodic.  Moreover, every ergodic measure is quasi-ergodic. This can be deduced obviously from Birkhoff ergodic theorem,
\begin {equation*}\label{Eq:erg-application}
\lim_{l\rightarrow +\infty}\frac 1l
  \sum_{n=0}^{l-1}\mu (f^{-n} (A)\cap B)=\mu (A)\mu (B)>0.
\end{equation*}
However,    quasi-invariant is not necessary to be truly invariant and quasi-ergodic is not necessary to be truly ergodic. Here we give some simple  examples.
It is easy to see $f^K (K\geq 2)$-invariant  (or ergodic) measures are always quasi-invariant but not necessarily invariant, for instance, the Dirac measure supported a point whose orbit is periodic with period larger than 1.  Given a infinite series  $\sum_{i=1}^{\infty}a_i=1 (a_i>0)$ and a sequence of probability measures $\mu_i (i\geq1)$ which is $f^{n_i}-$invariant for some large $n_i (\text{converging to }\infty\,\text{as }i \text{ goes to }\infty)$ but not invariant for $f^j$, $1\leq j\leq n_i$, define $\mu=\sum_{i=1}^{\infty}a_i\mu_i,$ then $\mu$ is still quasi-invariant but may be not $f^n$-invariant for any positive integer $n$. In particular, for example, if a system has infinite periodic orbits with different periods, then for every periodic orbit we   just take a point on it and  take $\mu_i$ to be
the Dirac measure on this chosen point and so $\mu=\sum_{i=1}^{\infty}a_i\mu_i$ is quasi-invariant but not $f^n$-invariant for any $n$.


For quasi-invariant or quasi-ergodic case, we have following.

\begin{Thm}\label{Thm:2015-quasi-invariant}Let $f\in \Diff^1 (M)$ and $\mu$ be a quasi-invariant measure. If there is an $f$-invariant set $\Delta$ with $\mu(\Delta)=1$, there is   be a $Df-$invariant quasi-dominated splitting $T_\Delta M=E\oplus F$ on
$\Delta$ and $\mu$ is  average-nonuniformly hyperbolic with respect to this splitting, then \\
   (1) the periodic points is dense in the support of $\mu;$\\
   (2) 
 there is horseshoe;\\
   (3) If further $\mu$ is quasi-ergodic, then there is a $\mu$ full-measured set
$\tilde{\Lambda}$ (corresponding to new established Liao-Pesin set as below) such that
for any $t>0,$ the set of periodic measures
is dense
in the set of all $f-$invariant measures supported on
$\tilde\Lambda$.

\end{Thm}

We will prove Theorem \ref{Thm:2015-quasi-invariant} in section \ref{proof-horseshoe-periodic-dense} and section \ref{subsection-densityperiodicmeas}.

\subsection{Livshitz Theorem}\label{section-Livshitz-introduction}

Now we state Livshitz Theorem in $C^1$ systems, which is  obtained in $\textit{C}^{1+\alpha}$ case (   \cite{Katok1}, see Theorem S.4.17 on Page 692, also see  \cite{KatokMondoz,BP}).

\begin{Thm}\label{Thm:LivThm}
 Let $f\in \Diff^1 (M)$ and   $\mu\in \mathcal{M}^{qldh}_{f} (M)$   be  a nonatomic   measure.   Let $\varphi:M\rightarrow \mathbb{R}$ be a H\"{o}lder continuous function such that for each  (hyperbolic) periodic point $z$ with $f^m (z)=z,$ we have $\sum_{i=0}^{m-1}\varphi (f^i (z))=0.$ Then there exists a Borel measurable function $\psi$ such that $$\varphi (x)=\psi (f (x))-\psi (x)$$ holds for $\mu$ almost every point $x.$
\end{Thm}

\begin{Rem} For  above system $f$ and  $C^\alpha$ cocycles, above kind of Livsic theorem should be true and moreover,    the  Lyapunov exponents  of $C^\alpha$ cocycles can be approximated by ones of periodic measures. The main  observation
 is that their proofs are just based on exponentially closing property.  One can see  \cite{Kal} (or    \cite{Dai, WS}) for
  details.
\end{Rem}

We will prove Theorem \ref{Thm:LivThm} in section \ref{proof-Livsic}.

\bigskip

Let $\phi\in C^0 (M).$ We say $\phi$  is a coboundary, if there is $h\in C^0 (M)$ such that $\phi=h\circ f-h.$
 Let $Cob (f)$ denote the space of all coboundary continuous functions.

\begin{Thm}\label{Thm:LivThm222222}
 Let $f\in \Diff^1 (M)$ and suppose that  $\mathcal{M}^{qldh}_{f} (M)\cap \m_e (M)=\m_e (M).$ Let $\varphi:M\rightarrow \mathbb{R}$ be a  continuous  (not necessarily H\"{o}lder continuous) function such that for each  (hyperbolic) periodic point $z$ with $f^m (z)=z,$ we have   $\sum_{i=0}^{m-1}\varphi (f^i (z))=0.$          Then $\varphi\in \overline{Cob (f)}.$
\end{Thm}

\begin{Rem}\label{Rem-good-example-allergodic-hypDom}
Recall that for the partially hyperbolic (not uniformly hyperbolic) systems introduced in  \cite{D2009}, Leplaideur et al proved in  \cite{LOR} for the central direction, all ergodic invariant measures only have negative exponents, with the exception of a Dirac measure supported on a saddle with positive exponent.   So above theorem can be applied in such systems (and all  its conjugations $h^{-1}fh$ where $h$ is a homeomorphism).
\end{Rem}

Theorem \ref{Thm:LivThm} needs a version of exponentially closing property and exponential closing is persisted for $C^\gamma$-conjugated ($\gamma>0$)  systems. However, Theorem \ref{Thm:LivThm222222} just needs closing lemma (not necessarily exponential) and this closing is persisted for all topologically-conjugated  systems.
We will prove Theorem \ref{Thm:LivThm222222} in section \ref{proof-Livsic}.

\subsection{Approximation of Horseshoes}

Now we state a result that the information of hyperbolic measure can be approximated by ones of horseshoes.

\begin{Thm}\label{Thm:Horshoe-Approximation}
 Let $f\in \Diff^1 (M)$ and   $\mu\in \mathcal{M}^{qldh}_{f} (M)$   be  a  nonatomic ergodic   measure.
  Then for any $\epsilon>0,$ any neighborhood $\mathcal{V}$ of $\mu$ in   weak$^*$ topology,
   there exists a horseshoe $H_\epsilon$ such that\\
 (1) $Haus (H_\epsilon,supp (\mu))<\epsilon;$\\
 (2) $h_{top} ( H_\epsilon)>h_\mu (f)-\epsilon;$\\
 (3) $\m_f (H_\epsilon)\subseteq \mathcal{V}$; \\
 (4) There is a dominated splitting on $H_\epsilon$:   $T_{H_\epsilon}=E\oplus F$ with $dim (E)=ind (\mu)$ such that for any $x\in H_\epsilon,$ $T_{H_\epsilon}=E\oplus F$  coincides the  extended-dominated Oseledec hyperbolic splitting;\\
 (5) there exists $C_\epsilon>0$ such that for any $x\in H_\epsilon,$ $n\geq 1,$ $$ \|D_xf^n|_{E (x)}\|\leq C_\epsilon exp (n (\lambda_s+\epsilon)) ,\,\,m (D_xf^{n}|_{F (x)}) \geq C^{-1}_\epsilon exp (n (\lambda_u-\epsilon)), $$ where $\lambda_s,\lambda_u$ denote the maximal negative Lyapunov exponent of $\mu$ and the minimal positive  Lyapunov exponent of $\mu$, respectively;

If further the  Oseledec  splitting of $\mu$ is dominated (i.e., every two distinct Oseledec bundles are dominated), the conclusions  (4) and  (5) can be stated better: \\
 (4')  If $\chi_1>\chi_2\cdots>\chi_m$ are the distinct Laypunov exponents of $\mu$, with multiplicities $n_1,\cdots,n_m\geq 1$, then there is a dominated splitting on $H_\epsilon$:   $T_{H_\epsilon}=E_1\oplus E_m$ with $dim (E_i)=n_i$ such that for any $x\in H_\epsilon,$ $T_{H_\epsilon}=E_1\oplus\cdots \oplus E_m$  coincides the  extended-dominated Oseledec   splitting;\\
 (5') there exists $C_\epsilon>0$ such that for any $x\in H_\epsilon,$ $n\geq 1,i=1,\cdots,m$ $$ C^{-1}_\epsilon exp (n (\chi_i-\epsilon))\leq  m (D_xf^n|_{E_i (x)} ) \leq \|D_xf^n|_{E_i (x)}\|\leq C_\epsilon exp (n (\chi_i+\epsilon)) , $$ where $\lambda_s,\lambda_u$ denote the maximal negative Lyapunov exponent of $\mu$ and the minimal positive  Lyapunov exponent of $\mu$, respectively, and moreover,\\
 (6) the Lyapunov exponents of $\mu$ can be $\epsilon-$approximated by ones of every ergodic measure  supported on $\m_f (H_\epsilon)$. In particular,\\
  (6.1) the Lyapunov exponents of $\mu$ can be $\epsilon-$approximated by ones of   hyperbolic  periodic measures;\\
  (6.2)
 the Lyapunov exponents of $\mu$  can be $\epsilon-$approximated by ones of  hyperbolic ergodic measures with positive metric entropy whose support are uniformly hyperbolic.

\end{Thm}


Recall that horseshoe has structural stability and every non-atomic ergodic measure supported on the horseshoe satisfies  the assumption Theorem \ref{Thm:Horshoe-Approximation}. So   hyperbolic measures with  (quasi-)limit-domination  are ``stable" or persistent under $C^1$ perturbations.

 In the $\textit{C}^{1+\alpha}$ case,  some similar statements can  be found in  \cite{Katok1,BP} and   \cite{ACW}. From  \cite{Katok1,BP,ACW} we know that stable manifold theorem and shadowing lemma are enough to prove  (1)  and  (3), but the technique of Lyapunov neighborhood  is important for other arguments. Here in our  present case, we do not have the technique of  (Lyapunov) regular neighborhood but fortunately we observe  that under the assumption of domination,    (1)  and  (3) are enough to imply other arguments.

\medskip


Moreover, we  state a following corollary which is obtained recently for $C^{1+\alpha}$ setting  \cite{LS2013}.

\begin{Thm}\label{cor1-Thm:Horshoe-Approximation}
 Let $f\in \Diff^1 (M)$ and   $\mu\in \mathcal{M}^{qldh}_{f} (M)$   be  a  nonatomic ergodic   measure.
Then $\mu$ is approximated by uniformly hyperbolic
sets in the sense that there exists a sequence $\Omega_n$ of compact, topologically
transitive, locally maximal, uniformly hyperbolic sets such that for any sequence $\mu_n$  of $f$-invariant ergodic probability measures with $supp (\mu_n)\subseteq \Omega_n,$
 we have $\mu_n\rightarrow \mu$ in the weak$^*$ topology.
\end{Thm}

 We say $f$ is Hausdorff-hyperbolic, if for any $\epsilon>0,$ there is a horseshoe $H_\epsilon$ such that
  $$Haus (H_\epsilon,M)<\epsilon.$$   Then by Theorem \ref{Thm:Horshoe-Approximation}  (1) we have

\begin{Cor}\label{Cor:1Hausdorff-hyperbolic-Horshoe-Approximation}
Let $f\in \Diff^1 (M)$ and   $\mu\in \mathcal{M}^{qldh}_{f} (M)$   be  a  nonatomic ergodic   measure with $supp (\mu)=M$. Then  $f$ is Hausdorff-hyperbolic.

\end{Cor}

By Theorem \ref{Thm:Horshoe-Approximation}  (2)  and  (5) we have

\begin{Cor}\label{Cor:1Horshoe-numbergrowth-periodic} Let $f\in \Diff^1 (M)$ and   $\mu\in \mathcal{M}^{qldh}_{f} (M)$   be  a  nonatomic ergodic   measure. Then
 there is some $ \zeta_0\in (0,1)$ such that for any $\zeta\in  (0,\zeta_0),$ we have  $$  h_{\mu} (f)\leq \lim_{C\rightarrow \infty} \limsup_{n\rightarrow \infty }\frac 1n \log \# P_n (f, \zeta, C) .$$
\end{Cor}

Obviously,
Corollary \ref{Cor:1Horshoe-numbergrowth-periodic} implies Theorem \ref{Thm:exponential-rate-per-pts}.
\smallskip

 Theorem \ref{Thm:Horshoe-Approximation}  (6.1) tells us that if  the  Oseledec  splitting of $\mu$ is dominated, then the Lyapunov
 exponents of $\mu$ can be $\epsilon-$approximated by ones of
  hyperbolic  periodic measures.
 For $C^{1+\alpha}$ diffeomorphism, similar results had been  proved in  \cite{WS} that for every ergodic hyperbolic measure $\mu$, the Lypunov exponents of $\mu$ can be approximated by ones of periodic measures. However, it is still unknown for $C^{1}$ case which only assume that the stable bundle dominates the unstable one.

\begin{Que}\label{Que:Lyapunov-approximation} Let $f\in \Diff^1 (M)$ and   $\mu\in \mathcal{M}^{qldh}_{f} (M)$   be  an ergodic measure. Then
 whether Lyapunov exponents of $\mu$ can be approximated by ones of periodic measures?
\end{Que}

For  Anosov case,  obviously every ergodic measure is hyperbolic and
its Osledec hyperbolic splitting  corresponds to the uniformly hyperbolic splitting so that it is dominated. However, it is also unknown   for the approximation of Lyapunov exponents if the system is just $C^1$:

\begin{Que}\label{Que:Laypunov-approxi-Anosov} Let $f\in \Diff^1 (M)$ be Anosov. Then for any ergodic measure $\mu$,
 whether Lyapunov exponents of $\mu$ can be approximated by ones of periodic measures?
\end{Que}



\subsection{Existence of Horseshoe, Hyperbolic Periodic Orbit in Partial Hyperbolic Systems}

\begin{Thm}\label{Thm-partial-hyperbolic}
Let $f:M\rightarrow M$ be a $C^1$ diffeomorphism on a   compact Riemanian manifold $M$ with a dominated splitting $TM=E\oplus F$. If one condition of following holds:

 \noindent {\bf  (A) }
If $E$ is quasi-conformal  and $F$ is uniformly expanding;

 \noindent {\bf  (A') }
If $dim (E)=1$ and $F$ is uniformly expanding;

 \noindent {\bf  (B) } If  $F$ is quasi-conformal  and  $E$ is uniformly contracting;

 \noindent {\bf  (B') } If $dim (F)=1$ and $E$ is uniformly contracting;

 \noindent {\bf  (C) } If   $E,F$ is quasi-conformal and $f$ has positive topological entropy;

 \noindent {\bf  (C') } If $dim (E)=dim (F)=1$ and $f$ has positive topological entropy;

Then  \\
 (1) There is  a horseshoe (in particular, there are infinite  hyperbolic periodic orbits);\\
   (2) $f$ is entropy-hyperbolic;\\
   (3)  there is some $ \zeta_0\in (0,1)$ such that for any $\zeta\in  (0,\zeta_0),$ we have  $$  h_{top} (f)=\lim_{C\rightarrow \infty} \limsup_{n\rightarrow \infty }\frac 1n \log \# P_n (f, \zeta, C) ,$$
  where \begin{eqnarray*}
  & &P_n (f,\zeta, C):=\{x\in M\,|\,f^nx=x
   \text{ and for any }l\geq 1, j=0,1,2,\cdots, n-1,\,\,\,\,\,\,\,\,\,
   \,\,\,\,\,\,\,\,\,\,\\\nonumber
   & & \,\,\,\,\,\,\,\,\,\,\,\,\,\,\,\,\,\,\,\,\,\,\,\,\,\,\,\,\,\,\,\,
   \,\,\,\,\,\,\,\,\,\,\,\,\,\,
  \,\max\{\|Df^l|_{E (f^jx)}\|,\|Df^{-l}|_{F (f^jx)}\|\}\leq C \zeta^l \}.
  \end{eqnarray*}

\end{Thm}

For the cases of  (A')  (B') and  (C'), the systems are all far from homoclinic tangencies so that by  \cite{LiaoVianaYang}
they  are entropy-expansive so that they have maximal entropy measures. However, for the cases of  (A)  (B) and  (C), it is still unknown the existence of maximal entropy measures.
 Let $\m_z (M)$ denotes the space of invariant measures with zero metric entropy. In a Baire space, a set is residual if it contains a countable intersection of dense open sets.

\begin{Thm}\label{Thm-partial-hyperbolic-zero-entropy-measure-dense}
Let $f:M\rightarrow M$ be the system as in Theorem \ref{Thm-partial-hyperbolic}. Then
$\m_z (M)$ is dense in $\m_f (M).$ In particular, for each case  of  (A')  (B')  (C'), $\m_z (M)$ is residual  in $\m_f (M).$

\end{Thm}

We will prove Theorem \ref{Thm-partial-hyperbolic-zero-entropy-measure-dense}   in section \ref{section-partialhyper}.

 The results of Theorem \ref{Thm-partial-hyperbolic} hold for  more general case than quasi-conformal case, since quasi-conformal case implies that there is only one  Lyapunov exponent in  the corresponding bundle for a.e. points.

\begin{Thm}\label{Thm-partial-hyperbolic-22222}
Let $f:M\rightarrow M$ be a $C^1$ diffeomorphism on a   compact Riemanian manifold $M$ with a dominated splitting $TM=E\oplus F$. If one condition of following holds:

 \noindent {\bf  (I) }
 Suppose that $F$ is uniformly expanding, i.e.  there exist $C>0$ and $0<\lambda<1$ such that $$ {\|Df^{-n}|_{F (x)}\|} \leq C \lambda^n, \forall x\in M,\,\, n\geq 1.$$ If there is some $a\in[0,-\frac{dim F}{dim E}\log\lambda)$ such that $E$ satisfies that
$$ \liminf_{n\rightarrow\infty} \frac1n\log\frac{\|Df^n|_{E (x)}\|}{m (Df^n|_{E (x)} )}\leq a.$$

 \noindent {\bf  (II) }Suppose that $E$ is uniformly contracting, i.e.  there exist $C>0$ and $0<\lambda<1$ such that $$ {\|Df^{n}|_{E (x)}\|} \leq C \lambda^n, \forall x\in M,\,\, n\geq 1.$$ If there is some $a\in[0,-\frac{dim E}{dim F}\log\lambda)$ such that $F$ satisfies that
$$ \liminf_{n\rightarrow\infty} \frac1n\log\frac{\|Df^n|_{F (x)}\|}{m (Df^n|_{F (x)} )}\leq a.$$
Then:
the results of Theorem \ref{Thm-partial-hyperbolic} hold.

\end{Thm}

Note that  (A)  (A') in Theorem \ref{Thm-partial-hyperbolic} imply  (I) and  (B),  (B') in Theorem \ref{Thm-partial-hyperbolic} imply   (II) (In fact, in these cases, $a=0$). So we only need to prove  (I) and  (II) of Theorem \ref{Thm-partial-hyperbolic-22222},  and  (C)  (C') in Theorem \ref{Thm-partial-hyperbolic}. Note that  (I) and  (II) are similar so that we only need to show one case, see section \ref{section-partialhyper}.



\begin{Thm}\label{Thm-partial-hyperbolic-333}
   Under the same assumptions as Theorem \ref{Thm-partial-hyperbolic-22222} (including the cases of   (A)  (A')  (B)  (B') in Theorem \ref{Thm-partial-hyperbolic}), we have: \\
   (1)
  Lebesgue measure is  average-nonuniformly hyperbolic;\\
(2) if  further Lebesgue measure is quasi-invariant, then    the periodic points form a dense subset of the whole manifold.
\end{Thm}

We will prove Theorem \ref{Thm-partial-hyperbolic-333} for system with condition (I) in section \ref{section-partialhyper} and the case (II) is similar.
Recall that from  \cite{ABV2000} and  \cite{BV2000}  SRB measures exist in $C^{2}$ partially hyperbolic systems, whose tangent bundle decomposes into two dominated bundles: one is uniformly expanding (or contracting) and another is  (average-)nonuniformly contracting (or expanding) on Lebesgue positive (or full) measure set (called mostly contracting or mostly expanding   \cite{ABV2000} and  \cite{BV2000}). Here for the system  in Theorem \ref{Thm-partial-hyperbolic-22222}
 which is further assumed $C^2$, mostly  contracting or mostly expanding   naturally holds on Lebesgue full measure set  for system $f^K$ ($K$ being large) and so it is possible to obtain SRB measures.




Moreover, we have another characterization for Lebesgue measure.  We say $f$ to be  {\it volume-non-expanding},
if Lebesgue a.e. $x$, $$\limsup_{n\rightarrow+\infty}\frac1n \log |det (Df^n)|\leq 0.$$

\begin{Thm}\label{Thm-partial-hyperbolic-444444444444444}
   Under the same assumptions as the cases of   (A)  (A')  (resp.,  (B)  (B'))  in Theorem \ref{Thm-partial-hyperbolic}:
 $f$ (resp., $f^{-1}$) is  volume-non-expanding.
\end{Thm}

We will prove Theorem \ref{Thm-partial-hyperbolic-444444444444444} in section \ref{section-exist-hyperbolic-srb-like}.

\bigskip

Above results are considered for partial hyperbolicity with two dominated bundles.
For the usual partial hyperbolicity, we have following result.

\begin{Thm}\label{Thm-partial-hyperbolic-usual}
Let $f:M\rightarrow M$ be a $C^1$ diffeomorphism on a   compact Riemanian manifold $M$ with a dominated splitting $TM=E^s\oplus E^c\oplus E^u$ (i.e., $E^s\prec E^c$ and $E^c\prec E^u$) and suppose that $E^s$ is uniformly contracting, $E^u$ is uniformly expanding and  $E^c$ is quasi-conformal. Then\\
  (1) either there exists at least one hyperbolic periodic orbit, or \\
   for any $\epsilon>0,$ there is some $C_\epsilon>0$ such that for any $n\geq 1, x\in M,$ $$ C^{-1}_\epsilon exp (-n\epsilon)\leq m (Df^n|_{E^c (x)})\leq \|Df^n|_{E^c (x)}\|\leq C_\epsilon exp (n\epsilon),$$  (in this case the Lyapunov exponents on the bundle $E^c$ of all points exist and equal to zero).\\
 (2) either there exists a horseshoe, or \\
   for Lebesgue a.e. $x\in M,$  the central Lyapunov exponents of $E^c$ at $x$ exist and equal to zero, that is
 $$\lim_{n\rightarrow +\infty}\frac1n \log\|Df^{ n}|_{E^c (x)}\|= \lim_{n\rightarrow +\infty}\frac1n \log m (Df^{ n}|_{E^c (x)})=0.$$
\end{Thm}

This theorem can be applied for all partially hyperbolic systems with one dimensional central bundle which form an open subset of $\Diff^1 (M).$
We will prove Theorem \ref{Thm-partial-hyperbolic-usual} in section \ref{proof-exitence-hyperbolic}. Remark that if one also consider $f^{-1}$, the consequence  (2) can be stated as follows: {\it   either there exists a horseshoe, or
   for Lebesgue a.e. $x\in M,$  the central Lyapunov exponents of $E^c$ at $x$ exist and equal to zero, that is    $$\lim_{n\rightarrow +\infty}\frac1n \log\|Df^{ \pm n}|_{E^c (x)}\|=  \lim_{n\rightarrow +\infty}\frac1n \log m  (Df^{\pm n}|_{E^c (x)})=0.$$}

Moreover, we can have a following theorem.

\begin{Thm}\label{Thm-partial-hyperbolic-usual-2}
Let $f:M\rightarrow M$ be a $C^1$ diffeomorphism on a   compact Riemanian manifold $M$ with a dominated splitting $TM=E^s\oplus E^c\oplus E^u$ (i.e., $E^s\prec E^c$ and $E^c\prec E^u$) where  $E^s$ is uniformly contracting and $E^u$ is uniformly expanding. Suppose for each fixed ergodic invariant measure with positive entropy, its   Lyapunov exponents in the direction $E^c$  are all non-zero with same sign (admitting coexistence of  two ergodic measures with  different sign Lyapunov exponents in the direction $E^c$). Then
the results of Theorem \ref{Thm-partial-hyperbolic} hold.  Moreover,  $\m_z (M)$ is  dense in $\m_f (M).$ In particular,  if $E^c$ is one-dimensional,  $\m_z (M)$ is residual in $\m_f (M).$

\end{Thm}

\begin{Rem} Similar results of Theorem \ref{Thm-partial-hyperbolic} hold provided that $E^c$ can be decomposed by a more fine dominated splitting $E^c_1\oplus \cdots E^c_k$ and for each fixed  $E^c_i $ and ergodic measure with positive entropy, its  Lyapunov exponents in the direction $E^c_i$   are all non-zero with same sign (admitting coexistence of  two ergodic measures with  different sign Lyapunov exponents in the direction $E^c_i$). In particular, all $E^c_i$ are one-dimensional, $\m_z (M)$ is residual in $\m_f (M).$

\end{Rem}

We will prove Theorem \ref{Thm-partial-hyperbolic-usual-2} in section \ref{section-partialhyper}.  Theorem \ref{Thm-partial-hyperbolic-usual-2}  can be applied for some systems which are not necessarily uniformly hyperbolic, for example,  partially hyperbolic   systems in  \cite{D2009,LOR} (and their  conjugated system  $h^{-1}fh$ where $h\in\Diff^1 (M)$).
Moreover, for  the partially hyperbolic systems of  \cite{D2009,LOR},  Lebesgue measure
 is  average-nonuniformly hyperbolic for $f$.

\begin{Thm}\label{Thm-partial-hyperbolic-usual-Lebegue-hyp}
Let $f:M\rightarrow M$ be a $C^1$ diffeomorphism on a   compact Riemanian manifold $M$ with a dominated splitting $TM=E^s\oplus E^c\oplus E^u$ (i.e., $E^s\prec E^c$ and $E^c\prec E^u$) where  $E^s$ is uniformly contracting and $E^u$ is uniformly expanding. Suppose for the direction $E^c$, every ergodic invariant measure with positive metric entropy has only negative Lyapunov exponents in
$E^c$ (resp., every ergodic invariant measure with positive metric entropy has only positive  Lyapunov exponents in
$E^c$).
Then
the results of Theorem \ref{Thm-partial-hyperbolic-usual-2} hold. Moreover, Lebesgue measure
is  average-nonuniformly hyperbolic for $f$ and $f$ (resp., $f^{-1}$)  is  volume-non-expanding. 
If  further Lebesgue measure is quasi-invariant, then   the periodic points form a dense subset of the whole manifold.

\end{Thm}




 By Theorem \ref{Thm-partial-hyperbolic-usual-2}, its results hold naturally under the assumptions of Theorem \ref{Thm-partial-hyperbolic-usual-Lebegue-hyp}, since the assumptions of Theorem \ref{Thm-partial-hyperbolic-usual-Lebegue-hyp} is stronger than   ones of Theorem \ref{Thm-partial-hyperbolic-usual-2}.
We will prove the left part of Theorem \ref{Thm-partial-hyperbolic-usual-Lebegue-hyp} in section \ref{section-partialhyper} and section \ref{section-exist-hyperbolic-srb-like}.

\bigskip

Observe that the assumption of Theorem \ref{Thm-partial-hyperbolic-22222} is an open condition, that is, there is a neighborhood such that every system in the neighborhood has similar condition. In other words, if $PH^1_{cc} (M)$  and $PH^1_{1} (M)$ denote the space of partially hyperbolic systems satisfying the assumption of Theorem \ref{Thm-partial-hyperbolic-22222} (which are {\it close to conformal}) and the space of  partially hyperbolic systems satisfying  (A') or  (B') in Theorem \ref{Thm-partial-hyperbolic} (i.e., one dimensional central bundle), respectively. Then $PH^1_{1} (M)\subseteq  PH^1_{cc} (M)$ are two open subsets of $PH^1  (M).$ Let $D^1 (M)$ be the space of all systems with a global dominated splitting which also is a open subset of $\Diff^1 (M)$.
Then we have following result for continuity of entropy function.

\begin{Thm}\label{Thm-partial-hyperbolic-22222-1111}
 (1) The entropy function $h: PH^1_{cc} (M)\rightarrow \mathbb{R}, f\mapsto h_{top} (f)$ is lower semi-continuous.\\
 (2) If $M$ is a surface, $h: D^1 (M)\rightarrow \mathbb{R}, f\mapsto h_{top} (f)$ is lower semi-continuous.\\
 (3) the following systems are also  lower semi-continuous points of entropy function $h: \Diff^1 (M)\rightarrow \mathbb{R}, f\mapsto h_{top} (f)$:\\
 (3.1) $f$ satisfies that every ergodic hyperbolic measure is in $\m_f^{qldh} (M).$ In particular, it includes the case of Theorem \ref{Thm-partial-hyperbolic-usual-2} and includes all  $C^1$ surface diffeomorphisms  with a global quasi-limit-dominated splitting. \\
 (3.2) $f\in D^1 (M)$ satisfies that two dominated bundles are quasi-conformal.
\end{Thm}

This theorem is not difficult to prove. Let us explain more precisely.  By above analysis of horseshoe and variational principle,
every system $f$ in above theorem satisfies entropy-hyperbolic. Then $$h_{top} (f)=\{h_{top} (f|_\Lambda)|\, \Lambda \text{ is a hyperbolic horseshoe  }\}.$$ By structural stability of horseshoe, lower semi-continuity follows.

\subsection{Some other recent known related  results }

For possible completeness of $C^1$ nonuniform hyperbolicity theory, we state several related known results including  connection of recurrent time and Lyapunov exponents, Pesin's entropy formula which builds the relation of metric entropy and Lyapunov exponents.

\subsubsection{Recurrence Time \& Lyapunov Exponents}

The first known result  is to connect ``Recurrence"   with  Lyapunov exponents (see  \cite{OliTian,STV,Oli}). Given $x\in M$ and $r>0$, denote the first return time of a ball $B (x,r)$ radius $r$ at $x$ by $$\tau (B (x,r)):= \min \{k>0 \,|\,\,f^k (B (x,r))\cap B (x,r)\neq \emptyset\}.$$ Then

\begin{Thm}\label{Thm:Recurrece-Lyaexp}  ( \cite{OliTian}) 
Let $f\in \Diff^1 (M)$ and   $\mu\in \mathcal{M}^{dh}_{f} (M)$   be  an ergodic measure.   Then  $\mu$  satisfies that for $\mu$ a.e. $x\in M$, $$\limsup_{r\rightarrow 0}\frac{\tau (B (x,r))}{-\log r}\leq \frac1{\lambda_u}-\frac1{\lambda_s},$$ where $\lambda_u,\,\lambda_s$ are the minimal positive Lyapunov exponent and maximal negative Lyapunov exponent of $\mu$, respectively.
\end{Thm}

\subsubsection{Pesin's Entropy Formula}
The second known result is Pesin entropy formula which was firstly obtained in $C^{1+\alpha}$ systems. Recall that an invariant measure $\mu$ satisfies Pesin entropy formula, if
$$h_{\mu} (f)=\int\sum_{\lambda_i (x)\geq0}\lambda_i (x)d\mu,$$ where $\lambda_1 (x)\geq\lambda_2 (x)\geq\cdots\geq\lambda_{dim\,M} (x)$ denote
the Lyapunov exponents of $\mu$ a.e. $x.$  Now we state $C^1$ Pesin's entropy formula for smooth measures (see  \cite{SunTianEntropy}) and its generalization for SRB-like measures (see  \cite{CCE}).

\begin{Thm}\label{PesFormula-Thm:1} ( \cite{SunTianEntropy})
Let $f\in \Diff^1 (M)$  preserve  an invariant
probability measure $\mu$ which is absolutely continuous relative to
Lebesgue measure.
 If there is an $f$-invariant measurable function $m (\cdot):M\rightarrow
\mathbb{N}$ such that  for
$\mu\,a.\,\,e.\,\,x\in M,$ there exists an $m (x)$-dominated splitting:
$T_{orb (x)}M=E_{orb (x)}\oplus F_{orb (x)}$, then
$$h_{\mu} (f)\geq\int \chi (x)d\mu,$$ where
$\chi (x)=\sum_{i=1}^{dim\,F (x)}\lambda_i (x).$ \\In particular, if
for $\mu\,\, a.\,e.\,\,x\in M$, $\lambda_{dim F (x)} (x) \geq  0 \geq \lambda_{dim F (x)+1} (x)$, then Pesin's entropy
formula holds.
\end{Thm}

 Consequently,  Pesin's entropy formula is valid on any  $f\in \Diff^1_m (M)$, assuming $f$ to be  Anosov, partially hyperbolic with one-dimensional center, far from tangency or be a generic system in $\Diff^1_m (M)$ (see \cite{SunTianEntropy}).
 \medskip

Now we state recent advance on entropy formula for the case of SRB-like measure, which always exists in any dynamics.
 Recall that $\mathcal{M} (M)$ denotes the space of all probability measures, and ${\mathcal M}_f (M) \subset {\mathcal M} (M) $  denotes the space of $f$-invariant probability measures. For a point $x\in M$ we consider  the following sequence  $$\Big \{\frac1n\sum_{j=0}^{n-1}\delta_{f^j (x)} \Big\}_{n\in\mathbb{N}} $$ where $\delta_y$ is the Dirac  probability measure supported at $y\in M.$ We define the nonempty and compact   set $p \omega (x)$ of  probability measures:
$$p \omega_f (x) = \Big\{\mu \in {\mathcal M} (M) : \ \ \exists \ n_i \rightarrow + \infty \mbox{ such that } \lim^*_{i \rightarrow + \infty} \ \frac{1}{n_i} \sum_{j= 0}^{n_i-1} \delta_{f^j (x)} = \mu \Big\}.$$
It is standard to check that $p \omega_f (x) \subset {\mathcal M}_f (M)$.

 \begin{Def} {\bf  (SRB-like measures)} \em
  \label{definitionSRB-like}

 A probability measure $\mu\in \mathcal{M}_f (M)$ is \em SRB-like \em    (or  observable or  pseudo-physical) if for any $\varepsilon>0$ the set $$A_\varepsilon (\mu)=\{x\in M\colon \ \mbox{dist} (p\omega_f (x),\mu)<\varepsilon\}$$ has positive Lebesgue measure. The set $A_\varepsilon (\mu)$ is called basin of $\varepsilon-$attraction of $\mu.$
 \end{Def}

We denote by $\mathcal{O}_{f}$ the set of all SRB-like measures for $f:M \mapsto M$. It is easy to see that every SRB-like measure for $f$ is $f$-invariant.

\begin{Thm}\label{PesFormula-Thm:2} ( \cite{CCE}) 
Let $f\in \Diff^1 (M)$ and there is a  global dominated splitting $T_MM=E\oplus F.$  Then for any SRB-like measure $\mu\in \mathcal{O}_{f},$
$$h_{\mu} (f)\geq\int \chi (x)d\mu,$$ where
$\chi (x)=\sum_{i=1}^{dim\,F (x)}\lambda_i (x)$ and $\lambda_1 (x)\geq\lambda_2 (x)\geq\cdots\geq\lambda_{dim\,M} (x)$ denote
the Lyapunov exponents of $\mu$ a.e. $x.$
\end{Thm}

Inspired by this result and the relation of domination  and limit-domination,  we ask  a following  question   for entropy formula in the case of global limit-dominated splitting.

\begin{Que}\label{Que-PesFormula:1} Let $f\in \Diff^1 (M)$ and there is a  global limit-dominated
splitting $T_MM=E\oplus F.$  Whether one has following result:  for any SRB-like measure $\mu,$
$$h_{\mu} (f)\geq\int \chi (x)d\mu,$$ where
$\chi (x)=\sum_{i=1}^{dim\,F (x)}\lambda_i (x)$ and $\lambda_1 (x)\geq\lambda_2 (x)\geq\cdots\geq\lambda_{dim\,M} (x)$ denote
the Lyapunov exponents of $\mu$ a.e. $x.$
\end{Que}

There are some examples for possible positive answer, see Example \ref{Ex:nonuni-hyp-sys}.










\bigskip

Moreover, there are some results on {\it large deviation} for $C^1$ non-uniformly hyperbolic systems, for example, see  \cite{QST,Paulo} etc. Here we omit the details.

\section{Pesin set, stable manifolds
\&  (exponential) shadowing lemma }

\subsection {Classical Pesin blocks, Pesin set}
In this subsection we give a quick review
concerning some notions and results of $C^{1+\alpha}$ Pesin theory for convenience to compare with our new Pesin set, called  Liao-Pesin  set.


Let $f: M\to M$ be a $C^{1+\alpha}$ diffeomorphism. We recall the classical  Pesin set \cite{P1} defined independently on measures  and corresponding  Katok's shadowing lemma.

\smallskip

\subsubsection{\bf Pesin set in $C^{1+\alpha}$ setting}
 Given $\lambda,\mu \gg \varepsilon>0$, and for all $k \in \mathbb{Z}^{+}$,
we define $\Lambda_{k}=\Lambda_{k} (\lambda,\mu;\varepsilon)$ to be
all points $x \in M$ for which there is a splitting
$T_{x}M=E_{x}^{s} \oplus E_{x}^{u}$ with invariant property
$D_{x}f^{m} (E_{x}^{s})=E_{f^{m}x}^{s}$ and
$D_{x}f^{m} (E_{x}^{u})=E_{f^{m}x}^{u}$ satisfying:

\begin{enumerate}
\item[ (a)] $\|Df^{n}|_{E_{f^{m}x}^{s}}\| \leq e^{\varepsilon
k}e^{- (\lambda-\varepsilon)n}e^{\varepsilon \mid m\mid},~\forall
m\in\mathbb{Z},~n\geq 1$;

\item[ (b)] $\|Df^{-n}|_{E_{f^{m}x}^{u}}\| \leq e^{\varepsilon
k}e^{- (\mu-\varepsilon)n}e^{\varepsilon \mid m\mid},~\forall m\in
\mathbb{Z}, ~n\geq 1$;

\item[ (c)] $\tan  (\angle (E_{f^{m}x}^{s},E_{f^{m}x}^{u})) \geq
e^{-\varepsilon k}e^{-\varepsilon \mid m\mid},~\forall m\in
\mathbb{Z}$.
\end{enumerate}
We set
$\Lambda=\Lambda (\lambda,\mu;\varepsilon)=\bigcup_{k=1}^{+\infty}
\Lambda_{k}$ and call $\Lambda$ a Pesin set.

According to Oseledec Theorem \cite{Os}, every ergodic hyperbolic measure $\mu$ has $s \,\, (s\leq
dim M)$ nonzero  Lyapunov exponents
$$\lambda_{1}< \cdot\cdot\cdot<\lambda_{r}<0<\lambda_{r+1}<\cdot\cdot\cdot< \lambda_{s}$$
with associated Oseledec splitting
$$T_{x}M=E_{x}^{1}\oplus\cdot\cdot\cdot\oplus E_{x}^{s},\,\,\,\,x\in L (\mu),$$
where we recall that $L (\mu)$ denotes an Oseledec basin of $\mu.$ If
we denote by $\lambda$ the absolute value of the largest negative
Lyapunov exponent $\lambda_r$ and $\mu$ the smallest positive
Lyapunov exponent $\lambda_{r+1}$, then for any $0<\varepsilon<\min\{\lambda,\,\mu\},$ one has $\mu$ full-measured Pesin set
$\Lambda=\Lambda (\lambda,\mu;\varepsilon)$  (see, for example, Proposition 4.2 in  \cite{P1}). And for any point $x\in L (\mu) \cap \Lambda,$ $E^s_x$ and $E^u_x$ coincide with $
E^1_x\oplus\cdot\cdot\cdot\oplus E^r_x$ and $E^{r+1}_x\oplus
\cdot\cdot\cdot\oplus E^s_x$ respectively.

The following statements are elementary properties of Pesin blocks (see  \cite{P1}):
\begin{enumerate}
\item[ (a)] $\Lambda_{1} \subseteq \Lambda_{2} \subseteq
\Lambda_{3}\subseteq \cdot\cdot\cdot$;

\item[ (b)] $f (\Lambda_{k}) \subseteq
\Lambda_{k+1},~f^{-1} (\Lambda_{k}) \subseteq \Lambda_{k+1}$;

\item[ (c)] $\Lambda_{k}$ is compact for $\forall\,\, k\geq 1$;

\item[ (d)] for $\forall\,\, k\geq 1$ £¬the splitting $x\to E_{x}^{u}\oplus
E_{x}^{s}$ depends continuously on  $x\in\Lambda_{k}$.
\end{enumerate}

\subsubsection {\bf Shadowing lemma}

We recall Katok's shadowing lemma \cite{P1} in this subsection. Let $ (\delta_{k})_{k=1}^{+
\infty}$ be a sequence of positive real numbers. Let
$ (x_{n})_{n=-\infty}^{+ \infty}$ be a sequence of points in
$\Lambda=\Lambda (\lambda, \mu, \varepsilon)$ for which there exists
a sequence $ (s_{n})_{n=-\infty}^{+ \infty}$ of positive integers
satisfying:

\begin{enumerate}
\item[ (a)] $x_{n}\in \Lambda _{s_{n}}, ~\forall n\in \mathbb{Z}$;

\item[ (b)] $\mid s_{n}-s_{n-1}\mid \leq 1, ~\forall n\in \mathbb{Z}$;

\item[ (c)] $   d (fx_{n},x_{n+1})\leq \delta_{s_{n}}, ~\forall n\in \mathbb{Z}$;
\end{enumerate}
then we call $ (x_{n})_{n=-\infty}^{+ \infty}$ a
$ (\delta_{k})_{k=1}^{+ \infty}$  pseudo-orbit. Given $\tau>0$, a
point $x\in M$ is a
 $\tau$-shadowing point for the $ (\delta_{k})_{k=1}^{+ \infty}$
 pseudo-orbit if $   d (f^{n} (x),x_{n+1})\leq \tau \varepsilon_{s_{n}},~ \forall n\in
\mathbb{Z}$, where $\varepsilon_{k}=\varepsilon_{0}e^{-\varepsilon
k}$ and $\varepsilon_0$ is a constant only dependent on the system itself.
\bigskip

\begin{Lem}\label{LemShadow}
 (Shadowing lemma) Let $f:M\rightarrow M$ be a
$C^{1+\alpha}$ diffeomorphism, with a non-empty Pesin set
$\Lambda=\Lambda (\lambda,\mu;\varepsilon)$ and fixed parameters,
$\lambda,\mu\gg\varepsilon>0$. For $\forall \tau >0$ there exists a
sequence $ (\delta_{k})_{k=1}^{+ \infty}$ such that for any
$ (\delta_{k})_{k=1}^{+ \infty}$ pseudo-orbit there exists a unique
$\tau$-shadowing
point.
\end{Lem}

\subsubsection {\bf Stable and unstable manifolds}

We recall stable manifold theorem (e.g.,  \cite{P1}) on Pesin set in $C^{1+\alpha}$ setting. Before that we recall the definition of  (local) stable manifold.

\begin{Def}\label{Def}
Given a non-empty Pesin set $\Lambda=\Lambda (\lambda,\mu;\varepsilon)$  (with
$\lambda,\mu\gg\varepsilon>0$) we shall define the  (local)  \emph{stable  (unstable) manifolds }through any point $x\in \Lambda$ by $$W^s_\delta (x)=\{y\in M |\,    d (f^{n} (x),f^{n} (y))\leq \delta e^{- (\lambda-\varepsilon)n},\,n\geq 0\}$$
 $$ (W^u_\delta (x)=\{y\in M |\,    d (f^{-n} (x),f^{-n} (y))\leq \delta e^{- (\mu-\varepsilon)n},\,n\geq 0\})$$ for some small $\delta>0.$
\end{Def}

\begin{Prop}\label{Prop-stable} (Stable Manifold Theorem)
Let $f:M\rightarrow M$ be a
$C^{1+\alpha}$ diffeomorphism and let  $\Lambda=\Lambda (\lambda,\mu;\varepsilon)$ a non-empty Pesin set  (with
$\lambda,\mu\gg\varepsilon>0$). There exists $\varepsilon_0>0$ such that for $x\in \Lambda (k\geq 1)$ and $\delta=\varepsilon_0e^{-\varepsilon k:}$:\\
 (a). $W^s_\delta (x),W^u_\delta (x) \text{are }C^1 \text{ submanifolds of M;}$\\
 (b). $T_xW^s_\delta=E^s_x,\,\,T_xW^u_\delta (x)=E^u_x.$
\end{Prop}

\subsection {New Pesin    Set, called Liao-Pesin   Set}

Parallel to $C^{1+\alpha}$ Pesin theory, we want to know whether Katok's shadowing lemma holds for $C^1$ systems. More precisely,

\begin{Que}\label{Que-LemShadow}
Let $f:M\rightarrow M$ be a
$C^{1}$ diffeomorphism. Is there  a non-empty Pesin set
$\Lambda
$ composed of a filtration of $\Lambda_k
$ 
such that

 (1) Katok's shadowing remains true?


 (2) Stable manifold theorem remains true?

\end{Que}

As noted in first section $C^1$ and $C^{1+\alpha}$ are different world so this generalization is difficult and  Pugh pointed out in
 \cite{Pugh} that the $\alpha$-H$\ddot{o}$lder condition for the
first derivative is necessary in the Pesin's stable manifold theorem. However, on one hand, inspired by Liao's idea of quasi-hyperbolic arc, we can give a {\it   partial but positive} answer to Katok's shadowing of Question \ref{Que-LemShadow} by  constructing a filtration of  new Pesin  blocks, called Liao-Pesin blocks and set. On the other hand, we also get stable manifolds on Liao-Pesin blocks.

Recall that a
so-called quasi-hyperbolic orbit segment  (Definition
\ref{def:quasihyp-orbit}) whose starting and ending points are near can be shadowed
by a periodic orbit from  Liao's closing lemma \cite{Liao79}. We are going to use this idea to realize our aim, that is, for an invariant measure, almost all orbit segments whose starting and ending points are near are quasi-hyperbolic and then can be shadowed by periodic points.  To guarantee any orbit segment in the
basin  of a hyperbolic measure to be quasi-hyperbolic  a condition called limit
domination is required.
More precisely, we can construct a filtration of new forms of Pesin blocks and Pesin
set such that  
all Pesin blocks have the same degree of mean
hyperbolicity  (Definition~\ref{def:degree-Pesinset}) and all
sufficiently long orbit segments with starting points and ending
points at the same block are of the same type of quasi-hyperbolicity
 (Definition~\ref{def:type-quasihyp}). These orbit segments satisfy the
conditions of the Liao  \cite{Liao79} closing lemma  so that they can be traced
by periodic orbits, which gives rise to the closing lemma. For obtaining
shadowing lemma in our setting, we apply a generalized idea  \cite{Gan} of
  Liao's closing lemma so that we can find a truth orbit to  trace any
pseudo-orbit consisting of finite or infinite orbit segments with starting and ending
points in a given  Pesin block. In particular, closing lemma deals with just one  orbit segment but shadowing lemma can deal with finite or infinite orbit segments so that shadowing lemma is much stronger.

\subsubsection{Establishment  of new Pesin set: Liao-Pesin set}

Now we start to introduce the definitions of our Liao-Pesin
blocks and set, and then state two theorems of shadowing and closing
lemma. Denote the minimal norm of an invertible linear map $A$ by
$m (A)=\|A^{-1}\|^{-1}$.

\begin{Def}\label{def:Pesinblock}
Given $K\in\mathbb{N},\,\,\zeta>0,$ and for all $k\in \mathbb{Z^+},$
we define $ \Lambda_k=\Lambda_k (K,\,\zeta)$ to be all points $x \in
M$ for which there is a splitting $T_{x}M=E (x)\oplus F (x)$ with the
invariance property $D_xf (E (x))=E (f (x))$ and $D_xf (F (x))=F (f (x))$
and satisfying:

$$ (a).\,\,\,\,\,\,\,\,\,\,\,\,\,\,\,\,\,\,\,\,\,\,\,\,\,\,\,\,\,\,\,\,
\frac{\log\|Df^r|_{E (x)}\|+\sum_{j=0}^{l-1}\log\|Df^{K}|_{E ({f^{jK+r} (x)})}\|}{lK+r}
   \leq-\zeta,\,\,\,\,\,\,\,\,\,\,\,\,\,\,\,\,\,\,\,\,\,\,\,\,\,\,\,\,\,\,\,\,\,\,\,\,
   \,\,\,\,\,\,\,\,\,\,\,\,\,\,\,\,\,\,\,\,\,\,\,\,\,\,\,\,\,\,\,\,\,\,\,\,\,\,\,\,\,
   \,\,\,\,\,\,\,\,\,\,\,\,\,\,\,\,\,\,\,\,\,\,\,\,\,\,\,$$$$\,\,\,\forall\,\,\,
l\,\geq \, k,\,\,\,\forall\,\,\,0\,\leq\, r\,\leq\,
K-1;\,\,\,\,\,\,\,\,\,\,\,\,\,\,\,\,\,\,\,\,\,\,\,\,\,\,\,\,\,\,\,\,\,\,\,\,\,\,\,$$

$$ (b).\,\,\,\,\,\,\,\,\,\,\,\,\,\,\,\,\,\,\,\,\,\,\,\,\,\,\,
\,\, \frac{\log\,m (Df^r|_{F (f^{-lK-r} (x))})+\sum_{j=-l}^{-1}\log
  m (Df^{K}|_{F ({f^{jK} (x)})})}{lK+r}\geq\zeta,\,\,\,\,\,\,\,\,\,\,\,\,\,\,\,\,\,\,
  \,\,\,\,\,\,\,\,\,\,\,\,\,\,\,\,\,\,\,\,\,\,\,\,\,\,
  \,\,\,\,\,\,\,\,\,\,\,\,\,\,\,\,\,\,\,\,\,\,\,
  \,\,\,\,\,\,\,\,$$$$\,\,\,\forall\,\,\,
l\,\geq \, k,\,\,\,\forall\,\,\,0\,\leq \,r\,\leq\,
K-1;\,\,\,\,\,\,\,\,\,\,\,\,\,\,\,\,\,\,\,\,\,\,\,\,\,\,\,\,\,\,\,\,\,\,\,\,\,\,\,$$

$$ (c).\,\,\,\,\,\,\,\,\,\,\,\,\,\,\,\,\,\,\,\,\,\,\,\,\,\,\,\,\,\,\,\,\,\,\,\,\,
\frac1{kK+r}\log\frac{\|Df^{kK+r}|_{E (x)}\|}{m (Df^{kK+r}|_{F (x)})}
\leq-2\zeta,\,\,\,\forall\,\,\,0\,\leq \,r\,\leq\,
K-1,\,\,\,\,\,\,\,\,\,\,\,\,\,\,\,\,\,\,\,\,\,\,\,\,\,\,\,\,\,\,\,\,\,\,\,\,\,\,\,\,\,\,\,\,\,\,\,
\,\,\,\,\,\,\,\,\,\,\,\,\,\,\,\,\,\,\,\,\,\,\,\,\,\,\,\,\,$$$$and
\,\,\,\,\,\,\,\,\,\,\,\,\,\,\,\,\,\,
\frac1K\log\frac{\|Df^K|_{E (f^{l} (x))}\|}{m (Df^K|_{F (f^{l} (x))})}\leq-2\zeta,\,\,\,\,
\forall\, \,\,l\,\geq \,kK.
\,\,\,\,\,\,\,\,\,\,\,\,\,\,\,\,\,\,\,\,\,\,\,\,\,\,\,\,\,\,\,\,\,\,\,\,\,\,\,\,\,\,\,\,\,\,\,\,\,\,\,\,\,\,$$
Denote by $\Lambda=\Lambda (K, \zeta)$ the maximal $f-$invariant
subset of $\bigcup_{k\geq 1}\Lambda_k$, meaning
$$\Lambda=\bigcap_{n\in\mathbb{Z}}f^n (\bigcup_{k\geq 1}\Lambda_k).$$ We call
$\Lambda$ a Liao-Pesin set and call $\Lambda_k (k\geq 1)$ Liao-Pesin blocks.
 (see Figure~\ref{pic:graph2.eps} and~\ref{pic:limitdomination} to
explain  (a),  (b) and  (c), respectively).
\end{Def}

 \setlength{\unitlength}{1mm}
  \begin{figure}[!htbp]
  \begin{center}
  \begin{picture} (180,50) (0,0)
  \put (0,0){\scalebox{1}[1]{\includegraphics[0,0][50,30]{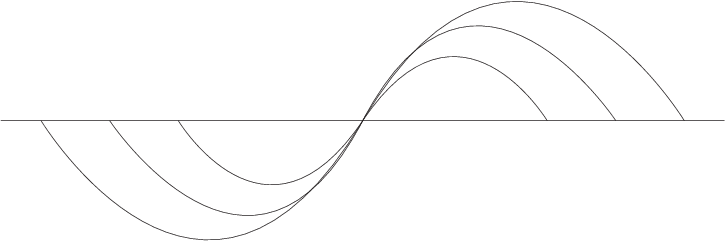}}}
  \put (67,27){${_0}$}\put (35,33){$_{_{_{{-kK-r}}}}$}\put (19,33){$_{_{_{{{- (k+1)K-r}}}}}$}
\put (1,33){$_{_{_{{{- (k+2)K-r}}}}}$}
\put (94,27){$_{_{{kK+r}}}$}\put (103,27){$_{_{{{ (k+1)K+r}}}}$}
\put (118,27){$_{_{{{ (k+2)K+r}}}}$}
  \put (8.8,20){${\cdots\cdots}$}\put (118,38){${\cdots\cdots}$}
  \put (93,50.1){$_{>}$}\put (85,45.9){$_{>}$}
  \put (80,41){$_{>}$}\put (37,10.2){$_{>}$}\put (46,14.55){$_{>}$}
  \put (52,19.53){$_{>}$}
  \end{picture}
  \caption{Graph to show  (a) and  (b).}
  \label{pic:graph2.eps}
  \end{center}
  \end{figure}
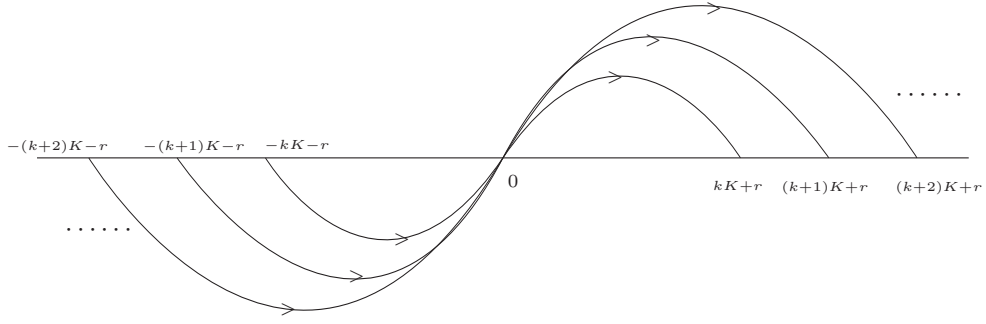
\setlength{\unitlength}{1mm}
  \begin{figure}[!htbp]
  \begin{center}
  \begin{picture} (180,32) (0,0)
  \put (0,0){\scalebox{1}[1]{\includegraphics[0,0][40,50]{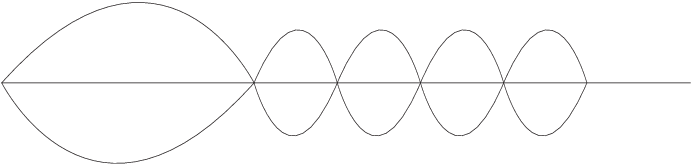}}}
  \put (4,14){$_{0}$}
\put (45,14){$_{_{kK+r}}$}
 \put (72,14){$_{_{ (k+2)K+r}}$}
  \put (58,14){$_{_{ (k+1)K+r}}$}
  \put (86,14){$_{_{ (k+3)K+r}}$}
  \put (100,14){$_{_{ (k+4)K+r}}$}
  \put (105,18.8){$\cdots\cdots$}
\put (28,0.8){$_{>}$}\put (28,27.4){$_{>}$}
\put (68.2,5.2){$_{>}$}\put (68.2,22.9){$_{>}$}
\put (82.2,5.2){$_{>}$}\put (82.2,22.9){$_{>}$}
\put (97.2,5.2){$_{>}$}\put (97.2,22.9){$_{>}$}
\put (55.2,5.2){$_{>}$}\put (55.2,22.9){$_{>}$}
    \end{picture}
  \caption{Graph to show  (c).}
  \label{pic:limitdomination}
  \end{center}
  \end{figure}
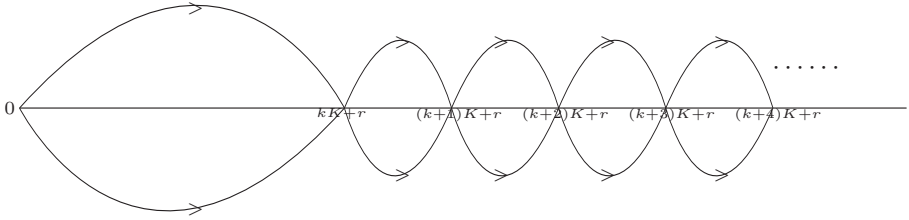

\begin {Rem}For  given   $K\in\mathbb{N}$ and $\zeta>0$, obviously
$\Lambda_k (K,\,\zeta)\subseteq\Lambda_{k+1} (K,\,\zeta),\,\,\forall\,
k \in \mathbb{N}$ and
$\Lambda (K,\,\zeta)\subseteq\Lambda (iK,\zeta),\,\,\forall\, i\in
\mathbb{N}.$ By sub-multiplication of norms, the conditions  (a) and  (b) imply that

$$ (a').\,\,\,\,\,\,\,\,\,\,\,\,\,\,\,\,\,\,\,\,\,\,\,\,\,\,\,\,\,\,\,\,
\frac{ \log\|Df^{n}|_{E (x)}\|}{n}
   \leq-\zeta,\,\,\,\,\,\,\,\,\,\,\,\,\forall\,\,\,\,\,\,
n\,\geq \, kK,\,\,\,  \,\,\,\,\,\,\,\,\,\,\,\, \,\,\,\,\,\,\,\,\,\,\,\,\,\,\,\,\,\,\,\,\,\,\,\,\,\,\,\,\,\,\,\,\,\,\,\,\,\,\,\,\,
   \,\,\,\,\,\,\,\,\,\,\,\,\,\,\,\,\,\,\,\,\,\,\,\,\,\,\,\,\,\,\,\,\,\,\,$$

$$ (b').\,\,\,\,\,\,\,\,\,\,\,\,\,\,\,\,\,\,\,\,\,\,\,\,\,\,\,
\,\, \frac{ \log
  m (Df^{n}|_{F (x)})}{n}\geq\zeta,\,\,\,\,\,\,\,\, \,\,\,\,\,\,\,\forall\,\,\,
n\,\geq \, kK .\,\,\,\,\,\,\,\,\,\,\,\,\,\,\,\,\,\,\,\,\,\,\,\,\,\,\,\,\,\,\,\,\,\,\,\,
  \,\,\,\,\,\,\,\,\,\,\,\,\,\,\,\,\,\,\,\,\,\,\,\,\,\,
  \,\,\,\,\,\,\,\,\,\,\,\,\,\,\,\,\,\,\,\,\,\,\,\,\,\,\,\,\,$$ Different with the construction of classical Pesin set, it is not necessarily required that $f^{\pm1}\Lambda_k\subseteq \Lambda_{k+1}.$ We will illustrate   Liao-Pesin  set more information in
Section \ref{section-New Pesin set}.
\end {Rem}
The definition of our Liao-Pesin set is based on a generalized
multiplicative ergodic theorem (Lemma
\ref{Lem:Generalized-Multi-erg-thm}) and limit domination
 (Definition \ref{def:limit-dom}). It enables us to realize
shadowing properties on nonempty Liao-Pesin blocks by using
Liao's closing and shadowing  lemma in a $C^1$ nonuniformly hyperbolic system with
limit domination.
Remark that for the case of flows, the constructed Liao-Pesin blocks in  \cite{SunTianV} display minor difference because there it needs to deal with continuous time.

\subsubsection{Basic properties of Liao-Pesin  set}\label{section-New Pesin set}

In this section we investigate more properties of Liao-Pesin set.
We recall a notion of Liao's quasi-hyperbolic orbit segment
 \cite{Gan,Liao79}.

\begin{Def}\label{def:quasihyp-orbit} Fix arbitrarily two constants $\zeta>0$
and $e\in\mathbb{Z}^+$ and consider an orbit segment
$$\{x,\, n\}:=\{ f^i (x)\,| \quad i=0,\,1,\,2,\,\cdots,\,n\},$$
where $x\in M$ and $n\in\mathbb{N}.$ We call $\{x,\, n\}$  a
$ (\zeta,\, e)$-quasi-hyperbolic
  orbit segment with respect to a
splitting
$$T_xM=E\oplus F,$$
if  
there is a partition
$$0=t_0<t_1<\cdots<t_m=n\,\,\,\,  (m\geq 1)$$
such that $   t_k-t_{k-1}\leq e$ and \\
 (1). $\,\,\,\,\,\, \,\frac 1{t_{k}}\Sigma_{j=1}^k
\log\|Df^{t_j-t_{j-1}}|_{ Df^{t_{j-1}} (E)}\|\leq -\zeta,$\\
 (2). $\,\,\,\,\,\,\, \frac 1{t_{m}-t_{k-1}}\Sigma_{j=k}^m
\log\,m (Df^{t_j-t_{j-1}}|_{ Df^{t_{j-1}} (F)})\geq  \zeta,$\\
 (3). $\,\,\,\,\,\,\,\frac
1{t_{k}-t_{k-1}}\log\frac{\|Df^{t_k-t_{k-1}}|_{ Df^{t_{k-1}} (E)}\|}{m (Df^{t_k-t_{k-1}}|_{ Df^{t_{k-1}} (F)})}
\leq -2\zeta,\,\,\, k=1,\, 2, \,\cdots,\, m.$
\end{Def}

We use the notion of quasi-hyperbolic orbit segment to introduce a
concept of \emph{type of quasi-hyperbolicity}.

\begin{Def}\label{def:type-quasihyp} Let $k,\,K\in\mathbb{N},\,\,\zeta>0$. For a  given
subset $\Delta$  (neither necessarily $f-$invariant  nor
$f^K$-invariant maybe), let $T_{x}M=E (x)\oplus F (x)\,\, (x\in
\Delta)$ be a $Df$-invariant splitting, meaning that it is invariant
on each orbit $orb (x, f)$  ($\Delta$ contains not necessarily the
whole orbit $orb (x, f)$) for
 $x\in \Delta.$ We say that $\Delta$ is a quasi-hyperbolic set of
$ (\zeta, (k+1)K)$-type, if any orbit segment  $\{x,n\}:=\{f^i (x)|
i=0,\,1,\,\cdots,n\}$ $ (n\geq2kK)$ with starting point and ending
point in $\Delta$ is $ (\zeta, (k+1)K)$-quasi-hyperbolic orbit segment
of $f$.
\end{Def}

Now we use the above notions to present more properties for Liao-Pesin
blocks and set.
\begin{Prop}\label{Prop:prop-Pesinblocks}
For given  $k,\,K\in\mathbb{N}$ and $\zeta>0$, the Liao-Pesin block
$\Lambda_k (K,\,\zeta)$ is  closed  and is a
 quasi-hyperbolic set of $ (\zeta, (k+1)K)$-type, and the splitting
$T_{x}M=E (x)\oplus  F (x)$ on  $\Lambda_k (K,\,\zeta)$ is continuous.
Further, the Liao-Pesin set $\Lambda (K,\,\zeta)$ is an average-nonuniformly
hyperbolic set with $ (K,\,\zeta)$-degree.

\end{Prop}

{\bf Proof}
 Given an orbit segment  $\{x,\,n\}:=\{f^i (x)\,|\,
i=0,\,1,\,\cdots,\,n\}$ $ (n\,\geq\,2kK)$ with starting point and
ending point in $\Lambda_k (K,\,\zeta)$, i.e., $x,\,
f^n (x)\in\Lambda_k (K,\,\zeta)$, we show that $\{x,\,n\}$ is
$ (\zeta,\, (k+1)K)$-quasi-hyperbolic orbit segment.

Write $n=lK+q,$ where $l\geq 2k$ and $0\leq q \leq K-1.$ Let
$m=l-2k+2$ and $$t_0=0,\,\,t_i= (k+i-1)K+q,
\,\,i=1,\,2,\,\cdots,\,m-1,\,\,t_m=lK+q.$$ Thus we get a partition
 $$t_0<\,t_1<\,t_2<\,\cdots<\,t_{m-1}<\,t_m.$$
In the partition, the starting subinterval has length of $kK+q$, the
ending subinterval has length of $kK$, and the rest  have length of
$K.$

Since $x\in\Lambda_k (K,\,\zeta)$, by  (a) in Definition
\ref{def:Pesinblock} while taking $r=q$ we have
$$ \frac 1{t_{i}}\Sigma_{j=1}^i \log\|Df^{t_j-t_{j-1}}|_{ Df^{t_{j-1}} (E (x))}\|$$$$\leq
\frac{\log\|Df^q|_{E (x)}\|+\sum_{j=0}^{k+i-2}\log\|Df^{K}|_{E ({f^{jK+q} (x)})}\|}{q+ (k+i-1)K}\leq
-\zeta,$$ $\,\,i=1,\, 2,\, \cdots,\, m-1,$ and
$$ \frac
1{t_{m}}\Sigma_{j=1}^m \log\|Df^{t_j-t_{j-1}}|_{ Df^{t_{j-1}} (E (x))}\|$$$$\leq
\frac{\log\|Df^q|_{E (x)}\|+\sum_{j=0}^{l-1}\log\|Df^{K}|_{E ({f^{jK+q} (x)})}\|}{q+lK}\leq
-\zeta,$$ which gives rise to  the first inequality in Definition
\ref{def:quasihyp-orbit}.

Notice that  $$l-k-i+2\geq l-k-m+2=k,\,\,i=1,\,2,\,\cdots,\,m.$$
Since
 $f^n (x)\in\Lambda_k (K,\,\zeta)$, by  (b) in the Definition
\ref{def:Pesinblock} while taking $r=0$ we have
$$ \frac 1{t_{m}-t_{i-1}}\Sigma_{j=i}^m \log\,m (Df^{t_j-t_{j-1}}|_{ Df^{t_{j-1}} (F (x))})$$$$\geq \sum_{j=- (l-k-i+2)}^{-1}\frac{\log
  m (Df^{K}|_{F ({f^{jK} (f^n (x))})})}{{ (l-k-i+2)K}}\geq\zeta,$$ $\,\,i= 2, \cdots,
  m,$ and while taking $r=q$ we have
$$ \frac 1{t_{m}-t_{0}}\Sigma_{j=1}^m \log\,m (Df^{t_j-t_{j-1}}|_{ Df^{t_{j-1}} (F (x))})$$$$\geq\frac{\log
m (Df^q|_{F (x)})+\sum_{j=-l}^{-1}{\log
  m (Df^{K}|_{F ({f^{jK} (f^n (x))})})}}{{lK+q}}\geq\zeta,$$
which gives rise to the  second inequality in Definition
\ref{def:quasihyp-orbit}.

Before continuing our proof we present  Figure
\ref{pic:quasi-hyperbolic-orbit} to illustrate the concepts.  We
denote $f^{t_i} (x)$ by $t_i$ and take $k=3,\,\,l=10,\,\,q=1,\,\,m=6$
and $ n=10K+1$ in the  Figure, and draw  a graph for the inequality
 (1) and  (2) of $ (\zeta,4K)= (\zeta,  (3+1)K)$ quasi-hyperbolic orbit
segment $\{x,10K+1\}$.
 \setlength{\unitlength}{1mm}
  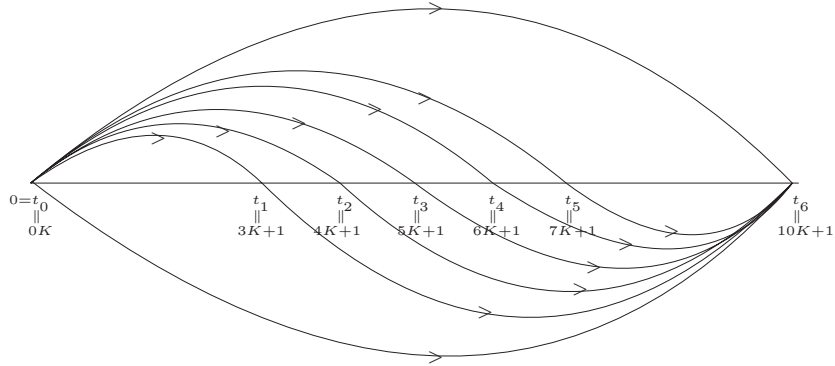
\begin{figure}[!htbp]
  \begin{center}
  \begin{picture} (120,58) (0,0)
  \put (0,0){\scalebox{1}[1]{\includegraphics[0,0][40,50]{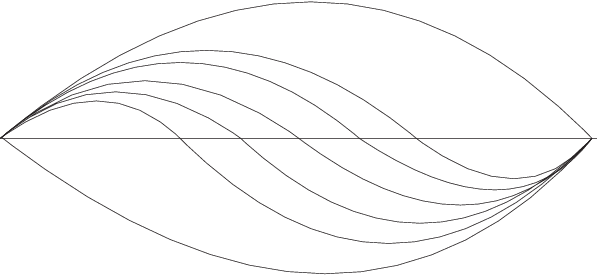}}}
  \put (5,30){$_{_{0=t_{_0}}}$}\put (8,28){$_{_{\|}}$}\put (7.5,26){$_{_{0K}}$}
  \put (37,30){$_{_{t_{_1}}}$} \put (37,28){$_{_\|}$}  \put (35,26){$_{_{3K+1}}$}
\put (48,30){$_{_{t_{_2}}}$}\put (48,28){$_{_{\|}}$}\put (45,26){$_{_{4K+1}}$}
\put (58,30){$_{_{t_{_3}}}$}\put (58,28){$_{_{\|}}$}\put (56,26){$_{_{5K+1}}$}
  \put (68,30){$_{_{t_{_4}}}$} \put (68,28){$_{_{\|}}$}\put (66,26){$_{_{6K+1}}$}
  \put (78,30){$_{_{t_{_{5}}}}$}\put (78,28){$_{_{\|}}$}\put (76,26){$_{_{7K+1}}$}
\put (108,30){$_{_{t_{_6}}}$}\put (108,28){$_{_{\|}}$}\put (106,26){$_{_{10K+1}}$}
\put (23.5,37.8){$_{>}$}
\put (32,38.7){$_{>}$}\put (42,39.83){$_{>}$}\put (52,41.6){$_{>}$}\put (58.6,43.1){$_{>}$}
\put (60,55){$_{>}$}\put (60,8.8){$_{>}$}\put (66.5,14.8){$_{>}$}\put (79,17.8){$_{>}$}
\put (80.8,20.7){$_{>}$}\put (85,23.8){$_{>}$}\put (91,25.5){$_{>}$}
  \end{picture}
  \caption{inequalities  (1),  (2) of $ (\zeta,4K)$ quasi-hyperbolic orbit segment.}
  \label{pic:quasi-hyperbolic-orbit}
  \end{center}
  \end{figure}

Now we continue our proof and verify  the third inequality in
Definition~\ref{def:quasihyp-orbit}.  Since
$x\in\Lambda_k (K,\,\zeta)$  and $t_1=kK+q$ , by  (c) in Definition
\ref{def:Pesinblock} while taking $r=q$ we have
$$\frac 1{t_{1}-t_{0}}\log\frac{\|Df^{t_1-t_{0}}|_{ Df^{t_{0}} (E (x))}\|}{m (Df^{t_1-t_{0}}|_{ Df^{t_{0}} (F (x))})}$$$$=
\frac1{kK+q}\log\frac{\|Df^{kK+q}|_{E (x)}\|}{m (Df^{kK+q}|_{F (x)})}\leq
-2\zeta,$$ while noting that $t_i\geq kK,\,i=1,\,2,\,\cdots,\,m-2$,
we have
$$\frac
1{t_{i}-t_{i-1}}\log\frac{\|Df^{t_i-t_{i-1}}|_{ Df^{t_{i-1}} (E (x))}\|}{m (Df^{t_i-t_{i-1}}|_{ Df^{t_{i-1}} (F (x))})}$$
$$\,\,\,\,\,\,\,\,\,\,\,\,\,\,\,\,\,\,\,\,\,\,\,\,\,\,\,=\frac
1{K}\log\frac{\|Df^{K}|_{ Df^{t_{i-1}} (E (x))}\|}{m (Df^{K}|_{ Df^{t_{i-1}} (F (x))})} \leq -2\zeta,\,\, i=2,\,\cdots,\, m-1,$$ and
while $t_m-t_{m-1}=kK$ and $t_{m-1}+jK\geq kK$, we have $$\frac
1{t_{m}-t_{m-1}}\log\frac{\|Df^{t_m-t_{m-1}}|_{ Df^{t_{m-1}} (E (x))}\|}{m (Df^{t_m-t_{m-1}}|_{ Df^{t_{m-1}} (F (x))})}$$$$\leq\frac
1{kK}\sum_{j=0}^{k-1}\log\frac{\|Df^{K}|_{ Df^{t_{m-1}+jK} (E (x))}\|}{m (Df^{K}|_{ Df^{t_{m-1}+jK} (F (x))})} \leq
-2\zeta,$$ which gives rise to the third inequality  in Definition
\ref{def:quasihyp-orbit}.

Before continuing our proof we also present  Figure
\ref{pic:dominationgraph} to give an explanation. We denote
$f^{t_i} (x)$ by $t_i$ and take $k=3,\,\,l=10,\,\,q=1,\,\,m=6$ and $
n=10K+1$ in the  Figure, and draw a graph for the inequality  (3) of
$ (\zeta,4K)= (\zeta,  (3+1)K)$ quasi-hyperbolic orbit segment
$\{x,10K+1\}$.

\setlength{\unitlength}{1mm}
  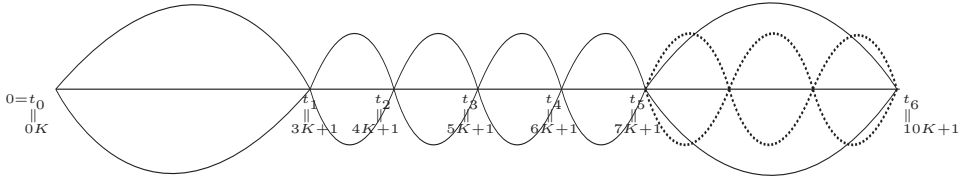
\begin{figure}[!htbp]
  \begin{center}
  \begin{picture} (180,30) (0,0)
  \put (0,0){\scalebox{1}[1]{\includegraphics[0,0][40,50]{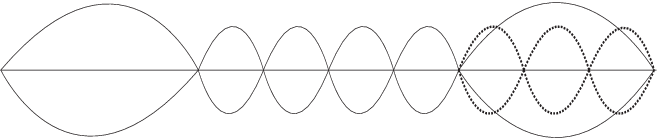}}}
  \put (2,12){$_{_{0=t_{_0}}}$}\put (5,10){$_{_{\|}}$}\put (4.5,8){$_{_{0K}}$}
  \put (41,12){$_{_{t_{_1}}}$} \put (41,10){$_{_\|}$}  \put (39.5,8){$_{_{3K+1}}$}
\put (62,12){$_{_{t_{_3}}}$}\put (62,10){$_{_{\|}}$}\put (60,8){$_{_{5K+1}}$}
\put (50.5,12){$_{_{t_{_2}}}$}\put (50.5,10){$_{_{\|}}$}\put (47.5,8){$_{_{4K+1}}$}
  \put (73,12){$_{_{t_{_4}}}$} \put (73,10){$_{_{\|}}$}\put (71,8){$_{_{6K+1}}$}
  \put (84,12){$_{_{t_{_{5}}}}$}\put (84,10){$_{_{\|}}$}\put (82,8){$_{_{7K+1}}$}
\put (120,12){$_{_{t_{_6}}}$}\put (120,10){$_{_{\|}}$}\put (120,8){$_{_{10K+1}}$}
  \end{picture}
  \caption{inequality  (3) of $ (\zeta,4K)$ quasi-hyperbolic orbit segment.}
  \label{pic:dominationgraph}
  \end{center}
  \end{figure}

Observe from the partition that $t_i-t_{i-1}\leq  (k+1)K.$ So
$\{x,n\}$ is a $ (\zeta, (k+1)K)$-quasi-hyperbolic orbit segment.
Thus, by Definition~\ref{def:type-quasihyp} $\Lambda_k (K,\,\zeta)$
is a
 quasi-hyperbolic set of $ (\zeta, (k+1)K)$-type.

Next we show that  $\Lambda_k (K,\,\zeta)$ is closed. Clearly the
three conditions in the definition of $\Lambda_k (K,\,\zeta)$ imply
that the splitting $T_{x}M=E (x)\oplus F (x)$ is unique. If $x\in M$
and $x_i\in \Lambda_k (K,\,\zeta) (i\geq1)$ is a convergent sequence
with $\lim_{i\rightarrow+\infty}x_i=x$, with the choice of $k$
fixed, then by a compactness argument we can choose a convergent
subsequence of the subspaces ${\xi} ({x_{i_j}})\rightarrow
{\xi'} ({x}),$ as $j\rightarrow+\infty$ where $\xi=E,\,\,F$. By
assumption $x_{i_j}\in \Lambda_k (K,\,\zeta),$ conditions in the
definition of $\Lambda_k (K,\,\zeta)$ are satisfied by $E ({x_{i_j}})$
and $F ({x_{i_j}})$. Letting $j\rightarrow +\infty$, then the three
conditions in the definition of $\Lambda_k (K,\,\zeta)$ hold for the
subbundles $E' ({x})$ and $F' ({x})$. By the uniqueness condition
above, $E' ({x})=E ({x})$ and $F' ({x})=F ({x})$. So
$x\in\Lambda_k (K,\,\zeta)$ and thus $\Lambda_k (K,\,\zeta)$ is
closed.

By the uniqueness condition above, there is only one possible limit
for $E ({x_{i_j}})$ and $F ({x_{i_j}})$. Thus the splitting $x\mapsto
E ({x})\oplus F ({x})$ is continuous on $\Lambda_k (K,\,\zeta)$. That
$\Lambda (K,\,\zeta)$ is a mean nonuniformly hyperbolic set with
$ (K,\,\zeta)$-degree is an easy consequence.

By  (a) and  (b) in Definition~\ref{def:Pesinblock}, it is easy to see
that the Pesin set $\Lambda (K,\,\zeta)$ is a  average-nonuniformly
hyperbolic set with $ (K,\,\zeta)$-degree.\hfill $\Box$

\bigskip

For the sets $\Lambda_k (K,\,\zeta),
k\geq1,$ we should observe that, in general, $\Lambda$ itself need
not necessarily be compact, nor  it is necessarily true that the
splitting $T_{x}M=E (x)\oplus  F (x)$ is continuous on
$\Lambda (K,\,\zeta).$



\medskip

\subsubsection{shadowing lemma and closing lemma}

To state shadowing lemma and closing lemma we need some notions.
Given $x\in M$ and $n\in\mathbb{N}$, let
   $$\{x,\,n\}:=
 \{f^j (x)\,|\,\,j=0,\,1,\,\cdots,\,n\}.$$ In other words, $\{x,\,n\}$
represents the orbit segment from $x$. For a sequence of points
$\{x_i\}_{i=-\infty}^{+\infty}$ in $M$
 and a sequence of positive integers
 $\{n_i\}_{i=-\infty}^{+\infty}$, we call $\{x_i,\,n_i\}_{i=-\infty}^{+\infty}$
 a $\delta$-pseudo-orbit,
 if $$   d (f^{n_i} (x_i),\,x_{i+1})<\delta$$ for all $i$.
  Given $\varepsilon>0,$ we call a point $x\,\in
M$  an $\varepsilon$-shadowing point for a pseudo-orbit
$\big{\{}x_i,\,n_i\big{\}}_{i=-\infty}^{+\infty},$  if
$   d\big{ (}f^{c_i+j} (x),f^j (x_i)\big{)}<\varepsilon$, $\forall\,\,
j=0,\,1,\,2,\,\cdots,\,n_i-1$ and $\forall\,\, i \in \mathbb Z$, where
$c_i$
 is defined as
  \begin {equation} \label{eq:respective-time-squens}c_i=\begin{cases}
 0,&\text{for }i=0\\
 \sum_{j=0}^{i-1}n_j,&\text{for }i>0\\
 -\sum_{j=i}^{-1}n_j,&\text{for }i<0.
\end{cases}
 \end {equation}
Now we state shadowing lemma and closing lemma on nonempty Pesin
blocks (Nonempty discussion will appear in next subsection).

\begin{Thm}\label{Thm:shadowinglem}  (Shadowing lemma) If $\Lambda_k (K,\,\zeta)\neq\emptyset$ for
some parameters  $k,\,K\in\mathbb{N}$ and $\zeta>0$, then $\Lambda_k (K,\,\zeta)$
satisfies the following shadowing property. For
$\forall\,\,\varepsilon>0,$ there exists $\delta>0$ such that if  a
$\delta$-pseudo-orbit $\{x_i,\,n_i\}_{i=-\infty}^{+\infty}$
satisfies $n_i\geq2kK $ and $x_i,f^{n_i} (x_i)\in
\Lambda_k (K,\,\zeta)$ for all $i$, then there exists a unique
$\varepsilon$-shadowing point $x\in M $ for
$\{x_i,n_i\}_{i=-\infty}^{+\infty}$. If further
$\{x_i,n_i\}_{i=-\infty}^{+\infty}$ is periodic, i.e., there
 exists an $m>0$ such that $x_{i+m}=x_i$ and $n_{i+m}=n_i$ for all i,
 then the shadowing point $x$ can be chosen to be periodic.
\end{Thm}

To  complete the proof of Theorem~\ref{Thm:shadowinglem}, we need
Liao's closing lemma  \cite{Liao79} and its generalization by Gan
 \cite{Gan}. Before that we need to a concept of quasi-hyperbolic
pseudo-orbit. Let $\zeta>0,\,\,e\in\mathbb{Z}^+,\,\,\delta>0.$
Given a sequence of orbit segments
$\{x_i,\,n_i\}_{i=-\infty}^{+\infty},$ we call
$\{x_i,\,n_i\}_{i=-\infty}^{+\infty}$ a  $ (\zeta,\,e
)$-quasi-hyperbolic $\delta$-pseudo-orbit with respect to splittings
$T_{x_i}M=E ({x_i})\oplus F ({x_i}),$ if 
$\{x_i,\,n_i\}_{i=-\infty}^{+\infty}$ is  a $\delta$-pseudo-orbit
and every orbit segment $\{x_i,n_i\}$ is a
$ (\zeta,\,e)$-quasi-hyperbolic  orbit segment with respect to
the $i$-th splitting  $T_{x_i}M=E  ({x_i})\oplus F ({x_i}).$ Now we
state Liao' closing lemma  \cite{Liao79} and its generalization, shadowing lemma
 \cite{Gan}.

\begin{Lem}\label{Lem:Liao-Gan-shadlem}  ( \cite{Gan,Liao79})
For any $\zeta>0$, $e\in\mathbb{Z}^+$, there exist $L>0,\,\, d_0
>0$ with following property for any $d\in (0,d_0]$. If for a
$ (\zeta,\,e)$-quasi-hyperbolic $d$-pseudo-orbit
$\{x_i,n_i\}_{i=-\infty}^{+\infty}$ with respect to $Df-$invariant
splittings $T_{x_i}M=E  ({x_i})\oplus F ({x_i})$, one has
\begin{eqnarray}\label{eqn-Bundle-nbhd} Df^{n_i} ( \xi ({x_i}))\cap T^\#M \subseteq U ( \xi ({x_{i+1})} \cap
T^\#M,d)\,\, (\xi=E,F), \text{ for all } i,
\end{eqnarray}  then there exists a unique
$Ld$-shadowing point $x\in M $ for
$\{x_i,n_i\}_{i=-\infty}^{+\infty}$. \\ Moreover, if
$\{x_i,n_i\}_{i=-\infty}^{+\infty}$ is periodic, i.e., there
 exists an $m>0$ such that $x_{i+m}=x_i$ and $n_{i+m}=n_i$ for all i,
 then the shadowing point $x$ can be chosen to be periodic with period $c_m$, where
$c_i$  is the same as  in  (\ref {eq:respective-time-squens}).

\end{Lem}

 In  \cite{Gan}, Lemma \ref{Lem:Liao-Gan-shadlem} is stated for the case $e=1$ and is assumed dominated splitting.  It     is stated  \cite{Liao79}   for  one orbit segment (That is, the particular case of Lemma \ref{Lem:Liao-Gan-shadlem} for $m=1$) in the case of flows. Remark that the method there is still suitable for any $e$ and the proof is almost similar. However, Lemma \ref{Lem:Liao-Gan-shadlem} is an important technique for present article and notice that its proof is very short, so that we will give a proof in section \ref{section-proof-shadowingliao} for completeness.



\bigskip

{\bf Proof of Theorem~\ref{Thm:shadowinglem}} For $\zeta>0$ and
$e= (k+1)K$, by Lemma~\ref{Lem:Liao-Gan-shadlem} there exist
$L=L (k,K,\zeta)>0,\,\, d_0=d_0 (k,K,\zeta)>0$ with the shadowing
property as stated in Lemma~\ref{Lem:Liao-Gan-shadlem} for any
$d\in (0,d_0]$.

Now we consider $d<\frac{\varepsilon}L.$ By Proposition
\ref{Prop:prop-Pesinblocks}, $\Lambda_k$ is closed, and the stable
subbundle and unstable subbundle in the Oseledec splitting are
 continuous restricted on   $\Lambda_k.$  Take  $\delta\in (0,d)$
 such that for all $x,y \in \Lambda_k$
 with $   d (x,y)<\delta$ we have $E^\xi_{x}\cap T^\#M \subseteq U (E^\xi_{y}
\cap T^\#M,d),\,\,\xi=s,u.$

For  a given $\delta$-pseudo-orbit
$\{x_i,n_i\}_{i=-\infty}^{+\infty}$ with  $x_i,f^{n_i} (x_i)\in
\Lambda_k$ and $n_i\geq2kK  (\forall i)$, we get by Proposition
\ref{Prop:prop-Pesinblocks}  that every orbit segment $\{x_i,n_i\}$
is $ (\zeta,\,e)$ quasi-hyperbolic orbit segment.  So
$\{x_i,n_i\}_{i=-\infty}^{+\infty}$ is a $ (\zeta,\,e)$
quasi-hyperbolic $d$-pseudo-orbit. That
$   d (f^{n_i} (x_i),x_{i+1})<\delta$ implies that
$$Df^{n_i} (E^\xi ({x_i}))\cap T^\#M \subseteq U (E^\xi ({x_{i+1})} \cap
T^\#M,d)\,\, (\xi=s,u)\text{for all } i.$$ By Lemma
\ref{Lem:Liao-Gan-shadlem} there exists a $Ld$-shadowing point $x\in
M $ for $\{x_i,n_i\}_{i=-\infty}^{+\infty}$ and thus $x$ is also a
$\varepsilon$-shadowing point. In particular, if
$\{x_i,n_i\}_{i=-\infty}^{+\infty}$ is periodic, i.e., there exists
an $m>0$ such that $x_{i+m}=x_i$ and $n_{i+m}=n_i$ for all $i$, then
the point $x$ can be chosen to be  periodic and thus we complete the
proof. \hfill $\Box$

Taking $m=1$ from Theorem~\ref{Thm:shadowinglem}, one deduces the
closing lemma.

\begin{Thm}\label{Thm:closinglem}  (Closing lemma) If $\Lambda_k (K,\,\zeta)\neq\emptyset$ for
some $k,\,K\in\mathbb{N}$ and $\zeta>0$, then $\Lambda_k (K,\,\zeta)$
satisfies closing property in the following sense. For $\forall
\varepsilon>0,$ there exists $\delta>0$ such that if for an orbit
segment $\{x,n\}$ with length $n\geq 2kK$, one has $x,f^{n} (x)\in
\Lambda_k (K,\,\zeta)$ and $   d (x,f^n (x))<\delta$, then there exists
a unique point $z=z (x)\in M $ satisfying:

 (1) $f^n (z)=z$;

 (2) $   d (f^i (x),\,\,f^i (z))<\varepsilon$, $i=0,\,1,\,\cdots,n-1$.
\end{Thm}

A direct application is to find periodic orbits.

\begin{Thm}\label{Thm:closinglem-find-periodic} Let $f\in \Diff^1 (M).$ Suppose that $\Lambda_k (K,\,\zeta)\neq\emptyset$ for
some $k,\,K\in\mathbb{N}$ and $\zeta>0$ and there exists some point $y\in M$ and a sequence of $\{n_i\}\uparrow\infty$  such that $f^{n_i} (y)\in \Lambda_k (K,\,\zeta).$ Then  the system $f$ has at least one periodic orbit.
\end{Thm}

Remark that the assumption of this result is a topological or analytical  condition independent of measures and so maybe there are other applications in future.

{\bf Proof.}
Given $\epsilon>0$, take $\delta>0$ as Theorem \ref{Thm:closinglem}. If necessary, take a subsequence of $\{n_i\}$, denoted by $\{m_i\},$ such that $m_{i+1}-m_i\geq 2kK.$ By the compactness of $M,$  we can  take large $i<j$ such that
$$d (f^{m_i} (y),f^{m_j} (y))<\delta.$$ Let $x=f^{m_i} (x)$ and $n=m_j-m_i.$ Then
$$x,f^n (x)\in \Lambda_k (K,\zeta),\,n\geq 2kK$$ and so by Theorem \ref{Thm:closinglem}, $f$ has at least one periodic orbit. \qed

\subsection{Improved Liao-Pesin blocks, stable manifolds and exponential shadowing}\label{StableMnfld-Clolem}

In this section all statements are independent on measures. The main aim of this section is preparing stable manifold theorem for proving Theorem \ref{Thm:ExpShadowing-measure}.

\subsubsection{Improved Liao-Pesin blocks}

If we add dominated condition in the definition of Liao-Pesin set:
$$ (d)\,\,\,\,\,\,\,\,\,\,\,\,\,\,\,\,\,\,\,\,\,\,\,\,\,\,\,\,\,\,\,\,\,\,\,\,\,\,\,\,
\frac1{kK}\log\frac{\|Df^{kK}|_{E (f^{l} (x))}\|}{m (Df^{kK}|_{F (f^{l} (x))})}\leq-2\zeta,\,\,\,\,\,\,
\forall\, \,\,l\in \mathbb{Z},\, \,\,\,\,\,\,\,\,\,\,\,\,\,\,\,\,
\,\,\,\,\,\,\,\,\,\,\,\,\,\,\,\,\,\,\,\,\,\,\,\,\,\,\,\,\,\,\,\,\,\,\,\,\,\,\,\,\,\,\,\,\,\,\,\,\,\,\,\,\,\,$$ (Remark that $K-$dominated splitting implies  (d)
for any  $\zeta\leq \frac1{2K}\log 2$,  and  (d) implies $kK-$dominated for large $k$),
then we can define a new subset of $\Lambda_k (K,\zeta),$
denoted by $\Lambda^{\#}_k=\Lambda^{\#}_k (K, \zeta)$.
Similarly, one can define $\Lambda^{\#}=\Lambda^{\#} (K, \zeta)$ the maximal $f-$invariant
subset of $\bigcup_{k\geq 1}\Lambda^{\#}_k$, meaning
$$\Lambda^{\#}=\bigcap_{n\in\mathbb{Z}}f^n (\bigcup_{k\geq 1}\Lambda^{\#}_k).$$

Note that $ \Lambda^{\#}_{k} (K,\zeta)\subseteq{\Lambda}_{k} (K,\zeta)$ and $\Lambda^{\#} (K,\zeta)\subseteq{\Lambda} (K,\zeta).$ Moreover, all  $\Lambda^{\#}_{k} (K,\zeta)$ $ (k\geq 1)$  also form an increasing sequence of closed subsets of $M$. Later we will establish stable manifold theorem  and exponential shadowing on $\Lambda^{\#}_{k}$.





\bigskip

Now let us define a set which bases on the Liao-Pesin block and the points nearby Liao-Pesin blocks.
Let $\sigma>0.$ Define $\Lambda^*_{k} (K,\zeta,\sigma)$ as the union set of $\Lambda^{\#}_{k} (K,\zeta)$ and $$\{z\in M\,\,|\,\,\exists \text{ pseudo-orbit } \{x_i,\,n_i\}_{i=-\infty}^{+\infty}\text{
satisfying }n_i\geq 2kK \text{ and } x_i,f^{n_i} (x_i)\in $$
$$
\Lambda^{\#}_{k} (K,\zeta)\text{ for all }i \text{ such that } z \text{ is a }
\sigma-\text{shadowing point } \text{ for }
\{x_i,n_i\}_{i=-\infty}^{+\infty} \}.$$

\medskip

\subsubsection {Dominated Splitting and  (Center) Stable Manifold}\label{DomiStab}

Let $I_1= (-1,1)$ be the open unit ball and
$I_{\epsilon}=\{x\in\mathbb{R}: |x|<\epsilon\}$ the open
$\epsilon-$ball. Denote by $Emb (I_1,M)$ the set of $C^1$
embeddings of $I_1$ in $M$, equipped with the  (uniform) $C^1$
topology. 
 The following lemma is taken from  \cite{HPS} about the
existence of center stable and unstable manifolds for a dominated
splitting.

\begin{Lem}\label{Lem:6} ( \cite{HPS})Let $\Delta$ be a closed
 $f-$invariant set $\&$ there is a dominated splitting $T_{\Delta}M=E\oplus F$
over $\Delta.$  There exist two continuous functions
$\theta^{cs}: \Delta\rightarrow Emb (I_1, M)$ and  $\theta^{cu}:
\Delta\rightarrow Emb (I_1, M)$ such that the following properties
hold:

 (a) $\theta^{cs} (x) (0)=x$ and  $\theta^{cu} (x) (0)=x$.

 (b) $T_xW_{\epsilon}^{cs} (x)=E (x)$ and
$T_xW_{\epsilon}^{cu} (x)=F (x)$, where
$W_{\epsilon}^{cs} (x)=\theta^{cs} (x)I_{\epsilon}$ and
$W_{\epsilon}^{cu} (x)=\theta^{cu} (x)I_{\epsilon}$.

 (c) For $0<\epsilon_1<1$ there exists $0<\epsilon_2<1$ such
that $$f (W_{\epsilon_2}^{cs} (x))\subseteq
W_{\epsilon_1}^{cs} (f (x))$$ and
$$f^{-1} (W_{\epsilon_2}^{cu} (x))\subseteq
W_{\epsilon_1}^{cu} (f^{-1} (x))$$ for all $x\in \Delta.$
\end{Lem}

The manifold
 $W^{cs}$ is called the  (local) center stable manifold and $W^{cu}$ the  (local) center unstable manifold.
They are only locally invariant. Generally, this kind of invariant
manifold is not unique. So when we use center invariant manifolds
for analysis, we will  (henceforth) fix a family in discs in Lemma
\ref{Lem:6}.
The following corollary is a direct consequence of the continuity of $W^{cs/cu} (x).$

\begin{Cor}\label{Cor:ProductStru} Let $\Delta$ be a closed
 $f-$invariant set $\&$ there is a dominated splitting $T_{\Delta}M=E\oplus F$
over $\Delta.$  Then there exists $\epsilon_0>0$ such that for any $0<\epsilon\leq\epsilon_0,$ there exists $\delta=\delta (\epsilon)$ such that for any $x,y\in\Delta,$ if $   d (x,y)<\delta,$ then $W_\epsilon^{cs} (x)$ and $W_\epsilon^{cu} (y)$ are transversal and have a unique intersection point.
\end{Cor}

Furthermore,   we also have such a corollary directly from condition  (c) in Lemma \ref{Lem:6} and the continuity of $W^{cs/cu} (x)$.

\begin{Cor}\label{Cor:MnfdInduction} Let $\Delta$ be a closed
 $f-$invariant set $\&$ there is a dominated splitting $T_{\Delta}M=E\oplus F$
over $\Delta.$  Then there exists $\epsilon_0>0$ such that for any $0<\epsilon\leq\epsilon_0,$ there exists $\delta=\delta (\epsilon)$ such that for any $x\in\Delta,$ \\
a) if $y\in W_\epsilon^{cs} (x)$ and $   d (f^{-1} (x),f^{-1} (y))<\delta$ then $f^{-1} (y)\in W_\epsilon^{cu} (f^{-1} (x));$\\
b) if $y\in W_\epsilon^{cu} (x)$ and $   d (f (x),f (y))<\delta$  then $f (y)\in W_\epsilon^{cu} (f (x)).$
\end{Cor}

Now let us consider whenever the center stable and unstable manifold
are the truth ones.

\begin{Prop}\label{Prop:StableManfld}Let $\Delta$ be a closed
 $f-$invariant set $\&$ there is a dominated splitting $T_{\Delta}M=E\oplus F$
over $\Delta.$   For $K\in \mathbb{N},$ $C>0$ and
$0<\lambda <\bar{\lambda}<1$, there exist
$\epsilon_0=\epsilon (K,C,\lambda,\bar{\lambda})>0$ and $\bar{C}=\bar{C} (K,C,\bar{\lambda})\geq 1$ such that if for
$x\in\Delta$, there is a two-sided sequence of integers $$\cdots<t_{-2}<t_{-1}<t_0=0<t_1<t_2<\cdots$$ with $t_i-t_{i-1}\leq K$ for all $i$ satisfying $$\prod_{j=0}^{i-1}\|Df^{t_{j+1}-t_{j}}|_{E ({f^{t_{j}} (x)})}\|
   \leq C \lambda^{t_i}  $$ and $$\prod_{j=0}^{-i+1}\|Df^{t_{j-1}-t_{j}}|_{F ({f^{t_{j}} (x)})}\|
   \leq C \lambda^{-t_{-i}}  $$$$(\text{in other words,} \prod_{j=-i}^{-1}m (Df^{t_{j+1}-t_{j}}|_{F ({f^{t_{j}} (x)})})
   \geq C \lambda^{t_{-i}}) $$ for all $i \in \mathbb{N},$ then for any $0<\epsilon\leq\epsilon_0$
$$diam (f^nW^{cs}_{\epsilon} (x))\rightarrow_n 0\,\,\,\,\,\,\,\,\,\text{and}\,\,\,\,\,\,\,\,\,
diam (f^{-n}W^{cu}_{\epsilon} (x))\rightarrow_n 0,$$ i.e., the center
stable manifold of size $\epsilon$ is in fact a stable manifold
and the center unstable manifold of size $\epsilon$ is in fact an
unstable manifold. Moreover,  for all $n\geq1,$ $$\frac {diam (f^nW^{cs}_{\epsilon} (x))}{diam (W^{cs}_{\epsilon} (x))}\leq \bar{C} \bar{\lambda}^n$$ and $$\frac {diam (f^{-n}W^{cu}_{\epsilon} (x))}{diam (W^{cu}_{\epsilon} (x))}\leq \bar{C} \bar{\lambda}^n.$$
\end{Prop}

{\bf Proof of Proposition \ref{Prop:StableManfld}} Let
$\Delta$ be a compact invariant set and
$T_{\Delta}M=E\oplus F$ be a $Df-$invariant dominated splitting over
$\Delta.$ By Lemma \ref{Lem:6} the center stable and unstable
manifolds exist. Here  fix a small constant $\tau\in (1,\frac{\bar{\lambda}}{\lambda}).$ One chooses
a small constant $\delta>0$ such that for any
$x\in\Delta,$ any point $y\in W_1^{cs} (x)$ with $   d (x,y)<\delta$ and
any vector $v\in T_yW_1^{cs} (x)$ and $1\leq j\leq L$ we have $$\|D_yf^j\cdot v\|\leq
\tau^j\,\|Df^j|_{E (x)}\|\cdot\|v\|.$$

Let $C_f>1$ be a bound
on the norm of the derivative $Df$. Define  $$\bar{C}=\max\{C_f^K\cdot C\cdot\bar{\lambda}^{-K},\,1\} .$$ Let $\epsilon_0>0$ small enough such that $$diam (W^{cs/cu}_{\epsilon_0} (y)<\frac{\delta}{\bar{C}}$$ for all $y\in \Delta.$
 Let $0<\epsilon\leq\epsilon_0.$ Write
$n=t_i+q$ where $0\leq q<t_{i+1}-t_i\leq K$.Since
$$\prod_{j=0}^{i-1}\|Df^{t_{j+1}-t_{j}}|_{E ({f^{t_{j}} (x)})}\|
   \leq C \lambda^{t_i},  $$ then one deduces inductively on $n$,
   \begin{eqnarray*}& &diam (f^n W^{cs}_{\epsilon} (x))\nonumber\\
   &\leq& C_f^q\cdot \Pi_{j=0}^{i-1} \tau^{ (t_{j+1}-t_j)}\,\|Df^{t_{j+1}-t_{j}}|_{E ({f^{t_{j}} (x)})}\|\cdot diam (W^{us}_{\epsilon} (x))\nonumber\\
&\leq& C_f^K\cdot  \tau^{t_i}\cdot C \lambda^{t_i}\cdot diam (W^{us}_{\epsilon} (x))\nonumber\\
&\leq&  \bar{C} \bar{\lambda}^n\cdot diam (W^{cs}_{\epsilon} (x)).\nonumber
\end{eqnarray*} Note that
$$diam (f^n W^{cs}_{\epsilon} (x))\leq  \bar{C} \bar{\lambda}^n\cdot \frac{\delta}{\bar{C}}=\bar{\lambda}^n\cdot\delta\leq \delta,$$ so that the
induction can go on and on. Moreover,
$$diam (f^n W^{cs}_{\epsilon} (x))\leq \bar{\lambda}^n\cdot\delta,$$ which goes to 0 as
$n\rightarrow+\infty.$

Similarly, one has get the result for the unstable manifold. \hfill$\Box$
\bigskip

\begin{Rem}\label{Rem:StaMnfd}  Similar discussions  in this subsection for $C^1$ \textit{surface} diffeomorphisms also appeared in  \cite{PujalsSam, Gan2}.
\end{Rem}

\subsubsection{\bf Uniform size of Stable Manifolds on Liao-Pesin blocks}

Note that we do not know whether improved  Liao-Pesin set $\Lambda^{\#}$ has dominated splitting. However,
From condition  (d) in the definition of Liao-Pesin block $\Lambda^{\#}_k$, the two bundles are dominated on invariant set $$\bigcup_{l\in\mathbb{Z}}f^l (\Lambda^{\#}_k)$$  and then on its closure, since dominated property can  always be extended on the closure even neighborhoods.  Thus, one can get the following corollary from Proposition \ref{Prop:StableManfld}, which shows that for every Pesin
block $\Lambda^{\#}_k$, the  (un)stable manifolds of all $x\in \Lambda^{\#}_k$
exist and have uniform length.

\begin{Prop}\label{prop:PesinStableManfld}  (Stable Manifold on $\Lambda^{\#}_{k}$) \\ For $K\in \mathbb{N},$
$\zeta>\bar{\zeta}>0$ and $k\geq 1$, there exist
$\epsilon_0=\epsilon (K,\zeta,\bar{\zeta}, k)>0$ and $\bar{C}=\bar{C} (K,k,\zeta,\bar{\zeta})\geq 1$  such that for any $0<\epsilon\leq\epsilon_0$ and for all
$x\in\Lambda^{\#}_{k}=\Lambda^{\#}_k (K,\zeta)$,
diam ($f^nW^{cs}_{\epsilon} (x))\rightarrow_n 0$ and
diam ($f^{-n}W^{cu}_{\epsilon} (x))\rightarrow_n 0$, i.e., the center
stable manifold of size $\epsilon$ is in fact a stable manifold
and the center unstable manifold of size $\epsilon$ is in fact an
unstable manifold. Moreover,  for all $n\geq1,$ $$\frac {diam (f^nW^{cs}_{\epsilon} (x))}{diam (W^{cs}_{\epsilon} (x))}\leq \bar{C} e^{-\bar{\zeta}n}$$ and $$\frac {diam (f^{-n}W^{cu}_{\epsilon} (x))}{diam (W^{cu}_{\epsilon} (x))}\leq \bar{C} e^{-\bar{\zeta}n}.$$
\end{Prop}

{\bf Proof of Proposition \ref{prop:PesinStableManfld}} Let
$\Delta:=\overline{\cup_{l\in\mathbb{Z}}f^l (\Lambda^{\#}_k)}.$ Clearly  $\Delta$ is a compact invariant set and
$T_{\Delta}M=E\oplus F$ be a $Df-$invariant dominated splitting over
$\Delta,$ due to condition  (d) in the definition of improved Liao-Pesin block $\Lambda^{\#}_k$.

Let $\lambda:=\lambda (\zeta)=e^{-\zeta}$, $\,\bar{\lambda}:=\bar{\lambda} (\bar{\zeta})=e^{-\bar{\zeta}}$ and define  $C:=C (\zeta,k,K)=\{{e^{\zeta}}\cdot C_f\}^{kK}$ where $C_f>1$ be a bound
on the norm of the derivative $Df$.   Then
 for $x\in \Lambda^{\#}_{k}$, by the first and second condition in the definition of
$\Lambda^{\#}_{k}$ one has
$$\Pi_{i=0}^{l-1}\|Df^K|_{E (f^iKx)}\|\leq C \lambda^{-lK}$$ and $$\Pi_{i=0}^{l-1}\|Df^{-K}|_{E (f^{-iK}x)}\|\leq C \lambda^{-lK}$$  for all
$l\geq 1.$ Clearly for all $x\in \Lambda^{\#}_{k},$ the sequence in Proposition \ref{Prop:StableManfld} can always be chosen $\{t_i=iK\}_{i=-\infty}^{+\infty}.$ So using Proposition \ref{Prop:StableManfld}, the stable and unstable manifolds on $\Lambda^{\#}_k$ can be chosen of uniform size only dependent on  $K,C,\lambda,\bar{\lambda}$ and thus only dependent on $K,\zeta,\bar{\zeta}, k.$ So we can complete the proof. \hfill$\Box$

\bigskip

Moreover, we can get uniform size for stable and unstable manifolds for all points in $\Lambda^*_{k} (K,\zeta,\sigma)$ if $\sigma$ is small enough.

\begin{Prop}\label{Prop:GeneralPesinStableManfld} (Stable Manifold on $\Lambda^{*}_{k}$)\\
  For $K\in \mathbb{N},$
$\zeta>\bar{\zeta}>0$ and $k\geq 1$, there exist
$\sigma=\sigma (K,\zeta,\bar{\zeta}, k)>0, \epsilon_0=\epsilon (K,\zeta,\bar{\zeta}, k)>0$ and $\bar{C}=\bar{C} (K,k,\zeta,\bar{\zeta})>0$ such that \\
 (1) For any $0<\epsilon\leq\epsilon_0$ and for all
$x\in\Lambda^*_{k} (K,\zeta,\sigma)$,
diam ($f^nW^{cs}_{\epsilon} (x))\rightarrow_n 0$ and
diam ($f^{-n}W^{cu}_{\epsilon} (x))\rightarrow_n 0$, i.e., the center
stable manifold of size $\epsilon$ is in fact a stable manifold
and the center unstable manifold of size $\epsilon$ is in fact an
unstable manifold. Moreover,  for all $n\geq1,$ $$\frac {diam (f^nW^{cs}_{\epsilon} (x))}{diam (W^{cs}_{\epsilon} (x))}\leq \bar{C} e^{-\bar{\zeta}n}$$ and $$\frac {diam (f^{-n}W^{cu}_{\epsilon} (x))}{diam (W^{cu}_{\epsilon} (x))}\leq \bar{C} e^{-\bar{\zeta}n}.$$\\
 (2) For any $0<\epsilon\leq\epsilon_0,$ there exists $\delta=\delta (\epsilon)$ such that for any $x,y\in\Lambda^*_{k} (K,\zeta,\sigma)$ if $   d (x,y)<\delta,$ then  $W_{\epsilon}^{cs} (x)$ and $W_{\epsilon}^{cu} (y)$  are transversal and have a unique intersection point, and thus the stable manifold and unstable manifold are transversal and have a intersection point.
\end{Prop}

{\bf Proof}  Let $\epsilon_0^1=\epsilon_0 (\zeta,\bar{\zeta},k,K)$ and $\bar{C}_1=\bar{C} (K,k,\zeta,\bar{\zeta})>0$ be the $\epsilon_0$ and $\bar{C}$ as in Proposition \ref{prop:PesinStableManfld} and thus the center stable and center unstable manifolds of all $x\in\Lambda^{\#}_{k} (K,\zeta)$ are the truth stable and unstable manifolds,  and  their sizes are  at least $\epsilon_0^1$.

 Let
$\Delta_1:=\overline{\cup_{l\in\mathbb{Z}}f^l (\Lambda^{\#}_k)}.$ Clearly  $\Delta_1$ is a compact invariant set and
$T_{\Delta_1}M=E\oplus F$ be a $Df-$invariant dominated splitting over
$\Delta_1,$ due to condition  (c) in the definition of Pesin block $\Lambda^{\#}_k$. Since dominated property can be extended to its neighborhoods  (see  \cite{BLV}), we can take an open neighborhood $U=B_{\Delta_1} (\tau)=\{x\,\,|\,   d (x,\Delta_1)<\tau\}$  (for some small $\tau>0$) of $\Delta_1$ such that the splitting is extended and dominated on $\Delta:=\cap_{n\in\mathbb{Z}}f^n (U)$. Take and fix $\zeta>\hat{\zeta}>\bar{\zeta}$. Let $\lambda:=\lambda (\hat{\zeta})=e^{-\hat{\zeta}}$, $\,\bar{\lambda}:=\bar{\lambda} (\bar{\zeta})=e^{-\bar{\zeta}}$ and define  $C:=C (\zeta,k,K)=\{{e^{\hat{\zeta}}}\cdot C_f\}^{kK}$ where $C_f>1$ be a bound
on the norm of the derivative $Df$.
 Let $\epsilon_0^2=\epsilon_0 (\lambda,\bar{\lambda},C, K)$ and $\bar{C}_2=\bar{C} (\bar{\lambda},C, K)>0$ be the $\epsilon_0$ and $\bar{C}$ as in Proposition \ref{Prop:StableManfld} for this $\Delta$.

  Note that dominated splitting is always continuous  (see  \cite{BLV}). So we can take $\gamma>0$ small enough such that if for any $x,y\in\Delta,$ if $   d (x,y)<\gamma,$ then for all $1\leq i\leq K,$ \begin{eqnarray}\label{Estimate-Hyp} e^{i (\hat{\zeta}-\zeta)}\leq\frac{\|Df^i|_E (x)\|}{\|Df^i|_E (y)\|},\,\,\frac{\|Df^i|_F (x)\|}{\|Df^i|_F (y)\|}\leq e^{i (\zeta-\hat{\zeta})}.
\end{eqnarray}
Take $\sigma=\min\{\tau,\gamma\}.$  If $z\in M$  and  there exists pseudo-orbit $\{x_i,\,n_i\}_{i=-\infty}^{+\infty}$
satisfying $n_i\geq 2kK$ and $x_i,f^{n_i} (x_i)\in
\Lambda^{\#}_{k} (K,\zeta) (\forall \,i)$ such that
$ z$  is a $\sigma$-shadowing point  for $
\{x_i,n_i\}_{i=-\infty}^{+\infty}.$  $   d\big{ (}f^{c_i+j} (z),f^j (x_i)\big{)}\leq\sigma\leq\tau$ implies that the  orbit of $z$ is contained in the $\tau-$neighborhood $U$ of $\Delta_1$, and thus $z\in\Delta.$  Since $x_i,f^{n_i} (x_i)\in\Lambda^{\#}_k (K,\,\zeta)$ and $   d\big{ (}f^{c_i+j} (z),f^j (x_i)\big{)}\leq\sigma\leq\gamma,$ it is easy to verify, using inequality  (\ref{Estimate-Hyp}),  that $z$ satisfies the condition as in Proposition \ref{Prop:StableManfld} for the above $\lambda$ and $C$.  So by the choice of $\epsilon_0^2$, the center stable and center unstable manifolds of all such points $z$ above are the truth stable and unstable manifolds,  and  their sizes are  at least $\epsilon_0^2.$

Using Corollary \ref{Cor:ProductStru} for the above set $\Delta$, there exists $0<\epsilon_0\leq\min\{\epsilon_0^1,\epsilon_0^2\}$ such that for any $0<\epsilon\leq\epsilon_0,$ there exists $\delta=\delta (\epsilon)$ such that for any $x,y\in\Delta,$ if $   d (x,y)<\delta,$ then $W_\epsilon^{cs} (x)$ and $W_\epsilon^{cu} (y)$ are transversal and have a unique intersection point. Clearly, $\bar{C}$ can be chosen $\max\{\bar{C_1}, \bar{C_2}\}.$

Since $\epsilon_0\leq\min\{\epsilon_0^1,\epsilon_0^2\}$, then for any $0<\epsilon\leq\epsilon_0,$ the center stable and center unstable manifolds of all $x\in\Lambda^*_{k} (K,\zeta,\sigma)$ are truth stable and unstable manifold, and their sizes are at least ${\epsilon}$.
 Note that $\Lambda^*_{k} (K,\zeta,\sigma) \subseteq\Delta.$ Thus  for any $x,y\in \Lambda^*_{k} (K,\zeta,\sigma) \subseteq\Delta,$ if $   d (x,y)<\delta,$ then the stable manifold $W_{\epsilon}^{cs} (x)$ and unstable manifold $W_{\epsilon}^{cu} (y)$  are transversal and have a unique intersection point. \qed

\begin{Rem}\label{Rem:ABCStableMnfd}The discussion of stable manifold for hyperbolic ergodic measures whose Oseledec splitting is dominated also
 appeared in  \cite{ABC} for the existence of stable manifold for  a.e. points. Here in our paper we
  give a clear and definite filtration of blocks, called Liao-Pesin blocks, such that every block  has stable manifold theorem and simultaneously has  (exponential) shadowing lemma. Here the definition is from topological viewpoint and independent of  invariant  measures.

\end{Rem}

Furthermore, we point out that from the above proof, if $\sigma$ is small enough, one can get a  strong relation of $\Lambda^*_{k}$ and $\Lambda^{\#}_{k}$ as follows.

\begin{Prop}\label{Prop:ShadowingPtIn}For $K\in \mathbb{N},$
$\zeta>\bar{\zeta}>0$ and $k\geq 1$, there  is a small enough number $\sigma>0$ and an integer $\bar{k}$ such that
$\Lambda^*_{k} (K,\zeta,\sigma)\subseteq\Lambda^{\#}_{\bar{k}} (K,\bar{\zeta}).$
\end{Prop}
{\bf Remark.} Note that the choice of $\bar{k}$ is only dependent on $k,K,\zeta,\bar{\zeta}$, though in particular  $\Lambda^*_{k} (K,\zeta,\sigma)$ contains  periodic points with arbitrarily large period. In other words, the choice of $\bar{k}$ is independent of the period of all periodic points in $\Lambda^*_{k} (K,\zeta,\sigma)$.
\medskip

{\bf Proof.} Let us use the same notions as in Proposition \ref{Prop:GeneralPesinStableManfld} and we only prove the first condition of Pesin block. From the proof of Proposition \ref{Prop:GeneralPesinStableManfld}, if $\sigma$ is small enough, using
inequality  (\ref{Estimate-Hyp}) and the definition of $\Lambda^*_{k} (K,\zeta,\sigma)$ we have
$$\frac{\log\|Df^r|_{E (f^{c_i} (z))}\|+\sum_{j=0}^{l-1}\log\|Df^{K}|_{E ({f^{jK+r} (f^{c_i} (z))})}\|}{lK+r}
   \leq-\hat{\zeta},$$ $\forall\,\,\,
k \leq l\,\leq [\frac{n_i}{K}] \, ,\,\,\,\forall\,\,\,0\,\leq\, r\,\leq\,
K-1.$
 This implies that if $c_i+kK\leq lK+r\leq c_{i+1}$ then $$\frac{\log\|Df^r|_{E (z))}\|+\sum_{j=0}^{l-1}\log\|Df^{K}|_{E ({f^{jK+r} (z))})}\|}{lK+r}
   \leq-\hat{\zeta}\leq-\bar{\zeta}.$$
To realize this inequality for $\forall\,\,\,
 l\,\geq \bar{k} \, ,\,\,\,\forall\,\,\,0\,\leq\, r\,\leq\,
K-1$ where $ \bar{k}$ is a constant,  we only need to consider $l$ satisfying $c_i< lK+r<c_i+kK.$ It is easy to see that
$$\frac{\log\|Df^r|_{E (z))}\|+\sum_{j=0}^{l-1}\log\|Df^{K}|_{E ({f^{jK+r} (z))})}\|}{lK+r}
   \leq\frac{- ([\frac{c_i-r}{K}]K+r)\hat{\zeta}+ (k+1)K C_f}{lK+r}$$
$$   \leq\frac{- ( (l-k-1)K+r)\hat{\zeta}+ (k+1)K C_f}{lK+r}.$$
Take $$\bar{k}=\max\{k,\,\,\frac{ (k+1)K (\hat{\zeta}+C_f)+K\bar{\zeta}}{\hat{\zeta}-\bar{\zeta}}\}.$$ Then for  $\forall\,\,\,
 l\,\geq \bar{k} \, ,\,\,\,\forall\,\,\,0\,\leq\, r\,\leq\,
K-1,$
$$\frac{\log\|Df^r|_{E (z))}\|+\sum_{j=0}^{l-1}\log\|Df^{K}|_{E ({f^{jK+r} (z))})}\|}{lK+r}
  \leq-\bar{\zeta}.$$\qed




\subsubsection{Exponential shadowing and closing lemma}

  Given $\theta>0$ and $\eta>0,$ we call a point $x\,\in
M$  an exponential $ (\eta,\theta)$-shadowing point for a pseudo-orbit
$\big{\{}x_i,\,n_i\big{\}}_{i=-\infty}^{+\infty},$  if
$$   d\big{ (}f^{c_i+j} (x),f^j (x_i)\big{)}<\eta\cdot e^{-\min\{j,n_i-j\}\theta},$$ $\forall\,\,
j=0,\,1,\,2,\,\cdots,\,n_i-1$ and $\forall\,\, i \in \mathbb Z$, where
$c_i$
 is defined as
  \begin {equation} \label{eq:exponential-respective-time-squens}c_i=\begin{cases}
 0,&\text{for }i=0\\
 \sum_{j=0}^{i-1}n_j,&\text{for }i>0\\
 -\sum_{j=i}^{-1}n_j,&\text{for }i<0.
\end{cases}
 \end {equation}


\medskip

We firstly state exponential shadowing for $C^{1+\alpha}$ systems which can be as a particular consequence of  Lemma \ref{LemShadow}.

\begin{Prop}\label{Prop-Shadow-expon}
 (Exponential shadowing lemma in $C^{1+\alpha}$ case)\\
 Let $f:M\rightarrow M$ be a
$C^{1+\alpha}$ diffeomorphism, with a non-empty Pesin block
$ {\Lambda}_k= {\Lambda}_k (\lambda,\mu;\varepsilon)$ and fixed parameters,
$\lambda,\mu\gg\varepsilon>0, k\geq 1$. For
$\forall\,\,\eta>0,$ there exists $\delta=\delta (k,\eta)>0$ such that if  a
$\delta$-pseudo-orbit $\{x_i,\,n_i\}_{i=-\infty}^{+\infty}$
satisfies $n_i\geq 1 $ and $x_i,f^{n_i} (x_i)\in
 {\Lambda}_k$ for all $i$, then there exists a unique exponential
$ (\eta,\varepsilon)$-shadowing point $x\in M $ for
$\{x_i,n_i\}_{i=-\infty}^{+\infty}$.\\
 If further
$\{x_i,n_i\}_{i=-\infty}^{+\infty}$ is periodic, i.e., there
 exists an integer $m>0$ such that $x_{i+m}=x_i$ and $n_{i+m}=n_i$ for all i,
 then the shadowing point $x$ can be chosen to be periodic.
\end{Prop}

{\bf Proof.} For given $\eta>0,$ take $\tau=\frac{\eta}{\varepsilon_0}$ in Lemma \ref{LemShadow}. Then there exists a
sequence $ (\delta_{k})_{k=1}^{+ \infty}$ such that for any
$ (\delta_{k})_{k=1}^{+ \infty}$ pseudo-orbit there exists a unique
$\tau$-shadowing
point. Here we take $\delta=\delta_{k+1}$ and consider  a
$\delta$-pseudo-orbit $\{x_i,\,n_i\}_{i=-\infty}^{+\infty}$
satisfies $n_i\geq 1 $ and $x_i,f^{n_i} (x_i)\in
\Lambda_k$ for all $i$.

Recall the property of Pesin blocks that $\Lambda_k, f^{\pm 1} (\Lambda_k)\subseteq \Lambda_{k+1}.$ Thus, if $u\in\Lambda_k$, then $f^{j} (u)\in \Lambda_{k+|j|},\,\forall\,i\in\mathbb{Z}.$ Note that $$x_i, f^{ n_i} (x_i)\in \Lambda_k$$
imply  $$\cdots ,  x_i\in  \Lambda_{k+1}, f (x_i)\in  \Lambda_{k+2},\cdots,   f^{ j} (x_i)\in  \Lambda_{\min\{k+j+1,k+n_i-j\}},$$$$ \cdots, f^{ n_i-1} (x_i)\in  \Lambda_{k+1}, x_{i+1}\in \Lambda_k,\cdots $$
and $d (f^{n_i} (x_i), x_{i+1})<\delta$  implies $$d (f^{n_i} (x_i), x_{i+1})<\delta_{k+1}$$ so that the above points form a
$ (\delta_{k})_{k=1}^{+ \infty}$ pseudo-orbit. Then there exists a unique
$\tau$-shadowing point $z$. It is easy to check  that $z$ is the needed  exponential $ (\eta,\varepsilon)$-shadowing point for
$\{x_i,n_i\}_{i=-\infty}^{+\infty}$. More precisely,  $$   d\big{ (}f^{c_i+j} (z),f^j (x_i)\big{)}<\tau
 \varepsilon_k=\tau\varepsilon_0e^{-\min\{k+j,k+n_i-j\}\varepsilon}<\eta\cdot e^{-\min\{j,n_i-j\}\varepsilon},$$ $\forall\,\,
j=0,\,1,\,2,\,\cdots,\,n_i-1$ and $\forall\,\, i \in \mathbb Z$.  \qed

\bigskip

Secondly, we state exponential shadowing in $C^1$ setting.
That is,  we show exponential shadowing property on Liao-Pesin blocks.

\begin{Prop}\label{Prop:ExpShadowinglem} (Exponential shadowing property on Liao-Pesin blocks) \\ Assume $\Lambda^{\#}_k (K,\,\zeta)\neq\emptyset$ for
some $k,\,K\in\mathbb{N}$ and $\zeta>0$. Then $\Lambda^{\#}_k (K,\,\zeta)$
satisfies exponential shadowing property  as follows: \\ there is $
\theta>0$ and  $T_k=T (k,K,\zeta)>0$ such that for
$\forall\,\,\eta>0,$ there exists $\delta=\delta (k, K, \zeta, \eta)>0$ such that if  a
$\delta$-pseudo-orbit $\{x_i,\,n_i\}_{i=-\infty}^{+\infty}$
satisfies $n_i\geq T_k $ and $x_i,f^{n_i} (x_i)\in
\Lambda^{\#}_k$ for all $i$, then there exists a  unique exponential
$ (\eta,\theta)$-shadowing point $x\in M $ for
$\{x_i,n_i\}_{i=-\infty}^{+\infty}$.

 If further
$\{x_i,n_i\}_{i=-\infty}^{+\infty}$ is periodic, i.e., there
 exists an integer $m>0$ such that $x_{i+m}=x_i$ and $n_{i+m}=n_i$ for all i,
 then the shadowing point $x$ can be chosen to be periodic. Moreover, the periodic orbit should be hyperbolic.
\end{Prop}

\begin{Rem}\label{Rem:shadowinglem}  In fact, the shadowing can be stated stronger to be  Lipschitz shadowing.  That is, there exists $L_k>0, \delta_k>0$ such that for any $0<\delta<\delta_k$ if a
$\delta$-pseudo-orbit $\{x_i,\,n_i\}_{i=-\infty}^{+\infty}$
satisfies $n_i\geq2kK $ and $x_i,f^{n_i} (x_i)\in
\Lambda^{\#}_k (K,\,\zeta)$ for all $i$,  then there exists
a
$L_k\delta$-shadowing point $x\in M $ for
$\{x_i,n_i\}_{i=-\infty}^{+\infty}$.
The main observation is that the used technique is Liao's shadowing lemma \cite{Liao79,Gan}
 for quasi-hyperbolic orbit segments and Liao's shadowing lemma can be Lipschitz.
\end{Rem}

\medskip

Using Theorem \ref{Thm:shadowinglem} and Proposition \ref{Prop:GeneralPesinStableManfld}, we start to prove exponential shadowing property.

{\bf Proof of Proposition \ref{Prop:ExpShadowinglem}} Let
$\bar{\zeta},$
$\sigma=\sigma (K,\zeta,\bar{\zeta}, k)>0, \epsilon_0=\epsilon (K,\zeta,\bar{\zeta}, k)>0$ and $\bar{C}=\bar{C} (K,k,\zeta,\bar{\zeta})>0$ be the numbers as in Proposition \ref{Prop:GeneralPesinStableManfld}.
Moreover, let $\Delta$ be the invariant set in the proof of Proposition \ref{Prop:GeneralPesinStableManfld}, which contains $\Lambda^*_{k} (K,\zeta,\sigma)$.

For any fixed $\eta\in (0,\sigma),$ take $\epsilon=\epsilon (\eta)>0$ small enough such that $$diam (W^{cs/cu}_{\epsilon} (y)<\frac{\eta}{2\bar{C}}$$ for all $y\in\Delta.$ By  Corollary \ref{Cor:MnfdInduction}, for this $\epsilon$ we can take $ b_1= b_1 (\epsilon)\in (0,\eta)$ such that for $x \in \Delta,\,$
 if  $y\in W_\epsilon^{cu} (x)$ and $   d (f (x),f (y))< b_1$ then $f (y)\in W_\epsilon^{cu} (f (x)).$

  For $ b_1$, take $\epsilon_*=\epsilon ( b_1)>0$ small enough such that $$diam (W^{cs/cu}_{\epsilon_*} (y)<\frac{b_1}{2\bar{C}}$$ for
  all $y\in\Delta.$ Moreover, by Proposition \ref{Prop:GeneralPesinStableManfld}, we can take a positive number $ b_2= b_2 (\epsilon_*)<\min\{\frac{ b_1}2,\sigma\}$ such that
for any $x,y\in\Lambda^*_{k} (K,\zeta,\sigma)$ if $   d (x,y)< b_2,$ then  $W_{\epsilon_*}^{cs} (x)$ and $W_{\epsilon_*}^{cu} (y)$  are transversal and have a unique intersection point.

 By Theorem \ref{Thm:shadowinglem}, for $ b_2$ there exists
$\delta=\delta ( b_2)>0$ such that if  a
$\delta$-pseudo-orbit $\{x_i,\,n_i\}_{i=-\infty}^{+\infty}$
satisfies $n_i\geq2kK $ and $x_i,f^{n_i} (x_i)\in
\Lambda^{\#}_k (K,\,\zeta)$ for all $i$, then there exists a
$ b_2$-shadowing point $x\in M $ for
$\{x_i,n_i\}_{i=-\infty}^{+\infty}$. If further
$\{x_i,n_i\}_{i=-\infty}^{+\infty}$ is periodic, i.e., there
 exists an $m>0$ such that $x_{i+m}=x_i$ and $n_{i+m}=n_i$ for all i,
 then the shadowing point $x$ can be chosen to be periodic. Since $b_2<\sigma$, the periodic point $x$ should be in $\Lambda^*_k.$
  By Proposition \ref{Prop:ShadowingPtIn}, $x$ is in  $\Lambda^{\#}_{\bar{k}}$  for some integer $\bar{k}.$ So the periodic point $x$ should be hyperbolic.

 Now let us prove that  the shadowing is  exponential. We only need to prove for orbit segment $\{x_0,n_0\}$, $\{x,n_0\}$ is exponential shadowing $\{x_0,n_0\}$, since the others are similar. That $  b_2<\sigma$ implies $x\in\Lambda^*_{k} (K,\zeta,\sigma).$
Notice that $   d (x_0,x)<   b_2$ and $x_0\in\Lambda^{\#}_k (K,\zeta)\subseteq\Lambda^*_{k} (K,\zeta,\sigma).$ So $W_{\epsilon_*}^{cs} (x_0)$ and $W_{\epsilon_*}^{cu} (x)$  are transversal and have a unique intersection point $y\in M$. Since $x_0\in\Lambda^*_{k} (K,\zeta,\sigma)$ and $y\in W_{\epsilon_*}^{cs} (x_0),$ by Proposition \ref{Prop:GeneralPesinStableManfld} we have
$$   d (f^j (x_0),f^j (y))\leq {diam (f^nW^{cs}_{\epsilon_*} (x_0))}\leq \bar{C} e^{-\bar{\zeta}j}{diam (W^{cs}_{\epsilon_*} (x_0))}\leq \frac{b_1}2 e^{-\bar{\zeta}j},$$ $\,j=0,1,\cdots,n_0.$
Combing this inequality with $   d (f^j (x_0),f^j (x))\leq  b_2<\frac{b_1}2,$  $\,j=0,1,\cdots,n_0,$
we have $$    d (f^j (x),f^j (y))\leq   d (f^j (x_0),f^j (x))+   d (f^j (x_0),f^j (y))\leq  b_1, \,j=0,1,\cdots,n_0.$$
From the choice of $b_1$ and $f^j (x)\in f^j (\Lambda^*_{k} (K,\zeta,\sigma))\subseteq f^j (\Delta)=\Delta$, by induction it is easy to get $f^j (y)\in W_\epsilon^{cu} (f^j (x))\,j=0,1,\cdots,n_0.$ Note that $f^{n_0} (x)$ is also in $ \Lambda^*_{k} (K,\zeta,\sigma)$. Using  $f^{n_0} (y)\in W_\epsilon^{cu} (f^{n_0} (x)),$ by Proposition \ref{Prop:GeneralPesinStableManfld} one can get $$   d (f^j (x),f^j (y))\leq {diam (f^{ (n_0-j)}W^{cu}_{\epsilon} (x_0))}\leq \bar{C} e^{-\bar{\zeta}{ (n_0-j)}}{diam (W^{cu}_{\epsilon} (x_0))}\leq $$ $\frac\eta2 e^{-\bar{\zeta}{ (n_0-j)}},\,j=0,1,\cdots,n_0. $
So
$$   d (f^j (x_0),f^j (x))\leq   d (f^j (x_0),f^j (y))+   d (f^j (x),f^j (y))$$
$$\leq  (\frac\eta2+\frac{b_1}2) \max\{e^{-\bar{\zeta}j},e^{-\bar{\zeta}{ (n_0-j)}}\}\leq \eta \max\{e^{-\bar{\zeta}j},e^{-\bar{\zeta}{ (n_0-j)}}\},\,j=0,1,\cdots,n_0. $$
Take $\theta=\bar{\zeta}$ and we complete the proof. \qed
\medskip

\begin{Rem}
Recently we notice an improved version of Liao's closing lemma by Dai  \cite{Dai2010}. This kind of closing is exponential closing for one quasi-hyperbolic orbit segment so that one can use his result (Theorem 2 of \cite{Dai2010}) to replace the role of Lemma \ref{Lem:Liao-Gan-shadlem} and then exponential closing should hold on $\Lambda_k,$ not just on its subset $\Lambda^{\#}_k$. Moreover, using the idea of  \cite{Dai2010} it should be straightforward to get exponential shadowing for finite and infinite quasi-hyperbolic orbit segments so that one can also get exponential shadowing on $\Lambda_k,$ not just on  $\Lambda^{\#}_k$. This method avoids the use of stable manifold. However, for dealing with the results of present paper, our above analysis on $\Lambda^{\#}_k$ are enough   so that we just give a remark here.

\end{Rem}

\smallskip

Remark that the statement of  Proposition \ref{Prop:ExpShadowinglem} is little weaker than    Proposition \ref{Prop-Shadow-expon} because in Proposition  \ref{Prop:ExpShadowinglem}  we require that the length of  every segment  $\{x_i,n_i\}$ is at least $T_k.$  So it is not convenient to find periodic orbit with small period.  But  in general, these small details do not infect the establishment of most Pesin theory.

On the other hand, exponential shadowing is just as one special  case of  Katok's shadowing (Lemma \ref{LemShadow}). It is  because exponential shadowing is just suitable for the orbit segments whose beginning and ending points are in the fixed Liao-Pesin block.  That is, exponential shadowing is a partial answer of Question \ref{Que-LemShadow}. Thus  a natural question aries:

{\bf Question.} {\it How about the shadowing for  orbit  segments whose beginning and ending points are from different Liao-Pesin blocks?}





\section {Existence of Liao-Pesin set 
}

What systems does Liao-Pesin set exist in? In this section we
answer this question in certain  $C^1$  systems.

\subsection{Average-nonuniform hyperbolicity and (limit-)dominated splitting}

Firstly let us give a topological condition to get Pesin set.
\begin{Thm}\label{Thm-topocondition-to-pesinset-0000000000}
 Let $\Delta$ be an $f-$invariant set
and $T_{\Delta}M=E\oplus F$ be a $Df-$invariant splitting on
$\Delta$.  If  $T_{\Delta}M=E\oplus F$
is  limit-dominated on $\Delta$ and $\Delta$ is an average-nonuniformly hyperbolic set
 corresponding to $T_{\Delta}M=E\oplus F$,
  then there is some $\zeta_0>0$  and $K\geq 1$ such that  for any $\chi\in (0, \zeta_0),$    $$\Delta\subseteq \Lambda (K,\chi).$$
  If further the splitting $T_{\Delta}M=E\oplus F$  is  dominated, then one can take $\zeta_0>0$  and $K\geq 1$ such that  for any $\chi\in (0, \zeta_0),$ $$\Delta\subseteq \Lambda^\# (K,\chi).$$

\end{Thm}
Suppose that the limit-domination holds for $S\geq 1$ and $\zeta_1>0$ in the sense that
$$\limsup_{l\rightarrow+\infty}  \frac1S \log\frac
{\|Df^S|_{E (f^{l} (x))}\|}{m (Df^S|_{F (f^{l} (x))})}\leq  -2\zeta_1,\,\,\forall
x \in \Delta.$$  And suppose  $\Delta$ is an average-nonuniformly hyperbolic set with $ (S',\zeta_2)$-degree
 corresponding to $T_{\Delta}M=E\oplus F$ for some $S'\geq 1$ and $\zeta_2>0$.
Take $\zeta_0=\min\{\zeta_1,\,\zeta_2\}.$ By sub-multiplication of norms, for $K=S\cdot S'$, it should satisfy that
$$\limsup_{l\rightarrow+\infty}  \frac1K \log\frac
{\|Df^K|_{E (f^{l} (x))}\|}{m (Df^K|_{F (f^{l} (x))})}\leq  -2\zeta_1,\,\,\forall
x \in \Delta ,$$  and $\Delta$ is an average-nonuniformly hyperbolic set with $ (K,\zeta_0)$-degree
 corresponding to $T_{\Delta}M=E\oplus F$.
So we can suppose  (limit-)domination and average-nonuniform hyperbolicity for common ``degrees" and then we can state Theorem \ref{Thm-topocondition-to-pesinset-0000000000} as follows.

\begin{Thm}\label{Thm-topocondition-to-pesinset}
 Let $\Delta$ be an $f-$invariant set
and $T_{\Delta}M=E\oplus F$ be a $Df-$invariant splitting on
$\Delta$. Let $K\in \mathbb{Z}^+$ and $\zeta_0>0.$  Suppose that $\Delta$ is an average-nonuniformly hyperbolic set with
$ (K,\,\zeta_0)$-degree corresponding to $T_{\Delta}M=E\oplus F$ and
suppose that   $T_{\Delta}M=E\oplus F$  is  limit-dominated on $\Delta$ in the sense that $$\limsup_{l\rightarrow+\infty}  \frac1K \log\frac
{\|Df^K|_{E (f^{l} (x))}\|}{m (Df^K|_{F (f^{l} (x))})}\leq  -2\zeta_0,\,\,\forall
x \in \Delta.$$
Then for any $\chi\in (0, \zeta_0),$  $$\Delta\subseteq \Lambda (K,\chi).$$
  If further the splitting $T_{\Delta}M=E\oplus F$  is  dominated in the sense that
  $$\frac1K \log\frac
{\|Df^K|_{E (x)}\|}{m (Df^K|_{F ( x)})}\leq  -2\zeta_0,\,\,\forall
x \in \Delta,$$  then   for any $\chi\in (0, \zeta_0),$  $$\Delta\subseteq \Lambda^\# (K,\chi).$$

\end{Thm}

{\bf Proof.} Fix    $\chi\in (0, \zeta_0).$    Note that the assumption of domination  implies the condition  (d) in the definition of Pesin set $\Lambda^\# (K,\chi).$ So we only need to show the first part: $\Delta\subseteq \Lambda (K,\chi).$ Before that we firstly state a basic fact.

\begin{Lem}\label{Lem-Lyap-suplim-relate-average} Let $\Delta\subseteq M$ be an invariant set and $E$ be a $Df$-invariant bundle over $\Delta$. Then for any $x\in \Delta,$ any $K\geq 1$,
$$\limsup_{n\rightarrow+\infty}\frac1n\log\|Df^n|_{E (x)}\|\leq
\limsup_{l\rightarrow+\infty}\frac1{lK}\sum_{i=0}^{l-1}\log \|Df^K|_{E (f^{iK} (x))}\|.$$
  The case $\liminf$ is also true.
  Similar inverse estimates  hold for the case of minimal norm.
\end{Lem}

{\bf Proof.} Write $n=lK+r,\,0\leq r<K.$
From the sub-additional multiplication of norms, $$\limsup_{n\rightarrow+\infty}\frac1n\log
{\|Df^n|_{E (x)}\|}$$ $$
 \leq \max_{ 0\leq r<K}\limsup_{l\rightarrow
+\infty}\frac{\sum_{j=0}^{l-1}\log\|Df^{K}|_{E ({f^{
jK} (x)})}\|+\log\|Df^r|_{E (f^{lK}x)}\|}{lK+r}$$
$$= \limsup_{l\rightarrow+\infty}\frac1{lK}\sum_{i=0}^{l-1}\log \|Df^K|_{E (f^{iK } (x))}\|.$$
\qed

 Given $x\in \Delta$, by  the average-nonuniformly hyperbolic assumption and invariance of $\Delta$,  for any $r\in\mathbb{Z},$ we have
$$\limsup_{l\rightarrow
+\infty}\sum_{j=0}^{l-1}\frac{\log\|Df^{K}|_{E ({f^{ jK} (f^{
r}x)})}\|}{lK}<-\chi,\,\,\,\liminf_{l\rightarrow
+\infty}\sum_{j=-l}^{-1}\frac{\log
m (Df^{K}|_{F ({f^{jK} (x)})})}{lK}>\chi.$$ By the boundary of $Df,Df^{-1},$ these imply    that for
$\forall\,\, 0\leq r\leq K-1,$
\begin{equation}\label{eq:get-condition-a-Pesset}
\begin{split}
&\,\,\,\,\, \limsup_{l\rightarrow
+\infty}\frac{\log\|Df^r|_{E (x)}\|+\sum_{j=0}^{l-1}\log\|Df^{K}|_{E ({f^{
jK+r} (x)})}\|}{lK+r} \\&\,\,\,\,\,=\limsup_{l\rightarrow
+\infty}\sum_{j=0}^{l-1}\frac{\log\|Df^{K}|_{E ({f^{ jK} (f^{
r}x)})}\|}{lK}<-\chi,
\end{split}
\end{equation}
\begin{equation}\label{eq:get-condition-b-Pesset}
\begin{split}
&\,\,\,\,\,\,\,\,\,\,\,\, \liminf_{l\rightarrow +\infty}\frac{\log
m (Df^r|_{F (f^{- (lK+r)} (x))})+\sum_{j=-l}^{-1}\log
m (Df^{K}|_{F ({f^{jK} (x)})})}{lK+r}\\&\,\,\,\,\,=\liminf_{l\rightarrow
+\infty}\sum_{j=-l}^{-1}\frac{\log
m (Df^{K}|_{F ({f^{jK} (x)})})}{lK}>\chi.
\end{split}
\end{equation}
By Lemma \ref{Lem-Lyap-suplim-relate-average},
\begin{eqnarray}\label{eq:get-condition-c-Pesset}
\begin{split}
&\,\,\,\,\,\limsup_{n\rightarrow+\infty}\frac1n\log\frac
{\|Df^n|_{E (x)}\|}{m (Df^n|_{F (x)})}\\
 &\,\,\,\,\,\leq \limsup_{n\rightarrow+\infty}\frac1n\log
{\|Df^n|_{E (x)}\|} -\liminf_{n\rightarrow+\infty}\frac1n\log {m (Df^n|_{F (x)})}\\
&\,\,\,\,\,\leq \limsup_{l\rightarrow
+\infty}\sum_{j=0}^{l-1}\frac{\log\|Df^{K}|_{E ({f^{ jK} ( x)})}\|}{lK}-\liminf_{l\rightarrow
+\infty}\sum_{j=-l}^{-1}\frac{\log
m (Df^{K}|_{F ({f^{jK} (x)})})}{lK} \\
&\,\,\,\,\,<   -2\chi.
\end{split}
\end{eqnarray}

By  (\ref{eq:get-condition-a-Pesset}) and
 (\ref{eq:get-condition-b-Pesset}) there exists $k_1=k_1 (x)$ such
that for all $l\geq k_1,\,\,\forall\, 0\leq r\leq K-1,$ one has
$$
\frac{\log\|Df^r|_{E (x)}\|+\sum_{j=0}^{l-1}\log\|Df^{K}|_{E ({f^{
jK+r} (x)})}\|}{lK+r}\leq-\chi,
$$
$$  \frac{\log m (Df^r|_{F (f^{- (lK+r)} (x))})+\sum_{j=-l}^{-1}\log
m (Df^{K}|_{F ({f^{jK} (x)})})}{lK+r}\geq\chi.$$
 Since the splitting $E\oplus F$ is
 limit-dominated, by  assumption  we can
choose $k_2=k_2 (x)$ such that for all $l\geq k_2K,$
$$\frac1K\log\frac{\|Df^K|_{E (f^{l} (x))}\|}{m (Df^K|_{F (f^{l} (x))})}\leq-2\chi.$$
Using the inequality  (\ref{eq:get-condition-c-Pesset}) there exists $k_3=k_3 (x)$ such that for all $l\geq k_3$
and $r=0,\,1,\,\cdots,K-1,$
$$\frac 1{lK+r}\log
\frac{\|Df^{lK+r}|_{E (x)}\|}{m (Df^{lK+r}|_{F (x)})}\leq-2
\chi.$$ Take $k\geq max\{k_1,k_2,k_3\}$ and then the three
conditions in Definition~\ref{def:Pesinblock} hold. Hence $x\in
\Lambda_k (K,\,\chi)\neq \emptyset$.  Recall that the Pesin set
$\Lambda (K,\,\chi)$ is the maximal $f$-invariant set in
$\bigcup_{k\geq1}\Lambda_k (K,\,\chi)$ and $\Delta$ is
$f$-invariant, one has $\Delta\subseteq\Lambda (K,\,\chi)$. \qed



\bigskip

If $\Delta$ is compact and the splitting is dominated (or just continuous), then the average-nonuniformly hyperbolic
assumption implies that $\Delta$ is uniformly hyperbolic. In fact the assumption can be weakened as nonzero Lyapunov exponents. Moreover,  it is enough to just assume on a subset with totally full measure (full measure for any invariant measure),   not necessarily  assuming for all points \cite{Cao}.

\begin{Thm}\label{Thm-Cao2003} Let $f:M\rightarrow M$ be a $C^1$ local diffeomorphism on a compact
manifold and let $\Lambda$ be a compact and $f-$invariant set.
Suppose that there exists a continuous $Df$-invariant splitting
$T_\Lambda M=E\oplus F$. If the Lyapunov exponents restricted on $E$
and $F$ of every f invariant probability measure are all negative
and positive respectively, then $\Lambda$ is uniformly hyperbolic.
\end{Thm}

Except the compactness of $\Lambda,$ the continuity of splitting is important and there are examples  \cite{CLR} that all Lyapunov exponents are far from zero but the system is not uniformly hyperbolic, which does not admit a continuous  (or dominated) splitting. However, it is interesting to ask how about the case that $T_\Lambda M=E\oplus F$ is just limit-dominated.


By Theorem \ref{Thm-topocondition-to-pesinset-0000000000} and Proposition \ref{Prop:ExpShadowinglem} and \ref{Prop:GeneralPesinStableManfld}, we can state shadowing for any probability measures without assumption of invariance or ergodicity.

\begin{Thm}\label{Thm:ExpShadowing-measure2015oct} Let $f\in \Diff^1 (M)$ and $\mu\in\m  (M).$  Suppose that there is  a $Df$-invariant dominated splitting $T_{\Delta}M=E\oplus F$ on an $f$-invariant set $\Delta$ with $\mu$ full measure and  $\mu$ is average-nonuniformly hyperbolic with respect to
 to this splitting.
   Then for each $\tau\in (0,1),$ there exist a compact set $\Lambda_\tau\subseteq M$, $\theta_\tau>0$ and $T_\tau \in \mathbb{N}$ such that $\mu (\Lambda_\tau)>1-\tau$ and following two properties hold:\\

 (i) \emph{ (Exponentially) Shadowing Lemma}: For
$\forall\,\,\eta>0,$ there exists $\delta=\delta (\tau,\eta)>0$ such that if  a
$\delta$-pseudo-orbit $\{x_i,\,n_i\}_{i=-\infty}^{+\infty}$
satisfies $n_i\geq T_\tau $ and $x_i,f^{n_i} (x_i)\in
\Lambda_\tau$ for all $i$, then there exists a unique  exponentially
$ (\eta,\theta_\tau)$-shadowing point $x\in M $ for
$\{x_i,n_i\}_{i=-\infty}^{+\infty}$. If further
$\{x_i,n_i\}_{i=-\infty}^{+\infty}$ is periodic, i.e., there
 exists an integer $m>0$ such that $x_{i+m}=x_i$ and $n_{i+m}=n_i$ for all i,
 then the shadowing point $x$ can be chosen to be periodic.\\

 (ii) \emph{Stable Manifold Theorem}:  There exists $\sigma>0$ such that all points in $\Lambda^*_\tau (\sigma)$  have uniform sizes of stable and unstable manifolds, where $\Lambda^*_\tau (\sigma)$  denotes the union of $\Lambda_\tau$ and the set $$\{z\in M\,\,|\,\,\exists \text{ pseudo-orbit } \{x_i,\,n_i\}_{i=-\infty}^{+\infty}\text{
satisfying }n_i\geq T_\tau \text{ and } x_i,f^{n_i} (x_i)\in
\Lambda_\tau$$
$$\text{ for all }i \text{ such that } z \text{ is a }
\sigma-\text{shadowing point } \text{ for }
\{x_i,n_i\}_{i=-\infty}^{+\infty} \}.$$ Moreover, if $x,y\in\Lambda^*_\tau (\sigma)$ are close enough, then the  (local) stable manifold at $x$ is transverse to the  (local) unstable manifold of $y.$

\end{Thm}

{\bf Proof of Theorem \ref{Thm:ExpShadowing-measure2015oct}}  By assumption we can take an $f$-invariant set $\Delta'\subseteq$ such that $T_{\Delta}M=E\oplus F$  is dominated on $\Delta'$, $\Delta'$  is average-nonuniformly hyperbolic with respect to this splitting and $\Delta'$  has $\mu$ full measure.  Fix $\tau>0.$ By Theorem  \ref{Thm-topocondition-to-pesinset-0000000000} we can take $\zeta>0$ small enough and $K$ large enough such that $\Delta\subseteq \Lambda^\# (K,\,\zeta) $ and then we can take large integer $k$ such that  $\mu(\Lambda^\#_k (K,\,\zeta))>1-\tau$.  Let $$\Lambda_\tau=\Lambda^\#_k (K,\,\zeta).$$
Then by Proposition \ref{Prop:GeneralPesinStableManfld} and Proposition \ref{Prop:ExpShadowinglem}, take $\theta_\tau=\theta (\zeta)>0,\,\,T_\tau=2kK \in \mathbb{N}$ and thus $\Lambda_\tau$ is the needed set.\qed


Remark that in fact $\theta_\tau$ is independent of $\tau,$ since   the choice of $\zeta$ is not necessarily small and can be a fixed positive number independent on the variation of $\Lambda_\tau$ and from Proposition \ref{Prop:ExpShadowinglem}  $\theta_\tau$ only depends on $\zeta.$

\bigskip

In particular, we state the exponentially closing lemma.

\begin{Thm}\label{Thm:Expclosinglem-measure2015oct}  (Exponentially Closing lemma) Let $f\in \Diff^1 (M)$ and $\mu\in\m  (M).$  Suppose that there is  a $Df$-invariant dominated splitting $T_{\Delta}M=E\oplus F$ on an $f$-invariant set $\Delta$ with $\mu$ full measure and  $\mu$ is average-nonuniformly hyperbolic with respect to
 to this splitting. Then for each $\tau>0,$ there exist a compact set $\Lambda_\tau\subseteq M$, $\theta_\tau>0$  and $T_\tau \in \mathbb{N}$ such that $\mu (\Lambda_\tau)>1-\tau$ and following two properties hold.\\
 (i) \emph{ (Exponentially) closing Lemma}:
For $\forall\,\,
\eta>0,$ there exists $\beta=\beta (\tau,\eta)>0$ such that if for an orbit
segment $\{x,n\}$ with length $n\geq T_\tau$, one has $x,f^{n} (x)\in
\Lambda_\tau$ and $   d (x,f^n (x))<\beta$, then there exists
a unique hyperbolic periodic point $z=z (x)\in M $ satisfying:

 (1) $f^n (z)=z$;

 (2) $   d (f^j (x),\,\,f^j (z))<\eta\cdot e^{-\min\{j,n-j\}\theta_\tau}$, $j=0,\,1,\,\cdots,n-1$.
\\
 (ii)  \emph{Stable Manifold Theorem}:  There exists $\sigma>0$ such that all points
in $\Lambda^*_\delta (\sigma)$  have uniform sizes of stable and unstable manifolds,
where $\Lambda^*_\tau (\sigma)$  denotes the union set of $\Lambda_\tau$ and the set of periodic points
nearby $\Lambda_\tau$ $$\{z\in M\,\,|\, \exists\, x, n\geq T_\tau \, \text{with}\, x, f^n (x)\in\Lambda_\tau
\, \text{s.t.}\,  f^n (z)=z,\,\,   d (f^i (x), f^i (z))\leq\sigma,0\leq i \leq n\}.$$ Moreover,
if $x,y\in\Lambda^*_\tau (\sigma)$  are close enough,
then the  (local) stable manifold $x$ is transverse to the  (local) unstable manifold of $y.$

\end{Thm}

By Theorem \ref{Thm-topocondition-to-pesinset} and average-nonuniform hyperbolicity of
Theorem \ref{Thm-partial-hyperbolic-333}, Theorem \ref{Thm-partial-hyperbolic-usual-Lebegue-hyp}, we have
\begin{Thm}\label{Thm-partial-hyperbolic-Lebesgue-Pesinset}
Under the same assumptions as Theorem \ref{Thm-partial-hyperbolic-22222} or the cases of   (A)  (A')  (B)  (B') in Theorem \ref{Thm-partial-hyperbolic} or Theorem \ref{Thm-partial-hyperbolic-usual-Lebegue-hyp}:\\
 there is some $\zeta_0>0$ and $K_0\geq 1$ such that for any $\zeta\in (0,\zeta_0)$ and any $K\geq K_0$,
 $ \Lambda^{\#} (K,\zeta)$ has Lebesgue full measure.

\end{Thm}
By Theorem \ref{Thm-partial-hyperbolic-Lebesgue-Pesinset},
Theorem \ref{Thm:ExpShadowing-measure2015oct} and \ref{Thm:Expclosinglem-measure2015oct} can be applicable to some partially hyperbolic systems for Lebesgue measure which is not necessarily invariant.

\subsection{Hyperbolic measure with  (quasi-)limit domination}\label{section-full Pesinset}
On the other hand, it is still unknown whether the assumption in Theorem \ref{Thm-topocondition-to-pesinset} can be replaced by non-zero Lyapunov exponents.
\begin{Que}\label{Que-topocondition-to-pesinset-Lyapunov}
 Let $\Delta$ be an $f-$invariant set
and $T_{\Delta}M=E\oplus F$ be a $Df-$invariant splitting on
$\Delta$. Let $K\in \mathbb{Z}^+$ and $\zeta_0>0.$  Suppose that  $T_{\Delta}M=E\oplus F$  is $K-$limit-dominated on $\Delta$ and $\Delta$ is a nonuniformly hyperbolic set with
$ (K,\,\zeta_0)$-degree corresponding to $T_{\Delta}M=E\oplus F$, that is, for any $x\in \Delta$
$$\limsup_{l\rightarrow
+\infty} \frac{\log\|Df^{l}|_{E (x)}\|}{l}\leq- \zeta_0,\,\,\,\,\,\,and$$
$$ \liminf_{l\rightarrow +\infty} \frac{\log
m (Df^{l}|_{F (f^{-l} (x))})}{l}\geq\zeta_0.$$
 Then whether    $$\Delta\subseteq \bigcup_{K\geq 1, \chi>0}\Lambda (K,\chi)?$$
  If further the splitting $T_{\Delta}M=E\oplus F$  is $K-$dominated, whether one has  $$\Delta\subseteq  \bigcup_{K\geq 1, \chi>0} \Lambda^\# (K,\chi)?$$

\end{Que}




We will give a positive answer in the sense of probabilistic perspective that if $\m_f (\Delta)\neq \emptyset,$ then for any $\mu\in\m_f (\Delta)$ $\mu$ a.e.  point in $\Delta$ should
 be in the Liao-Pesin set.   That is,

 \begin{Thm}\label{Thm-topocondition-to-pesinset-Lyapunov}
 Let $\Delta$ be an $f-$invariant set
and $T_{\Delta}M=E\oplus F$ be a $Df-$invariant splitting on
$\Delta$. Let $K\in \mathbb{Z}^+$ and $\zeta_0>0.$  Suppose that  $T_{\Delta}M=E\oplus F$  is $K-$limit-dominated on $\Delta$ and $\Delta$ is a nonuniformly hyperbolic set with
$ (K,\,\zeta_0)$-degree corresponding to $T_{\Delta}M=E\oplus F$, that is, for any $x\in \Delta$
$$\limsup_{l\rightarrow
+\infty} \frac{\log\|Df^{l}|_{E (x)}\|}{l}\leq- \zeta_0,\,\,\,\,\,\,and$$
$$ \liminf_{l\rightarrow +\infty} \frac{\log
m (Df^{l}|_{F (f^{-l} (x))})}{l}\geq\zeta_0.$$ If $\m_f (\Delta)\neq \emptyset,$ then there is $\Delta'\subseteq \Delta$ such that $$\Delta'\subseteq \bigcup_{K\geq 1, \chi>0}\Lambda (K,\chi)$$
and for any $\mu\in\m_f (\Delta)$, $\mu (\Delta')=1.$  If further the splitting $T_{\Delta}M=E\oplus F$  is dominated, then there is $\Delta'\subseteq \Delta$ such that $$\Delta'\subseteq \bigcup_{K\geq 1, \chi>0}\Lambda^\# (K,\chi)$$
and for any $\mu\in\m_f (\Delta)$, $\mu (\Delta')=1.$
\end{Thm}

 Obviously for any $\mu\in\m_f (\Delta)$,  all the Lyapunov exponents of $\mu$
  are far from zero and $\mu$ has limit-domination or domination.
  Thus, we only need to show that   for any hyperbolic measure $\mu$ with limit-domination,
  $\mu$ a.e. points are in the Liao-Pesin set.
  Now for nonuniformly hyperbolic systems with limit  (quasi-)domination,  let us  state and show the existence of Pesin set with full
measure or  with measure arbitrarily close to 1.


\begin{Thm}\label{Thm:fullPesinset}
Let $f\in \Diff^1 (M)$ and $\mu\in\m_f^{qldh} (M).$
Then
$$\mu (\cup_{\zeta>0}\cup_{K\in\mathbb{N}}\Lambda (K,\,\zeta))=\mu (\cup_{\zeta>0}\cup_{K\in\mathbb{N}}\Lambda^{\#} (K,\,\zeta))=1.$$ If further all the Lyapunov exponents of $\mu$ are far from zero, then there exists $\zeta_0>0$ such that for all $\zeta\in (0,\zeta_0)$ $$\mu (\cup_{K\in\mathbb{N}}\Lambda (K,\,\zeta))=\mu (\cup_{K\in\mathbb{N}}\Lambda^{\#} (K,\,\zeta))=1.$$ Moreover, if $\mu$ is ergodic, then there exists $K_0\in\mathbb{N},$ $\zeta_0>0$ such that for all $K\geq K_0,\,\,\zeta\in (0,\zeta_0),$ $$\mu (\Lambda (K,\,\zeta))=\mu (\Lambda^{\#} (K,\,\zeta))=1.$$
\end{Thm}

To prove this theorem we need a lemma as follows.
\begin{Lem}\label{Lem:Generalized-Multi-erg-thm}{\bf  (Generalized Oseledec Multiplicative Ergodic Theorem or Sub-additional Ergodic Theorem)}\\
Let $f\in \Diff^1 (M)$ preserve an  invariant probability
measure $\mu$, and $E\subseteq TM$ be a $Df-$invariant subbundle
defined over an $f-$invariant set with $\mu$ full measure. Let
$\lambda_{E}^+ (x)$ be the maximal Lyapunov exponent in $E (x)$ of the
measure $\mu.$ Then, for any $\tau>0,\varepsilon>0$ there exist an invariant set $B_{\tau,\varepsilon}$ with $\mu (B_{\tau,\varepsilon})\geq 1-\tau$ and an
integer $K_{\tau,\varepsilon}$ such that for  every point
$x\in B_{\tau,\varepsilon}$ and any $K\geq K_{\tau,\varepsilon}$, the Birkhoff averages
$$\frac1{lK}\sum_{j=0}^{l-1}\log\|Df^K|_{E (f^{jK} (x))}\|$$ converge
towards a number contained in
$[\lambda_{E}^+ (x),\lambda_{E}^+ (x)+\varepsilon)$, when
$l\rightarrow+\infty.$ In particular, for $\mu$ a.e. $x$, we have the following limits exist and
$$\lim_{K\rightarrow +\infty}\lim_{l\rightarrow +\infty}\frac1{lK}\sum_{j=0}^{l-1}\log\|Df^K|_{E (f^{jK} (x))}\|=\lambda_{E}^+ (x).$$
Moreover, one should have same result for the Birkhoff averages
$$
\frac1{lK}\sum_{j=1}^{l }\log\|Df^K|_{E (f^{-jK} (x))}\| (\text{ or written by }  \frac1{lK}\sum_{j=-l}^{-1 }\log\|Df^K|_{E (f^{jK} (x))}\|).$$

For the  case of minimal Lyapunov exponent,  similar results work.
 \end{Lem}

{\bf Proof.}
 The ergodic case firstly appeared in  \cite{ABC}. Here we prove the version for general invariant (not necessarily ergodic) measures. In particular, the ergodicity is useful in the  proof of  \cite{ABC} and  it is needed to consider $f^N-$ergodic measure for large $N$.   Here we choose another way which is not necessary to consider $f^N-$invariant or $f^N-$ergodic   measures.
Now we start to prove.

By Sub-additional Ergodic Theorem \cite{Walters} (or Oseledec Multiplicative Ergodic Theorem), one has
\begin{eqnarray*} &\lim_{l\rightarrow
+\infty}\int|\frac{\log\|Df^{l}|_{E (x)}\|}{l}-\lambda^+_{E} (x)|d\mu=\int\lim_{l\rightarrow
+\infty}|\frac{\log\|Df^{l}|_{E (x)}\|}{l}-\lambda^+_{E} (x)|d\mu=0.
   \end{eqnarray*}
For convenience to write,  later  we   always assume that following limits exist since they always exist for $\mu$ a.e. $x$.
 For any $\varepsilon>0$ and
$\tau\in (0,1)$, we can take large  $L$ such that  $\int h (x)d\mu<\frac12\tau\varepsilon$  where
$$h (x):=|\frac{\log\|Df^{{L}}|_{E (x)}\|}{{L}}-\lambda^+_E (x)|.
   $$
   Since $\mu$ is $f-$invariant,
$\frac1l\sum_{i=0}^{l-1}h (f^{i} (x))$ converge a.e. to
integrable nonnegative functions $h^* (x).$  Also
$h^* (f (x))=h^* (x)$ a.e. and $\int
h^* (x)d\mu=\int
h (x)d\mu<\frac12\tau\varepsilon.$

We
set $B=\{x\in M|\,h^* (x)<\frac12\varepsilon\}$ and  claim
$\mu (B)\geq 1-{\tau}$.  In fact, if $\mu (B)< 1-\tau$, then $\int
h^* (x)d\mu\geq\int_{M\setminus B} h^* (x)d\mu\geq
\frac12\varepsilon\mu ({M\setminus B})>\frac12\tau\varepsilon,$ a
contradiction to the inequality $\int
h^* (x)d\mu<\frac12\tau\varepsilon.$ So if $x\in B,$ then
$$\lim_{l\rightarrow
+\infty}\frac1l\sum_{i=0}^{l-1}|\frac{\log\|Df^{{L}}|_{E (^{i} (x))}\|}{{L}}-\lambda^+_E (x)|<\frac\varepsilon2.$$
This implies that for all $x\in B,$
\begin{eqnarray}\label{eq-average-estimate-Lya}
\lim_{l\rightarrow
+\infty}\frac1{lL}\sum_{i=0}^{l-1}\log\|Df^{{L}}|_{E (f^{i} (x))}\|<\lambda^+_E (x)+\frac\varepsilon2.
\end{eqnarray}
Remark that in the case that $\mu$ is $f$-ergodic,   this estimate is more easy to get for a set with full measure.
More precisely, by sub-additional ergodic principle, there is some $L$ large enough such that $$\int\frac{\log\|Df^{{L}}|_{E (x)}\|}{{L}}d\mu\leq\lambda^+_E (x)+\frac{\varepsilon}2.$$
By Birkhoff ergodic theorem and ergodicity of $\mu$, for $ \mu$  a.e. $x$, $$\lim_{l\rightarrow
+\infty}\frac1{lL}\sum_{i=0}^{l-1}\log\|Df^{{L}}|_{E (f^{i} (x))}\|=\int\frac{\log\|Df^{{L}}|_{E (x)}\|}{{L}}d\mu
<\lambda^+_E (x)+\frac\varepsilon2.$$

\smallskip

 Before continuing the proof, we need a lemma as follows which is also useful for other results below.

\begin{Lem}\label{Lem-norm-estimate}
Let $C_f=\max\{\log\|Df\|,\log\|Df^{-1}\|,0\}, $  let $\Delta\subseteq M$ be an invariant set and $E$ be a $Df$-invariant bundle over $\Delta$.
Suppose  $L\geq 1,$   $K\geq 2L.$ 
Then for any $x\in \Delta$,    \\
  (1) Take $r$ such that $rL\geq K> (r-1)L,$  for any $0\leq p<L$, one has   $$\log\|Df^K|_{E (x)}\|\leq 4L\cdot
 C_f+\sum_{s=0}^{r-1}\log\|Df^{L}|_{E (f^{p+sL} (x))}\|.$$
  (2)
$$L \log\|Df^K|_{E (x)}\|\leq 5L^2 \cdot
C_f+\sum_{t=0}^{K-1}\log\|Df^{L}|_{E (f^{t } (x))}\|. $$
  (2') $$\limsup_{n\rightarrow+\infty}\frac1n\log\|Df^n|_{E (x)}\|\leq \limsup_{l\rightarrow+\infty}\frac1{lK}\sum_{i=0}^{l-1}\log\|Df^K|_{E (f^{  iK} (x))}\|$$$$\leq
\frac{5L\cdot
C_f}K+\limsup_{l\rightarrow+\infty}\frac1{lL}\sum_{i=0}^{l-1}\log\|Df^{L}|_{E (f^{i} (x))}\|.$$
\end{Lem}
{\bf Proof.}

Take $r$ such that $\, rL\geq K> (r-1)L$.  For any $0\leq p<L$, one decomposes the orbit segment  of length $K+L$ of
$x$ as $$ (x, f (x),\cdots,f^{p-1} (x)), $$$$
\,\,\,\,\,\,\,\, (f^p (x),\cdots,f^{p+ (r-1)L-1} (x)),$$$$\,\,\,\,\,\,\,
 (f^{p+ (r-1)L} (x),\cdots,f^{K+L-1} (x)).$$    One
deduces that \begin{eqnarray*}
& &\|Df^K|_{E (x)}\|\leq\|Df^{K+L}|_{E (x)}\|\|Df^{-L}|_{E (f^{K+L}x)}\|\leq e^{ L\cdot
C_f}\|Df^{K+L}|_{E (x)}\|\\
&\leq &e^{ L\cdot
C_f}\|Df^p|_{E (x)}\|\cdot\,\,\, (\|Df^{L}|_{E (f^p (x))}\|\cdots\|Df^{L}|_{E (f^{p+ (r-2)L} (x))}\|)\,\,\,\,
\\
& &\cdot\|Df^{K+L- (p+ (r-1)L)}|_{E (f^{p+ (r-1)L} (x))}\|\\
&\leq& e^{3L\cdot
C_f}\cdot (\|Df^{L}|_{E (f^p (x))}\|\cdots\|Df^{L}|_{E (f^{p+ (r-1)L} (x))}\|)
\cdot\|Df^{L}|_{E (f^{p+ (r-1)L} (x))}\|^{-1} \\
&\leq& e^{4L\cdot
C_f}\cdot (\|Df^{L}|_{E (f^p (x))}\|\cdots\|Df^{L}|_{E (f^{p+ (r-1)L} (x))}\|).
 \end{eqnarray*}
 Hence,   $$\log\|Df^K|_{E (x)}\|\leq 4L\cdot
C_f+\sum_{s=0}^{r-1}\log\|Df^{L}|_{E (f^{p+sL} (x))}\|.$$
 So summing the inequalities for $p=0,1,\cdots,L-1,$ one has  $$L \log\|Df^K|_{E (x)}\|\leq 4 L^2\cdot
C_f+\sum_{p=0}^{L-1}\sum_{s=0}^{r-1}\log\|Df^{L}|_{E (f^{p+sL} (x))}\|$$
$$=4L^2\cdot
C_f+\sum_{t=0}^{Lr-1}\log\|Df^{L}|_{E (f^{t } (x))}\| \leq 5L^2 \cdot
C_f+\sum_{t=0}^{K-1}\log\|Df^{L}|_{E (f^{t } (x))}\|. $$
These imply  (1) and  (2) hold for all $x\in \Delta$.

For any $y\in Orb (x),$  by  (2) one has  $$  \log\|Df^K|_{E (y)}\|\leq  5L \cdot
C_f+\frac1L\sum_{j=0}^{K-1}\log\|Df^{L}|_{E (f^{j } (y))}\|.$$
 Then summing the inequalities for $y=x,f^K (x),\cdots,f^{ (l-1)K} (x),$
$$ \sum_{i=0}^{l-1}\log \|Df^K|_{E (f^{iK} (x))}\|\leq   {5lL\cdot
C_f} + \sum_{i=0}^{l-1}\sum_{j=0}^{K-1}\frac1L\log\|Df^{L}|_{E (f^{j+iK } (x))}\|$$
$$={5lL\cdot
C_f} + \frac1L\sum_{t=0}^{lK-1} \log\|Df^{L}|_{E (f^{t } (x))}\|.$$ It follows that
$$\limsup_{l\rightarrow+\infty}\frac1{lK}\sum_{i=0}^{l-1}\log \|Df^K|_{E (f^{iK} (x))}\|\leq \frac{5L\cdot
C_f}K+\limsup_{l\rightarrow+\infty}\frac1{lKL}\sum_{i=0}^{lK-1}\log\|Df^{L}|_{E (f^{i} (x))}\| $$
$$\leq \frac{5L\cdot
C_f}K+\limsup_{l\rightarrow+\infty}\frac1{lL}\sum_{i=0}^{l-1}\log\|Df^{L}|_{E (f^{i} (x))}\|.$$

 By sub-multiplications of the norms (write $n=lK+r,\,0\leq r<K$), we have $$\log\|Df^n|_{E (x)}\|\leq \sum_{i=0}^{l-1}\log\|Df^K|_{E (f^{  iK} (x))}\|+ (K-1)C_f
 $$ Then
$$\limsup_{l\rightarrow+\infty}\frac1{lK}\sum_{i=0}^{l-1}\log\|Df^K|_{E (f^{  iK} (x))}\|\geq
\limsup_{n\rightarrow+\infty}\frac1n\log\|Df^n|_{E (x)}\|.$$
 (The case of $\liminf $ is also true). This ends the proof of  (2').\qed


 Now we continue to prove  Lemma \ref{Lem:Generalized-Multi-erg-thm}. By Lemma \ref{Lem-norm-estimate}  (2') and sub-additional ergodic theorem,  take and fix
$K_{\tau,\varepsilon}\geq \max\{2L, \frac{10L\cdot C_f}{\varepsilon}\}$ and then by  (\ref{eq-average-estimate-Lya}) for all $K\geq K_{\tau,\varepsilon}$, one
gets for $\mu$ a.e. $x\in B$,
$$\lambda_{E}^+ (x)=\lim_{n\rightarrow+\infty}\frac1n\log\|Df^n|_{E (x)}\|\leq\lim_{l\rightarrow+\infty}\frac1{lK}\sum_{i=0}^{l-1}\log\|Df^K|_{E (f^{iK} (x))}\|<
\lambda_{E}^+ (x)+\varepsilon.$$
 So we complete the proof.\hfill $\Box$

\bigskip

 Moreover,    by same methods, similar results should hold for continuous linear cocycles on homeomorphisms and general sub-additional sequence of integrable functions. Here we just state the case of cocycles, and the case of sub-additional functions is left for readers.   More precisely,  let $f$ be an invertible map of a compact metric space $X$ and let $A:X\rightarrow GL (m,\mathbb{R}) (m\geq 1)$  be a  measurable  function. One main object of interest is the asymptotic behavior
of the products of $A$ along the orbits of the
transformation $f$, called cocycle induced from $A$: for $n>0$
 $$A (x,n):=A (f^{n-1} (x))\cdots A (f (x))A (x), $$ and $$A (x,-n):=A (f^{-n} (x))^{-1}\cdots A (f^{-2} (x))^{-1} A (f^{-1} (x))^{-1}=A (f^{-n}x,n)^{-1}.$$  For any $x\in X$ and any $v\in \mathbb{R}^m,$ define the Lyapunov exponent of vector $v$ at $x$,  $$\lambda^A (x,v):=\lim_{n\rightarrow \infty}\frac1n{\log\|A (x,n)v\|},
   $$ if the limit exists. If the above limit exists for all $v\neq 0$, we say $x$ to be {\it Lyapunov-regular.}

 \begin{Thm}\label{thm-generalizeOseledec}
 Let $f:X\rightarrow X$ be a measure-preserving invertible map of a Lebesgue space $ (X,\mu)$ and let $A:X\rightarrow GL (n,\mathbb{R})$ be a bounded measurable matrix functions over $X$.  Then there exists  a set $Y\subseteq X$ such that $\mu (X \setminus Y)=0$ and for each $x\in Y:$ there exists a decomposition of $$\mathbb{R}^n=\oplus_{i=1}^{k (x)}H_i (x)$$
  that is invariant under the linear extension of $f$ determined by $A$. The Lyapunov exponents $\chi_1 (x)<\chi_2 (x)\cdots<\chi_{k (x)} (x)$ exist and are $f$-invariant and
   $$\lim_{K\rightarrow \pm\infty}\lim_{l\rightarrow+\infty}\frac1{lK}\sum_{i=0}^{l-1}\log\|A (f^{iK} (x),K)|_{H (f^{\pm  iK} (x))}\|$$
   $$=\lim_{K\rightarrow \pm\infty}\lim_{l\rightarrow+\infty}\frac1{lK}\sum_{i=0}^{l-1}\log m (A (f^{iK} (x),K)|_{H (f^{\pm  iK} (x))})$$
  $$=\lim_{m\rightarrow \pm \infty} \log \|A (x,m)|_{H_i (x)}\|=\lim_{m\rightarrow \pm \infty} \log m (A (x,m)|_{H_i (x)})=\chi_i.$$
 \end{Thm}


\bigskip

Now we start to prove Theorem \ref{Thm:fullPesinset}.

{\bf Proof of Theorem \ref{Thm:fullPesinset}}
Let
$\lambda_s (x),\lambda_u (x)$ denote the maximal and minimal Lyapunov exponents of $E^s (x)$ and $E^u (x)$ respectively for $\mu$ a.e. $x$.
Notice that $$\mu\big{ (}\,\cup_{\zeta>0}\{x\in M \,\,|\,\,\lambda_s (x)<-2\zeta,\lambda_u (x)>2\zeta\}\,\big{)}=1$$ and by Corollary \ref{cor-hyp-quasilimitdomination},
$$\mu\big{ (}\,\cup_{\zeta>0}\cup_{S_0\geq 1}\{x\in M \,\,|\,\,\frac1S\log\frac
{\|Df^S|_{E^s (y)}\|}{m (Df^S|_{E^u (y)})}\leq-2\zeta,\,\,\forall y\in Orb (x),\,\,S\geq S_0\}\,\,\big{)}=1.$$
So for any $\tau\in (0,1),$ we can take small $\zeta>0$, large $S_0$ and an $f$-invariant set $\Delta_1$ with $\mu (\Delta_1)>1-\frac\tau2$ such that if $x\in\Delta_1$ then
 $\lambda_s (x)<-2\zeta,\,\,\lambda_u (x)>2\zeta$ and $$\frac1S\log\frac
{\|Df^S|_{E^s (y)}\|}{m (Df^S|_{E^u (y)})}\leq-2\zeta,\,\,\forall y\in Orb (x)\,\,S\geq S_0.$$

On the other hand, for $\varepsilon=\zeta$ in Lemma \ref{Lem:Generalized-Multi-erg-thm} we can take an invariant set $\Delta_2=B_{\frac\tau2,\varepsilon}$ with $\mu (\Delta_2)>1-\frac\tau2$ such that $\Delta_2$ satisfies the result of Lemma \ref{Lem:Generalized-Multi-erg-thm} for both bundles $E^s (x)$ and $E^u (x)$.
More precisely, there is
$K_{\frac\tau2,\varepsilon}$ such that for  every point
$x\in \Delta_2$ and any $K\geq K_{\frac\tau2,\varepsilon}$,  the Birkhoff averages satisfy
$$\lim_{l\rightarrow+\infty}\frac1{lK}\sum_{j=0}^{l-1}\log\|Df^K|_{E^s (f^{jK} (x))}\|
\leq \lambda_{s} (x)+\varepsilon$$  and simultaneously $$\lim_{l\rightarrow
+\infty}\frac1{lK}\sum_{j=-l}^{-1}{\log
m (Df^{K}|_{E^u ({f^{jK} (x)})})}
\geq \lambda_{u} (x)-\varepsilon.$$

Take $\Delta_\tau=\Delta_1\cap\Delta_2$ and $K_\tau=\max\{S_0,\,K_{\frac\tau2,\varepsilon}\}$. Clearly $\mu (\Delta_\tau)>1-\tau$. Then
for  every point
$x\in \Delta_\tau$ and any $K\geq K_\tau$, we have
$$\lim_{l\rightarrow+\infty}\frac1{lK}\sum_{j=0}^{l-1}\log\|Df^K|_{E^s (f^{jK} (x))}\|\leq \lambda_{s} (x)+\varepsilon\leq - \zeta ,$$   $$\lim_{l\rightarrow
+\infty}\frac1{lK}\sum_{j=-l}^{-1}{\log
m (Df^{K}|_{E^u ({f^{jK} (x)})})}\geq \lambda_{u} (x)-\varepsilon\geq  \zeta$$ and
$$\frac1L\log\frac
{\|Df^L|_{E^s (y)}\|}{m (Df^L|_{E^u (y)})}\leq-2\zeta,\,\,\forall y\in Orb (x),\,\,L\geq K.$$
 By Theorem \ref{Thm-topocondition-to-pesinset}, $\Delta_\tau\subseteq \Lambda^\# (K,\zeta)$ (in fact, here  (c) and  (d) in the definition of Liao-Pesin set are obvious from above last equality). By the arbitrary choice of small $\tau$ we complete the proof.
\qed

\bigskip

 In particular if further the hyperbolic measure is  ergodic, we have
\begin{Thm}\label{Thm:fullPesinset-ergodic}
Let $f\in \Diff^1 (M)$   and $\mu\in\m_f^{qldh} (M)$ be ergodic. Then there is some $\beta>0$ and $K_0\in\mathbb{Z}^+$ such that for any $0<\zeta<\beta$,  any  $K\geq K_0$,  Pesin set
$\Lambda (K,\,\zeta)$ is of $\mu$ full measure.
\end{Thm}

\begin{Rem}\label{Rem:Pesinblocks-for-surfacediff}  Every ergodic measure supported on hyperbolic sets obviously satisfies the assumption of Theorem \ref{Thm:fullPesinset}. And every ergodic measure $\mu$ of partially hyperbolic diffeomorphisms satisfies the assumption of Theorem \ref{Thm:fullPesinset}, provided that the Lyapunov exponents of $\mu$ in the center bundle are all positive (or all negative). In particular, if the dimension of the center bundle of a partially hyperbolic diffeomorphism is one, then every hyperbolic ergodic invariant measure naturally satisfies the assumption of Theorem \ref{Thm:fullPesinset}.  For a nonuniformly
hyperbolic system on a {\it surface} with  limit
domination (or domination), since the subbudles are both one dimensional, one can
take $K_0\geq1$ and get a  Pesin set $\Lambda (K,\,\zeta)$ of full
measure for any $K\geq K_0$ according to the proof of
Theorem~\ref{Thm:fullPesinset}. 
In this case  (a) and  (b) in definition of Liao-Pesin set can be written simply: \\
$$ (a).\,\,\,\,\,\,\,\,\,\,\,\,\,\,\,\,\,\,\,\,\,\,\,\,\,\,\,\,\,\,\,\,\,\,\,
 \frac 1{l}\log\|Df^{l}|_{E (x)}\|
   \leq-\zeta,\,\,\,\,\,\,\, \,\forall\,\,\,
l\,\geq \, kK;\,\,\,\,\,\,\,\,\,\,\,\,\,\,\,\,\,\,\,\,\,\,\,\,\,\,\,\,\,\,\,\,\,\,\,\,\,\,\,$$
$$ (b).\,\,\,\,\,\,\,\,\,\,\,\,\,\,\,\,\,\,\,\,\,\,\,\,\,\,\,\,\,\,
\,\, \frac1l \log
  m (Df^{l}|_{F (f^{-l}x)}) \geq\zeta,\, \,\,\,\forall\,\,\,
l\,\geq \, kK.\,\,\,\,\,\,\,\,\,\,\,\,\,\,\,\,\,\,\,\,\,\,\,\,\,\,\,\,\,\,\,\,\,\,\,\,\,\,\,$$

\end{Rem}
\bigskip

By Theorem \ref{Thm:fullPesinset} and Proposition \ref{Prop:ExpShadowinglem} and \ref{Prop:GeneralPesinStableManfld}, we have

\begin{Thm}\label{Thm:ExpShadowing-measure} Let $f\in \Diff^1 (M)$ and $\mu\in\m_f^{qldh} (M).$ For each $\tau\in (0,1),$ there exist a compact set $\Lambda_\tau\subseteq M$, $\theta_\tau>0$ and $T_\tau \in \mathbb{N}$ such that $\mu (\Lambda_\tau)>1-\tau$ and following two properties hold:\\

 (i) \emph{ (Exponentially) Shadowing Lemma}: For
$\forall\,\,\eta>0,$ there exists $\delta=\delta (\tau,\eta)>0$ such that if  a
$\delta$-pseudo-orbit $\{x_i,\,n_i\}_{i=-\infty}^{+\infty}$
satisfies $n_i\geq T_\tau $ and $x_i,f^{n_i} (x_i)\in
\Lambda_\tau$ for all $i$, then there exists a unique  exponentially
$ (\eta,\theta_\tau)$-shadowing point $x\in M $ for
$\{x_i,n_i\}_{i=-\infty}^{+\infty}$. If further
$\{x_i,n_i\}_{i=-\infty}^{+\infty}$ is periodic, i.e., there
 exists an integer $m>0$ such that $x_{i+m}=x_i$ and $n_{i+m}=n_i$ for all i,
 then the shadowing point $x$ can be chosen to be periodic.\\

 (ii) \emph{Stable Manifold Theorem}:  There exists $\sigma>0$ such that all points in $\Lambda^*_\tau (\sigma)$  have uniform sizes of stable and unstable manifolds, where $\Lambda^*_\tau (\sigma)$  denotes the union of $\Lambda_\tau$ and the set $$\{z\in M\,\,|\,\,\exists \text{ pseudo-orbit } \{x_i,\,n_i\}_{i=-\infty}^{+\infty}\text{
satisfying }n_i\geq T_\tau \text{ and } x_i,f^{n_i} (x_i)\in
\Lambda_\tau$$
$$\text{ for all }i \text{ such that } z \text{ is a }
\sigma-\text{shadowing point } \text{ for }
\{x_i,n_i\}_{i=-\infty}^{+\infty} \}.$$ Moreover, if $x,y\in\Lambda^*_\tau (\sigma)$ are close enough, then the  (local) stable manifold at $x$ is transverse to the  (local) unstable manifold of $y.$

In particular, if
the hyperbolic invariant measure $\mu$ is also \emph{ergodic} (or all Lyapunov exponents far away from zero), then $\theta_\tau$ can be chosen independent of  $\tau.$

\end{Thm}

\begin{Rem}\label{Rem:ExpShadowing}
 Note that for every hyperbolic  basic set of Axiom A systems, Anosov shadowing lemma and the local product
 structure yield the exponential shadowing property. In the $C^{1+\alpha}$ non-uniformly hyperbolic case,
 exponential shadowing property
 is a particular case of Katok Shadowing lemma  \cite{P1} (see Proposition \ref{Prop-Shadow-expon} for a detailed proof).
\end{Rem}


{\bf Proof of Theorem \ref{Thm:ExpShadowing-measure}} Fix $\tau>0.$ By Theorem \ref{Thm:fullPesinset} we can take $\zeta>0$ small enough and $K$ large enough such that $\mu (\Lambda^\# (K,\,\zeta))>1-\tau$ and then we can take large integer $k$ such that  $\mu (\Lambda^\#_k (K,\,\zeta))>1-\tau$.  Let $$\Lambda_\tau=\Lambda^\#_k (K,\,\zeta).$$
Then by Proposition \ref{Prop:GeneralPesinStableManfld} and Proposition \ref{Prop:ExpShadowinglem}, take $\theta_\tau=\theta (\zeta)>0,\,\,T_\tau=2kK \in \mathbb{N}$ and thus $\Lambda_\tau$ is the needed set.


If all Lyapunov exponents are far away from zero for a.e. points (ergodic hyperbolic measure is a particular case), the choice of $\zeta$ is not necessarily small and can be a fixed positive number independent on the variation of $\Lambda_\tau$ and from Proposition \ref{Prop:ExpShadowinglem}  $\theta_\tau$ only depends on $\zeta.$ \qed

\bigskip

In particular, we state the exponentially closing lemma.

\begin{Thm}\label{Thm:Expclosinglem-measure}  (Exponentially Closing lemma) Let $f\in \Diff^1 (M)$ and $\mu\in\m_f^{qldh} (M).$   For each $\tau>0,$ there exist a compact set $\Lambda_\tau\subseteq M$, $\theta_\tau>0$  and $T_\tau \in \mathbb{N}$ such that $\mu (\Lambda_\tau)>1-\tau$ and following two properties hold.\\
 (i) \emph{ (Exponentially) closing Lemma}:
For $\forall\,\,
\eta>0,$ there exists $\beta=\beta (\tau,\eta)>0$ such that if for an orbit
segment $\{x,n\}$ with length $n\geq T_\tau$, one has $x,f^{n} (x)\in
\Lambda_\tau$ and $   d (x,f^n (x))<\beta$, then there exists
a unique hyperbolic periodic point $z=z (x)\in M $ satisfying:

 (1) $f^n (z)=z$;

 (2) $   d (f^j (x),\,\,f^j (z))<\eta\cdot e^{-\min\{j,n-j\}\theta_\tau}$, $j=0,\,1,\,\cdots,n-1$.
\\
 (ii)  \emph{Stable Manifold Theorem}:  There exists $\sigma>0$ such that all points
in $\Lambda^*_\delta (\sigma)$  have uniform sizes of stable and unstable manifolds,
where $\Lambda^*_\tau (\sigma)$  denotes the union set of $\Lambda_\tau$ and the set of periodic points
nearby $\Lambda_\tau$ $$\{z\in M\,\,|\, \exists\, x, n\geq T_\tau \, \text{with}\, x, f^n (x)\in\Lambda_\tau
\, \text{s.t.}\,  f^n (z)=z,\,\,   d (f^i (x), f^i (z))\leq\sigma,0\leq i \leq n\}.$$ Moreover,
if $x,y\in\Lambda^*_\tau (\sigma)$  are close enough,
then the  (local) stable manifold $x$ is transverse to the  (local) unstable manifold of $y.$

In particular, if further
the hyperbolic invariant measure $\mu$ is   \emph{ergodic} (or all Lyapunov exponents of $\mu$ a.e. points  are far away from zero), then $\theta_\tau$ can be chosen independent of  $\tau.$
\end{Thm}

\section{Hyperbolic Measures, Partial  Hyperbolicity and Lebesgue measure}
\subsection{Existence of hyperbolic measures}\label{proof-exitence-hyperbolic}

Now we state a basic fact for systems with  (limit-)dominated splitting.

\begin{Thm}\label{Thm-partial-hyperbolic-hyperbolicmeasure}
Let $f:M\rightarrow M$ be a $C^1$ diffeomorphism on a   compact Riemanian manifold $M$.  Let $\Delta\subseteq M$ be an  $f-$invariant set (not necessarily compact)  and $T_{\Delta}M=E\oplus F$ be
a  $Df-$invariant splitting on $\Delta$.  Suppose $dim (E (x))$ is constant on $\Delta$, denoted by $dim (E).$\\
 (1) If 
  $T_{\Delta}M=E\oplus F$ is  (quasi-)limit-dominated, then $$\{\mu\in \mathcal{M}^{h}_{f} (\Delta)| \,ind (\mu)=dim E\}=\{\mu\in \mathcal{M}^{ldh}_{f} (\Delta)| \,ind (\mu)=dim E\}$$$$ (\text{ or } =\{\mu\in \mathcal{M}^{qldh}_{f} (\Delta)| \,ind (\mu)=dim E\}).$$\\
  (2) If  $T_{\Delta}M=E\oplus F$ is   (quasi-)dominated, then$$\{\mu\in \mathcal{M}^{h}_{f} (\Delta)| \,ind (\mu)=dim E\}=\{\mu\in \mathcal{M}^{qdh}_{f} (\Delta)| \,ind (\mu)=dim E\}$$$$ (\text{ or }  =\{\mu\in \mathcal{M}^{dh}_{f} (\Delta)| \,ind (\mu)=dim E\}).$$


\end{Thm}

Let us recall a result for the existence of hyperbolic measures in $C^1$ generic systems.

\begin{Thm}\label{Thmgeneric-measure-hyperbolic} ( \cite{ABC}) Let $\Lambda$ be an isolated non-trivial transitive set of
a $C^1$ generic diffeomorphism $f\in \Diff (M)$ and let $T_\Lambda=E_1\oplus E_2\cdots E_m$ be a finest dominated splitting (that is, every $E_i$ can not be decomposed into two dominated subbundles).
 Then generic measures $\mu\in \m_f (\Lambda)$ are ergodic, hyperbolic with support $\Lambda$ and their Oseledec splitting all coincide  with
$T_\Lambda=E_1\oplus E_2\cdots E_m$. \\In particular, this implies $ \m_f^{dh} (\Lambda)\cap \m_e (\Lambda)$ contains a dense $G_\delta$ subset of $\m_f (\Lambda).$

\end{Thm}

Note that in a partially hyperbolic diffeomorphism, the corresponding splitting is always dominated. Moreover, if the Lyapunov exponents of some ergodic measure $\mu$ in the center bundle are all positive (or all negative), then $\mu$ is hyperbolic measure and its Oseledec's hyperbolic splitting $T_{x}M=E^s (x)\oplus
E^u (x)$ is dominated.  Here we state a sufficient condition to show the existence of ergodic hyperbolic measure with domination  in partially hyperbolic diffeomorphisms (or partially hyperbolic invariant sets).

\begin{Thm}\label{Thm:HypMeas-paritialhyper}
Let $f\in \Diff^1 (M)$  be a partially hyperbolic diffeomorphism with a corresponding global $Df$ invariant splitting $T_{M}M=E^s\oplus E^c
\oplus E^u$. If there exists a point $x_0\in M$ such that $$\liminf_{n\rightarrow +\infty}\frac1n\sum_{i=0}^{n-1}\log{\|Df|_{E^c (f^i (x_0))}\|}<0$$$$\text{  (or respectively, }\limsup_{n\rightarrow +\infty}\frac1n\sum_{i=0}^{n-1}\log{m (Df|_{E^c (f^i (x_0))})}>0),$$ then there exists at least one  ergodic  measure  $\mu$ such that $\mu$ is hyperbolic and its Oseledec's hyperbolic splitting $T_{x}M=E^s (x)\oplus
E^u (x)$ coincides with $ (E^s (x)\oplus E^c (x))
\oplus E^u (x)$ (or respectively, $E^s (x)\oplus  (E^c (x)
\oplus E^u (x))$) and thus is dominated. This implies $ \m_f^{dh} (M)\cap \m_e (M)\neq \emptyset.$
\end{Thm}

It is also allowable without stable bundle  (or unstable bundle) in the assumption.
This theorem suggests that for partially hyperbolic systems,  average-nonuniform hyperbolicity of {\it one point} in the central direction can determine the existence of ergodic hyperbolic measures with dominated splitting and thus by Theorem \ref{Thm:Horshoe-PosiStaMnfd}, this assumption is a method to prove the existence of hyperbolic periodic orbits.

{\bf Proof.}
Denote by
$V_f (x)$ the set of accumulation measures of time averages
$${\delta (x)}^N=\frac 1{N}{\sum_{j=0}^{N-1}\delta (f^jx)},$$ where $\delta (x)$ denotes the Dirac measure at $x$.
 Then $V_f (x)$ is a nonempty, closed and connected subset of the space of $f$ invariant measures. Recall that dominated splitting is always continuous \cite{BLV}. Thus the center bundle  $E^c$ is continuous and the function $\log\|Df|_{E^c (x)}\|$ is continuous. By assumption,
there exists  (at least) one invariant measure $\nu\in V_f (x_0)$ such that
$$\int\log\|Df|_{E^c (x)}\|d\nu=\liminf_{n\rightarrow +\infty}\frac1n\sum_{i=0}^{n-1}\log{\|Df|_{E^c (f^i (x_0))}\|}<0.$$ By ergodic decomposition theorem, there exists at least one ergodic  measure  $\mu$ such that $$\int\log\|Df|_{E^c (y)}\|d\mu<0.$$ Note that
$\|Df^n|_{E^c (x)}\|\leq\prod_{i=0}^{n-1}\|Df|_{E^c (f^i (x))}\|.$ Then by Subadditive Ergodic Theorem and Birkhoff Ergodic Theorem (see  \cite{Walters}),
 $\mu$ a.e. $x$ satisfies
$$\lim_{n\rightarrow +\infty}\frac1n\log\|Df^n|_{E^c (x)}\|\leq\lim_{n\rightarrow +\infty}\frac1n\sum_{i=0}^{n-1}\log{\|Df|_{E^c (f^i (x))}\|}$$$$=\int\log\|Df|_{E^c (y)}\|d\mu<0.$$
This implies that all Lyapunov exponents of   $\mu$ in the center bundle are  negative and
thus $\mu$ is the required  hyperbolic ergodic measure.\qed

\bigskip

Moreover, we can replace the assumption of above theorem more weaker.
\begin{Thm}\label{Thm:HypMeas-paritialhyper-2222}
Let $f\in \Diff^1 (M)$  be a partially hyperbolic diffeomorphism with a corresponding global $Df$ invariant splitting $T_{M}M=E^s\oplus E^c
\oplus E^u$. If there exists a point $x_0\in M$ and $S\geq 1$ such that at least one following case happens: \\
 (1) $$\liminf_{n\rightarrow +\infty}\frac1n\sum_{i=0}^{n-1}\log{\|Df^S|_{E^c (f^i (x_0))}\|}<0$$$$\text{  (or respectively, }\limsup_{n\rightarrow +\infty}\frac1n\sum_{i=0}^{n-1}\log{m (Df^S|_{E^c (f^i (x_0))})}>0);$$
 (2) $$\liminf_{n\rightarrow +\infty}\frac1n\sum_{i=0}^{n-1}\log{\|Df^S|_{E^c (f^{iS} (x_0))}\|}<0$$$$\text{  (or respectively, }\limsup_{n\rightarrow +\infty}\frac1n\sum_{i=0}^{n-1}\log{m (Df^S|_{E^c (f^{iS} (x_0))})}>0),$$
then there exists at least one  ergodic  measure  $\mu$ such that $\mu$ is hyperbolic and its Oseledec's hyperbolic splitting $T_{x}M=E^s (x)\oplus
E^u (x)$ coincides with $ (E^s (x)\oplus E^c (x))
\oplus E^u (x)$ (or respectively, $E^s (x)\oplus  (E^c (x)
\oplus E^u (x))$) and thus is dominated. This implies $ \m_f^{dh} (M)\cap \m_e (M)\neq \emptyset.$
\end{Thm}

{\bf Proof.}  Replacing  the function $\log\|Df|_{E^c (x)}\|$ by $\log\|Df^S|_{E^c (x)}\|$ in the proof of above theorem, there exists at least one ergodic  measure  $\mu$ such that $$\int\log\|Df^S|_{E^c (x)}\|d\mu<0.$$ We will prove this $\mu$ is needed for case  (1).


 Take $L=S$ and take $K$ large enough such that $$\frac{5L\cdot
C_f}K+\frac1L\int\log\|Df^L|_{E^c (x)}\|d\mu<0.$$ Recall that $\mu$ is $f$-ergodic.  Then by Subadditive Ergodic Theorem and Birkhoff Ergodic Theorem (see  \cite{Walters}),
 $\mu$ a.e. $x$ the following limits exist and by Lemma \ref{Lem-norm-estimate} one has
$$\lim_{n\rightarrow +\infty}\frac1n\log\|Df^n|_{E^c (x)}\|\leq
\frac{5L\cdot
C_f}K+\lim_{l\rightarrow+\infty}\frac1{lL}\sum_{i=0}^{l-1}\log\|Df^{L}|_{E^c (f^{i} (x))}\|$$$$=\frac{5L\cdot
C_f}K+\frac1L\int\log\|Df^L|_{E^c (x)}\|d\mu<0.$$
This implies that all Lyapunov exponents of   $\mu$ in the center bundle are  negative and
thus $\mu$ is the required  hyperbolic ergodic measure. This ends the proof for case  (1) of Theorem \ref{Thm:HypMeas-paritialhyper-2222}.

Now we start to consider the case  (2) of Theorem \ref{Thm:HypMeas-paritialhyper-2222}. Replace $f$ of Theorem \ref{Thm:HypMeas-paritialhyper} by $f^S$. One can have $f^S$-ergodic hyperbolic measure $\nu$ such that the Lyapunov exponents of $\nu$ in the central bundle are all negative.  More precisely, by the proof of Theorem \ref{Thm:HypMeas-paritialhyper}, this $\nu$ also satisfies that there is an $f^S$-invariant set $B$ with $\nu (B)=1$ such that for any $x\in B$,
$$\lim_{n\rightarrow +\infty}\frac1n\log\|Df^{nS}|_{E^c (x)}\| $$$$\leq\lim_{n\rightarrow +\infty}\frac1n\sum_{i=0}^{n-1}\log{\|Df^S|_{E^c (f^{iS} (x))}\|}=\int\log\|Df|_{E^c (y)}\|d\nu<0.$$
Take $\mu=\frac1S (\nu+f_*\nu+\cdots+f^{S-1}_*\nu) $ and $\Gamma=B\cup f B\cdots f^{S-1}B.$ Then $\mu$ is $f$-ergodic   invariant and  $f\Gamma=\Gamma$ with $\mu (\Gamma)=1.$ By sub-additional ergodic theorem, there is $\Gamma'\subseteq \Gamma$ with $\mu (\Gamma')=1$  such that for any $y\in\Gamma'$, the following limit exists   $$\lim_{n\rightarrow +\infty}\frac1n\log\|Df^{n}|_{E^c (y)}\|. $$ Moreover, for any  $y\in\Gamma'$, there is some $0\leq p\leq S-1$
and $x\in B$ such that $f^p (x)=y.$  Note that $$C_f^{-S}\|Df^{n+p}|_{E^c (x)}\|\leq \|Df^n|_{E^c (y)}\|\leq C_f^S\|Df^{n+p}|_{E^c (x)}\|.$$ Thus   $\lim_{n\rightarrow +\infty}\frac1n\log\|Df^{n}|_{E^c (y)}\|$ $$=\lim_{n\rightarrow +\infty}\frac1n\log\|Df^{n}|_{E^c (x)}\|=\frac1S\lim_{n\rightarrow +\infty}\frac1n\log\|Df^{nS}|_{E^c (x)}\| <0.$$
That is, the Lyapunov exponents of ergodic $\mu$ in central direction is negative. We complete the proof.\qed

\bigskip

However, it is still unknown for the following case.
\begin{Que}\label{Que:HypMeas-paritialhyper}
Let $f\in \Diff^1 (M)$  be a partially hyperbolic diffeomorphism with a corresponding global $Df$ invariant splitting $T_{M}M=E^s\oplus E^c
\oplus E^u$. If  there exists a point $x_0\in M$ such that $$\liminf_{n\rightarrow +\infty}\frac1n \log{\|Df^n|_{E^c (x_0)}\|}<0$$$$\text{  (or respectively, }\limsup_{n\rightarrow +\infty}\frac1n \log{m (Df^n|_{E^c ( x_0) })}>0),$$ then whether $ \m_f^{dh} (M)\neq \emptyset?$
\end{Que}

It is true for  $dim (E^c)=1$ (or more general, conformal case), because in this case the condition in Question \ref{Que:HypMeas-paritialhyper} is equivalent to the condition in Theorem \ref{Thm:HypMeas-paritialhyper}.
It is also true when  $E^c$ is a quasi-conformal bundle,  see Theorem  \ref{Thm-partial-hyperbolic-usual}  (1).
However, for non-quasi-conformal case, it is still unknown. Moreover, remark that if the condition in Question \ref{Que:HypMeas-paritialhyper} holds for a set with Lebesgue positive measure  and suppose $f\in C^{1+\alpha}$, it is still unknown whether there is SRB measure?  (see  \cite{ADL2014}).

\bigskip

On the other hand,   replacing the partial hyperbolicity by dominated splitting, it is also unknown whether the average-nonuniform hyperbolicity of a point can permit the existence of hyperbolic ergodic measures with dominated splitting.

\begin{Que}\label{Que:HypMeas-replace-paritialhyper-by Dominated}
Let $f\in \Diff^1 (M)$  be a   diffeomorphism with a  global dominated splitting $T_{M}M=E
\oplus F$. If  there exists a point $x_0\in M$ such that  $$\liminf_{n\rightarrow +\infty}\frac1n\sum_{i=0}^{n-1}\log{\|Df|_{E  (f^i (x_0))}\|}<0 (\text{ or }\limsup_{n\rightarrow +\infty}\frac1n\sum_{i=0}^{n-1}\log{\|Df|_{E  (f^i (x_0))}\|}<0),$$ and  $$  \limsup_{n\rightarrow +\infty}\frac1n\sum_{i=0}^{n-1}\log{m (Df|_{F (f^i (x_0))})}>0 (\text{ or }\liminf_{n\rightarrow +\infty}\frac1n\sum_{i=0}^{n-1}\log{m (Df|_{F (f^i (x_0))})}>0 ) ,$$ then whether $ \m_f^{dh} (M)\neq \emptyset?$
\end{Que}

In the case of $$ \limsup_{n\rightarrow +\infty}\frac1n\sum_{i=0}^{n-1}\log{\|Df|_{E  (f^i (x_0))}\|}<0,\,\,\,\,\,\liminf_{n\rightarrow +\infty}\frac1n\sum_{i=0}^{n-1}\log{m (Df|_{F (f^i (x_0))})}>0,$$
every  $\mu\in V_f (x_0)$ satisfies that  $\int \log\|Df|_{E (x)}\|d\mu<0,\,\int \log\|Df|_{F (x)}\|d\mu>0.$
However, from these it is not sure one can get a hyperbolic measure  (by using Ergodic Decomposition theorem).  It is feasible when $V_f (x_0)\cap\m_e (M)\neq \emptyset$.
\bigskip

Now we start to prove Theorem  \ref{Thm-partial-hyperbolic-usual}.

\medskip

{\bf Proof of Theorem  \ref{Thm-partial-hyperbolic-usual} }       

 (1) Note that for any ergodic measure, there is only one  (different) Lyapunov exponent  in the quasi-conformal central bundles. Thus, if there is an ergodic measure $\mu$ such that all the Lyapunov exponents  in the central bundle are non-zero, then $\mu$ is hyperbolic and its Oseledec hyperbolic splitting is dominated. By Theorem \ref{Thm:Horshoe-PosiStaMnfd}, there exists at least one hyperbolic periodic orbit. So if we face the inverse  case, for any  ergodic measure $\mu$, its Lyapunov exponents  in the central bundle are zero.
Before continuing the proof, we need to recall a known result of  \cite{Herz} .

\begin{Lem}\label{Lem:1} (Proposition 3.4 in  \cite{Herz}) Let $f : X \rightarrow X$ be a continuous map of a compact metric space.
Let $a_n : X \rightarrow R, n \geq 0,$ be a sequence of continuous functions such that
\begin{eqnarray}\label{Lem1Eq1}
a_{n+k} (x) \leq a_n (f^k (x)) + a_k (x) \text{ for every  } x \in X, n, k \geq 0.
\end{eqnarray}
and such that there is a sequence of continuous functions $b_n : X \rightarrow R, n \geq 0,$ satisfying
\begin{eqnarray}\label{Lem1Eq2}
 a_n (x) \leq a_n (f^k (x)) + a_k (x) + b_k (f^n (x)) \text{ for every }  x \in X, n, k \geq 0.
\end{eqnarray}
If $$\inf\frac1n\int_Xa_n (x) d\mu<0$$ for every ergodic $f$-invariant measure, then there is $N > 0$
such that for any $n\geq N,$ $a_n (x) < 0$ for every $x \in X.$

\end{Lem}

Fix $\epsilon>0. $
Define functions for  $x\in M$ $$a_n (x):=\log {\|Df^n|_{E^c (x)}\|}-n\epsilon.$$ Then $\inf\frac1n\int_M a_n (x) d\mu<0$ holds for ergodic invariant measure $\mu$. Recall that $\|AB\|\leq \|A\| \|B\|$. Then it is easy to see that $a_n$ satisfy  (\ref{Lem1Eq1}) of Lemma \ref{Lem:1}. Taking into account  (\ref{Lem1Eq1}) we see that  (\ref{Lem1Eq2}) holds once $a_n (x) \leq a_{n+k} (x) + b_k (f^n (x)).$ This is easily verified
for  $b_k (x):=\log {\| (Df^{k}|_{E^c (x)})^{-1}\|} $
since $$ {\|Df^{n}|_{E^c (x)}\|} \leq {\|Df^{n+k}|_{E^c (x)}\|} \times {\| (Df^{k}|_{E^c (f^n (x))})^{-1}\|}.$$
Recall that  $E^c$ is a continuous splitting so  that $a_n (x),b_n (x)$ are continuous functions.
Then all assumptions of Lemma \ref{Lem:1} are satisfied. So there is $N > 0$
such that $a_N (x) < 0$ for every $x \in M.$ This implies that for any $x \in M,$
$$\|Df^N|_{E^c (x)}\|<e^{N\epsilon}.$$ Equivalently, it is not difficult  to  see that  there is $C^{ (1)}_\epsilon>0$  such that for any $x \in M,n\geq 1$,
$$\|Df^n|_{E^c (x)}\|\leq C^{ (1)}_\epsilon e^{n\epsilon}.$$

Similarly, by defining $ a_n (x):=-\log m (Df^n|_{E^c (x)})+n\epsilon,$ we can get that   there is ${C^{ (2)}_\epsilon}>0$ such that  for any $x \in M,n\geq 1$,  $$m (Df^n|_{E^c (x)}\|\geq \frac 1{C^{ (2)}_\epsilon} e^{-n\epsilon}.$$
Let $C_\epsilon=\max\{C^{ (1)}_\epsilon,\,C^{ (2)}_\epsilon\}$, then we complete the proof of  (1).

 (2)
Let $$\lambda_{max}^c (x)=\limsup_{n\rightarrow \infty}\frac1n\log\|Df^n|_{E^c (x)}\|,\,\lambda_{min}^c (x)=\liminf_{n\rightarrow \infty} \frac1n\log m (Df^n|_{E^c (x)}).$$
By quasi-conformal condition, for any invariant measure $\mu$, $\mu$ a.e. $x$, the above limits exist and $\lambda_{max}^c (x)=\lambda_{min}^c (x).$
 Suppose that the set  $\{x\in M|\,\, \lambda_{max}^c (x)=\lambda_{min}^c (x)=0 \}$ does not have Lebesgue full measure.  We aim to show that  $f$ has a horseshoe.


Firstly we need to prove that
\begin{Lem}\label{lem-srb-integralhyperbolic-positiveentropy}
Suppose that the set  $\{x\in M|\,\, \lambda_{max}^c (x)=\lambda_{min}^c (x)=0 \}$ does not have Lebesgue full measure.  Then
there is some SRB-like measure $\mu$
  such that $$\int \lambda_{max}^c (x)d\mu=\int \lambda_{min}^c (x)d\mu\neq 0.$$
\end{Lem}

{\bf Proof.}  By contradiction, suppose that for any SRB-like measure $\mu$, one has
  \begin{eqnarray}\label{eq-srbLyapu==0}
  \int \lambda_{max}^c (x)d\mu=\int \lambda_{min}^c (x)d\mu= 0.
   \end{eqnarray}
   Before continuing the proof we   want to show following proposition which is also useful for the proof of  Theorem \ref{Thm-partial-hyperbolic-333}. Recall that $\mathcal{O}_{f}$ denotes  the set of all SRB-like measures for $f:M \mapsto M$.


\begin{Prop}\label{Prop-srb-integrablehyp-Lebesguehyp} Let $E$ be a $Df$-invariant continuous bundle and $q\in\mathbb{R}$. Let $$\lambda_{max}^E (x)=\limsup_{n\rightarrow \infty}\frac1n\log\|Df^n|_{E (x)}\|,\,\lambda_{min}^E (x)=\liminf_{n\rightarrow \infty} \frac1n\log m (Df^n|_{E (x)}).$$ If for any SRB-like measure   $\mu\in \mathcal{O}_f$, $$\int\lambda_{max}^E (x)d\mu  <q  (\,\text{ resp.,} \int\lambda_{min}^E (x)d\mu  >q),$$ then there exists $\zeta<q$ (resp., $\zeta>q$), $K_0\geq 1$ and an invariant set $\Theta$ with Lebesgue full measure such that for any $x\in \Theta$ and  any $K\geq K_0,$ $$\limsup_{l\rightarrow
+\infty}\sum_{j=0}^{l-1}\frac{\log\|Df^{K}|_{E ({f^{
j K} (x)})}\|}{lK }\leq \zeta,\,\, (\text{ resp., }\liminf_{l\rightarrow
+\infty}\sum_{j=0}^{l-1}\frac{\log m (Df^{K}|_{E ({f^{
j K} (x)})})}{lK }\geq \zeta).$$

\end{Prop}

{\bf Proof.} In this proposition, we   consider the derivative case, in fact one can state a similar result for sub-additional sequence as Lemma \ref{Lem:1}. Here we just prove the derivative case.


 We will divide the proof into two steps.

{\bf Step 1.} We first prove following lemma.

\begin{Lem}\label{lem-central-srb-Lyap-estimate} There exists $L\geq 1$  and $\epsilon>0$ such that  for 
any $\mu\in \mathcal{O}_f$, $$\int\frac1{L}\log\|Df^{L}|_{E (x)}\|d\mu <q-\epsilon<0.$$
\end{Lem}

{\bf Proof.}
Let's recall some basic  results related with SRB-like measures.
We call basin of attraction $A (\mathcal{K})$ of any nonempty weak$^*$ compact subset $\mathcal{K}$ of probabilities, to $$A (\mathcal{K}):=\{x\in M: p\omega_f (x)\subseteq \mathcal{K}\}.$$  We need a following theorem, which is  a   reformulation of the main results of  \cite{CE}:

\begin{Thm}\label{SRB-like} {\bf  (  \cite{CE})}

 The set $\mathcal{O}_{f}$ of all SRB-like measures for $f$ is the minimal weak$^*$ compact subset of $\mathcal{M} (M)$
 whose basin of  attraction has total Lebesgue measure.

 \em

 In other words: $\mathcal{O}_{f}$ is nonempty and weak$^*$ compact, and the minimal nonempty weak$^*$ compact set that contains, for  Lebegue almost all the initial states $x \in M$, the limits of the convergent subsequences of $\{\frac1n\sum_{j=0}^{n-1}\delta_{f^j (x)}\}_{n\in\mathbb{N}}.$
\end{Thm}

\bigskip

For any $\mu\in \mathcal{O}_f$, by assumption 
 one can choose  some $K (\mu)\geq 1$ and $\epsilon_{\mu}>0$ such that $$\int\frac1{K (\mu)}\log\|Df^{K (\mu)}|_{E (x)}\|d\mu<  q-\epsilon_\mu .$$
By continuity of the bundle $E$, the function $\|Df^{L}|_{E (x)}\|$ is continuous. Thus we can take a neighborhood $B (\mu)\subseteq \m_f (M)$ such that for any $\nu\in B (\mu),$   $$\int\frac1{K (\mu)}\log\|Df^{K (\mu)}|_{E (x)}\|d\nu< q -\epsilon_\mu.$$
Since $\mathcal{O}_f$ is compact, one can take a finite $\mu_1,\cdots,\mu_m\in \mathcal{O}_f$ such that  $\mathcal{O}_f\subseteq \cup_{i=1}^m B (\mu_i)$. Take $L=\Pi_{i=1}^m K (\mu_i) $ and $\epsilon=\min\{\epsilon_{\mu_i}\}$. Then  for any $\mu\in \mathcal{O}_f$, there is $1\leq i\leq m$ such that  $\mu\in B (\mu_i)$ and  so by the sub-multiplication of norms and invariance of $\mu$, $$\int\frac1{L}\log\|Df^{L}|_{E (x)}\|d\mu\leq \int\frac1{L}\sum_{j=0}^{\frac{L}{K (\mu_i)}}\log\|Df^{K (\mu_i)}|_{E (f^{jK (\mu_i)}x)} \|d\mu$$$$=\int\frac1{K (\mu_i)}\log\|Df^{K (\mu_i)}|_{E (x)}\|d\mu< q -\epsilon .$$ \qed


 {\bf Step 2. Complete the proof of Proposition \ref{Prop-srb-integrablehyp-Lebesguehyp}.}

 By Lemma \ref{SRB-like}, for Lebesgue a.e. $x$, $pw_f (x)\subseteq \mathcal{O}_f.$ Fix such a point $x$.  By weak$^*$ topology and the continuity of $\|Df^{L}|_{E (y)}\|$, we have
$$ \limsup_{l\rightarrow
+\infty}\sum_{j=0}^{l-1}\frac{\log\|Df^{L}|_{E ({f^{
j } (x)})}\|}{lL}\leq\sup_{\mu\in \mathcal{O}_f } \int\frac1{L}\log\|Df^{L}|_{E (y)}\|d\mu< q -\epsilon.$$

Let     $K\geq 2L. $ 
Then  by Lemma \ref{Lem-norm-estimate}   (2'),
$$\limsup_{l\rightarrow+\infty}\frac1{lK}\sum_{i=0}^{l-1}\log \|Df^K|_{E (f^{iK} (x))}\|\leq   \frac{5L\cdot
C_f}K+\limsup_{l\rightarrow+\infty}\frac1{lL}\sum_{i=0}^{l-1}\log\|Df^{L}|_{E (f^{i} (x))}\|.$$
   Taking
$K_{0}\geq \max\{2L, \frac{10L\cdot C_f}{\varepsilon}\}$ and then for all $K\geq K_{0}$,
$$\limsup_{l\rightarrow+\infty}\frac1{lK}\sum_{i=0}^{l-1}\log\|Df^K|_{E (f^{iK} (x))}\|<q +
\frac12\epsilon-\epsilon=q-\frac12\epsilon.$$ Take $\zeta=q-\frac12\epsilon$ we complete the proof.\hfill $\Box$

Now we continue our proof.  Let $E$ of Proposition \ref{Prop-srb-integrablehyp-Lebesguehyp} be the central bundle $E^c$, $q>0$. Then by  (\ref{eq-srbLyapu==0}) the assumptions in Proposition \ref{Prop-srb-integrablehyp-Lebesguehyp} hold and then there is $\zeta<q,$ $K_0\geq 1$ and an invariant set $\Theta$ with Lebesgue full measure such that for any $x\in \Theta$ and  any $K\geq K_0,$ $$\limsup_{l\rightarrow
+\infty}\sum_{j=0}^{l-1}\frac{\log\|Df^{K}|_{E ({f^{
j K} (x)})}\|}{lK }\leq \zeta.$$ This implies that $$\lambda_{max}^c (x)= \limsup_{n\rightarrow
+\infty} \frac{\log\|Df^{n}|_{E ( x)}\|}{n }\leq \limsup_{l\rightarrow
+\infty}\sum_{j=0}^{l-1}\frac{\log\|Df^{K}|_{E ({f^{
j K} (x)})}\|}{lK }\leq\zeta< q.$$ By arbitrariness of $q$,  for any $x\in \Theta$,  $$\lambda_{max}^c (x)\leq 0.$$

Similarly, letting $E$ of Proposition \ref{Prop-srb-integrablehyp-Lebesguehyp} be the central bundle $E^c$, $q<0$ and then by  (\ref{eq-srbLyapu==0}) one has that Lebesgue a.e. $x$, $$\lambda_{min}^c (x)\geq 0.$$
This contradicts the assumption that the set  $\{x\in M|\,\, \lambda_{max}^c (x)=\lambda_{min}^c (x)=0 \}$ does not have Lebesgue full measure.
This ends the proof of Lemma \ref{lem-srb-integralhyperbolic-positiveentropy}.\qed

To complete the proof of Theorem  \ref{Thm-partial-hyperbolic-usual}  (2), by Lemma \ref{lem-srb-integralhyperbolic-positiveentropy} we only need to show that

\begin{Lem}\label{lem-srb-imply-horseshoe}
Suppose that  there is some SRB-like measure $\mu$
  such that $$\int \lambda_{max}^c (x)d\mu=\int \lambda_{min}^c (x)d\mu\neq 0.$$ Then $f$ has a horseshoe. \\  More precisely, we get a horseshoe  with large entropy as follows: from uniform expanding of $E^u,$ we can assume that there exist $C>0$ and $0<\lambda<1$ such that $$ {\|Df^{-n}|_{E^u (x)}\|} \leq C \lambda^n, \forall x\in M,\,\, n\geq 1.$$ Then for any $\omega>0,$ there is a horseshoe with positive topological entropy larger than $-dim (E^u)\log\lambda-\omega.$
\end{Lem}

{\bf Proof.}
  Recall that $\mu$ a.e. $x$, $\lambda_{max}^c (x) =  \lambda_{min}^c (x)$, denoted by $\lambda^c (x)$ for convenience. Then $\int \lambda^c (x)d\mu\neq 0.$ Let $\Delta_1=\{x|\lambda^c (x)\text{ exists and }\neq 0\}$,  then $t:=\mu (\Delta_1)>0.$ Define
   $\mu_1=\mu_{\Delta_1},$ we will show that $h_{\mu_1} (f)>0.$ Recall that the system is partially hyperbolic, by Theorem \ref{PesFormula-Thm:2} we have   $$h_\mu (f)\geq \int \sum_{i=1}^{dim\,E^u (x)}\lambda_i (x) d\mu>0,$$ where $\lambda_1 (x)\geq \lambda_2 (x)\geq\cdots \geq \lambda_{dim (M)} (x)$ denote the Lyapunov exponents of   $x$.

   If $t=1$, then $h_{\mu_1} (f)=h_{\mu} (f)\geq \int \sum_{i=1}^{dim\,E^u (x)}\lambda_i (x) d\mu=  \int \sum_{i=1}^{dim\,E^u (x)}\lambda_i (x) d\mu_1\geq -dim (E^u)\log\lambda>0.$ Otherwise,
      $\mu\neq \mu_1$ and $t= \mu_1 (\Delta_1)< 1.$ Define  $\Delta_2=\{x|\lambda^c (x)\text{ exists and }= 0\}$ and
   $ \mu_2=\mu_{\Delta_2}.$ By Ruelle's inequality and the definition of $\Delta_2$ and $\mu_2$, $h_{\mu_2} (f)\leq \int \sum_{\lambda_i (x)>0}\lambda_i (x) d\mu_2= \int \sum_{i=1}^{dim\,E^u (x)}\lambda_i (x) d\mu_2.$ By affine property of metric entropy, $$h_\mu (f)=th_{\mu_1} (f)+ (1-t)h_{\mu_2} (f) \leq th_{\mu_1} (f)+  (1-t)\int \sum_{i=1}^{dim\,E^u (x)}\lambda_i (x) d\mu_2.$$ However, from above analysis we know that $$ h_\mu (f)\geq \int \sum_{i=1}^{dim\,E^u (x)}\lambda_i (x) d\mu=t\int \sum_{i=1}^{dim\,E^u (x)}\lambda_i (x) d\mu_1 + (1-t)\int \sum_{i=1}^{dim\,E^u (x)}\lambda_i (x) d\mu_2.$$ So $$h_{\mu_1} (f)\geq \int \sum_{i=1}^{dim\,E^u (x)}\lambda_i (x) d\mu_1\geq -dim (E^u)\log\lambda>0.$$

   By Ergodic Decomposition theorem of metric entropy, there is at least one ergodic component of $\mu_1$, denoted by $\nu,$ has positive entropy $>\max\{-dim (E^u)\log\lambda-\frac12\omega,\,0\}$. By definition of $\Delta_1$ and $\mu_1$, the ergodic $\nu$ can be taken with $\nu (\Delta_1)=1.$ This $\nu$ is a hyperbolic measure with positive metric entropy and its Oseledec hyperbolic splitting is dominated. In other words, $\nu\in \mathcal{M}^{dh}_{f} (M)$ and is non-atomic. By Theorem \ref{Thm:Horshoe-Approximation}  (2), $f$ has a horseshoe with topological entropy larger than $h_\nu (f)-\frac12\omega>-dim (E^u)\log\lambda-\omega.$
  \qed

\bigskip

By Lemma \ref{lem-srb-integralhyperbolic-positiveentropy} and Lemma \ref{lem-srb-imply-horseshoe}, we state Theorem  \ref{Thm-partial-hyperbolic-usual}  (2) in another way.

\begin{Thm}\label{Thm-Que:HypMeas-paritialhyper-2222222}
Let $f\in \Diff^1 (M)$  be a partially hyperbolic diffeomorphism with a corresponding global $Df$ invariant splitting $T_{M}M=E^s\oplus E^c
\oplus E^u$. Suppose that  $E^c$  is quasi-conformal. If  there exists a  Lebesgue positive measure set $H$ such that each point $x_0\in H$ satisfies that $$\liminf_{n\rightarrow +\infty}\frac1n \log{\|Df^n|_{E^c (x_0)}\|} (=\liminf_{n\rightarrow +\infty}\frac1n \log{m (Df^n|_{E^c ( x_0) })})\neq 0$$$$\text{ or respectively, }\limsup_{n\rightarrow +\infty}\frac1n \log{\|Df^n|_{E^c (x_0)}\|} (=\limsup_{n\rightarrow +\infty}\frac1n \log{m (Df^n|_{E^c ( x_0) })})\neq 0,$$ then  $ \m_f^{dh} (M)\cap\m_f^n (M)\neq \emptyset$ and then  by Theorem \ref{Thm:Horshoe-PosiStaMnfd} there is horseshoe.
\end{Thm}


Furthermore, we want to ask   a following question for  the case of non-conformal central bundle.

\begin{Que}\label{Que:HypMeas-paritialhyper-2222222}
Let $f\in \Diff^1 (M)$  be a partially hyperbolic diffeomorphism with a corresponding global $Df$ invariant splitting $T_{M}M=E^s\oplus E^c
\oplus E^u$. If  there exists a  Lebesgue positive measure set $H$ such that each point $x_0\in H$ satisfies that $$\liminf_{n\rightarrow +\infty}\frac1n \log{\|Df^n|_{E^c (x_0)}\|}<0\,\text{  (or  }\limsup_{n\rightarrow +\infty}\frac1n \log{\|Df^n|_{E^c ( x_0) }\|}<0),$$$$\text{ or }\limsup_{n\rightarrow +\infty}\frac1n \log{m (Df^n|_{E^c ( x_0) })}>0\,\text{  (or  }\liminf_{n\rightarrow +\infty}\frac1n \log{m (Df^n|_{E^c ( x_0) })}>0),$$ then whether $ \m_f^{dh} (M)\cap\m_f^n (M)\neq \emptyset$  (which implies the existence of horseshoe by Theorem \ref{Thm:Horshoe-PosiStaMnfd})?

In other words, whether the system satisfies that  either there is some horseshoe, or Lebesgue a.e. $x$,  $$\liminf_{n\rightarrow +\infty}\frac1n \log{\|Df^n|_{E^c (x)}\|}\geq 0 \geq  \limsup_{n\rightarrow +\infty}\frac1n \log{m (Df^n|_{E^c ( x) })}?$$
\end{Que}

It is easy to see $H$ can be extended to the set of  $\Delta:=\cup_{n\in \mathbb{Z}} f^n H$. So if Lebesgue measure is invariant (i.e., the system is volume-preserving), then
$\mu=Leb|_\Delta$ is one  required measure.

For the non-conformal case, it is also unknown for the cases replacing the above assumption by
$$\liminf_{n\rightarrow +\infty}\frac1n\sum_{i=0}^{n-1}\log{\|Df|_{E  (f^i (x_0))}\|}<0 (\text{ or }\limsup_{n\rightarrow +\infty}\frac1n\sum_{i=0}^{n-1}\log{\|Df|_{E (f^i (x_0))}\|}<0),$$$$\text{ or }\,\,\,\,\,\, \limsup_{n\rightarrow +\infty}\frac1n\sum_{i=0}^{n-1}\log{m (Df|_{F (f^i (x_0))})}>0 (\text{ or }\liminf_{n\rightarrow +\infty}\frac1n\sum_{i=0}^{n-1}\log{m (Df|_{F (f^i (x_0))})}>0 ) .$$ By Theorem \ref{Thm:HypMeas-paritialhyper}, we only know that $ \m_f^{dh} (M) \neq \emptyset$.
 We   give a partial answer as follows. We say a continuous $Df-$invariant bundle $G\subseteq TM$ to be {\it non-contracting} ( resp., {\it non-expanding}), if for any $x\in M$,
 $$\limsup_{n\rightarrow +\infty}\frac1n \log |det  (Df^n|_{G (x)})|\geq 0\,\,\text{ ( resp., }\liminf_{n\rightarrow +\infty}\frac1n \log |det  (Df^n|_{G (x)})|\leq 0).$$

 \begin{Thm}\label{Thm-Que:HypMeas-paritialhyper-2222222}
Let $f\in \Diff^1 (M)$  be a partially hyperbolic diffeomorphism with a corresponding global $Df$ invariant splitting $T_{M}M=E^s\oplus E^c
\oplus E^u$.  Suppose that If $E^c\oplus E^u$ is non-contracting  (resp., $E^s\oplus E^c$ is non-expanding), there exists a  Lebesgue positive measure set $H$ such that each point $x_0\in H$ satisfies that $$\liminf_{n\rightarrow +\infty}\frac1n\sum_{i=0}^{n-1}\log{\|Df|_{E  (f^i (x_0))}\|}<0 \,\text{  ( resp.,  }\limsup_{n\rightarrow +\infty}\frac1n\sum_{i=0}^{n-1}\log{m (Df|_{F (f^i (x_0))})}>0),$$ then   $ \m_f^{dh} (M)\cap\m_f^n (M)\neq \emptyset$ and then there is horseshoe by Theorem \ref{Thm:Horshoe-PosiStaMnfd}.
\end{Thm}

{\bf Proof.} Let $$\lambda_{max}^c (x)=\limsup_{n\rightarrow \infty}\frac1n\log\|Df^n|_{E^c (x)}\|,\,\lambda_{min}^c (x)=\liminf_{n\rightarrow \infty} \frac1n\log m (Df^n|_{E^c (x)}).$$
For any invariant measure $\mu$, $\mu$ a.e. $x$, the above limits exist and $\lambda_{max}^c (x)\geq\lambda_{min}^c (x).$   By assumption and Lemma \ref{SRB-like}, there is $x_0\in H$ such that $pw_f (x_0)\subseteq \mathcal{O}_f$ and
$$\liminf_{n\rightarrow +\infty}\frac1n\sum_{i=0}^{n-1}\log{\|Df|_{E  (f^i (x_0))}\|}<0.$$
It follows that there is at least one SRB-like measure $\mu\in pw_f (x_0)\subseteq \mathcal{O}_f$ such that $$ \int \log{\|Df|_{E^c (x)}\|}d\mu=\liminf_{n\rightarrow +\infty}\frac1n\sum_{i=0}^{n-1}\log{\|Df|_{E  (f^i (x_0))}\|}<0.$$ It follows that $$ \int \lambda_{max}^c (x)d\mu\leq \int \log{\|Df|_{E^c (x)}\|}d\mu <0.$$

 Let $\Delta_1=\{x|\lambda^c_{max} (x)\text{ exists and }< 0\}$,  then $t:=\mu (\Delta_1)>0.$ Define
   $\mu_1=\mu_{\Delta_1},$ we will show that $h_{\mu_1} (f)>0.$ Recall that the system is partially hyperbolic, by Theorem \ref{PesFormula-Thm:2} we have   $$h_\mu (f)\geq \int \sum_{i=1}^{dim\,E^u (x)}\lambda_i (x) d\mu>0,$$ where $\lambda_1 (x)\geq \lambda_2 (x)\geq\cdots \geq \lambda_{dim (M)} (x)$ denote the Lyapunov exponents of   $x$.

   If $t=1$, then $h_{\mu_1} (f)=h_{\mu} (f)\geq \int \sum_{i=1}^{dim\,E^u (x)}\lambda_i (x) d\mu=  \int \sum_{i=1}^{dim\,E^u (x)}\lambda_i (x) d\mu_1\geq -dim (E^u)\log\lambda>0.$ Otherwise,
      $\mu\neq \mu_1$ and $t= \mu_1 (\Delta_1)< 1.$ Define  $\Delta_2=\{x|\lambda^c_{max} (x)\text{ exists and }\geq  0\}$ and
   $ \mu_2=\mu_{\Delta_2}.$ By Ruelle's inequality and the definition of $\Delta_2$ and $\mu_2$, $$h_{\mu_2} (f)\leq \int \sum_{\lambda_i (x)>0}\lambda_i (x) d\mu_2 \leq  \int \sum_{i=1}^{dim\,E^u (x)}\lambda_i (x) d\mu_2+dim (E^c)\int \lambda_{max}^c (x)d\mu_2.$$ By affine property of metric entropy,
   \begin{eqnarray*}& &h_\mu (f)=th_{\mu_1} (f)+ (1-t)h_{\mu_2} (f)
   \\ &\leq& th_{\mu_1} (f)+  (1-t) (\int \sum_{i=1}^{dim\,E^u (x)}\lambda_i (x) d\mu_2+dim (E^c)\int \lambda_{max}^c (x)d\mu_2)\\& =&th_{\mu_1} (f)+  (1-t)\int \sum_{i=1}^{dim\,E^u (x)}\lambda_i (x) d\mu_2 \\&+&dim (E^c)\int \lambda_{max}^c (x)d\mu-tdim (E^c)\int \lambda_{max}^c (x)d\mu_1.\\& <&th_{\mu_1} (f)+  (1-t)\int \sum_{i=1}^{dim\,E^u (x)}\lambda_i (x) d\mu_2 -tdim (E^c)\int \lambda_{max}^c (x)d\mu_1. \end{eqnarray*}
    However, from above analysis we know that $$ h_\mu (f)\geq \int \sum_{i=1}^{dim\,E^u (x)}\lambda_i (x) d\mu=t\int \sum_{i=1}^{dim\,E^u (x)}\lambda_i (x) d\mu_1 + (1-t)\int \sum_{i=1}^{dim\,E^u (x)}\lambda_i (x) d\mu_2.$$ So $$h_{\mu_1} (f)>\int ( \sum_{i=1}^{dim\,E^u (x)}\lambda_i (x)+dim (E^c)  \lambda_{max}^c (x)) d\mu_1 $$$$\geq  \int \limsup_{n\rightarrow +\infty}\frac1n \log |det  (Df^n|_{E^c (x)\oplus E^u (x)}) d\mu_1 \geq 0.$$

   By Ergodic Decomposition theorem of metric entropy, there is at least one ergodic component of $\mu_1$, denoted by $\nu,$ has positive entropy. By definition of $\Delta_1$ and $\mu_1$, the ergodic $\nu$ can be taken with $\nu (\Delta_1)=1.$ This $\nu$ is a hyperbolic measure with positive metric entropy and its Oseledec hyperbolic splitting is dominated. In other words, $\nu\in \mathcal{M}^{dh}_{f} (M)$ and is non-atomic. By   Theorem \ref{Thm:Horshoe-PosiStaMnfd} or Theorem \ref{Thm:Horshoe-Approximation}  (2), $f$ has a horseshoe with positive topological entropy.
  \qed

\subsection{Some Partially Hyperbolic Systems}\label{section-partialhyper}

Following theorem is to show the existence of hyperbolic measures with dominated splitting in some partially hyperbolic systems.

\begin{Thm}\label{Thm-SUBSEction-partial-hyperbolic-22222}

 (1) Under the assumption  of Theorem \ref{Thm-partial-hyperbolic}, every ergodic measure with positive metric entropy is hyperbolic with a dominated Oseledec hyperbolic splitting.  In other words,
$$ \mathcal{M}_{e} (M)\cap\m_f^+ (M)=\mathcal{M}_{e} (M)\cap\mathcal{M}^{dh}_{f} (M).$$
 (2) Under the assumption  of Theorem \ref{Thm-partial-hyperbolic-22222}, every ergodic measure $\mu$ with
metric entropy satisfying  $$h_\mu (f)  > {dimE}\cdot a, $$
 is hyperbolic with a dominated Oseledec hyperbolic splitting.

\end{Thm}

{\bf Proof.}  (1) We only need to consider the cases  (C) and  (C'), since others are particular cases  of   (2) for $a=0$.   Note that  (C') is a particular case of  (C) so that we only need to prove  (C).

Let $\mu$ be an ergodic measure with positive metric entropy. By Ruelle's inequality \cite{Ru} and quasi-conformal condition, $\mu$ a.e. $x$, the Lyapunov exponents of $E$ are negative and the Lyapunov exponents of $F$ are positive. This implies that   the Oseledec hyperbolic splitting of $\mu$ coincides with the given dominated splitting.

\medskip

 (2) Let $$\lambda_{max}^E (x)=\limsup_{n\rightarrow \infty}\frac1n\log\|Df^n|_{E (x)}\|,\,\lambda_{min}^E (x)=\liminf_{n\rightarrow \infty} \frac1n\log m (Df^n|_{E (x)}).$$ For $\mu$ a.e. $x\in M,$ the $\limsup,\liminf$ can be written by $\lim$ so that by assumption for $\mu$ a.e. $x\in M,$  $$ \lambda_{max}^E (x)\leq \lambda_{min}^E (x)+a .$$ By Ruelle's inequality \cite{Ru} for $f^{-1}$ and ergodicity of $\mu$,  $\mu$ a.e.  $x$,  $$ dimE\cdot a  < h_\mu (f)=h_\mu (f^{-1})\leq dim E  (-\lambda_{min}^E (x))\leq dim E (-\lambda_{max}^E (x)+a).$$ This implies $\mu$ a.e $x$,
$$\lambda_{max}^E (x)  <0.$$
So the Lyapunov exponents of $\mu$ in the bundle $E$ are all negative so that the Oseledec hyperbolic splitting of $\mu$ coincides with the given dominated splitting.\qed

\bigskip

Now we will prove Theorem \ref{Thm-partial-hyperbolic}, Theorem \ref{Thm-partial-hyperbolic-22222} and Theorem \ref{Thm-partial-hyperbolic-333}.  For the last two theorems, we only need to consider the system   (I)  of Theorem \ref{Thm-partial-hyperbolic-22222}, since the system  (II) is similar (just needing to consider   $f^{-1}$).

\bigskip

{\bf Proof of Theorem \ref{Thm-partial-hyperbolic}.}

We only need to consider the cases  (C) and  (C'), since others are particular cases  of  Theorem \ref{Thm-partial-hyperbolic-22222} for $a=0$.   Note that  (C') is a particular case of  (C) so that we only need to prove  (C).

By Variational Principle, there exists  ergodic $\mu$ with positive entropy arbitrarily close the topological entropy. By Theorem \ref{Thm-SUBSEction-partial-hyperbolic-22222} $\mu$ is a hyperbolic ergodic measure with domination. By Theorem \ref{Thm:Horshoe-Approximation}, one can get  (1) and  (2). For  (3), they are obvious from   Corollary \ref{Cor:1Horshoe-numbergrowth-periodic}. \qed

\bigskip

Remark that   the  main idea  for  proof  is to find hyperbolic measures with domination which have metric entropy arbitrarily close to the topological entropy. Similar observation will appear  in the proof of Theorem \ref{Thm-partial-hyperbolic-22222}.

\bigskip

{\bf Proof of  (I)   in Theorem \ref{Thm-partial-hyperbolic-22222}.}
Similar as the proof of  (C) in Theorem \ref{Thm-partial-hyperbolic}, we only need to construct an ergodic $\mu$
such that its metric entropy is arbitrarily close to the topological entropy and $\mu$ is hyperbolic with domination. From a recent a result of  \cite{c-tian} that for a $C^1$ partial hyperbolic system, if $F$ is uniformly expanding, i.e.  there exist $C>0$ and $0<\lambda<1$ such that $$ {\|Df^{-n}|_{F (x)}\|} \leq C \lambda^n, \forall x\in M,\,\, n\geq 1,$$ then
\begin{Lem}\label{lem-c-tian} (  \cite{c-tian} ) Every SRB-like measure for $f$ has   positive entropy larger or equal to $\mbox{\em dim} (F)\log \lambda^{-1}$. In particular, $$h_{top} (f)\geq \mbox{\em dim} (F)\log \lambda^{-1}.$$
\end{Lem}
Fix $\epsilon\in (0,-\frac{dimF}{dimE}\log\lambda-a).$ By Variational Principle, there exists  ergodic $\mu$ such that $$h_\mu (f)>h_{top} (f)-\epsilon\geq -dimF \log\lambda-\epsilon>dim E\cdot a.$$  By Theorem \ref{Thm-SUBSEction-partial-hyperbolic-22222} $\mu$ is a hyperbolic ergodic measure with domination. By Theorem \ref{Thm:Horshoe-Approximation}, one can get  (1) and  (2). For  (3), it is  obvious from  Corollary \ref{Cor:1Horshoe-numbergrowth-periodic}.  We complete the proof. \qed

\bigskip
{\bf Proof of   Theorem \ref{Thm-partial-hyperbolic-zero-entropy-measure-dense}}
Firstly we give a simple lemma.
\begin{Lem}\label{lem-zeroentropy2015Oct} If every ergodic measure can be approximated by invariant measures with zero entropy, then $\m_z (M)$ is dense in $\m_f (M)$. If further the entropy map $\mu\mapsto h_\mu (f)$ is upper semi-continuous, $\m_z (M)$ is residual in $\m_f (M)$.
\end{Lem}
{\bf Proof} By Ergodic Decomposition Theorem, every invariant measure $\mu$ can be approximated by convex sum of ergodic measures and so $\mu$ also can be approximated by convex sum of invariant measures with zero entropy. Note  that the set of all convex sum of invariant measures with zero entropy   are contained in $\m_z (M)$ and thus $\m_z (M)$ is dense in $\m_f (M)$.

If  the entropy map $\mu\mapsto h_\mu (f)$ is upper semi-continuous, then the sets $$U_n=\{\nu\in \m_f (M)|\,h_\nu (f)<\frac1n\}$$ are all open subsets of $\m_f (M).$ It is obvious that $\m_z (M)=\cap _{n\geq 1}U_n.$ Thus $\m_z (M)$ is residual in $\m_f (M)$.\qed

Under the assumption  of Theorem \ref{Thm-partial-hyperbolic}, by Theorem \ref{Thm-SUBSEction-partial-hyperbolic-22222} every ergodic measure with positive metric entropy is hyperbolic with a dominated Oseledec hyperbolic splitting.    By Theorem \ref{Thm:density-per-meas22222}, every ergodic measure with positive metric entropy can be approximated by periodic measures. This implies that every ergodic measure can be approximated by invariant measures with zero metric entropy. By lemma \ref{lem-zeroentropy2015Oct}, $\m_z (M)$ is dense in $\m_f (M)$.


 In particular, for each case of  (A')  (B')  (C'), $f$ is far from homoclinic tangency and by  \cite{LiaoVianaYang}  $f$ is entropy-expansive so that the entropy map $\mu\mapsto h_\mu (f)$ is upper semi-continuous.  By lemma \ref{lem-zeroentropy2015Oct}, $\m_z (M)$ is residual in $\m_f (M)$. Now we complete the proof of Theorem \ref{Thm-partial-hyperbolic-zero-entropy-measure-dense}. \qed

\bigskip

{\bf Proof of   Theorem \ref{Thm-partial-hyperbolic-333} for system   (I)   in Theorem \ref{Thm-partial-hyperbolic-22222}.}

For the direction of $F$, there exist $C>0$ and $0<\lambda<1$ such that $$ {\|Df^{-n}|_{F (x)}\|} \leq C \lambda^n, \forall x\in M,\,\, n\geq 1.$$ Fix $\bar{\lambda}\in (\lambda,1).$ Take $K_1\geq 1$ large enough such that  for any $n\geq K_1$
$$C\lambda^n\leq \bar{\lambda}^n.$$ Thus, for any $x\in M$ and any $K\geq K_1,$ $${\|Df^{-K}|_{F (x)}\|} \leq\bar{\lambda}^K. $$ This implies that  (Lebesgue) every $x\in M$ and any $K\geq K_1,$
$$  \liminf_{l\rightarrow +\infty}\sum_{j=-l}^{-1}\frac{\log
m (Df^{K}|_{F ({f^{ \pm jK} (x)})})}{lK}$$$$=-\limsup_{n\rightarrow \infty} \frac1{nK} \sum_{j=0}^{n-1} \log {\|Df^{-K}|_{F (f^{\pm jK}x)}\|} \geq -\log \bar{\lambda}>0. $$
So we only needs to prove the direction $E$.

By Lemma \ref{lem-c-tian}, every SRB-like measure $\mu$ for $f$ has   positive entropy larger or equal to $\mbox{\em dim} (F)\log \lambda^{-1}$. Let $$\lambda_{max}^E (x)=\limsup_{n\rightarrow \infty}\frac1n\log\|Df^n|_{E (x)}\|,\,\lambda_{min}^E (x)=\liminf_{n\rightarrow \infty} \frac1n\log m (Df^n|_{E (x)}).$$ For $\mu$ a.e. $x\in M,$ the $\limsup,\liminf$ can be written by $\lim$ so that by assumption for $\mu$ a.e. $x\in M,$  $$ \lambda_{max}^E (x)\leq \lambda_{min}^E (x)+a .$$ By invariance of $\mu$ and Ruelle's inequality \cite{Ru} for $f^{-1}$,  $$-dimF \log\lambda \leq h_\mu (f)=h_\mu (f^{-1})$$
$$\leq \int dim E  (-\lambda_{min}^E (x))d\mu\leq \int dim E (-\lambda_{max}^E (x)+a)d\mu.$$ This implies
\begin{eqnarray}\label{eq-srbLyapunov}
\int\lambda_{max}^E (x)d\mu \leq \frac{dimF}{dimE}\log\lambda +a <0.
\end{eqnarray}

By  (\ref{eq-srbLyapunov}),  the result (1) of Theorem \ref{Thm-partial-hyperbolic-333} can be deduced from Proposition \ref{Prop-srb-integrablehyp-Lebesguehyp} for $q=0$.  Combining Theorem \ref{Thm:2015-quasi-invariant} with (1), the result (2) is obtained. \qed

\medskip

 Theorem \ref{Thm-partial-hyperbolic-333} studies quasi-invariant measure with average-nonuniform  hyperbolicity. However,   if average-nonuniform  hyperbolicity is replaced by nonuniform hyperbolicity, it is still unknown whether there is similar result.  The reason is that we do not know whether there is some Liao-Pesin set with full measure for a quasi-invariant measure. More precisely, we state a following question.

\begin{Que} Let $f\in \Diff^1 (M)$ and  $TM=E\oplus F$ be a $Df-$invariant splitting on
$M$.  Suppose that   $TM=E\oplus F$
is  limit-dominated on $M$ and $\mu$   is quasi-invariant (or quasi-ergodic). If $\mu$
    is   nonuniformly hyperbolic with respect to this splitting, that is, there is a $\mu$ full measure set $H$ such that each point $x_0\in H$ satisfies that $$\limsup_{n\rightarrow +\infty}\frac1n \log{\|Df^n|_{E (x_0)}\|}<0,$$$$\text{ and }\liminf_{n\rightarrow +\infty}\frac1n \log{m (Df^n|_{F ( x_0) })}>0.$$
 Then whether there is some Liao-Pesin set $\Lambda (K,\zeta) $ with $\mu$ positive measure  (or even full measure)?

\end{Que}

In this question, the bundles $E$ and $F$ are general bundles, not required one-dimensional or quasi-conformal. The main obstruction is that it is not known whether for quasi-invariant measure, nonuniform hyperbolicity implies average-nonuniform  hyperbolicity, as proved in Lemma \ref{Lem:Generalized-Multi-erg-thm}.

\bigskip


{\bf Proof of Theorem \ref{Thm-partial-hyperbolic-usual-Lebegue-hyp} for
average-nonuniform hyperbolicity and density of periodic points.}
 Let $E=E^s\oplus E^c$ and $F= E^u$. Then by uniform expanding of $E^u,$ $F$ is a average-nonuniformly expanding bundle. So we only need to consider $E$.

Let $$A=\{x|\,pw (x) \text{ contains only one measure which is  ergodic   with positive entropy}\}$$ and $$B=\{x|\,pw (x) \text{ contains only one measure which is  ergodic   with zero  entropy}\}.$$ They are two disjoint invariant subsets of $M$. By Birkhoff Ergodic theorem for any invariant measure $A\cup B$ has full measure and moreover, for every ergodic measure $\mu$ with positive entropy, $\mu (A)=1.$ By   Ergodic Decomposition theorem and assumption,  every invariant measure $\nu$ supported on $A$ has positive metric entropy and for $\nu$ a.e. point, its Lyapunov exponents of central direction are all negative.  However, every invariant measure supported on $B$ has zero metric entropy.

 Let $\mu$ be a SRB-like measure,   by Lemma \ref{lem-c-tian}, $h_\mu (f)>0.$ We will prove $\mu (A)=1,\mu (B)=0.$ That is,
 \begin{Lem}\label{lem-srb-meas-positiveentropysigular-zeroentropy}
 If there is $\theta\in[0,1]$ and two invariant measures $\mu_1$ and $\mu_2$ such that
$\mu=\theta\mu_1+ (1-\theta)\mu_2$ and $\mu_1 (A)=1,\,\,h_{\mu_2} (f)=0$, then $\theta=1.$
 \end{Lem}

 {\bf Proof.} Othewise,  $\theta< 1$. By Theorem \ref{PesFormula-Thm:2},
$$h_{\mu} (f)\geq\int \chi (x)d\mu,$$ where
$\chi (x)=\sum_{i=1}^{dim\,F (x)}\lambda_i (x)$ and $\lambda_1 (x)\geq\lambda_2 (x)\geq\cdots\geq\lambda_{dim\,M} (x)$ denote
the Lyapunov exponents of $\mu$ a.e. $x.$  Then by the uniform expanding of $F$,
$$h_{\mu} (f)\geq\int \chi (x)d\mu>\theta\int \chi (x)d\mu_1.$$

On the other hand, $\mu_1 (A)=1$ implies that for $\mu_1$ a.e. point, its Lyapunov exponents of central direction are all negative. By Ruelle's inequality \cite{Ru} and the assumption, $$h_{\mu} (f)=\theta h_{\mu_1} (f)\leq \theta \int \sum_{\lambda_i (x)>0}\lambda_i (x)d\mu_1=\theta\int \chi (x)d\mu_1.$$ This is a contradiction. \qed

Now we know that for every SRB-like measure, a.e. points only have negative Lyapunov exponents in $E^c$ and then by domination $E^s\oplus_\prec E^c$, a.e. points only have negative Lyapunov exponents in $E^c$. By Proposition \ref{Prop-srb-integrablehyp-Lebesguehyp}, we complete the proof of average-nonuniform hyperbolicity.

Combining Theorem \ref{Thm:2015-quasi-invariant} with average-nonuniform hyperbolicity, if Lebesgue is quasi-invariant, then the periodic points form a dense subset in the whole space.\qed


\medskip


{\bf Proof of Theorem \ref{Thm-partial-hyperbolic-usual-2}.}
For the consequences of Theorem \ref{Thm-partial-hyperbolic} under the assumptions of  Theorem \ref{Thm-partial-hyperbolic-usual-2}, the proofs   are simple.
From Lemma \ref{lem-c-tian}, the topological entropy is positive. By  variational principle, one can  find ergodic (not periodic) measures  which have positive metric entropy arbitrarily close to the topological entropy. By assumption, these measures are hyperbolic  with domination. By Theorem \ref{Thm:Horshoe-Approximation}, one can get results (1) and  (2). For  (3), it is obvious from  Corollary \ref{Cor:1Horshoe-numbergrowth-periodic}. Now we start to discuss measures with zero entropy.

By assumption and Theorem \ref{Thm:density-per-meas22222}, every ergodic measure can be approximated by periodic measures.
Since periodic measures have zero entropy, by lemma \ref{lem-zeroentropy2015Oct}  $\m_z (M)$ is dense in $\m_f (M)$.

If $E^c$ is one-dimensional, then $f$ is far from homoclinic tangency and by  \cite{LiaoVianaYang}  $f$ is entropy-expansive so that the entropy map $\mu\mapsto h_\mu (f)$ is upper semi-continuous. By lemma \ref{lem-zeroentropy2015Oct}, we  complete the proof of Theorem \ref{Thm-partial-hyperbolic-usual-2}. \qed



\subsection{Existence of Hyperbolic SRB-like Measure}\label{section-exist-hyperbolic-srb-like}

Let $\m_f^+ (M)$ denotes the space of all invariant measures with positive entropy and let $\m_f^{PF}$ denotes the space of all invariant measures satisfying Pesin's entropy formula (i.e., $h_\mu (f)=\int \sum_{\lambda_i (x)>0}\lambda_i (x)d\mu$).
By Lemma \ref{lem-srb-meas-positiveentropysigular-zeroentropy}, Lemma \ref{lem-c-tian} and Theorem \ref{PesFormula-Thm:2}, we know that

\begin{Thm}
Under the same assumptions as Theorem \ref{Thm-partial-hyperbolic-usual-Lebegue-hyp},  $\mathcal{O}_f\subseteq \m_f^{dh} (M)\cap \m_f^n (M)\cap \m_f^+ (M)\cap \m_f^{PF}.$
It follows that by  Theorem  \ref{Thm:Horshoe-PosiStaMnfd} $\overline{\cup_{\mu \in \mathcal{O}_f} supp (\mu)}\subseteq \overline{Per (f)}$.
\end{Thm}

{\bf Proof of Theorem \ref{Thm-partial-hyperbolic-usual-Lebegue-hyp} for the result on volume-non-expanding. }
  By Lemma \ref{lem-srb-meas-positiveentropysigular-zeroentropy} and Theorem \ref{PesFormula-Thm:2}, for every SRB-like measure $\mu,$ $\mu$ a.e. $x$, the Lyapunov exponents of $E^s\oplus E^c$ are all negative and $h_{\mu} (f)=\int \log |det  (D_xf|_{E^u})| d\mu$. By Ruelle's inequality for $f^{-1},$
 $$h_{\mu} (f)\leq \int \sum_{ -\lambda_i (x)>0}-\lambda_i (x) d\mu= -\int \log |det  (D_xf|_{E^s\oplus E^c})| d\mu ,$$
where $\lambda_1 (x)\geq \lambda_2 (x)\geq\cdots \geq \lambda_{dim (M)} (x)$ denote the Lyapunov exponents of   $x$.
 This implies that for any SRB-like measure $\mu,$  $$ \int \log |det D_xf | d \mu \leq 0.$$ By Lemma \ref{SRB-like}, Lebesgue a.e. $x$,
$$\lim_{n\rightarrow+\infty}\frac1n \log |det (D_xf^n)| \leq \sup_{\mu \in \mathcal{O}_f } \int \log |det D_xf | d \mu  \leq 0.$$
\qed


For existence of hyperbolic SRB-like measure, we also have
\begin{Thm}\label{Thm-partial-hyperbolic-hyperbolic-srblike-existence}
Under the same assumptions as     (A)  (A')  (resp.,  (B)  (B')) in Theorem \ref{Thm-partial-hyperbolic}
   $\mathcal{O}_f\subseteq \m_f^{dh} (M)\cap \m_f^n (M) \cap \m_f^+ (M)\cap \m_f^{PF}$  (resp., $\mathcal{O}_{f^{-1}}\subseteq \m_f^{dh} (M)\cap \m_f^n (M) \cap \m_f^+ (M)\cap \m_f^{PF}).$
It follows that by  Theorem  \ref{Thm:Horshoe-PosiStaMnfd} $\overline{\cup_{\mu \in \mathcal{O}_f} supp (\mu)}\subseteq \overline{Per (f)}$ ( resp., $\overline{\cup_{\mu \in \mathcal{O}_{f^{-1}}} supp (\mu)}\subseteq \overline{Per (f)}$).
\end{Thm}

{\bf Proof.} We only need to consider $ (A), (A')$, since $ (B), (B')$ are similar.

Note that $ (A')$ is a particular case  of $ (A )$. We only need to consider $ (A).$ Let $\mu$ be a SRB-like measure of $f$.  Recall that the system is partially hyperbolic, by Theorem \ref{PesFormula-Thm:2} we have   $$h_\mu (f)\geq \int \sum_{i=1}^{dim\,F (x)}\lambda_i (x) d\mu>0,$$ where $\lambda_1 (x)\geq \lambda_2 (x)\geq\cdots \geq \lambda_{dim (M)} (x)$ denote the Lyapunov exponents of   $x$.
 Let $$\lambda_{max}^E (x)=\limsup_{n\rightarrow \infty}\frac1n\log\|Df^n|_{E (x)}\|,\,\lambda_{min}^E (x)=\liminf_{n\rightarrow \infty} \frac1n\log m (Df^n|_{E (x)}).$$ For $\mu$ a.e. $x\in M,$ the $\limsup,\liminf$ can be written by $\lim$ so that by quasi-conformal assumption for $\mu$ a.e. $x\in M,$  $\lambda_{max}^E (x)= \lambda_{min}^E (x) ,$ denoted by $\lambda^E (x)$ for convenience.
  Let $\Delta_1=\{x|\lambda^E  (x)\text{ exists and }< 0\}$, firstly we show that $t:=\mu (\Delta_1)>0.$ Otherwise, $\mu$ a.e. points, there are no negative Lyapunov expoents. By Ruelle's inequality for $f^{-1}$ and invariance of $\mu,$ $$h_\mu (f)=h_\mu (f^{-1})\leq\int \sum_{ -\lambda_i (x)>0}-\lambda_i (x) d\mu=0.$$ It contradicts $h_\mu (f)>0.$
 Now  we start  to prove $t:=\mu (\Delta_1)=1.$
By contradiction,  $0<t< 1.$
 Define
   $\mu_1=\mu_{\Delta_1},$
   define  $\Delta_2=\{x|\lambda^E (x) \text{ exists and }\geq  0\}$ and
   $ \mu_2=\mu_{\Delta_2}.$ By Ruelle's inequality for $f^{-1}$, $$h_{\mu_2} (f)\leq\int \sum_{ -\lambda_i (x)>0}-\lambda_i (x) d\mu_2=0. $$
   So $$h_\mu (f)=th_{\mu_1} (f)+ (1-t)h_{\mu_2} (f)$$$$ = th_{\mu_1} (f) \leq t\int \sum_{ \lambda_i (x)>0}\lambda_i (x) d\mu_1\leq t\int \sum_{i=1}^{dim\,F (x)}\lambda_i (x) d\mu_1.$$

      However, by above analysis we know  that  $$h_\mu (f)\geq \int \sum_{i=1}^{dim\,F (x)}\lambda_i (x) d\mu$$$$=t\int \sum_{i=1}^{dim\,F (x)}\lambda_i (x) d\mu_1 + (1-t)\int \sum_{i=1}^{dim\,F (x)}\lambda_i (x) d\mu_2>t\int \sum_{i=1}^{dim\,F (x)}\lambda_i (x) d\mu_1.$$  It is a contradiction.

      Now we know that $\mu$ a.e. $x$, the Lyapunov exponents of bundle $E (x)$ are all negative (that is, $\mu$ is hyperbolic) and $$h_\mu (f)\geq \int \sum_{i=1}^{dim\,F (x)}\lambda_i (x) d\mu.$$
      So by Ruelle's inequality and uniform  expanding of bundle $F$, the inverse inequality is also true and so Pesin's entropy formula holds.

      Now we only need to prove that $\mu$ is non-atomic. Otherwise,  there is a point $x\in M$ such that $\mu\{x\}>0,$ then $x$ should be a periodic orbit. Denote the measure supported on the orbit of $x$ by $\nu.$ Then $h_\nu (f)=0.$ Since $h_\mu (f)>0,$ then $\mu\neq \nu.$ So $\theta:=\mu (M\setminus Orb (x))\in (0,1).$ Let $\omega=\mu|_{M\setminus Orb (x)}$, then  by Ruelle's inequality, $$h_\mu (f)=\theta h_{\omega} (f)+ (1-\theta) h_{\nu} (f) =\theta h_{\omega} (f)\leq \theta\int \sum_{i=1}^{dim\,F (x)}\lambda_i (x) d\omega.$$ It contradicts
$$h_\mu (f)\geq \int \sum_{i=1}^{dim\,F (x)}\lambda_i (x) d\mu$$$$=\theta\int \sum_{i=1}^{dim\,F (x)}\lambda_i (x) d\omega+ (1-\theta)\int \sum_{i=1}^{dim\,F (x)}\lambda_i (x) d\nu> \theta\int \sum_{i=1}^{dim\,F (x)}\lambda_i (x) d\omega .$$ \qed

\bigskip

{\bf Proof of  Theorem \ref{Thm-partial-hyperbolic-444444444444444}. }
  By  Theorem \ref{Thm-partial-hyperbolic-hyperbolic-srblike-existence}, for every SRB-like measure $\mu,$ $\mu$ a.e. $x$, the Lyapunov exponents of $E $ are all negative and $h_{\mu} (f)=\int \log |det  (D_xf|_{F})| d\mu$. By Ruelle's inequality for $f^{-1},$
 $$h_{\mu} (f)\leq \int \sum_{ -\lambda_i (x)>0}-\lambda_i (x) d\mu= -\int \log |det  (D_xf|_{E })| d\mu ,$$
where $\lambda_1 (x)\geq \lambda_2 (x)\geq\cdots \geq \lambda_{dim (M)} (x)$ denote the Lyapunov exponents of   $x$.
 This implies that for any SRB-like measure $\mu,$  $$ \int \log |det D_xf | d \mu \leq 0.$$ By Lemma \ref{SRB-like}, Lebesgue a.e. $x$,
$$\limsup_{n\rightarrow+\infty}\frac1n \log |det (D_xf^n)| =\sup_{\mu \in pw_f (x) } \int \log |det D_xf | d \mu \leq \sup_{\mu \in \mathcal{O}_f } \int \log |det D_xf | d \mu  \leq 0.$$ \qed

Under the assumption  of Theorem \ref{Thm-partial-hyperbolic-22222}, it is still unknown whether every SRB-like measure    is hyperbolic, even though we know that, for example, in the case  (I) for any SRB-like measure the integrable maximal Lyapunov exponents of $E$ is always negative.
\subsection{Stably ergodic diffeomorphims}

It was proved in  \cite{BFP} that in the space of stably ergodic volume-preserving systems, there is an open and dense subset such that every system in this subset is non-uniformly hyperbolic: volume measure is a hyperbolic measure and admits a dominated splitting $T M = E^-
\oplus E^+$,
where   $E^-$
   (resp. $E^+$) coincides a.e. with the sum of the Oseledets spaces corresponding to  negative (resp. positive)
Lyapunov exponents. Moreover, it is known that the space of stably ergodic   systems  contains all Anosov diffeomorphisms  and many partially hyperbolic ones ( see e.g.  \cite{GPS}). It is
not true that every stably ergodic diffeomorphism can be approximated by a partially hyperbolic system, see  \cite{Ta2004,BV2000}.

Thus, there are some  non-uniformly hyperbolic systems with domination which is not Anosov and not partially hyperbolic.

\section{Existence of Horseshoe, Density of Hyperbolic Periodic Orbits \& Liv$\check{\text{s}}$ic  Theorem}


In this section we use Theorem \ref{Thm:Expclosinglem-measure} to prove Theorem \ref{Thm:Horshoe-PosiStaMnfd} and \ref{Thm:LivThm} (resp. one can use Theorem \ref{Thm:Expclosinglem-measure2015oct} to prove Theorem \ref{Thm:2015-quasi-invariant} (1) and (2)). Note that different to Katok's  (shadowing and) closing lemma, the length of every pseud-orbit segment in the present paper must be larger than $T_\tau$  and thus we must take care when we use  (shadowing and) closing lemma. Firstly, we prove Theorem \ref{Thm:Horshoe-PosiStaMnfd}.

\medskip
\subsection{Density of periodic orbits}\label{proof-horseshoe-periodic-dense}
{\bf Proof of Theorem \ref{Thm:Horshoe-PosiStaMnfd}}
The  proof is an adaption of  \cite{Uga}.
Let $x\in supp (\mu)$ and let $\eta>0$ be small enough. There exists a positive number $\tau$ small enough  such that $\mu (B (x,\eta/2)\cap\Lambda_\tau)>0.$   By Theorem \ref{Thm:Expclosinglem-measure}, we can take $\sigma>0$ such that the stable and unstable manifolds of  $\Lambda^*_\tau (\sigma)$  defined in Theorem \ref{Thm:Expclosinglem-measure} have the same uniform size.  Let  $\delta>0$  be the number  such that if $x,y\in\Lambda^*_\tau (\sigma)$ are $\delta$ close, then their stable and unstable manifolds  are transversal. Let $T_\tau$ be the number as in Theorem \ref{Thm:Expclosinglem-measure} which is only dependent on $\Lambda^*_\tau.$
Now let $r>0$ be an arbitrary small number, with $r<\min\{\eta/2,\sigma,\delta/3\}$.

  Pick a set $B\subseteq B (x,\eta/2)\cap\Lambda_\tau$   of diameter less than  $\beta=\beta (\tau,r)$  (as in Theorem \ref{Thm:Expclosinglem-measure}) less than $r$ and of positive measure. Let $x_1\in B$ be a recurrent point  (by the Poincare Recurrence Theorem), and $n (x_1)\geq T_\tau$ a positive integer such that $f^{n (x_1)} (x_1)\in B$. Since $   d (x_1,f^{n (x_1)} (x_1))<\beta$, in applying Theorem \ref{Thm:Expclosinglem-measure} we obtain that there exists a periodic point $z_1$ of period $n (x_1)$ such that $   d (f^i (x_1),f^i (z_1))<r<\min\{\eta/2,\sigma\},\,0\leq i\leq n (x_1)-1.$ This implies that $z_1\in\Lambda^*_\tau (\sigma)$ and $   d (x,z_1)<\eta.$

  We will prove that $\mu$-almost every point of $B$ belongs to $\overline{W^u (z_1)}$. For that, let $y\in B$ be a Borel density point, and let $\tau'>0$ be small enough such that $\tau'<r/2$, $   d (z_1,y)>\tau'$ and $\mu (B\cap B (y,\tau'/2))>0$. Pick a set $\tilde{B}\subseteq B\cap B (y,\tau'/2)$ of diameter less than $\tilde{\beta}=\beta (k,\tau'/2)$  (as in Theorem \ref{Thm:Expclosinglem-measure}) and of positive measure. Let $x_2\in \tilde{B}$  be a recurrent point, and $n (x_2)\geq T_\tau$ a positive integer such that  $f^{n (x_2)} (x_2)\in B$. Since $   d (x_2,f^{n (x_2)} (x_2))<\tilde{\beta}$. Since $   d (x_2,f^{n (x_2)} (x_2))<\tilde{\beta}$, in applying Theorem \ref{Thm:Expclosinglem-measure} we obtain that there exists a periodic point $z_2$ of period $n (x_2)$ such that $   d (f^i (x_2),f^i (z_2))<\tau'/2<r<\sigma,\,0\leq i\leq n (x_2)-1.$ This implies that $z_2\in\Lambda^*_\tau (\sigma)$ and $   d (y,z_2)<   d (y,x_2)+   d (x_2,z_2)<\tau',$ thus $z_1\neq z_2.$

   Note that $$   d (z_1,z_2)<   d (x_1,z_1)+   d (x_1,y)+   d (y,z_2)<r+r+\tau'<3r<\delta.$$ By the choice of $\delta$,   the stable  (unstable) manifold $W^{s/u} (z_1)$ and unstable (stable) manifold $W^{u/s} (z_2)$  are transversal.  So $z_1$ has homoclinic point and hence there exists horseshoe.

   Let $w$ be a transverse intersection point on $W^u (z_1)$ and $W^s (z_2).$
  By the Inclination Lemma (see e.g.  \cite{Katok1}), we have that $W^u (z_1)$ accumulates locally on $W^u (z_2)$. Indeed, let $N$ be a common period for $z_1$ and $z_2.$ The points $z_1$ and $z_2$ are hyperbolic fixed points for $f^N$ with  $w$ a  transverse intersection point of $W^u (z_1)$ and $W^s (z_2).$ The images under $f^N$ of a ball around $w$ in $W^u (z_1)=f^NW^u (z_1)$ accumulate on $W^u (z_2)$. Hence $$   d  (W^u (z_1),y)\leq   d  (W^u (z_1),x_2)+   d  (x_2,y)\leq\tau'$$ for all $\tau'>0$. This implies that $y\in\overline{W^u (z_1)}$, for all Borel density points $y\in B$. Therefore $\mu (\overline{W^u (z_1)})\geq \mu (B)>0.$

  Since the  hyperbolic periodic point $z_1$ is such that $   d (x,z_1)<\eta$, where $x$ is an arbitrary point of $supp (\mu)$ and $\eta>0$ is arbitrarily small, we obtain that the support of the measure $\mu$ is contained in the closure of all hyperbolic periodic points with $\mu (\overline{W^u (z_1)})>0.$ \qed

 \bigskip

 {\bf Proof of Theorem \ref{Thm:2015-quasi-invariant} (1) and (2).} One can follow the proof of Theorem \ref{Thm:Horshoe-PosiStaMnfd} to give the proof, by replacing Theorem \ref{Thm:Expclosinglem-measure} by Theorem \ref{Thm:Expclosinglem-measure2015oct} and replacing invariance of the measure by quasi-invariant.

\subsection{Livshitz Theorem}\label{proof-Livsic}
{\bf Proof of Theorem \ref{Thm:LivThm}}
 The  proof is an adaption of  \cite{Katok2,KatokMondoz,BP} and not too long,  here we only prove the ergodic case similar as Theorem 15.3.1 in  \cite{BP}. Then there is $x\in supp (\mu)$ such that $supp (\mu)=\overline{orb (x)}$ (In fact this property holds for $\mu$ a.e. $x,$ see Theorem 5.15 in  \cite{Walters}).

  Given $\tau>0$ and let $\Lambda_\tau$ be as in Theorem \ref{Thm:Expclosinglem-measure}. We may assume that $supp (\mu)=\overline{orb (x)}$ for some $x\in\Lambda_\tau$ and that for every $\tau>0$, the intersection $orb (x)\cap\Lambda_\tau$ is dense in $\Lambda_\tau$. Let $T_\tau>0$ is the number only dependent on $\tau$ in Theorem \ref{Thm:Expclosinglem-measure}.
Set $\psi=0$ and
$$\psi (f^n (x))=\sum_{k=0}^{n-1} \varphi (f^k (x))$$
for each $n\in \mathbb{N}.$ To extend $\psi$ continuously to $\Lambda_\tau$
we show that $\psi$ is uniformly continuous on $orb (x)\cap\Lambda_\tau$.
For any $\eta>0$, let $\beta=\beta (\tau,\eta)>0$ be as in Theorem \ref{Thm:Expclosinglem-measure}.

 Firstly, we consider $n_2>n_1$ satisfying $n_2-n_1\geq T_\tau,\,\,f^{n_1} (x),f^{n_2} (x)\in\Lambda_\tau$ and $   d (f^{n_1} (x),\,\,f^{n_2} (x))<\beta,$ then by Theorem \ref{Thm:Expclosinglem-measure}, there exists a hyperbolic periodic point $z$ with period $n_2-n_1$ such that
 $$   d (f^j (f^{n_1} (x)),\,\,f^j (z))<\eta\cdot e^{-\min\{j,n_2-n_1-j\}\theta_\tau},$$ where $j=0,\,1,\,\cdots,n_2-n_1-1$ and $\theta_\tau>0$ is the number only dependent on $\tau$ in Theorem \ref{Thm:Expclosinglem-measure}. Since $\varphi$ is H$\ddot{o}$lder continuous, $$|\varphi (x)-\varphi (y)|\leq C    d (x,y)^\kappa$$ for some $C>0$ and $0<\kappa\leq 1$, and we obtain
 $$|\varphi (f^j (f^{n_1} (x))-\varphi (f^j (z))|\leq C\eta^\kappa\cdot e^{-\kappa\min\{j,n_2-n_1-j\}\theta_\tau}.$$
 Thus, \begin{eqnarray*}
 & &|\psi (f^{n_1} (x))-\psi (f^{n_2} (x))|\\
 &=& |\sum_{j=n_1}^{n_2-n_1-1}  (\varphi (f^j (f^{n_1} (x))-\varphi (f^j (z))) + \sum_{j=n_1}^{n_2-n_1-1} \varphi (f^j (z))|\\
 &\leq&C\eta^\kappa\cdot\sum_{j=n_1}^{n_2-n_1-1} \max\{e^{-j\kappa\theta_\tau},e^{- (n_2-n_1-j)\kappa\theta_\tau}\}+|\sum_{j=n_1}^{n_2-n_1-1} \varphi (f^j (z))|\\
  &\leq& C \eta^\kappa N
 \end{eqnarray*}
 for some constant $N>0$ independent of $n_1$ and $n_2.$

  Secondly, we consider $n_2>n_1$ satisfying $f^{n_1} (x),f^{n_2} (x)\in\Lambda_\tau$ and $$   d (f^{n_1} (x),\,\,f^{n_2} (x))<\frac\beta2.$$ Since $orb (x)\cap\Lambda_\tau$ is dense in $\Lambda_\tau$, then
  we can take $n_3$ large enough such that $n_3-n_1> n_3-n_2 \geq T_\tau$, $f^{n_3} (x)\in \Lambda_\tau$ and $   d (f^{n_3} (x),\,\,f^{n_1} (x))<\frac\beta2$. So $   d (f^{n_3} (x),\,\,f^{n_2} (x))\leq d (f^{n_3} (x),\,\,f^{n_1} (x))+   d (f^{n_1} (x),f^{n_2} (x))<\beta.$
  Thus by the above discussion, we have $$|\psi (f^{n_3} (x))-\psi (f^{n_1} (x))|\leq C \eta^\kappa N$$ and
  $$|\psi (f^{n_3} (x))-\psi (f^{n_2} (x))|\leq C \eta^\kappa N.$$ Therefore,
   \begin{eqnarray*}& &|\psi (f^{n_1} (x))-\psi (f^{n_2} (x))|\\
   &\leq&|\psi (f^{n_3} (x))-\psi (f^{n_1} (x))|+|\psi (f^{n_3} (x))-\psi (f^{n_2} (x))|\\
   &\leq& 2C \eta^\kappa N.
   \end{eqnarray*}
  This shows that we can extend $\psi$ continuously to $\Lambda_\tau.$ Now we extend $\psi$ to $\cup_{i=0}^\infty f^i (\Lambda_\tau)$ as follows. If $y\in f (\Lambda_\tau)\setminus \Lambda_\tau,$ then let $$\psi (y)=\psi (f^{-1} (y))+ \varphi (f^{-1} (y)),$$ and we use the same expression for all
 $y\in \cup_{i=0}^{n+1}f^i (\Lambda_\tau)\setminus \cup_{i=0}^n f^i (\Lambda_\tau).$ Since $$\mu (\cup_{i=0}^\infty f^i (\Lambda_\tau))=1,$$ the function $\psi$ is defined almost everywhere and clearly satisfies the needed condition.
\qed

\medskip
Remark that exponentially closing plays the crucial role  in Theorem \ref{Thm:LivThm} but the stable manifold theorem is not used  in the proof.

\bigskip

Now we start to prove Theorem \ref{Thm:LivThm222222}.

{\bf Proof of Theorem \ref{Thm:LivThm222222}.} By assumption and Theorem \ref{Thm:density-per-meas22222}, every ergodic measure can be approximated by periodic ones.
By the information assumed on periodic orbits and weak$^*$ topology. Every ergodic $\mu$ satisfies $\int \varphi d\mu=0.$ Fix $\epsilon>0.$ Take  $a_n (x)=\sum_{j=0}^{n-1}\varphi (f^jx)-n\epsilon$
 and $b_n=-\sum_{j=0}^{n-1}\varphi (f^jx)+n\epsilon.$ It is easy to check the assumptions in Lemma \ref{Lem:1} (In fact the inequality of sub-additional property is equality). So there is some  $N_1\geq 1$  such that for any $x\in M$ and $n\geq N_1,$ $$\sum_{j=0}^{n-1}\varphi (f^jx)\leq n\epsilon. $$
Similarly, taking $a_n (x)=-\sum_{j=0}^{n-1}\varphi (f^jx)-n\epsilon$ to get that there is some   $N_2\geq 1$  such that for any $x\in M$ and $n\geq N_2,$ $$-\sum_{j=0}^{n-1}\varphi (f^jx)\leq n\epsilon. $$
So $\frac1n\sum_{j=0}^{n-1}\varphi (f^jx)$ converges uniformly to zero.

Define for $n\geq 1 $  and $x\in M$
$$h_n (x)=\frac1n\sum_{j=1}^{n-1} (n-j)\varphi (f^{j-1}x).$$
We can verify that $$\varphi (x)-\frac1n\sum_{j=0}^{n-1}\varphi (f^jx)=h_n (x)-h_n (fx).$$ So $\varphi\in \overline{Cob (f)}.$\qed

\subsection{Approximation of Horseshoes}\label{proof-horshoe-approximation}

{\bf Proof of Theorem \ref{Thm:Horshoe-Approximation} } 
 The proofs of  (1) and  (3) can be as an adaption  of  \cite{BP,ACW}. For others in our setting, without the technique of Lyapunov neigborhood,   we can use  (1) and  (3) to prove.

 For ergodic measure, the splitting should be dominated on $supp (\mu).$
Since dominated property can be extended to its neighborhoods  (see  \cite{BLV}), we can take an open neighborhood $U=B_{supp (\mu)} (\tau_0)=\{x\,\,|\,\gamma (x,supp (\mu))<\tau_0\}$ (for some small $\tau_0>0$) of $supp (\mu)$ such that the splitting is extended and dominated on $\Delta_0:=\cap_{n\in\mathbb{Z}}f^n (U)$.
Without loss of generality, we can suppose that $\epsilon<\tau_0.$ Take $\epsilon'\in (0,\epsilon)$ small enough such that $$h_\mu (f)-\epsilon<\frac1{1+\epsilon'} ({h_\mu (f)-3\epsilon'}).$$

 Assume
$\{\varphi_i\}_{i=1}^{\infty}$ is the dense subset of $C (M)$ giving
the weak$^*$ topology, that is,
$$D (\mu,\,m)=\sum_{i=1}^{\infty}\frac{|\int \varphi_id\mu-\int \varphi_idm|}{2^i\|\varphi_i\|}$$
for  $\mu,m\in\mathcal{M} (f)$. Take $ \varsigma>0$  such that
\begin{eqnarray}\label{eq-measure-nbhd}
\{\nu\,\,| D (\mu,\nu)\leq 3\varsigma\}\subset \mathcal{V}.
\end{eqnarray}
 Then there exists
integer $T>0$ such that
$$\sum_{i=T+1}^{\infty}\frac{|\int \varphi_id\mu-\int
\varphi_idm|}{2^i\|\varphi_i\|}\leq\sum_{i=T+1}^{\infty}\frac{1}{2^{i-1}}< \varsigma.$$
Now we denote $\phi_i=\varphi_i/\|\varphi_i\|$ for
$i=1,\cdots,T$.

{\bf Step 1.} {\it Choice of Separated Set.}

We use Katok's definition of metric entropy ( see
 \cite{Katok2}). For $x,y\in M$ and  $l\in \mathbb{N}$, let
$$d_l (x,y)=\max_{0\leq i\leq l-1}d (f^i (x),\,f^i (y)).$$ For
$\gamma>0,\,\delta\in  (0,1)$, let $N_l (\gamma,\,\delta)$ be the
minimal number of $\gamma-$balls $B_l (x,\,\gamma)$ in the
$d_l-$metric, which cover a set of measure at least $1-\delta$. Then
$$h_{\mu} (f)=\lim_{\gamma\rightarrow0}\liminf_{l\rightarrow\infty}\frac{\log N_l (\gamma,\,\delta)}{l}
=\lim_{\gamma\rightarrow0}\limsup_{n\rightarrow\infty}\frac{\log
N_l (\gamma,\,\delta)}{l}.$$
Here we can fix a $\delta\in  (0,1)$.

 For above $\epsilon, \epsilon',  \varsigma$, $\delta,$
and functions
$\phi_1,\,\phi_2,\,\cdots,\,\,\phi_T,\,$
take $\gamma\in (0,\min\{\frac{\epsilon}2,\varsigma\})$ such that if
$d (x,y)<\gamma$ then
$$|\phi_i (x)-\phi_i (y)|<\varsigma (i=1,2,\cdots,T)$$ and
$$h_{\mu} (f)-{\epsilon'}<\liminf_{l\rightarrow+\infty}\frac1l\log N_l (\gamma,\delta)
\leq\limsup_{l\rightarrow+\infty}\frac1l\log
N_l (\gamma,\delta).$$

From Theorem \ref{Thm:ExpShadowing-measure}, take $\Delta=\Lambda_\delta\cap supp (\mu)$  such that
$\mu (\Delta)>1-\delta $ and  one has exponential  shadowing lemma on $\Delta.$ More precisely, there is  $
\theta>0$ (here $\mu$ is ergodic so that $\theta>0$ is independent of $\delta$) and  $T_\delta>0$ such that following  holds: for above
$\,\frac{\gamma}{2}>0,$ there exists     $\tau\in (0,\frac{\gamma}2)$ such that if  a
$\tau$-pseudo-orbit $\{x_i,\,n_i\}_{i=-\infty}^{+\infty}$
satisfies $n_i\geq T_k $ and $x_i,f^{n_i} (x_i)\in
\Delta$ for all $i$, then there exists a
$\frac{\gamma}{2}$-shadowing point $x\in M $ for
$\{x_i,n_i\}_{i=-\infty}^{+\infty}$ (here it is not necessarily the shadowing  to be exponential for $\theta$).

Take   a
 finite open cover $\{B (a_i, \tau)\}_{i=1}^{Q}$ of $\mbox{supp} (\mu)$
 where $a_i\in \mbox{supp} (\mu).$

 Let $\xi$ be a finite partition of $M$
with $\mbox{diam}\,\xi <\tau$ and $\xi>\{\Delta,M\setminus
\Delta\}$. Consider the set\begin{eqnarray*}
\Delta_{s}&=&\big{\{}x\in \Delta_k:
\{x,f (x),\cdots, f^{m-1} (x)\}\cap B (a_i,
\gamma)\neq \varnothing\\[2mm]
&& (i=1,2,\cdots,Q),\,\,f^m (x)\in \xi (x)~ \mbox{for some} ~
m\in[l, (1+\epsilon')l)\\[2mm]
 &&\mbox{and}~
 |\frac1{l}\sum_{j=0}^{l-1}\phi_i (f^jx)-\int\phi_id\mu|<\varsigma (i=1,2,\cdots,T),
 ~\mbox{for}~ l\geq s\big{\}}.
 \end{eqnarray*}
By ergodicity of $\mu$ and Birkhoff's Ergodic theorem,
$\mu (\Delta_{s})\rightarrow \mu (\Delta)\, (
as\,\,s\rightarrow+\infty)$. Take sufficiently large $s$ such that
$\mu (\Delta_{s})>1-\delta.$ Let $E_l\subseteq \Delta_{s}$ be an
$ (l,\gamma)-$separated set of maximal cardinality; in other
words, the cover by $\gamma-$balls in the $d_l-$metric centered
at points in $E_l$ is a minimal cover $\Delta_{s}.$ Then we have
$$h_{\mu} (f)-{\epsilon'}<\liminf_{l\rightarrow+\infty}\frac1l\log \sharp
\,E_l \leq\limsup_{l\rightarrow+\infty}\frac1l\log \sharp \,E_l.$$
So there exists $L>2kK$ such that for all $l\geq L$,
$$h_{\mu} (f)-{\epsilon'}<\frac1l\log \sharp \,E_l, $$
and $$Q<e^{l\epsilon'},\,\,l\epsilon'<e^{l\epsilon'} .$$

For $m=l,l+1,\cdots,[ (1+\epsilon')l]-1$, set $$R_m:=\{x\mid\, x\in E_l, f^m (x)\in \Delta_{s} \text{ but } \{f^l (x),f^{l+1} (x),\cdots,f^{m-1}\}\cap \Delta_{s}=\emptyset\}.$$ and thus
$$\sum_{m=l}^{[ (1+\epsilon')\,l]-1} \sharp R_m= \sharp E_l
\leq l\epsilon'\,\,\max_{l\leq m\leq [ (1+{\epsilon'})\,l]-1} \sharp
R_m.$$ Therefore,   for every $l\geq L$ (here we fix a $l$) one can
find an
 $N=m_l$  ( $l\leq m_l\leq[ (1+{\epsilon'})l]-1$)\, satisfying $$e^{-l\epsilon'}\sharp
E_l\leq \frac1{l\epsilon'} \sharp E_l\leq \sharp R_{N}.$$ Hence,
$$\frac 1 {1+{\epsilon'}} (h_{\mu} (f)-2{\epsilon'})<\frac 1 {1+{\epsilon}} ({- \epsilon' }+ \frac 1
{l}\log\sharp E_l)\leq \frac1{N} \log \sharp\,R_{N}.$$

Now take $i\in\{ 1,2,\cdots, Q\}$ such that $R_N\cap B (a_i,\tau)$ carries maximal cardinality.
Let $B=B (a_i,\tau)$ and $Y=R_N \cap B (a_i,\tau)$, then  $$f^N (Y)\subseteq \Delta\cap B.$$ Moreover,
\begin{eqnarray}\label{eq-entropy-estimate}
h_\mu (f)-\epsilon <\frac 1 {1+{\epsilon'}} (h_{\mu} (f)-3{\epsilon'})<\frac 1 {1+{\epsilon'}} (h_{\mu} (f)-2{\epsilon'}-\log Q)
\end{eqnarray}
$$\leq \frac1{N} \log  (\frac1Q\sharp\,R_{N})\leq \frac1{N} \log \sharp\,Y.$$

{\bf Step 2.} {\it Construction of  Invariant Compact Set $\Lambda$.}

 Consider the shift $X=Y^{\mathbb{Z}}$ over the alphabet $Y$. Note that $diam (Y)\leq diam (B)\leq \tau.$ So for each $ y= (y_n)\in X ,$ there is a corresponding $\tau-$pseudo-orbit:
 $$\cdots,y_0,\cdots,f^{N-1} (y_0),y_1,\cdots,f^{N-1} (y_1),\cdots.$$
By Theorem \ref{Thm:shadowinglem}, for each such pseudo-orbit there is a   point $\pi (y)$ whose orbit $\frac\gamma2-$shadows.
 Thus, we can obtain an  $f^N$ invariant set $\Lambda_0$ which is the union of all orbits that shadow the elements of $X$  (In fact $$\Lambda_0= \cap_{n=1}^\infty \cup_{\{y_0,y_1,\cdots,y_{n-1}\}\in Y^n} \cap_{j=0}^{n-1} f^{-jN}\overline{B_N (y_j,\frac \gamma 2)} ,$$ where $ Y^n$ is the $n-$product space of $Y\times\cdots\times Y$). Define $$\Lambda=\Lambda_0\cup f (\Lambda_0)\cdots\cup f^{N-1}  (\Lambda_0). $$
Then $\Lambda$ is a transitive compact invariant set.

Since for any $y,y'\in Y,$ there is $j\in[0,N)$ such that $d (f^jy,f^jy')>\gamma$ and since the shadowing is $\frac\gamma2-$shadowing so that if $x$ and $x'$ $ \frac\gamma2 -$shadows distinct elements of $X$, then $x$ and $x'$ are distinct as well. In particular, distinct periodic elements of $X$ correspond to different shadowing periodic points of $\Lambda$, then   $$\# P_{nN} (f)\cap \Lambda\geq  (\# Y)^n.$$ So \begin{eqnarray}\label{eq-periodicgrowth-for-constructed-horseshoe}
\limsup_{n\rightarrow \infty }\frac 1n \log \# P_n (f)\cap \Lambda \geq \frac1{N} \log \sharp\,Y .
\end{eqnarray}

From the choice of $Y$, for any $y_j\in Y,$ one has $$\max_{1\leq i\leq T} |\frac 1N \sum_{m=1}^N\phi_i (f^m (y_j))-\int \phi_i d\mu|<\varsigma.$$
By the choice $\gamma$, if for some $x\in M$, $d_N (x,y_j)<\gamma$, then
 $$\max_{1\leq i\leq T} |\frac 1N \sum_{m=1}^N\phi_i (f^m (y_j))-\frac 1N \sum_{m=1}^N\phi_i (f^m (x))|<\varsigma.$$ It follows that a point $x$ $ (\frac{\gamma}2,\theta)$-shadows a pseudo-orbit in $X$, then
$$\limsup_{n\rightarrow\infty}\max_{1\leq i\leq T} |\frac 1n \sum_{m=1}^n\phi_i (f^m (x))-\int \phi_i d\mu|<2\varsigma.$$ This implies for any ergodic $\nu\in\m_f (\Lambda),$ $D (\mu,\nu)\leq 3\varsigma$ and so by Ergodic Decomposition Theorem for any invariant measure $\nu\in\m_f (\Lambda),$  $D (\mu,\nu)\leq3\varsigma.$ By    (\ref{eq-measure-nbhd}),  we have  the conclusion  (3).

By construction $\Lambda\subseteq B (supp (\mu),\gamma)\subseteq B (supp (\mu),\epsilon)$. Recall that $Y\subseteq \Delta_s$ so that for any $y\in Y$, $$ \{y,\cdots,f^{N-1}\}\cap B (a_i,\tau)\neq\emptyset, i=1,2,\cdots,Q.$$
This follows that $$ supp (\mu)\subseteq B (Y\cup\cdots\cup f^{N-1}Y,\tau) \subseteq B (\Lambda, \tau+\frac{\gamma}2)\subseteq B (\Lambda, \gamma)\subseteq B (\Lambda, \epsilon).$$
So the conclusion  (1) holds. Moreover, $\epsilon<\tau_0$ implies that  $\Lambda$ should be contained in $\Delta_0$ so that $\Lambda$ carries  dominated splitting (conclusion  (4)) which coincides with the extended  (continuous) dominated splitting $$T_{\Delta_0}M=E^s\oplus E^u. $$

\medskip

{\bf Step 3.} {\it Uniform Hyperbolicity of $\Lambda$.}

Firstly let us state a basic fact.

\begin{Prop}\label{Prop-semi-continuity-Lyapunov} Let $\Gamma$ be an $f$ invariant compact set and $E\subset TM$ be a continuous $Df$ invariant subbundle on $\Gamma$. Let $$\lambda_{max} (E,x)=\limsup_{n\rightarrow\infty}\frac1n \log \|Df^n|_{E (x)}\|,\,\,\lambda_{min} (E,x)=\liminf_{n\rightarrow\infty}\frac1n \log m (Df^n|_{E (x)}).$$ Then
the function $$ \Psi_{max}: \m_f (\Gamma)\rightarrow\mathbb{R}, \mu\mapsto \int \lambda_{max} (E,x) d\mu$$ are upper semi-continuous, and respectively,
the function $$ \Psi_{min}: \m_f (\Gamma)\rightarrow\mathbb{R}, \mu\mapsto \int \lambda_{min} (E,x) d\mu$$ are lower semi-continuous.

\end{Prop}

{\bf Proof.} We only prove the first one and the second one is similar.
By sub-additional ergodic theorem, $$\Psi_{max} (\mu)= \inf_{n\geq 1 }\int \frac1n \log \|Df^n|_{E (x)}\|d\mu.$$ It follows the upper semi-continuity. \qed

\medskip

If necessary, take $\varsigma$ small enough such that all the invariant measures supported on $\Lambda$ are close to $\mu$  enough  in weak$^*$ topology.  By Proposition \ref{Prop-semi-continuity-Lyapunov}, this follows that for any $ \nu\in\m_f (\Gamma),$
$$ \int \lambda_{max} (E^s,x) d\nu\leq \int \lambda_{max} (E^s,x) d\mu+\frac{\epsilon}2\leq \lambda_s+\frac{\epsilon}2 ,$$
$$\,  \int \lambda_{min} (E^s,x) d\nu\geq\int \lambda_{min} (E^u,x) d\mu -\frac{\epsilon}2=\lambda_u-\frac{\epsilon}2.$$  where
$\lambda_s,\lambda_u$ denote the maximal negative Lyapunov exponent of $\mu$ and the minimal positive  Lyapunov exponent of $\mu$, respectively.  This implies all invariant measures of $\Gamma$ are hyperbolic so that by Theorem \ref{Thm-Cao2003} $\Gamma$ is uniformly hyperbolic. By the argument of Lemma \ref{Lem:1}, we can get that there exists $C_\epsilon>0$ such that the conclusion  (5) holds.

If further the Oseledec splitting of $\mu$ is dominated, then  (4') is similar as above by extending the dominated splitting to the neighborhood. Moreover,  one can use Proposition \ref{Prop-semi-continuity-Lyapunov} for every Oseledec bundle $E_i$ and by the argument of Lemma \ref{Lem:1}, it is similar as the proof of  (5) to get  (5'). And it is easy to see that  (6) is
 a particular case of  (5').

 Now we start  to prove  (2).  By uniform hyperbolicity of $\Lambda,$ from  (\ref{eq-growth-period-hyperbolic})
 $$ h_{top} (f|_\Lambda)\geq\limsup_{n\rightarrow \infty }\frac 1n \log \# P_n (f)\cap \Lambda .$$
 Together with  (\ref{eq-entropy-estimate})  and  (\ref{eq-periodicgrowth-for-constructed-horseshoe}), this ends the proof of conclusion  (2).

 {\bf Step 4.} {\it Horseshoe  $H_\epsilon.$}

 If $\Lambda$ is a horseshoe, then it is the needed $H_\epsilon.$ Otherwise, one just needs to apply following proposition of  \cite{ACW}.

 \begin{Prop} Let $\Lambda$ be a hyperbolic set of $C^1$ diffeomorphism $f$. Then for any $\delta>0$, there is a horseshoe $K$ that is $\delta-$close to $\Lambda$ in the Hausdorff topology and satisfies $h_{top} (f|_K)\geq h_{top} (f|_\Lambda)-\delta.$

 \end{Prop} This is  Proposition 8.7 of  \cite{ACW}, here we omit the details. Now we complete the proof of  of Theorem \ref{Thm:Horshoe-Approximation}. \qed

\section{Distribution of Periodic Measures \& other possible applications}

Let $\mu\in \m  (M).$ For a Liao-Pesin set $\Lambda (K,\,\zeta)$, define $\tilde \Lambda_k (K,\,\zeta)=\textrm{supp} (
\mu|_{\Lambda_{k} (K,\,\zeta)})$ and $\tilde
\Lambda (K,\,\zeta)=\bigcup_{k=1}^\infty\tilde \Lambda_k (K,\,\zeta).$
For convenience, we say they are {\it $\mu$-Liao-Pesin blocks} and {\it $\mu$-Liao-Pesin set}, respectively.
 Similarly, for improved Liao-Pesin set $  \Lambda^\# (K,\zeta)$, one can define improved {\it $\mu$-Liao-Pesin block} $\tilde \Lambda^\#_k$ and {\it $\mu$-Liao-Pesin set} $\tilde \Lambda^\#$. They are subsets of Liao-Pesin blocks  and Liao-Pesin set and they depend on the given measure.


\subsection{Non-uniform Specification}


\begin{Lem}\label{Lem:weak-specification}Let $\mu\in \m_f (M).$ If $\mu$ is  ergodic, then for every $\mu$ positive-measured set $\tilde
\Lambda_k (K,\,\zeta)$, one has the  specification property as follows:

 For any $\varepsilon>0, $ there exist
 $X_{\tau}=X_{\tau} (k,\varepsilon)> 0 (\tau=1,2)$ such that  for a given sequence of  points $x_1,x_2,\cdots,x_N\in\,
 \widetilde\Lambda_k (N\in\mathbb{N})$ and  a
sequence of positive numbers ${n_1,n_2,\cdots,n_N}$, one has
$n_i\geq 2kK $ and $f^{n_i}x_i \in\,\widetilde\Lambda_k$ for
$i=1,\,2,\,\cdots,\,N,$ then there exist  a periodic point $z  \in
M$, a positive number
$p\in[\sum_{i=1}^{N}n_i+NX_{1},\,\,\sum_{i=1}^{N}n_i+NX_{2}],$
 and a sequence of nonnegative numbers $c_0=0,c_1,\cdots,c_{N-1}$, such that\\
 (1) $f^pz=z;$\\
 (2)
$d (f^{c_{i-1}+j} (z),f^j (x_i))<\varepsilon,\,\forall\,j=0,\,1,\,\cdots,n_i,
\,i=1,\,2,\,\cdots,\,N$;\\
 (3) $Orbit (z)\subseteq B (supp (\mu),\varepsilon )$.
\end{Lem}

{\bf Proof of  Lemma~\ref{Lem:weak-specification}}  Given  Pesin
block $\Lambda_k=\Lambda_k (K,\,\zeta)$ and $\varepsilon>0$, by
Theorem~\ref{Thm:shadowinglem} there exists $\delta>0$ such that for
any periodic $\delta$-pseudo-orbit
$\{x_i,\,n_i\}_{i=-\infty}^{+\infty}$ with $n_i\geq2kK $ and
$x_i,f^{n_i} (x_i)\in \Lambda_k (K,\,\zeta)$ $\forall \,\,i$, there is
a $\varepsilon$-shadowing periodic point $x\in M $ for
$\{x_i,n_i\}_{i=-\infty}^{+\infty}$.

Choose and fix for $\widetilde\Lambda_k$ a finite cover
$\alpha=\{V_1,V_2,\cdots,V_{r_0}\}$ by nonempty open balls $V_i$ in
$M$ such that $diam (U_i) <\delta$ and $\mu (U_i)>0$ where
$U_i=V_i\cap \widetilde\Lambda_k$, $i=1,\,2,\,\cdots,\,r_0$. Since
$\mu$ is $f$ ergodic, by Birkhoff ergodic theorem we have
\begin {equation}\label{Eq:erg-application}
\lim_{l\rightarrow +\infty}\frac 1l
  \sum_{n=0}^{l-1}\mu (f^{-n} (U_i)\cap U_j)=\mu (U_i)\mu (U_j)>0.
\end{equation}
Take
\begin {equation}\label{Eq:min-transfer-time} X_{i,j}=\min\{n\in \mathbb{N}\,\,|\,\, n\geq
 2kK, \,\,\mu (f^{-n} (U_i)\cap
  U_j)>0\}.
\end{equation}
By   (\ref{Eq:erg-application}), $1\leq X_{i,j}<+\infty$. Let
$$X_{1}=\min_{1\leq i, j\leq
r_0}X_{i,j}>0,\,\,\,\,\,\,X_{2}=\max_{1\leq i, j\leq
r_0}X_{i,j}>0.$$

Now let us consider a given  sequence of points
$x_1,\,x_2,\,\cdots,\,x_N\in\,\widetilde\Lambda_k,$ and a sequence
of positive numbers ${n_1,n_2,\cdots,n_N}$ satisfying $n_i\geq 2kK$
and $f^{n_i}x_i \in\,\widetilde\Lambda_k.$ Fix
$U_{i_0},U_{i_1}\in\alpha$ so that $$x_i\in U_{i_0},f^{n_i}x_i\in
U_{i_1},i=1,\,2,\,\cdots,\,N.$$ Take $y_i\in U_{i_1}$ by
 (\ref{Eq:min-transfer-time}) such that $f^{X_{ (i+1)_0,i_1}}y_i\in
U_{ (i+1)_0}$ for $i=1,\,2,\,\cdots,\,N-1$ and choose $y_N\in
U_{N_1}$ such that $f^{X_{ (N+1)_0,N_1}}y_N\in U_{1_0}$. Thus we get
a periodic
 $\delta$-pseudo-orbit in $supp (\mu)\subseteq M$:
$$\{f^t (x_1)\}_{t=0}^{n_1}\,\,\,\bigcup\,\,\,
\{f^t (y_1)\}_{t=0}^{X_{2_0,1_1}}\,\,\,\bigcup\,\,\,
\{f^t (x_2)\}_{t=0}^{n_2}\,\,\,\bigcup\,\,\,\{f^t (y_2)\}_{t=0}^{X_{3_0,2_1}}$$$$\,\,\,\bigcup\,\,\,
\cdots\cdots\,\,\,\bigcup\,\,\,\{f^t (x_N)\}_{t=0}^{n_N}
\,\,\,\bigcup\,\,\,\{f^t (y_N)\}_{t=0}^{X_{ (N+1)_0,N_1}}$$ satisfying $$x_i,\,\,
f^{n_i} (x_i),\,\,y_i,\,\,f^{X_{ (i+1)_0,i_1}}y_i\in
\widetilde\Lambda_k\subseteq\Lambda_k\cap supp (\mu)\,\,\, (\forall i).$$ Hence by
Theorem~\ref{Thm:shadowinglem} there exists a periodic point $z\in
M$ with period $p=\sum_{i=1}^{N} (n_i+X_{ (i+1)_0,i_1})$
$\varepsilon$-shadowing the above sequence. This implies the orbit $Orbit (z)$ lies in the $\varepsilon-$neigborhood of above pseudo-orbit so that  $$Orbit (z)\subseteq B (supp (\mu),\varepsilon )$$ and
$$   d (f^{c_{i-1}+j} (z),f^j (x_i))<\varepsilon,\,\,\forall\,j=0,\,1,\,\cdots,n_i,
\,\,\,i=1,\,2,\,\cdots,\,N,$$
 where  $$c_i=\begin{cases}
 \,\,0,&\text{for }\,\,i=0\\
 \,\,\sum_{j=1}^{i}[n_j+X_{ (j+1)_0,j_1}],\,\,&\text{for}\,\,i=1,\,2,\,\cdots,\,N.\\
\end{cases}
 $$ Clearly
$p\in[\sum_{i=1}^{N}n_i+NX_{1},\,\,\sum_{i=1}^{N}n_i+NX_{2}].$
 \hfill $\Box$

\bigskip

\begin{Rem}\label{Rem-NS}  In particular, if  $\mu$ is a mixing hyperbolic measure, we can replace inequality  (\ref{Eq:min-transfer-time}) by \begin {equation}\label{Eq:Remark}
\lim_{n\rightarrow +\infty}\mu (f^{-n} (U_i)\cap U_j)=\mu (U_i)\mu (U_j)>0.
\end{equation}
 Then  by   (\ref{Eq:Remark}) we can
take a finite integer
$$ X_{i,\,j}=\max\{n\in \mathbb{N}\,\,|\,\, n\geq
 1, \,\,\mu (f^{-n} (U_i)\cap
  U_j)=0\}+1.
$$ Let
$$M_{k}=
\max_{1\leq i, j\leq r_k}X_{i,\,j}.$$
Then for any $ N\geq M_k$ there exist $y\in U_j$  such that $f^N (y)\in U_i.$
So we can follow the above proof and then the non-uniform specification can be stronger:
 If for a given sequence of  points $x_1,x_2,\cdots,x_N\in\,
 \widetilde\Lambda_k (N\in\mathbb{N})$ and  a
sequence of positive numbers ${n_1,n_2,\cdots,n_N}$, one has
$n_i\geq 2kK $ and $f^{n_i}x_i \in\,\widetilde\Lambda_k$ for
$i=1,\,2,\,\cdots,\,N,$ then for any sequence of nonnegative numbers $c_0=0,c_1,\cdots,c_{N-1}$ with $c_i-c_{i-1}\geq n_{i-1}+M_k$ and $p\geq \sum_{i=1}^Nn_i+NM_k,$  there exists  a periodic point $z  \in
M$   such that\\
 (1) $f^pz=z;$\\
 (2)
$   d (f^{c_{i-1}+j} (z),f^j (x_i))<\varepsilon,\,\forall\,j=0,\,1,\,\cdots,n_i,
\,i=1,\,2,\,\cdots,\,N$;\\
 (3) $Orbit (z)\subseteq B (supp (\mu),\varepsilon )$.


\end{Rem}

\subsection{Density of Periodic Measures}\label{subsection-densityperiodicmeas}

{\bf Proof of Theorem~\ref{Thm:density-per-meas}} To deduce the
density property of periodic measures, the first two statements in weaker specification
property of Lemma~\ref{Lem:weak-specification} is enough. One can
prove this theorem word by word by same method as  employed in  \cite{Hir} or  \cite{LLS}, just also considering $Orbit (z)\subseteq B (supp (\mu),\varepsilon )$. Here we omit the
details. \hfill $\Box$
\bigskip

 {\bf Proof of Theorem \ref{Thm:2015-quasi-invariant} (3).} One can follow the proof of Theorem \ref{Thm:Horshoe-PosiStaMnfd} to give the proof, by replacing Theorem \ref{Thm:Expclosinglem-measure} by Theorem \ref{Thm:Expclosinglem-measure2015oct} and replacing ergodicity of the measure by quasi-ergodic.

 \subsection{Newhouse's theorem on maximal entropy measure} 

One potential application of topological definition of Pesin set, independent of measures, is possibly useful to find maximal entropy measure in $C^{1}$ diffeomorphisms. Recall that   for $C^{1+\alpha}$ diffeomorphisms, Newhouse  \cite{Newhouse} gives
 necessary and sufficient conditions for the existence of maximal measure.

Let $\{\Lambda_l:=\Lambda_k(\lambda,\mu;\epsilon)\}_{l\geq 1}$ be the classical Pesin set.  Let
$\bar{\varepsilon}=(\varepsilon_1,\varepsilon_2,\cdots,)$ be a nonincreasing sequence of positive real numbers which approach zero. 
Let
$$\m_{\lambda,\mu;\epsilon,\bar{\varepsilon}}=\{\nu\in\m_{f}(M)\,:\,\nu( {\Lambda}_l)\geq1-\varepsilon_l,\,l=1,2,...\}.$$
Since each $ \Lambda_l$ is compact, the map $\nu\to
\nu( \Lambda_l)$ is upper-semicontinuous. Hence,
$\m_{\lambda,\mu;\epsilon,\bar\varepsilon}$ is a closed convex subset of $\m_{f}(\Lambda)$.  Let $V\subseteq \m_f(M)$ be a subset. If $P$ is any finite set of $M$, we say $P$ is related to $V$ if the discrete measure $\frac 1{\text{card} P}\sum_{x\in P} \delta_x$ is in $V.$ Since $V$ consists of invariant measures, if $P$ is related to $V$, then $P$ must be an invariant set. Analogously, we say that a  peridoc point $p$, or its orbit, is related to $V$ is the discrete measure uniformly distributed on its orbit is in $V$.

 \begin{Thm}\label{Newhouse}
Suppose $f\in \Diff^{1+\alpha} (M)$ and $h_{top}(f)>0$. A necessary and sufficient condition for the existence of a hyperbolic measure $\nu$  with $h_{\nu}( f) = h_{top}( f)$ is that there exist $\lambda>0,\mu>0, \epsilon>0$ with $\epsilon\ll \min\{\lambda,\mu\}$,   a sequence $\bar{\varepsilon}$  and a sequence $p_1, p_2, \cdots$ of periodic points related to $\m_{\lambda,\mu;\epsilon,\bar{\varepsilon}}$ such that   $\lim_{n\rightarrow \infty}1/n\log \text{card}\{ p_i: per(p_i) = n\} = h_{top}(f).$

 \end{Thm}

 A natural question is for the case of $C^{1}$ diffeomorphisms.  Similar as $\m_{\lambda,\mu;\epsilon,\bar{\varepsilon}}$, we can define
  $\m_{\zeta,K,\bar{\varepsilon}}$ by replacing $\Lambda_k(\lambda,\mu;\epsilon)$ by Liao-Pesin blocks $\Lambda_k(\zeta,K)$.

 \begin{Que}\label{QueNewhouse}
Suppose $f\in \Diff^{1} (M)$ and $h_{top}(f)>0$. Whether the following two conditions are equivalent? \\
 (1) there is a hyperbolic measure $\nu$ with (limit-)dominated Oseledec's hyperbolic splitting such that  $h_{\nu}( f) = h_{top}( f)$; \\
  (2) there exist $\zeta>0,$ integer $ K\geq 1$,   a sequence $\bar{\varepsilon}$  and a sequence $p_1, p_2, \cdots$ of periodic points related to $\m_{\zeta,K,\bar{\varepsilon}}$ such that   $\lim_{n\rightarrow \infty}1/n\log \text{card}\{ p_i: per(p_i) = n\} = h_{top}(f).$

 \end{Que}

Newhouse's proof \cite{Newhouse} of Theorem \ref{Newhouse} is mainly based on the construction of Pesin blocks, closing lemma and Lyapunov neighborhood over Pesin blocks. Here we have constructed  Liao-Pesin blocks $\Lambda_k(\zeta,K)$ and realize closing lemma over it but it is still unknown how to overcome the role of Lyapunov neighborhood, since Lyapunov neighborhood is a technique for $C^{1+\alpha}$ case. So up to now Question \ref{QueNewhouse} is an open problem.

\section{Appendix A: Limit dominated splitting, dominated splitting and Lyapunov exponents}\label{section-limitdomination}

In this section we mainly discuss certain properties of limit
domination itself, and the relations between  limit domination and
domination in topological perspective. Moreover, we point out that
the limit domination is closely related with the gap of \textit{mean}
expanding Lyapunov exponent on the unstable bundle and \textit{mean}
contracting Lyapunov exponent on the stable bundle.

\subsection{Limit-dominated and dominated splitting}

Similar as the equivalent definition of dominated splitting, we have following equivalent statements for limit domination.
Let $$\alpha = \max_{x\in M}\,\,
\log\frac{\|Df_x\|}{m (Df_x)}\,.$$ Clearly $\alpha\geq 0.$

\begin{Prop}\label{Prop:property-limitdom-itself-1} Given  an $f$
invariant set $\Delta$, assume that there is a limit-dominated
 $Df$-invariant splitting $T_{x}M=E (x)\oplus
 F (x)$ on $\Delta$.  Then the following properties hold.

 (1) there exists  $\bar{\lambda}\in (0,1)$ and  $K_0>0$ such that for any $K\geq K_0,$
$$\limsup_{l\rightarrow+\infty}  \frac
{\|Df^K|_{E (f^{l} (x))}\|}{m (Df^K|_{F (f^{l} (x))})}\leq  \bar{\lambda}^K,\,\,\forall
x \in \Delta.$$

  (2)there exists  $\bar{\lambda}\in (0,1)$ and  $\bar{C}>0$ such that for any $K\geq 1,$
$$\limsup_{l\rightarrow+\infty}  \frac
{\|Df^K|_{E (f^{l} (x))}\|}{m (Df^K|_{F (f^{l} (x))})}\leq  \bar{C }\bar{\lambda}^K,\,\,\forall
x \in \Delta.$$

\end{Prop}

{\bf Proof}
 By assumption,  there
exists $S\in \mathbb{Z}^+,\,\lambda\in (0,1)$  such that
$$\limsup_{l\rightarrow+\infty}  \frac
{\|Df^S|_{E (f^{lS} (x))}\|}{m (Df^S|_{F (f^{lS} (x))})}\leq  \lambda^S,\,\,\forall
x \in \Delta.$$

 (1) Let $K=kS+q.$ Since  $T_{x}M=E (x)\oplus  F (x)$ is
$ (S,\lambda)$-limit-dominated and  $$\frac
{\|Df^K|_{E (f^{l} (x))}\|}{m (Df^K|_{F (f^{l} (x))})}\leq
 (\prod_{i=0}^{k-1}\frac
{\|Df^S|_{E (f^{l+iS} (x))}\|}{m (Df^S|_{F (f^{l+iS} (x))})})\times\frac
{\|Df^q|_{E (f^{l+kS} (x))}\|}{m (Df^q|_{F (f^{l+kS} (x))})}$$
$$\,\,\,\,\,\,\leq (\prod_{i=0}^{k-1}\frac
{\|Df^S|_{E (f^{l+iS} (x))}\|}{m (Df^S|_{F (f^{l+iS} (x))})})\times
e^{q\alpha},$$ by assumption and   (\ref{def-equal-limit-domination}) one has
$$
\limsup_{l\rightarrow+\infty} \frac
{\|Df^K|_{E (f^{l} (x))}\|}{m (Df^K|_{F (f^{l} (x))})}\leq\lambda^{kS}e^{S\alpha},\,\,
\forall x \in \Delta.
$$ Take and fix $K_0> S\frac{\log e^\alpha-\log \bar{\lambda}}{\log \bar{\lambda}-\log  {\lambda}}.$ Then for any $K\geq K_0,$ $$
\limsup_{l\rightarrow+\infty} \frac
{\|Df^K|_{E (f^{l} (x))}\|}{m (Df^K|_{F (f^{l} (x))})}\leq\bar{\lambda}^{K},\,\,
\forall x \in \Delta.
$$

  (2) Take same numbers as in  (1). Let $\bar{C}=\frac{e^{K_0\alpha}}{\bar{\lambda}^{K_0}}$
  Then $\bar{C}\geq 1$ and  for any $K\geq 1,$ $$
\limsup_{l\rightarrow+\infty} \frac
{\|Df^K|_{E (f^{l} (x))}\|}{m (Df^K|_{F (f^{l} (x))})}\leq \max\{ e^{K\alpha}, \bar{\lambda}^{K}\}\leq \bar{C}\bar{\lambda}^{K} ,\,\,
\forall x \in \Delta.
$$ \qed

Similar to the case of domination whose splitting is unique if one
fixes the dimensions of the subbundles  (see  \cite{BLV}), the
limit-dominated splitting is also unique if one fixes the dimensions
of the subbundles.  The
following proposition points out more properties that limit
domination has in common with domination. For two subbundles $E (x)$ and $F (x)$, define the angle $$\angle(E (x),F (x))= \inf\{\|u-v\|: \,u\in E (x), v\in
F (x),\|u\|=1\,\text{or} \,\|v\|=1\}.$$
Let $$Ang (E (x),F (x))=\inf\{\|u-v\|:\|u\|=\|v\|=1,\,u\in E (x), v\in
F (x)\}.$$ Remark that $\angle(E(x),F (x))\leq Ang (E (x),F (x))\leq \sqrt{2}\angle(E(x),F (x)).$

 \begin{Prop}\label{Prop:property-limitdom-itself-2} Let us make the same assumptions as in Proposition
\ref{Prop:property-limitdom-itself-1}.    Then  there exists $e_1>0$ such that
$\liminf_{n\rightarrow+\infty}\angle (E (f^{n} (x)),F (f^{n} (x)))\geq
e_1,\,\,\,\forall x\in \Delta.$
\end{Prop}
{\bf Proof.} By equivalence of $\angle (\cdot,\cdot)$ and $Ang (\cdot,\cdot)$, we only need to prove the result for $Ang (\cdot,\cdot).$
 By assumption,  there
exists $S\in \mathbb{Z}^+,\,\lambda\in (0,1)$  such that
$$\limsup_{l\rightarrow+\infty}  \frac
{\|Df^S|_{E (f^{lS} (x))}\|}{m (Df^S|_{F (f^{lS} (x))})}\leq  \lambda^S,\,\,\forall
x \in \Delta.$$

Let $$c=\min_{x\in M}\,m (D_xf^S),\,\,\,\,\,\,\,C=\max_{x\in
M}\|D_xf^S\|.$$ Clearly $c,\,\,C \in  (0, +\infty).$ By continuity of
the tangent bundle $T_xM$, there exists real number $e_0>0$ such
that if $\|u-v\|<e_0,\,\|u\|=\|v\|=1,\,u,\,v\in T_xM$ then
$$
\big{|}\|D_xf^S (u)\|-\|D_xf^S (v)\|\big{|}<c (1-e^{-\lambda}),\,\forall
x \in M,$$ which implies
\begin {equation}\label{eq:anti-dominated}
\frac{\|D_xf^S (u)\|}{\|D_xf^S (v)\|}\geq
\frac{\|D_xf^S (v)\|-c (1-e^{-\lambda} )}{\|D_xf^S (v)\|}\geq
e^{-\lambda}.
\end{equation}

 Since  $T_{x}M=E (x)\oplus  F (x)$ is
$ (S,\lambda)$-limit-dominated, for any $x\in \Delta,$ there exists
an integer $N (x)\geq 1$ such that for  $n\geq N (x),$
\begin {equation}\label{eq:dominated}
\log\frac
{\|Df^S|_{E (f^{nS} (x))}\|}{m (Df^S|_{F (f^{nS} (x))})}<-2\lambda+\lambda=-\lambda.
\end{equation}
For $n\geq N (x)$ and two vectors $u\in E (f^{nS} (x)), v\in
F (f^{nS} (x))$ with $\|u\|=\|v\|=1,$ we claim that $\|u-v\|\geq e_0.$
Otherwise, it holds that $\|u-v\|< e_0,$ and thus by the inequality
 (\ref{eq:anti-dominated}),
$$\frac{\|Df^S|_{E (f^{nS} (x))}\|}{m (Df^S|_{F (f^{nS} (x))})}\geq
\frac{\|D_xf^S (u)\|}{\|D_xf^S (v)\|}\geq e^{-\lambda},$$ which
contradicts  (\ref{eq:dominated}). So, we have
$$Ang (E (f^{nS} (x)),F (f^{nS} (x)))$$$$=\inf\{\|u-v\|:\|u\|=\|v\|=1,\,u\in
E (f^{nS}x), v\in F (f^{nS}x)\}\geq e_0,$$ for all $n\geq N (x)$ and
therefore
$\liminf_{n\rightarrow+\infty}Ang (E (f^{nS} (x)),F (f^{nS} (x)))\geq e_0,\,\, x\in \Delta.$
By the invariance of $\Delta,$ we have $\liminf_{n\rightarrow+\infty}Ang (E (f^{n} (x)),F (f^{n} (x)))\geq
e_0,$ $\forall x\in \Delta.$
\hfill $\Box$

\bigskip
Clearly,  domination implies limit
domination and the later is  weaker. The following proposition
focuses on the inverse implication, provided that the space is compact and the splitting is continuous.

\begin{Prop}\label{Prop:relation-dom-limdom}
Let us make the same assumptions as in Proposition
\ref{Prop:property-limitdom-itself-1}. If further $\Delta$ is compact
and  $T_{x}M=E (x)\oplus F (x)$  is continuous on $\Delta $,
then $T_{x}M=E (x)\oplus  F (x)$ is   dominated on $\Delta$.

\end{Prop}

 {\bf Proof.} One can use Lemma \ref{Lem:1} to prove but here we prefer to give a direct proof.

 Since  $T_{x}M=E (x)\oplus  F (x)$ is
limit-dominated, there
exists $S\in \mathbb{Z}^+,\,\lambda\in (0,1)$  such that
$$\limsup_{l\rightarrow+\infty}  \frac
{\|Df^S|_{E (f^{lS} (x))}\|}{m (Df^S|_{F (f^{lS} (x))})}\leq  \lambda^S,\,\,\forall
x \in \Delta.$$
Take $\chi=-\frac12 \log \lambda^S.$ For any $\delta>0$, there exists an
integer $l (x)\geq 1$ such that for all $l\geq l (x),$ $$\log\frac
{\|Df^S|_{E (f^{lS} (x))}\|}{m (Df^S|_{F (f^{lS} (x))})}\leq-2\chi+\delta,\,\,x\in\Delta,
$$
which implies
$$\frac1{l-l (x)}\sum_{i=l (x)}^{l-1}\log\frac
{\|Df^S|_{E (f^{iS} (x))}\|}{m (Df^S|_{F (f^{iS} (x))})}\leq
-2\chi+\delta,\,\,\,\forall\,\, l>l (x).$$
 Letting
$l\rightarrow +\infty$ we have $$\limsup_{l\rightarrow
+\infty}\frac1{l}\sum_{i=0}^{l-1}\log\frac
{\|Df^S|_{E (f^{iS} (x))}\|}{m (Df^S|_{F (f^{iS} (x))})}\,\,\,\,\,\,\,\,\,\,\,\,\,\,\,\,\,\,\,\,\,\,\,\,\,\,\,\,\,\,\,\,\,$$
$$=\limsup_{l\rightarrow
+\infty}\frac1{l-l (x)}\sum_{i=l (x)}^{l-1}\log\frac
{\|Df^S|_{E (f^{iS} (x))}\|}{m (Df^S|_{F (f^{iS} (x))})}\leq
-2\chi+\delta.$$ Letting $\delta\rightarrow 0$ one has
$$\limsup_{l\rightarrow
+\infty}\frac1{l}\sum_{i=0}^{l-1}\log\frac
{\|Df^S|_{E (f^{iS} (x))}\|}{m (Df^S|_{F (f^{iS} (x))})}\leq -2\chi.$$
Thus we can take $n (x)\geq 1$  such that
$$\frac1{n (x)}\sum_{i=0}^{n (x)-1}\log\frac
{\|Df^S|_{E (f^{iS} (x))}\|}{m (Df^S|_{F (f^{iS} (x))})}<-2\chi+\chi=-\chi,\,\,x\in
\Delta.
$$

Since $TM=E\oplus F$ is continuous on  $\Delta$, there exists a
neighborhood $V_x$ of $x$ such that for every $y \in V_x$ one has
$$\frac1{n (x)}\sum_{i=0}^{n (x)-1}\log\frac
{\|Df^S|_{E (f^{iS} (y))}\|}{m (Df^S|_{F (f^{iS} (y))})}<-\chi.
$$
We take a finite cover  $\{V_{x_1},...,V_{x_q}\}$  for  the compact
$\Delta$ and let   $ N = \max\{n (x_1)\,,\,...\,,\,n (x_q)\}.$ Let
$$\gamma = max_{x\in \Delta} |\log\frac
{\|Df^S|_{E (x)}\|}{m (Df^S|_{F (x)})}|.$$ Then $\gamma <\infty$
because of the continuity of splittings $E$ and $F$ and the
compactness of $\Delta$. We define inductively a sequence  $N_k :
\Delta \rightarrow\mathbb{N}$ by
$$ N_0 (x) = 0,\,\,N_1 (x) = min\{n (x_i) : x \in V_{x_i} ,i =
1,...,q\},$$$$N_{k+1} (x) = N_k (x)+N_1 (f^{N_k (x)S}  (x)),\,\,k\geq
1.$$ Thus, for all $x\in \Delta$ and $n$, there exists $k$ such that
$N_k (x)\leq n< N_{k+1} (x)$. Hence
$$\sum_{i=0}^{n-1}\log\frac
{\|Df^S|_{E (f^{iS} (x))}\|}{m (Df^S|_{F (f^{iS} (x))})}
<-N_k (x)\chi+ (n-N_k (x))\gamma\leq-n\chi+N (\chi+\gamma).$$
By taking $n_0=[2+\frac{N (\chi+\gamma)}\chi]+1$ where $[a]$
denotes the maximal integer less than or equal to $a$, we then have
$$\log\frac
{\|Df^{nS}|_{E (x)}\|}{m (Df^{nS}|_{F (x)})}\leq\sum_{i=0}^{n-1}\log\frac
{\|Df^S|_{E (f^{iS} (x))}\|}{m (Df^S|_{F (f^{iS} (x))})}\leq-2\chi$$
for all $x\in \Delta$ and $n\geq n_0$, which deduces from  (\ref{eq-domination-equvialent2}) that
$T_xM=E (x)\oplus F (x)$ is  dominated on $\Delta$.
\hfill $\Box$

\bigskip

However, if we do not assume $\Delta$ to be compact and the splitting to be continuous, it is unknown the above proposition.  Here we only have a following result.
 Let $\omega (\Delta)= \bigcup_{x\in \Delta} \omega (x) .$ It is an invariant set.

\begin{Prop}\label{prop-limitdom-imply-dom-omegaset}
 Let us make the same assumptions as in Proposition
\ref{Prop:property-limitdom-itself-1}.    Then
there exists $n_0>0$ such that for every $n\geq n_0,$
$T_{x}M=E (x)\oplus  F (x)$ is $ (nS,\lambda)$ dominated on the closure of $\omega (\Delta)$.

\end{Prop}

{\bf Proof.}
 By assumption,  there
exists $S\in \mathbb{Z}^+,\,\lambda\in (0,1)$  such that
$$\limsup_{l\rightarrow+\infty}  \frac
{\|Df^S|_{E (f^{l} (x))}\|}{m (Df^S|_{F (f^{l} (x))})}\leq  \lambda^S,\,\,\forall
x \in \Delta.$$
Take $\hat{\lambda}\in  (\lambda, 1)$. Define
 $$\Delta_N:=\{x\in \Delta|\,\frac
{\|Df^S|_{E (f^{l} (x))}\|}{m (Df^S|_{F (f^{l} (x))})}\leq  \hat{\lambda}^S,\, \forall \,l\,\,\, \geq N\}.$$ Then
$\Delta_N\subseteq \Delta_{N+1}.$ $\Delta=\cup_{N\geq 1}\Delta_N.$ Let $\omega (\Delta_N)=\cup_{x\in \Delta_N} \omega (x)$. Then $\omega (\Delta)=\cup{N\geq 1}\omega (\Delta_N).$
 Similar as domination, the splitting on $\Delta_N$ should be unique and can be extended on the closure. This implies that there is an extended splitting $E\oplus F$ on $\omega (\Delta_N)$ such that for any $y\in\omega (\Delta_N), $ $$\frac
{\|Df^S|_{E (y)}\|}{m (Df^S|_{F (y)})}\leq  \hat{\lambda}^S.$$
By the arbitrariness of $N$, the extended splitting  $E\oplus F$ on $\omega (\Delta$ satisfies that
 for any $y\in\omega (\Delta), $ $$\frac
{\|Df^S|_{E (y)}\|}{m (Df^S|_{F (y)})}\leq  \hat{\lambda}^S.$$
This is the general domination so that one has same estimate on the closure of $\omega (\Delta)$. \qed

 \bigskip

A point $x$ is recurrent
if there is $n_i\uparrow+\infty$ such that $\lim_{i\rightarrow
+\infty}f^{n_i} (x)=x$. Denote the set of all recurrent points by $Rec (f).$ It is well-known that for any invariant measure, $Rec (f)$ is of full
measure for any invariant measure.  Let $\mu$ be an invariant measure with a limit-dominated splitting $E\oplus F$ on some invariant set $\Delta$ with $\mu$ full measure. Then by Proposition \ref{prop-limitdom-imply-dom-omegaset} this splitting can be extended as a dominated splitting on the closure of $\omega (\Delta).$ Note that $Rec (f)\cap\Delta\subseteq \omega (\Delta)$ is of $\mu$ full measure so that  $\omega (\Delta)$ and its closure are   of $\mu$ full measure. So, in particular,  the splitting can be extended as a dominated splitting on support of $\mu.$ So

\begin{Cor}\label{cor-hyp-domination-same-limitdom000000} $  \m_f^{ldh} (M)=\m_f^{dh} (M)$.

\end{Cor}




Moreover, we have
\begin{Cor}\label{cor-hyp-domination-same-limitdom}     $  \m_f^{qdh} (M)=\m_f^{qldh} (M).$

\end{Cor}

{\bf Proof.} We only need to the direction $\supseteq.$  Let $\mu\in \m_f^{qldh} (M)$ and define
$$\Delta_{S,\lambda}=\{x|\,\,\limsup_{l\rightarrow+\infty}  \frac
{\|Df^S|_{E^u (f^{l} (x))}\|}{m (Df^S|_{E^u (f^{l} (x))})}\leq  \lambda^S\}.$$
Then $\mu (\cup_{\lambda \in (0,1)}\cup_{S\geq 1} \Delta_{S,\lambda})=1.$ Let $\mu_{\lambda,S}=\mu|_{ \Delta_{S,\lambda}}.$ Then $\mu_{\lambda,S}$ has limit-dominated hyperbolic Oseledec splitting. By Corollary \ref{cor-hyp-domination-same-limitdom000000}, $\mu_{\lambda,S}$ has dominated hyperbolic Oseledec splitting.
Thus, by equivalent statements of domination  (\ref{eq-domination-equvialent1}), there is some $\zeta>0$ and $S_0\geq 1$ such that
 $$\mu_{\lambda,S}\big{ (}\, \{x\in M \,\,|\,\,\frac1L\log\frac
{\|Df^L|_{E^s (y)}\|}{m (Df^L|_{E^u (y)})}\leq-3\zeta,\,\,\forall y\in Orb (x),\,\,L\geq S_0\}\,\,\big{)}=1.$$
By arbitrariness of $\lambda,S,$
 $$\mu\big{ (}\,\cup_{\zeta>0}\cup_{S_0\geq 1}\{x\in M \,\,|\,\,\frac1L\log\frac
{\|Df^L|_{E^s (y)}\|}{m (Df^L|_{E^u (y)})}\leq-3\zeta,\,\,\forall y\in Orb (x),\,\,L\geq S_0\}\,\,\big{)}=1.$$


In particular, we sate a following corollary for hyperbolic measures with quasi-limit-domination.

\begin{Cor}\label{cor-hyp-quasilimitdomination} Let $\mu\in \m_f^{qldh} (M).$  Then its Oseledec hyperbolic splitting $T_xM=E^s (x)\oplus E^u (x)$ satisfies that
$$\mu\big{ (}\,\cup_{\zeta>0}\cup_{S_0\geq 1}\{x\in M \,\,|\,\,\frac1S\log\frac
{\|Df^S|_{E^s (y)}\|}{m (Df^S|_{E^u (y)})}\leq-3\zeta,\,\,\forall y\in Orb (x),\,\,S\geq S_0\}\,\,\big{)}=1.$$

\end{Cor}

\bigskip
Though there are no differences between limit domination and domination in the probabilistic perspective, it is unknown whether it is same from the topological perspective. Limit domination is
weaker in the geometric or topological  perspective, because limit domination only require that $E$ can dominate $F$ for  large enough positive iterate  of the orbit (e.g., see
Example~\ref{Ex:nonuni-hyp-sys}). This is the general case. Now we consider one particular case. Assume that there is  a global limit-dominated splitting  which coincides with the hyperbolic splitting in  an Axiom A system, it is still unknown whether the splitting is a global dominated splitting so that it is Anosov. That is,

\begin{Que}  Let $f\in \Diff^1 (M)$ be of Axiom A.   Suppose there is a limit-dominated
 $Df$-invariant splitting $T_{x}M=E (x)\oplus
 F (x)$ on $M$. If this splitting  restricted on the non-wandering set $\Omega (f)$ coincides with the corresponding stable bundle and unstable bundle of $f$, whether $f$ is Anosov?

\end{Que}

Up to now, it is easy to know  that
\begin{Thm}  Let $f\in \Diff^1 (M)$ be of Axiom A.   Suppose there is a  dominated (or just continuous)
 $Df$-invariant splitting $T_{x}M=E (x)\oplus
 F (x)$ on $M$. If this splitting  restricted on the non-wandering set $\Omega (f)$ coincides with the corresponding stable bundle and unstable bundle of $f$, then  $f$ is Anosov.

\end{Thm}

{\bf Proof.}  For Axiom A system, since every invariant measure should be supported on the non-wandering set $\Omega (f)$, then every invariant measure should be hyperbolic.  Recall the result of  \cite{Cao} that

\begin{Lem}\label{lem-Cao2003}  Let $f\in \Diff^1 (M)$ and suppose there is a  dominated  (or just continuous)
 $Df$-invariant splitting $T_{\Lambda}M=E (x)\oplus
 F (x)$ on a compact invariant set $\Lambda\subseteq M$. If for every invariant measure, its Lyapunov exponents restricted on $E$ are negative and its   Lyapunov exponents restricted on $F$ are positive, then $\Lambda$ is   uniformly hyperbolic.
\end{Lem}

For our case the set  $\Lambda$ in Lemma \ref{lem-Cao2003} is $M$ so that $M$ should be uniformly hyperbolic. This means $f$ to be Anosov. \qed

\bigskip

Following same idea,  one should have a following result.
\begin{Thm} Let $f\in \Diff^1 (M)$ and suppose there is a  limit-dominated
 $Df$-invariant splitting $T_{\Lambda}M=E (x)\oplus
 F (x)$ on a  compact invariant set $\Lambda\subseteq M$. If for every invariant measure, its Lyapunov exponents restricted on $E$ are negative and its   Lyapunov exponents restricted on $F$ are positive, then $\overline{\omega (\Lambda)}\subseteq \Lambda$ is   uniformly hyperbolic.
\end{Thm}
{\bf Proof.} One one hand, from Proposition \ref{prop-limitdom-imply-dom-omegaset}, the splitting $T_{x}M=E (x)\oplus
 F (x)$ is  dominated on $\overline{\omega (\Lambda)}$ so that it is continuous.
 On the other hand, notice that $\m_f (\overline{\omega (\Lambda)})=\m_f (\Lambda)$ so that every invariant measure supported on
 $\overline{\omega (\Lambda)}$ is hyperbolic. By Lemma \ref{lem-Cao2003}, $\overline{\omega (\Lambda)}$ should be uniformly hyperbolic.\qed

 \bigskip

To further illustrate the differences between limit domination and domination, we
construct a simple example as follows.

\begin{Ex}\label{Ex:nonuni-hyp-sys} Let $g$ be a $C^r  (r\geq1)$ increasing function on $[0,1]$,
 satisfying: $$g (0)=0,\,\,g' (0)=\frac12,\,\,g (1)=1,\,\,g' (1)=\frac12,\,\,g (\frac12)=\frac12,
 \,\,g' (\frac12)=\frac{3+\sqrt[]{5}}2,\,\,\,and$$
 $$g (x)<x,\,\,for\,\,all\,\,
 x\in (0,\frac12),\,\,g (x)>x,\,\,for\,\,all\,\,x\in (\frac12,1).$$
 \setlength{\unitlength}{1mm}
  \begin{figure}[!htbp]
  \begin{center}
  \begin{picture} (80,60) (0,0)
  \put (0,0){\scalebox{1}[1]{\includegraphics[0,0][30,40]{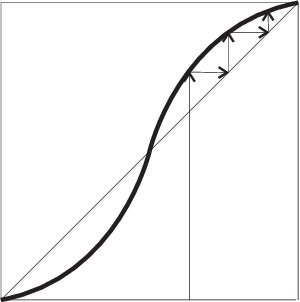}}}
  \put (6.6,54){$_{1}$}
\put (6.6,2.8){$_{0}$}\put (57.9,2.8){$_{1}$}\put (39.2,2.8){$_{x_0}$}
  \end{picture}
  \caption{Graph of the function $g$.}
  \label{pic:graph-of-g}
  \end{center}
  \end{figure}
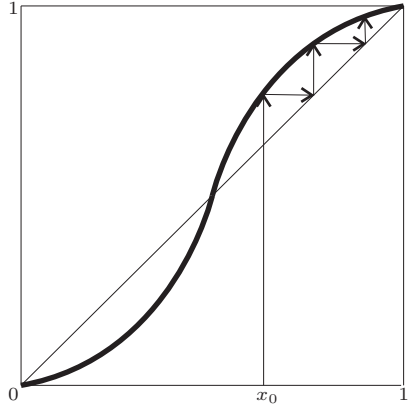
And let $h: T^2\rightarrow T^2$ be the hyperbolic Torus automorphism
$$ (y,z)\mapsto (2y+z,y+z),\,\, y,z\in S^1=\mathbb{R}/\mathbb{Z}.
 $$ Define $f=g\times h:T^3\rightarrow
T^3$. Clearly,  $$Df (x,y,z)= \left (
\begin {array}{ccc}
g' (x)&0 &0\\
0&2&1\\
0&1&1
\end {array}
\right).
$$

There exists naturally a continuous splitting $TT^3=E_1\oplus
E_2\oplus E_3$, where $E_2$ and $E_3$ are from  the hyperbolic Torus
automorphism $h$ and $E_1$ is $g-$invariant. The forward  Lyapunov
exponent of $E_1$ is $\log\frac12$ over $T^3-\{\frac12\}\times T^2$
 (we only use the forward  Lyapunov exponent for the stable subbundle
in construction of our Liao-Pesin  set in section 3) and the Lyapunov
exponents of $ E_2\oplus E_3$ over $\{0\}\times T^2$ are
$\,\,\log\frac{3-\sqrt[]{5}}2,\,\,\log\frac{3+\sqrt[]{5}}2$
respectively. Set $E^s=E_1\oplus E_2$ and $E^u=E_3$, then $E^s\oplus
E^u$ construct a continuous $Df$-invariant splitting of $TT^3$ over
the whole space $T^3$ and the two subbundles are
 limit-dominated over
$T^3-\{\frac12\}\times T^2$  (hint:
$\log\frac{\frac12}{\frac{3+\sqrt[]{5}}2}=-2\log ({3+\sqrt[]{5}})^{\frac12}$).
However, the limit-dominated property can not be extended to the
whole space $T^3$ because the splitting over every point
$ (\frac12,y,z)\,\, (y,z \in S^1)$ does not have limit
domination (hint:
$\log\frac{\frac{3+\sqrt[]{5}}2}{\frac{3+\sqrt[]{5}}2}=0$), even
though the splitting is continuous on the whole space $T^3$. In this
example, the maximal $f-$invariant set admitting domination is
$\{0\}\times T^2$. Moreover, every point $ (0,p)$ is hyperbolic
periodic point where $p$ is a hyperbolic periodic point for $h$.
Denote by $\delta_0$ the point measure at point $0\in S^1$ and
denote by $m$ the Lebesgue measure on $T^2$, then the product
measure $\mu=\delta_0\times m$ is a hyperbolic ergodic measure of
the diffeomorphism $f=g\times h$ with three nonzero Lyapunov exponents
$-\log2,\,\,-\log\frac{3+\sqrt[]{5}}2,\,\,\log\frac{3+\sqrt[]{5}}2$. Moreover, this $\mu$ is the unique SRB-like measure (or SRB, but not absolutely continuous to Lebesgue measure)  and it is easy to see Pesin's entropy formula holds for   $\mu$. So this example gives a positive direction for Question \ref{Que-PesFormula:1}.
 \hfill
$\Box$
\end{Ex}
\begin{Rem}\label{Rem:Example-Pesinset} In Example~\ref{Ex:nonuni-hyp-sys} the Oseledec basin of hyperbolic ergodic
measure $\mu=\delta_0\times m$ is
$L (\mu)=\textrm{supp} (\mu)=\{0\}\times T^2$. Taking $K=1$ and
$0<\zeta<\log2$, our  Pesin set $\Lambda (1,\zeta)$ is
$T^3\setminus\{\frac12\}\times T^2$ and every Pesin block
$\Lambda_k (1,\zeta)$ is $[0,a_k]\bigcup[b_k,1]\times T^2$ for some
$a_k\in (0,\frac12)$ and $b_k\in (\frac12,1).$  Clearly
$\Lambda (K,\,\zeta)- (L (\mu)\bigcup \textrm{supp} (\mu))$ is not empty
and  has Lebesgue full measure.
\end{Rem}

\subsection{Some relations on Lyapunov exponents \&  (Limit-)Dominated Splitting}

Let $\Delta$ be an $f-$invariant set
and $T_{\Delta}M=E\oplus F$ be a $Df-$invariant splitting on
$\Delta$. For $K\geq 1,\,x \in \Delta$ define
$$ \lambda (K,x):=\limsup_{l\rightarrow+\infty}\frac1K\log\frac
{\|Df^K|_{E (f^{lK} (x))}\|}{m (Df^K|_{F (f^{lK} (x))})}.
$$
It is easy to see $T_{\Delta}M=E\oplus F$ is a limit-dominated (or quasi-limit-dominated) splitting on $\Delta$ $\Leftrightarrow$ there is some $K\geq 1$ such that $\sup_{x\in\Delta}\lambda (K,x)<0 (\text{ or for any }x\in \Delta, \text{ there exists } K (x)\geq 1,\,\, \lambda (K (x),x)<0).$

Given an $f$-invariant measure $\nu$ with $\nu (\Delta)=1$ (or just considering $f^K$-invariant measure if necessary), by sub-additional ergodic theorem,
the following limits
exist for $\nu$ a.e. $x$:
$$\lim_{l\rightarrow
+\infty}\frac 1{lK}{\log\|Df^{lK}|_{E (x)}\|}=\lim_{l\rightarrow
+\infty}\frac 1l{\log\|Df^{l}|_{E (x)}\|},$$$$
\lim_{l\rightarrow +\infty}\frac 1{lK}{\log m (Df^{lK}|_{F  (x)})}=\lim_{l\rightarrow +\infty}\frac 1l{\log m (Df^{l}|_{F  (x)})},$$
 and
$$\lim_{l\rightarrow +\infty}\frac
1{lK}\sum_{j=0}^{l-1}{\log\|Df^{K}|_{E ({f^{ jK} (x)})}\|},\,\,\,\,
\lim_{l\rightarrow +\infty}\frac 1{lK}\sum_{j=0}^{l-1}{\log
m (Df^{K}|_{F  ({f^{ jK} (x)})})},$$ which are denoted respectively by
$\lambda_E (x)$, $\lambda_F (x)$ and
 $\lambda^E (K,x)$, $\lambda^F (K,x)$.
   Clearly the functions $\lambda^E (K,x)$, $\lambda_E (x)$, $\lambda_F (x)$,  $\lambda^F (K,x)$ and
$\lambda (K,x)$  are $f^K$-invariant (if $\nu$ is $f^K$-ergodic, then
$\lambda (K,x),\lambda^E (K,x),\lambda^F (K,x)$ are constants $\nu\,\, a.
e.\,\,x $)  and
 $ \lambda_E (x)\leq\lambda^E (K,x)$ and $ \lambda_F (x)\geq\lambda^F (K,x).$
 By Lemma \ref{Lem:Generalized-Multi-erg-thm},  $\nu-a.\,e.\,\,x$ satisfies that
 $$\lambda_E (x)=\lim_{K\rightarrow \infty}\lambda^E (K,x)\text{ and } \lambda_F (x)=\lim_{K\rightarrow \infty}\lambda^F (K,x).$$

  If $\lambda^E (K,x)<0$ and $\lambda^F (K,x)>0$ hold for $\nu-a.\,e.\,\,x$, we say,
for convenience,  that $E$ is mean contracting, $F $ is mean
expanding and $ (f^K, \nu)$ is mean hyperbolic.

The concept of
limit  domination relates in a  more natural way to  mean
expansion and mean contraction by the following proposition. Given an $f$-invariant measure $\mu,$  we have

\begin{Prop}\label{Prop:relation-limdom-Lyaexp}
 (1) For $\nu$ a.e.  $x $,
$\lambda (K,x)\geq\lambda^E (K,x)-\lambda^F (K,x)\geq\lambda_E (x)-\lambda_F (x)$.

 (2)
 For $\nu$ a.e.  $x $, if the limit
$$\lim_{l\rightarrow+\infty}\log\frac
{\|Df^K|_{E (f^{lK} (x))}\|}{m (Df^K|_{F  (f^{lK} (x))})}$$ exists,
then $\lambda (K,x)=\lambda^E (K,x)-\lambda^F (K,x).$
\end{Prop}

{\bf Proof}  (1) By the definitions of $\lambda^E (K,x)$ and
$\lambda^F (K,x)$ one has
$$\lim_{l\rightarrow
+\infty}\frac1{l}\sum_{i=0}^{l-1}\log\frac
{\|Df^K|_{E (f^{iK} (x))}\|}{m (Df^K|_{F  (f^{iK} (x))})}=\lambda^E (K,x)-\lambda^F (K,x).$$
For  $\varepsilon>0$, by definition of $\lambda (x)$ there exists  a
positive integer $l (x)\geq1$ such that
 $$\log\frac
{\|Df^K|_{E (f^{lK} (x))}\|}{m (Df^K|_{F  (f^{lK} (x))})}\leq
\lambda (x)+\varepsilon, \,\,\,\,  \forall \,\, l\geq l (x),$$ which
implies
$$\frac1{l-l (x)}\sum_{i=l (x)}^{l-1}\log\frac
{\|Df^K|_{E (f^{iK} (x))}\|}{m (Df^K|_{F  (f^{iK} (x))})}\leq
\lambda (x)+\varepsilon,\,\,\,\forall \,\,l>l (x).$$  Letting
$l\rightarrow +\infty$ we have
$$\lim_{l\rightarrow +\infty}\frac1{l}\sum_{i=0}^{l-1}\log\frac
{\|Df^K|_{E (f^{iK} (x))}\|}{m (Df^K|_{F  (f^{iK} (x))})}$$$$=\lim_{l\rightarrow
+\infty}\frac1{l-l (x)}\sum_{i=l (x)}^{l-1}\log\frac
{\|Df^K|_{E (f^{iK} (x))}\|}{m (Df^K|_{F  (f^{iK} (x))})}\leq
\lambda (x)+\varepsilon.$$ Letting $\varepsilon\rightarrow 0$,
$$\lambda^E (K,x)-\lambda^F (K,x)\leq
\lambda (x).$$

 (2) Since the limit exists, one can do same work as  the proof of   (1) for another direction. \hfill $\Box$
\bigskip

Proposition \ref{Prop:relation-limdom-Lyaexp}  (1) tells us that
 if  there is an invariant measure $\mu$ and $K\geq 1$ such that for $\mu$ a.e. $x$, $$0>\lambda_E (x)>\lambda (K,x) (\text{ or } 0>-\lambda_F (x)>\lambda (K,x)),$$
then $\mu$ is hyperbolic with quasi-limit-dominated splitting.

By Proposition \ref{Prop:relation-limdom-Lyaexp} (2), we have a following proposition when the system naturally has
 (quasi-)limit-domination. For convenience to state, we introduce a concept. Let $T_{orb (x)}M=E^s\oplus E^u$ be a $Df$-invariant splitting. We say $T_{orb (x)}M=E \oplus F$ to be {\it compatible}, if  there is some $K (x)\geq 1$ such that $K\geq K_0,$
the limit $$\lim_{l\rightarrow+\infty}\log\frac
{\|Df^K|_{E (f^{lK} (x))}\|}{m (Df^K|_{F  (f^{lK} (x))})}$$ exists.

\begin{Prop}\label{prop-natural-limit-dominated}
  Let $f\in \Diff^1 (M)$ and   $\mu\in \mathcal{M}^{h}_{f} (M)$. If  for $\mu$ a.e $x$,  the Oseledec hyperbolic splitting of $\mu$ at $x$ $T_{orb (x)}M=E^s\oplus E^u$ is compatible,  then this splitting is $\mu$-quasi-limit-dominated.
In other words, $\mu \in \mathcal{M}^{qldh}_{f} (M).$
\end{Prop}

{\bf Proof.} By Lemma \ref{Lem:Generalized-Multi-erg-thm} and the hyperbolicity of $\mu$, $\mu$ a.e. $x$, there is some $K (x)\geq 1$ such that for any $K\geq K (x),$ $\lambda^{E^s} (K,x)<0,\lambda^{E^u} (K,x)>0$  (close to $\lambda_{E^s} (x),\lambda_{E^u} (x)$ enough).
By Proposition \ref{Prop:relation-limdom-Lyaexp} (2), $\lambda (K,x)<0.$ We complete the proof. \qed

\smallskip

In particular, we can state this proposition in another way.

\begin{Prop}\label{prop-natural-limit-dominated-2}
  Let $f\in \Diff^1 (M)$,  $\Delta$ be an $f-$invariant set
and $T_{\Delta}M=E\oplus F$ be a $Df-$invariant splitting on
$\Delta$. Suppose $dim E (x)$ is constant  and   for every  $x\in \Delta$,  the splitting $T_{orb (x)}M=E \oplus F$ is compatible,  then  $$\{\mu\in \mathcal{M}^{h}_{f} (\Delta)| \,ind (\mu)=dim E\}=\{\mu\in \mathcal{M}^{qldh}_{f} (\Delta)| \,ind (\mu)=dim E\}.$$
\end{Prop}

Remark that for quasi-invariant measures, such kind of results are unknown.

\section{Appendix B: Proof of Lemma \ref{Lem:Liao-Gan-shadlem} on  Liao's shadowing}\label{section-proof-shadowingliao}

Before proving let us state some useful lemmas.

Let $\lambda \in(0,1).$ A pair of sequence $\{a_i,b_i\}_{i=1}^n$ of positive numbers is called $\lambda$-hyperbolic if $a_k\leq \lambda$ and $b_k\geq\lambda^{-1}$ for  $k=1,2,\cdots,n.$ A pair of sequences
$\{a_i,b_i\}_{i=1}^n$ of positive numbers is called $\lambda$-quasi-hyperbolic if following three conditions are satisfied:
(1) $\Pi_{j=1}^k a_j\leq\lambda^k,$  (2)  $\Pi_{j=k}^n  b_j\geq\lambda^{k-n-1},$ (3) $\frac{b_k}{a_k}\geq\lambda^{-2}$, for $k=1,2,\cdots,n.$

\begin{Lem}\label{lem-Gan-sequence-quasihyperbolic}
Let $e>0$ be an integer and $\zeta>0$.  Take $\lambda=e^{-\zeta}.$ 
If orbit segment $\{x,n\}$ is $(\zeta,e)$-quasi-hyperbolic with respect to $T_xM=E\oplus F$ and a partition $0=t_0<t_1<\cdots<t_m=n$ with $t_k-t_{k-1}\leq e$, then $\{a_i,b_i \}_{i=1}^m$ is a $\lambda$-quasi-hyperbolic pair of sequence where $a_j=\|Df^{t_j-t_{j-1}}|_{ Df^{t_{j-1}} (E)}\|,\,b_j=m (Df^{t_j-t_{j-1}}|_{ Df^{t_{j-1}} (F)})$.
\end{Lem}

{\bf Proof} Recall that $(\zeta,e)$-quasi-hyperbolic means that \\ (1). $\,\,\,\,\,\, \,\frac 1{t_{k}}\Sigma_{j=1}^k
\log a_j=\frac 1{t_{k}}\Sigma_{j=1}^k
\log\|Df^{t_j-t_{j-1}}|_{ Df^{t_{j-1}} (E)}\|\leq -\zeta,$\\
 (2). $\,\,\,\,\,\,\, \frac 1{t_{m}-t_{k-1}}\Sigma_{j=k}^m
\log\,b_j=\frac 1{t_{m}-t_{k-1}}\Sigma_{j=k}^m
\log\,m (Df^{t_j-t_{j-1}}|_{ Df^{t_{j-1}} (F)})\geq  \zeta,$\\
 (3). $\,\,\,\,\,\,\,\frac
1{t_{k}-t_{k-1}}\log\frac{a_k}{b_k}=\frac
1{t_{k}-t_{k-1}}\log\frac{\|Df^{t_k-t_{k-1}}|_{ Df^{t_{k-1}} (E)}\|}{m (Df^{t_k-t_{k-1}}|_{ Df^{t_{k-1}} (F)})}
\leq -2\zeta,\,\,\, k=1,\, 2, \,\cdots,\, m.$\\
Note that $t_k\geq k$ and $t_m-t_{k-1}\geq m-k+1.$ Fix  one $k=1,2,\cdots,n$. By (1),  $\Pi_{j=1}^k a_j\leq e^{-\zeta t_k}\leq e^{-\zeta k}=\lambda^k.$ Similarly, by (2) one has $\Pi_{j=k}^n  b_j\geq\lambda^{k-m-1}.$ By (3) $\frac{b_k}{a_k}\geq e^{2\zeta(t_k-t_{k-1})}\geq\lambda^{-2},$ since $t_k-t_{k-1}\geq  1.$\qed

\bigskip

A sequence $\{c_i\}^n_{
i=1}$ of positive numbers is called a balance sequence if $\Pi_{j=1}^k c_j\leq 1,  k = 1, 2, \cdots , n -1, $ and $\Pi_{j=1}^n c_j = 1.$
A  balance sequence $\{c_i\}^n_{
i=1}$ is called well-adapted
to a $\lambda$-quasi-hyperbolic sequence pair$\{a_i,b_i\}_{i=1}^n$  if $\{a_i/c_i, b_i/c_i\}_{i=1}^n$ is $\lambda$-hyperbolic.

\begin{Lem}\label{lem-Gan-balance}(\cite{Gan})
Let $\lambda \in(0,1).$ Any $\lambda$-quasi-hyperbolic $\{a_i,b_i\}_{i=1}^n$ of positive numbers has a well-adapted sequence $\{c_i\}^n_{
i=1}$.
\end{Lem}

Let $(X,\|\cdot\|)$ be a Banach space. For any $\eta>0$ denote by $X(\eta)$   the closed ball in $X$ with radius $\eta,$ that is,
$X(\eta):=\{v\in X: \|v\|\leq \eta\}.$ Let $X$ be a direct sum of two closed subspacws $E$ and $F$.
For   $E $ and $F$, recall that the angle  between $E$ and $F$ is defined as $$\angle (E ,F )=\inf\{\|u-v\|:\,u\in E , v\in
F,\, \|u\| \text{ or } \|v\|=1\}.$$ Remark that $0< \angle (E ,F )\leq 1$, since $E,F$ are closed. Let
$$Ang (E ,F )=\inf\{\|u-v\|:\|u\|=\|v\|=1,\,u\in E , v\in
F \}.$$ Remark that $\angle (E ,F )\leq  Ang (E ,F )\leq \sqrt{2}\angle (E ,F ).$

 \begin{Lem}\label{Lem:property-dominated-imply angle2015} (1) Let   $S>0$ be an integer and $\lambda>0.$     Then  there exists $\alpha_S>0$ such that if the splitting $T_yM=E\oplus F$ satisfies $$\frac
{\|Df^S|_{E (y)}\|}{m (Df^S|_{F (y)})}\leq e^{-\lambda S} ,$$ then
$\angle (E (y),F (y))\geq
\alpha_S.$

(2)Let $e>0$ be an integer and $\zeta>0$, there exists $\alpha>0$ such that if orbit segment $\{x,n\}$ is $(\zeta,e)$-quasi-hyperbolic with respect to $T_xM=E\oplus F$ and a partition $0=t_0<t_1<\cdots<t_m=n$ with $t_k-t_{k-1}\leq e$, then $\angle (Df^{t_j}(E) ,Df^{t_j}(F) )\geq \alpha$ holds for any $j=0,1,2,\cdots,m.$
\end{Lem}

{\bf Proof.} By the equivalent relation of  $\angle (E ,F )$ and $ Ang (E ,F )$. We only need to prove the conclusion for $ Ang (E ,F )$.

 (1) Let $$c=\min_{x\in M}\,m (D_xf^S),\,\,\,\,\,\,\,C=\max_{x\in
M}\|D_xf^S\|.$$ Clearly $c,\,\,C \in  (0, +\infty).$ By continuity of
the tangent bundle $T_xM$, there exists real number $\beta_S>0$ such
that if $\|u-v\|<\beta_S,\,\|u\|=\|v\|=1,\,u,\,v\in T_xM$ then
$$
\big{|}\|D_xf^S (u)\|-\|D_xf^S (v)\|\big{|}<c (1-e^{-\lambda}),\,\forall
x \in M,$$ which implies
\begin {equation}\label{eq:anti-dominated-2015}
\frac{\|D_xf^S (u)\|}{\|D_xf^S (v)\|}>
\frac{\|D_xf^S (v)\|-c (1-e^{-\lambda} )}{\|D_xf^S (v)\|}\geq
e^{-\lambda}.
\end{equation}

Suppose point  $y$ satisfies
\begin {equation}\label{eq:dominated-2015}
\frac
{\|Df^S|_{E (y)}\|}{m (Df^S|_{F (y)})}\leq e^{-\lambda S}.
\end{equation}
Then for two vectors $u\in E (y), v\in
F (y)$ with $\|u\|=\|v\|=1,$ we claim that $\|u-v\|\geq \beta_S.$
Otherwise, it holds that $\|u-v\|< \beta_S,$ and thus by the inequality
 (\ref{eq:anti-dominated-2015}),
$$\frac{\|Df^S|_{E (y)}\|}{m (Df^S|_{F (y)})}\geq
\frac{\|D_yf^S (u)\|}{\|D_yf^S (v)\|}\geq e^{-\lambda}\geq e^{-\lambda S},$$ which
contradicts  (\ref{eq:dominated-2015}). So, we have
$$Ang (E (y),F (y))=\inf\{\|u-v\|:\|u\|=\|v\|=1,\,u\in
E (y), v\in F (y)\}\geq \beta_S.$$

(2) Let $e>0$ be an integer and $\zeta>0$, take  $\beta=\min_{1\leq S\leq e}\{\beta_S\}>0$. Then  if orbit segment $\{x,n\}$ is $(\zeta,e)$-quasi-hyperbolic with respect to $T_xM=E\oplus F$ and a partition $0=t_0<t_1<\cdots<t_m=n$, then for any $j=0,1,2,\cdots,m,$ by the third condition from definition of  quasi-hyperbolic (taking $S=t_k-t_{k-1}$ and $y=f^{t_{k-1}}(x)$), one has $$\frac
{\|Df^S|_{E (y)}\|}{m (Df^S|_{F (y)})}\leq e^{-2\lambda S}\leq e^{-\lambda S}.$$ So  $Ang (Df^{t_j}(E) ,Df^{t_j}(F) )\geq \beta_S\geq \beta,$ since $S\leq e.$   \qed

\medskip

In the following, $X_i$ will denote a $W$-dimensional Euclidean space for any integer
$i$. Assume that $X_i$ has a direct sum decomposition $X_i = E_i \oplus F_i.$
Let $Y = \prod_{i=-\infty}^{+\infty}X_i.$ Endowed with the supremum norm $\|v\|=\sup\{\|v_i\|\}\,(v =
(v_i))$, $Y$ is a Banach space. Let us consider the mapping $ \Phi: Y \rightarrow Y$ with the
following form $(\Phi v)_{i+1} = \Phi_i(v_i)$ where $\Phi_i : X_i \rightarrow X_{i+1}.$  Now we recall a theorem on fixed point from \cite{Gan}.

\begin{Lem}\label{lem-Gan-fixedpoint}
For any $\mu\in(0,1),\epsilon>0,\alpha\in(0,1]$ satisfying $\epsilon_1=\frac{2\epsilon (1+\mu)}{\alpha^2(1-\mu)}<1$ and $\eta>0$, denote by $R=R(\mu,\epsilon,\alpha)=\frac{1+\mu}{\alpha(1-\mu)(1-\epsilon_1)},$ $H=2R,\, d_0=\frac\eta H,$  $d\in(0,d_0].$ If $ \Phi: Y(\eta) \rightarrow Y$ has the form
$\Phi_i=L_i+\phi_i:X_i(\eta)\rightarrow X_{i+1}$ and $L_i$ has the block form  $ L_i= \left (
\begin {array}{cc}
A_i&B_i  \\
C_i&D_i
\end {array}
\right) 
$ with respect to the decomposition $X_i=E_i\oplus F_i,$ $\angle(E_i,F_i)\geq \alpha,$ $\max\{\|A_i\|,\,\|D_i^{-1}\|\}\leq \mu,$ $\max\{\|B_i\|,\,\|C_i\|\}\leq\epsilon$, $Lip\phi\leq\frac1H,$  $\|\phi_i(0)\|\leq d$,
 then $\Phi$ has a unique fixed point $ v\in Y(\eta)$ and $\|v\|\leq Hd.$
\end{Lem}

\begin{Thm}\label{Thm-Corollary-Gan-fixedpoint}
For any $K>1,\mu\in(0,1),\epsilon>0,\alpha\in(0,1]$ satisfying $\epsilon_1=\frac{2\epsilon (1+\mu)}{\alpha^2(1-\mu)}<1$ and $\eta>0$, denote by $R=R(\mu,\epsilon,\alpha)=\frac{1+\mu}{\alpha(1-\mu)(1-\epsilon_1)},$ $H=2R,\, d_0=\frac\eta H,$  $d\in(0,d_0].$ If $ \Phi: Y(\eta) \rightarrow Y$ has the form
$\Phi_i=L_i+\phi_i:X_i(\eta)\rightarrow X_{i+1}$,  $L_i$ has the block form  $ L_i= \left (
\begin {array}{cc}
A_i&B_i  \\
C_i&D_i
\end {array}
\right) 
$ with respect to the decomposition $X_i=E_i\oplus F_i,$ $\angle(E_i,F_i)\geq \alpha,$ $\frac1K\leq\|A_i\|,\,\,\|D_i^{-1}\|^{-1}\leq K,$ $\max\{\|B_i\|,\,\|C_i\|\}\leq\frac{\epsilon}K$, $Lip\phi\leq\frac1{KH},$  $\|\phi_i(0)\|\leq d$ and if there is a strictly increasing two-sided sequence $\{N_j\}_{j=-\infty}^{+\infty}$ such that  for each $j$, $\{\|A_i\|,\,\|D_i^{-1}\|^{-1}\}_{i=N_j}^{N_{j+1}-1}$ is $\mu$-quasi-hyperbolic and $\phi_i(0)=0$ for $i=N_j,\cdots,N_{j+1}-2$,
 then $\Phi$ has a unique fixed point $ v\in Y(\eta)$ and $\|v\|\leq Hd.$
\end{Thm}

{\bf Proof.} By Lemma \ref{lem-Gan-balance},  for each $j$, $\{\|A_i\|,\,\|D_i^{-1}\|^{-1}\}_{i=N_j}^{N_{j+1}-1}$ has  a well-adapted sequence $\{h_i\}_{i=N_j}^{N_{j+1}-1}$ such that 
$\{\frac{\|A_i\|}{h_i},\,\frac{\|D_i^{-1}\|^{-1}}{h_i}\}_{i=N_j}^{N_{j+1}-1}$ is $\mu$-hyperbolic, that is,
$\Pi_{i=N_j}^{k} h_i\leq 1,  k = N_j,N_j+1,  \cdots , N_{j+1}-2, $   $\Pi_{i=N_j}^{N_{j+1}-1} h_i = 1,$  and 
$ \frac{\|A_i\|}{h_i}\leq \mu\leq 1, \frac{\|D_i^{-1}\|^{-1}}{h_i}\geq \mu^{-1}\geq 1. $
 It follows that $h_i\geq \|A_i\|\geq \frac1K, \,h_i\leq \|D_i^{-1}\|^{-1}\leq K.$

Let $g_k=\Pi_{i=N_j}^{k} h_k,\, \tilde{L}_k=h_k^{-1} L_k,\,\tilde\phi_k(v)=g_k^{-1}\phi_k(g_{k-1}v)$ (note that $g_{N_{j}-1}=1$ and $g_k\leq 1$), and 
 $\tilde\Phi_k=\tilde L_k+\tilde\phi_k.$ Write $\tilde  L_k= \left (
\begin {array}{cc}
\tilde A_k&\tilde B_k  \\
\tilde C_k&\tilde D_k
\end {array}
\right)
.$  Denote by $\Psi_k=\Phi_k\circ\cdots\circ\Phi_{N_i}$ and $\tilde\Psi_k=\tilde\Phi_k\circ\cdots\circ\tilde\Phi_{N_i}$.  
 Then we have $ \tilde\Psi_k=g_k^{-1}\Phi_k,$  
  and in particular  $\tilde\Psi_{N_{j+1}-1}= \Phi_{N_{j+1}-1},$ since $g_{N_{j+1}-1}=1$.  Note that 
 $ \|\tilde A_k\|=h^{-1}_k\|A_k\|\leq \mu,\,\|\tilde D_k^{-1}\|^{-1}=h^{-1}_k\|D_k^{-1}\|^{-1}\geq \mu^{-1},$  $\max\{\|\tilde B_k\|,\,\|\tilde C_k\|\}=h^{-1}_k\max\{\|\tilde B_k\|,\,\|\tilde C_k\|\}\leq K\max\{\|B_k\|,\,\|C_k\|\}\leq {\epsilon} $,
   $Lip\tilde\phi_k=g_k^{-1}Lip\phi_k g_{k-1}=h^{-1}_kLip\phi_k\leq K\frac1{KH}=\frac1H,$  $\tilde\phi_k(0)=g_k^{-1}\phi_k(0)=0      $ for $k=N_j,\cdots,N_{j+1}-2$ and $\tilde\phi_k(0)=g_k^{-1}\phi_k(0)  =\phi_k(0)    $ for $k=N_{j+1}-1$ since $g_k=1.$ Then we can apply Lemma \ref{lem-Gan-fixedpoint} for $ \tilde\Phi=\{\Phi_k\}: Y(\eta) \rightarrow Y$ to obtain a unique fixed point $\tilde v=\{\tilde v_k\}$ of $ \tilde\Phi$
 and $\|\tilde v\|\leq Hd.$  Let $v_{N_j}=\tilde v_{N_j}$ and for $N_j<k<N_{j+1}$ define $v_{k}=\Phi_{k-1} (v_{k-1})$ inductively. 
 Now we only need to prove that $v$ is a fixed point of $\Phi$ and $\|v\|\leq Hd.$

Since $v_k=\Phi_{k-1} (v_{k-1})=\Psi_{k-1}(v_{N_j})=g_{k-1}\tilde\Psi_{k-1}(v_{N_j})=g_{k-1}\tilde\Psi_{k-1}(\tilde v_{N_j})=g_{k-1}\tilde\Phi_{k-1}(\tilde v_{k-1})=g_{k-1} \tilde v_k, $ one has $\|v_k\|\leq \|\tilde v_k\|\leq Hd.$  Since $v_{N_{j+1}}=\tilde v_{N_{j+1}}=\tilde\Psi_{N_{j+1}-1}(\tilde v_{N_{j}})=\tilde\Psi_{N_{j+1}-1}(v_{N_{j}})= \Psi_{N_{j+1}-1}(v_{N_{j}})=\Phi_{N_{j+1}-1}(v_{N_{j+1}-1})$, $v$ is a fixed point of $\Phi.$ \qed

\begin{Lem}\label{lem-Gan-bycontinuity-bundledistanceclose-normclose}(\cite{Gan}) Let $g\in \Diff^1(M).$
For any $r>1,\gamma>0,\epsilon>0,\alpha\in(0,1]$, there exists $\eta>0$ such that if $x,y\in M$, $T_xM=E_x\oplus F_x$, $T_yM=E_y\oplus F_y$, $\angle(E_x,F_x)\geq\alpha,$ $\angle(E_y,F_y)\geq\alpha,$
$Dg ( \xi ({x}))\cap T^\#M \subseteq U ( \xi ({y )} \cap
T^\#M,\eta)\,\, (\xi=E,F),$ then $\Phi=exp^{-1}_y\circ g \circ exp_x: T_xM(\eta)\rightarrow T_yM$ can be written as $\Phi=L+\phi,$ where $L$ has the form $ L = \left (\begin {array}{cc}
A &B   \\
C &D
\end {array} \right) $ with respect to splittings $T_xM=E_x\oplus F_x$, $T_yM=E_y\oplus F_y$, and $\frac1r\leq \frac{\|A\|}{\|Dg|_{E_x}\|}\leq r,$ $\frac1r\leq \frac{\|D^{-1}\|^{-1}}{m(Dg|_{F_x})}\leq r,$
$\|B\|, \|C\|\leq \epsilon,$ $Lip\phi\leq \gamma.$ Here the norms are induced by the Riemannian metric on $T_xM$ and $T_yM$.
\end{Lem}

For a fixed integer $e\geq1$, let $g=f,f^2,\cdots,f^e$, by Lemma \ref{lem-Gan-bycontinuity-bundledistanceclose-normclose} we can choose a small suitable $\eta>0$ such that  Lemma \ref{lem-Gan-bycontinuity-bundledistanceclose-normclose} holds simultaneously for $f,f^2,\cdots,f^e$. That is,

\begin{Cor}\label{Cor-bundledistanceclose-normclose}
 Let $e\geq1$ be a fixed integer  and $f\in \Diff^1(M).$
For any $r>1,\gamma>0,\theta>0,\alpha\in(0,1]$, there exists $\eta>0$ such that if $x,y\in M$, $T_xM=E_x\oplus F_x$, $T_yM=E_y\oplus F_y$, $\angle(E_x,F_x)\geq\alpha,$ $\angle(E_y,F_y)\geq\alpha,$
$Df^i ( \xi ({x}))\cap T^\#M \subseteq U ( \xi ({y )} \cap
T^\#M,\eta)\,\, (\xi=E,F) $ for some $1\leq i\leq e,$ then $\Phi=exp^{-1}_y\circ f^i \circ exp_x: T_xM(\eta)\rightarrow T_yM$ can be written as $\Phi=L+\phi,$ where $L$ has the form $ L = \left (\begin {array}{cc}
A &B   \\
C &D
\end {array} \right) $ with respect to splittings $T_xM=E_x\oplus F_x$, $T_yM=E_y\oplus F_y$, and $\frac1r\leq \frac{\|A\|}{\|Df^i|_{E_x}\|}\leq r,$ $\frac1r\leq \frac{\|D^{-1}\|^{-1}}{m(Df^i|_{F_x})}\leq r,$
$\|B\|, \|C\|\leq \theta,$ $Lip\phi\leq \gamma.$
\end{Cor}

Now we start to prove
Lemma \ref{Lem:Liao-Gan-shadlem}.

{\bf Proof of Lemma \ref{Lem:Liao-Gan-shadlem}} Let $e\geq1$ be a fixed integer and $\zeta>0.$  Take $\lambda=e^{-\zeta}.$ Denote $K_0=\max_{1\leq i\leq e}\sup_{x\in M}\{\|D_xf^i\|,\|(D_xf^i)^{-1}\|\}$.
 Let $K=2K_0.$  By Lemma
\ref{Lem:property-dominated-imply angle2015} there exists $\alpha>0$ such that  if orbit segment $\{x,n\}$ is $(\zeta,e)$-quasi-hyperbolic with respect to $T_xM=E\oplus F$ and a partition $0=t_0<t_1<\cdots<t_m=n$ with $t_k-t_{k-1}\leq e$, then $\angle (Df^{t_j}(E) ,Df^{t_j}(F) )\geq \alpha$ holds for any $j=0,1,2,\cdots,m-1.$

Let $\mu=\min\{\frac{1+\lambda}2,2\lambda\}$, $r=\frac\mu\lambda$. Then $\mu\in(0,1),\,1<r\leq 2.$ Let $\epsilon>0$ small enough such that $\epsilon_1=\frac{2\epsilon (1+\mu)}{\alpha^2(1-\mu)}<1$.
Denote by $R=R(\mu,\epsilon,\alpha)=\frac{1+\mu}{\alpha(1-\mu)(1-\epsilon_1)},$ $H=2R, $ $\theta=\frac\epsilon K,\,\gamma=\frac1{KH}.$  For above $r>1$ $\alpha,\theta,\gamma$
 take $\eta>0$ small enough such that Corollary \ref{Cor-bundledistanceclose-normclose} holds. 
 Let  $L=KH$ and $d_0=\frac\eta H,$
 fix $d\in(0,d_0].$

Let $\{x_j,n_j\}^{\infty}_{j=-\infty}$ be a $(\zeta,e)$-quasi-hyperbolic $d$-pseudoorbit. We may suppose that for each $j,$ $\{x_j,n_j\} $ is a $(\zeta,e)$-quasi-hyperbolic with respect to  a partition 
  $0=t^{(j)}_0<t^{(j)}_1<\cdots<t^{(j)}_{m_j}=n_j$ with $t^{(j)}_{i+1}-t^{(j)}_i\leq e$ for some $m_j\leq n_j.$ Let $N_j$
 be defined as
  \begin {equation} \label{eq:respective-time-squens-2015}N_j=\begin{cases}
 0,&\text{for }j=0\\
 \sum_{k=0}^{j-1}m_k,&\text{for }j>0\\
 -\sum_{k=j}^{-1}m_k,&\text{for }j<0.
\end{cases}
 \end {equation}

  Let $y_k=f^{t^{(j)}_{k-N_j}}x_j$ and $l_{k}={t^{(j)}_{k+1-N_j}}-{t^{(j)}_{k-N_j}}$ if $N_j\leq k<N_{j+1}$. 
  Note that $1\leq l_{k}\leq e$.  Denote by $X_k=T_{y_k}M,E_k=E_{y_k},\,F_k=F_{y_k},
\, Y = \prod_{i=-\infty}^{+\infty}X_i.$   Let
 $\Phi_k=exp^{-1}_{y_{k+1}}\circ f^{l_{k}} \circ exp_{y_k}: T_{y_{k}}M(\eta)\rightarrow T_{y_{k+1}}M$.  Then by Corollary \ref{Cor-bundledistanceclose-normclose},   $ \Phi=\{\Phi_k\}: Y(\eta) \rightarrow Y$ has the form
$\Phi_k=L_k+\phi_k:X_k(\eta)\rightarrow X_{k+1}$,  $L_k$ has the block form  $ L_k= \left (
\begin {array}{cc}
A_k&B_k  \\
C_k&D_k
\end {array}
\right)
$ with respect to the decomposition $X_k=E_k\oplus F_k,$ $\angle(E_k,F_k)\geq \alpha,$ $\frac1r \|Df^{l_{k}}|_{E_k}\| \leq\|A_k\|\leq r \|Df^{l_{k}}|_{E_k}\| ,$ 
$\frac1r   m(Df^{l_{k}}|_{F_k})^{-1} \leq\|D_k^{-1}\|\leq r m(Df^{l_{k}}|_{F_k})^{-1} ,$ 
$\max\{\|B_k\|,\,\|C_k\|\}\leq\theta=\frac{\epsilon}K$, $Lip\phi\leq\gamma=\frac1{KH}.$ 
Note that
$\frac1K\leq\frac1{rK_0}\leq \|A_k\|\leq rK_0\leq K$ and $\frac1K\leq\frac1{rK_0}\leq \|D_k^{-1}\|\leq rK_0\leq K.$

From  Lemma \ref{lem-Gan-sequence-quasihyperbolic} we know for each $j$, $\{\|Df^{l_{k}}|_{E_k}\|,\,m(Df^{l_{k}}|_{F_k})\}_{k=N_j}^{N_{j+1}-1}$ is $\lambda$-quasi-hyperbolic. Then
for each $j$, $\{\|A_k\|,\,\|D_k^{-1}\|^{-1}\}_{k=N_j}^{N_{j+1}-1}$ is $\mu$-quasi-hyperbolic. From $d$-pseudoorbit orbit we have $\|\phi_i(0)\|\leq d $  and $\phi_i(0)=0$ for $i=N_j,\cdots,N_{j+1}-2$.

 Thus we can use 
 Theorem \ref{Thm-Corollary-Gan-fixedpoint} to get   $\Phi$ has a unique fixed point $ v\in Y(\eta)$ and $\|v\|\leq Hd.$
 Let $z=exp_{y_0}v_0.$ Then $z$ $Hd$-shadows $\{y_k\}$, that is, if $N_j\leq k<N_{j+1}$, then $d(f^{t^{(j)}_k}(x_j),f^{c_j+t^{(j)}_k}(z))\leq Hd,$  where the notation $c_j$ is from (\ref{eq:respective-time-squens}). 
  Notice that $\sup_{1\leq l\leq e} d(f^lx,f^ly)\leq K_0d(x,y)\leq K d(x,y)$ and recall the partition $0=t^{(j)}_0<t^{(j)}_1<\cdots<t^{(j)}_{m_j}=n_j$  satisfies $t^{(j)}_{k+1}-t^{(k)}_i\leq e$.  Then $\forall\,\,
i=0,\,1,\,2,\,\cdots,\,n_j-1$ and $\forall\,\, j \in \mathbb Z$, there is $t^{(j)}_k$ and $l\leq e$  such that $i=t^{(j)}_k+l.$ Thus $ d\big{ (}f^{c_j+i} (x),f^i (x_j)\big{)} \leq K d\big{ (}f^{c_j+t^{(j)}_k} (x),f^{t^{(j)}_k} (x_j)\big{)}\leq KHd=Ld. $

Now suppose 
$\{x_i,n_i\}_{i=-\infty}^{+\infty}$ is periodic, i.e., there
 exists an $m>0$ such that $x_{i+m}=x_i$ and $n_{i+m}=n_i$ for all $i$. Define $w$ by $w_k=v_{N_m+k}$. Then $w,v$ are fixed points of $\Phi$ in  $Y(\eta)$. By uniqueness, $v=w$ and so $f^{c_m}z=z.$ \qed

\section*{Acknowlegements} The research of W. Sun is supported by National Natural Science
Foundation  (No. 11231001).  The research of X. Tian was  supported by National Natural Science Foundation of China (grant no. 11301088) and  Specialized
  Research Fund for the Doctoral Program of Higher Education (No.  20130071120026).


\bigskip
\end{document}